\documentclass[12pt,oneside]{amsart}
\usepackage{stmaryrd}
\usepackage{lmodern}
\usepackage{hyperref}
\usepackage{amsrefs}
\usepackage{amssymb}
\usepackage{color}
\usepackage{tikz-cd}
\usepackage{mathrsfs}
\usepackage{comment}
\usepackage{soul,xcolor}
\setstcolor{red}

\makeindex
\newtheorem{theorem}{Theorem}[section]
\newtheorem{lemma}[theorem]{Lemma}
\newtheorem{proposition}[theorem]{Proposition}

\newtheorem{condition}[theorem]{Condition}

\theoremstyle{definition}
\newtheorem{definition}[theorem]{Definition}
\newtheorem{example}[theorem]{Example}
\newtheorem{remark}[theorem]{Remark}
\numberwithin{equation}{section}

\newcommand {\mat} [1] {\left[\begin{array}{#1}}
\newcommand {\rix}     {\end{array}\right]}

\def\div{{\rm div}}
\def\grad{{\rm grad}}

\def\Tan{{\rm Tan}}

 \def\Dom{{\rm Dom}}
\def\Lip{{\rm Lip}}
\def\loc{{\rm loc}}
\def\per{{\rm per}}
\def\CBB{{\rm CBB}}
\def\USC{{\rm USC}}
\def\LSC{{\rm LSC}}
\def\SCC{{\rm SCC}}

\def\supp{{\rm supp}}

\def\opt{{\rm opt}}
\def\Str{{\rm Str}}

\def\Id{{\rm Id}}
\def\LHS{{\rm LHS}}
\def\RHS{{\rm RHS}}

\def\N{{\mathbb{N}}}
\def\P{{\mathbb{P}}}

\def\R{{\mathbb{R}}}
\def\T{{\mathbb{T}}}
\def\Z{{\mathbb{Z}}}

\newcommand{\sfX}{{\mathsf X}}
\newcommand{\sfY}{{\mathsf Y}}

\newcommand{\sfH}{{\mathsf H}}
\newcommand{\sfU}{{\mathsf U}}

\newcommand{\sfD}{{\mathsf D}}

\newcommand{\sfG}{{\mathsf G}}
\newcommand{\sfL}{{\mathsf L}}
\newcommand{\sfP}{{\mathsf P}}
\newcommand{\sfS}{{\mathsf S}}
\newcommand{\sfd}{{\mathsf d}}
\newcommand{\sfdist}{{\mathsf{dist}}}
\newcommand{\sfm}{{\mathsf m}}
\newcommand{\sfp}{{\mathsf p}}

\newcommand{\bmu}{{\boldsymbol \mu}}
\newcommand{{\bnu}}{{\boldsymbol \nu}}
\newcommand{{\btau}}{{\boldsymbol \tau}}
\newcommand{\bsigma}{{\boldsymbol \sigma}}

\newcommand{\bpartial}{{\boldsymbol \partial}}
\newcommand{{\bxi}}{{\boldsymbol \xi}}
\newcommand{\bpi}{{\boldsymbol \pi}}

\newcommand{\bm}{{\boldsymbol m}}

\newcommand{\bM}{{\boldsymbol M}}
\newcommand{\bN}{{\boldsymbol N}} 
\newcommand{\bH}{{\bf H}} 
\newcommand{\bx}{{\bf x}} 
\newcommand{\by}{{\bf y}}

\begin{document}
 
\title{On a Hamilton-Jacobi PDE theory for hydrodynamic limit of action minimizing collective dynamics}
\def \non{{\nonumber}}

\author{Jin Feng}
\address{Mathematics Department\\
University of Kansas \\
Lawrence, KS 66045, USA.
}


\date{\today}


 \thanks{During initial phase of this research, the author was supported in part by 
LABEX MILYON (ANR-10-LABX-0070) of Universit{\'e} de Lyon, 
within the program ``Investissements d'Avenir" (ANR-11-IDEX-0007) operated 
 by the French National Research Agency (ANR). The author thanks 
 Albert Fathi for consistent helps and encouragements, Luigi Ambrosio 
 and Andrea Mondino for references and useful discussions. 
 He also thanks Toshio Mikami, Mark Peletier, Giuseppe Savar\'e,  Lawrence C. Evans, Wei Cheng, Nicola Gigli, Zhao Dong and Philippe Thieullen for helpful conversations on related mathematical issues over the years.}

\begin{abstract}
We establish multi-scale convergence theory for a class of Hamilton-Jacobi PDEs in space of probability measures. They arise from context of hydrodynamic limit of $N$-particle deterministic action minimizing (global) Lagrangian dynamics.  

From a Lagrangian point of view, this can also be viewed as a limit result on two scale convergence of action minimizing  probability-measure-valued paths. However, we focus on the Hamiltonian formulation here mostly. We derive and study convergence of the associated abstract but scalar Hamilton-Jacobi equations, defined in space of probability measures. There is an infinite dimensional singular averaging structure within these equations. 
We develop an indirect variational approach to apply finite dimensional weak K.A.M. theory to such infinite dimensional setting here. With a weakly interacting particle assumption, the averaging step only involves that of individual particles, which is implicitly but rigorously treated using the weak K.A.M. theory. Consequently, we can  close the above mentioned averaging step by identifying limiting Hamiltonian, and  arrive at a rigorous convergence result on solutions of the nonlinear PDEs in space of probability measures.
     
In technical development parts of the paper, we devise new viscosity solution techniques regarding projection of equations with a submetry structure in state space,  multi-scale convergence for certain abstract Hamilton-Jacobi equations in metric spaces, as well as comparison principles for equations in space of probability measures.  The space of probability measure we consider is a special case of Alexandrov metric space with curvature bounded from below. Since some results are better explained in such metric space setting, we also develop some techniques in the general settings which are of independent interests.  
\end{abstract}

\maketitle
\tableofcontents
\clearpage

\section{Introduction}\label{intro}
In a 2020 Oberwalfach online workshop, jointly with Toshio Mikami, the author informally presented a multi-scale convergence theory for a class of Hamilton-Jacobi PDEs in space of probability measures~\cite{FM20}. We proposed using such a theory to understand hydrodynamic limit behaviors of certain deterministic {\em global action minimizing} collective dynamics. Such an approach is variational by nature. The purpose of this work is to provide details of this program. In order to highlight main ideas by addressing difficult issues one step at a time, in this paper, we only consider situations where the microscopic Hamiltonian particles are weakly interacting through a mean field. Strongly interacting cases are of course more interesting. But they will involve additional thermodynamic type variables, such as various forms of free energies. Hence additional ideas and techniques are needed, we leave them for future explorations.

\subsection{Models of collective dynamics and their limits, formal derivations}
We give an overview on our program at the heuristic level. First, we go through some direct and formal calculations. Then, we describe the program and its relation with existing approaches.
\subsubsection{A particle model}
Let $\sfH:=\sfH(q,p): \R^d \times \R^d \mapsto \R$ be a given function. We define an associated Hamiltonian ODE system 
\begin{align*}
\begin{cases}
 \dot{q}  = \nabla_p \sfH(q,p), \\
  \dot{p}  = - \nabla_q \sfH(q,p).
\end{cases}
\end{align*}
Such system describes movements of a single particle with position variable $q$ and generalized momentum variable $p$. We are interested in collective behaviors of $N$ such independently acting particles. For such purpose, we introduce another Hamiltonian function 
  \begin{align}\label{HN}
 \sfH_N({\bf q},{\bf p}) :=    \frac1N \sum_{i=1}^N   \sfH(q_i,p_i),
 \end{align}
where the $({\bf q, p}) := \big( (q_1,p_1), \ldots, (q_N, p_N)\big) \in (\R^d \times \R^d)^N$. We consider phase space $(\R^d \times \R^d)^N$ as a normed vector space with a weighted norm
\begin{align}\label{Nnorm}
| ({\bf q},{\bf p}) |_N^2: = \frac1N \sum_{i=1}^N (|q_i |^2 + |p_i |^2).
\end{align} 
Let $\nabla_{N, {\bf p}}, \nabla_{N,{\bf q}}$ denote gradients for functions in the Hilbert space $\big((\R^d \times \R^d)^N, |\cdot|_N\big)$. Note that this is different than
the $\nabla, \nabla_{p}, \nabla_{q}$ we just used a few lines earlier, which are gradients when $\R^d$ or $(\R^d \times \R^d)^N$ is endowed with the usual (un-weighted) Euclidean norm $|\cdot|$. The conversion relation is that (identified as vectors in $(\R^d \times \R^d)^N$)
\begin{align}\label{defNabN}
\nabla_N f({\bf q,p}) = N \nabla f({\bf q,p}).  
\end{align}
With these notations and relations clarified, we can write the $N$-particle level Hamiltonian dynamic for the $\sfH_N$ as  
\begin{align}\label{HNODE}
\begin{cases}
\dot{q}_i = \nabla_{N,p_i} \sfH_N({\bf q,p}) = \nabla_{p_i} \sfH(q_i,p_i), \\
\dot{p}_i  = - \nabla_{N,q_i} \sfH_N({\bf q,p})=  - \nabla_{q_i} \sfH(q_i, p_i).  
\end{cases}
\end{align} 
The special case where $\sfH(q,p) := \frac12 |p|^2 - U(q)$ deserves particular mention as the corresponding dynamic describes Newton's law in classical mechanics. 

If we introduce additional pair-wise interaction terms, then with proper rescaling on the time-space variables, one can, in a certain regime, derive the Boltzmann equation in physicist's  hand-waving ways (e.g. Chapters 3.3 and 3.4 in Kardar~\cite{Kard07}). Such derivation requires a physical ansatz (molecular chaos) whose rigorous derivation does not exist so far, and has been a ground for controversies.  
We skip the Boltzmann's kinetic limit and go directly to the next one -- a hydrodynamic limit. In addition to the particle number, we also introduce anther scaling parameter $\epsilon: = \epsilon_N \to 0$ as $N \to \infty$. We  speed up time and scale down space to introduce macroscopic variables:
\begin{align*}
{\bf x}(t):={\bf x}_\epsilon(t) = \epsilon{\bf q}(\frac{t}{\epsilon}), \quad {\bf P}(t):={\bf P}_\epsilon(t) ={\bf p}(\frac{t}{\epsilon}).
\end{align*}
Then
\begin{align*}
{\bf \dot{x}}(t) = {\bf \dot{q}}(\frac{t}{\epsilon}), \quad {\bf \dot{P}}(t) = \frac{1}{\epsilon}{\bf \dot{p}}(\frac{t}{\epsilon}),
\end{align*} 
and the Hamiltonian dynamic \eqref{HNODE} can be re-written as a closed system
\begin{align}\label{Hflow}
\begin{cases}
\dot{x}_i(t)  = \nabla_{2} \sfH\big( \frac{x_i}{\epsilon}, P_i\big),   \\
\dot{P}_i(t)  =    - \frac{1}{\epsilon}  \nabla_{1} \sfH\big(\frac{x_i}{\epsilon}, P_i \big), 
\end{cases}
\end{align}
where $\nabla_1 \sfH(q,p) := \nabla_q \sfH(q,p)$ 
and $\nabla_2 \sfH(q,p) := \nabla_p \sfH(q,p)$. A singular perturbation structure now emerges. 

We consider a situation where macro-scale particles don't escape to infinity too quickly in the macroscopic scale. We add a confinement potential  $U \in \R^d \mapsto \R$ by replacing single particle level Hamiltonian from $\sfH$ into  $\sfH_{U}(q,p):=\sfH(q,p) - U(\epsilon q)$. We also introduce pairwise weak interaction modeled through a smooth even function $V \in C^2(\R^d;\R)$ (hence $\nabla V(0) =0$ in particular): 
for given ${\bf q} :=(q_1, \ldots, q_N)$, we define
\footnote{Note that we are considering a (mathematically easier) scaling where the perturbative term 
$U(x):= U(\epsilon q)$ is small in the microscopic $({\bf q,p})$-level but non-ignorable in the macroscopic $({\bf x,P})$-level, similarly for the interaction terms $V$. }
\begin{align*}
\sfH_{U,V}(q,p; {\bf q}) := \sfH(q,p) - U(\epsilon q) 
 - \frac1N \sum_{j=1}^N V\big(\epsilon(q-q_j)\big), 
 \quad \forall (q,p) \in \R^d \times \R^d.
\end{align*}
Equivalently, we can also directly start with a rescaled Hamiltonian in the Hilbert space $\big( (\R^d \times \R^d)^N; | \cdot |_N\big)$ by
\begin{align}\label{defHN}
  H_N \big( {\bf x}, {\bf P} \big) &:= \frac1N \sum_{i=1}^N \Big( \sfH\big( \frac{x_i}{\epsilon}, P_i \big)  - U(x_i) - \frac1N \sum_{j=1}^N V( x_i - x_j )\Big) \\
 & = \frac1N \sum_{i=1}^N \sfH_{U,V}( \frac{x_i}{\epsilon}, P_i; \frac{{\bf x}}{\epsilon}).
 \nonumber
\end{align}

Letting $\epsilon:=\epsilon_N \to 0$ at appropriate speed as $N \to \infty$. A classical problem in this context is to understand how point-particle level Hamiltonian dynamics given by the 
$H_N({\bf x}, {\bf P})$ converge to a continuum level dynamic modeled by conservation laws in continuum mechanics. We approach this issue from an indirect manner.
Instead of considering solution to Hamiltonian ODE given by the $\sfH_N$ converges to system of partial differential equations of the Euler equation type, we consider convergence of respective Lagrangian actions. Actions are scalar quantities. In continuum, the Lagrangian action is defined over probability-measure-valued curves. In classical particle mechanics, the method of generating functions characterizes this action using Hamilton-Jacobi partial differential equations.  In continuum level, these equations become defined in space of probability measures. In summary, we are lead to a mathematical problem on convergence of a class of multi-scale Hamilton-Jacobi PDEs in space of probability measures. Within the class of (global) action minimizing Lagrangian dynamics, the associated Hamiltonian operators giving the PDEs enjoy a nonlinear version of the maximum principle. We can devise abstract (derivative free) viscosity solution techniques, based upon maximum principle, to study such type of Hamilton-Jacobi equations in space of probability measure, including convergence of solutions. This is precisely what we do in this paper. We summarize main results and explain some of the new technical tools we develop for achieving these results in Section~\ref{Sec1:main}.
 
\subsubsection{One scale mean-field limit}\label{Mfield}
We explain why Hamilton-Jacobi equations in space of probability measures are natural in our context.  To simplify, we do not yet consider the multi-scale hydrodynamic limit in this subsection. We step back to consider a single scale mean-field limit for now. That is, we choose $\epsilon =1$ fixed while $N \to \infty$. To simplify, we even take the $U=0= V$.

Our problem has an obvious symmetry -- the $N$-particle Hamiltonian $\sfH_N$ is invariant under permutation on particle indices.  We consider the position vector ${\bf q}$ modulo permutation of particle index as a physical observable variable. To formalize, let $\sfG_N$ denote the permutation group on $N$ indices, it acts on ${\bf q} \in (\R^d)^N$ by
\begin{align*}
 \pi {\bf q} := (q_{\pi(1)}, \ldots, q_{\pi(N)}), \quad \forall \pi \in \sfG_N, \quad {\bf q} \in (\R^d)^N.
\end{align*}
We use quotient space $(\R^d)^N / \sfG_N$ as the physical space of observables with a quotient metric
\begin{align}\label{dNstar}
\sfd^*_N({\bf q}^*, \hat{\bf q}^*) :=\inf_{\substack{{\bf q} \in {\bf q}^*\\\hat{\bf q} \in \hat{\bf q}^*}}| {\bf q} - \pi {\hat{\bf q}}|_N := \inf_{\pi \in \sfG_N} \sqrt{\frac1N \sum_{i=1}^N |q_i - \hat{q}_{\pi(i)}|^2}.
\end{align}
Note that such space is non-smooth. Singularities arise whenever particles collide (i.e. 
\begin{align*}
 q_{i_1} = q_{i_2} = \ldots = q_{i_l}, \qquad \exists \quad i_1 \neq i_2 \neq \ldots \neq i_{i_l}).
\end{align*}
 This singular space can be identified with space of empirical probability measures  
 \begin{align}\label{defXN}
\sfX_N:=\Big\{  \sigma(dq) := \frac1N \sum_{i=1}^N \delta_{q_i}(dq) : q_i \in \R^d \Big\}.
\end{align}
Let $\sfX:=\mathcal P_2(\R^d)$ be the space of all probability measures on $\R^d$ with finite second moment. We denote $\sfd$ the Wasserstein order-2 metric on this space
(e.g. Ambrosio, Gigli and Sav\'are~\cite{AGS08}). Then $\sfX_N$ is a closed subspace in the Wasserstein metric space $(\sfX, \sfd)$. 
If we identify ${\bf q}^*, \hat{\bf q}^*$ with empirical probability measures $\sigma, \hat{\sigma}$ respectively, then by a result on optimization of linear functional over convex set, and by Birkhoff's theorem on doubly stochastic matrix, we have
\begin{align}\label{isoddN}
 \sfd (\sigma, \hat{\sigma}) = \sfd^*_N({\bf q}^*, \hat{\bf q}^*).
\end{align}
For a short proof, see for instance the second example on page 5 of Villani~\cite{Villani03}. The identity map of $N$-point empirical probability measure to the Wasserstein space, therefore, is an isometric embedding map
\begin{align}\label{embeta}
\sfp_N: \sfX_N \mapsto \sfX.
\end{align}  
In summary, the quotient metric space $\big( (\R^d)^N /\sfG_N, \sfd_N^* \big)$ can be identified isometrically with the $(\sfX_N, \sfd)$ which is isometrically embedded into $(\sfX, \sfd)$. It is useful to keep the following facts in mind:  the space $(\sfX, \sfd)$ is a complete separable metric space, it is a geodesic Alexandrov space with curvature bound below by zero (e.g. Ambrosio, Gigli and Savar\'e~\cite{AGS08} and Villani~\cite{Villani09}), we have in current literature a well-developed modern theory on optimal mass transport (\cites{AGS08,Villani09}) and more generally, a first order calculus theory on Alexandrov metric space\cite{AKP19}. We recall that the $\sfd$ topologizes a version of weak convergence of probability measures which is narrow convergence plus convergences of moments up to the second order.  

We introduce Hamiltonian operator
\begin{align*}
 \sfH_N f({\bf q}) := \sfH_N\big({\bf q}, \nabla_{N} f({\bf q}) \big). 
\end{align*}
Next, we shall rewrite Hamiltonian operators in the $\sigma$-coordinate. We consider a class of smooth test functions 
\begin{align}\label{Poly}
f(\sigma):= \psi(\langle \vec{\varphi}, \sigma\rangle) &:=  \psi(\langle \varphi_1, \sigma\rangle, \ldots, \langle \varphi_K, \sigma\rangle) \\
& = \psi\Big( \frac1N \sum_{i=1}^N \varphi_1(q_i), \ldots, \frac1N \sum_{i=1}^N \varphi_K(q_i) \Big) =:f({\bf q}), \quad \forall \varphi_i \in C^1(\R^d);  \nonumber
\end{align}
and denote
\begin{align*}
\frac{\delta f}{\delta \sigma} = \sum_{k=1}^K \partial_k \psi(\langle \vec{\varphi}, \sigma\rangle) \varphi_k.
\end{align*}
We have
\begin{align*}
\nabla_{N, q_i} f({\bf q}) = N \nabla_{q_i} f({\bf q})= \sum_{k=1}^K \partial_k \psi\big(\langle \vec{\varphi}, \sigma\rangle\big) \nabla \varphi_k(q_i) =   \nabla \frac{\delta f}{\delta \sigma}(q_i).
\end{align*}
Therefore
\begin{align}\label{HPN}
 H_N f(\sigma)  := \sfH_N\big( {\bf q},  \nabla_N f({\bf q}) \big)  =  \int_{\R^d} \sfH\big(q, \nabla \frac{\delta f}{\delta \sigma} (q)\big)  \sigma (dq).
\end{align}
These $H_N$s converges to 
\begin{align}\label{Hmean}
 H f(\sigma) :=  \int_{\R^d \times \R^d} \sfH(q, P) {\boldsymbol \nu}(dq, dP), 
\end{align}
where the
\begin{align*}
 \bnu (dq, dP) := \delta_{\nabla \frac{\delta f}{\delta \sigma} (q)}(dP)\sigma(dq) =  \grad_\sigma f \in \mathcal P(\R^d \times \R^d).
\end{align*}
The last equality above follows by abstract first order calculus result in Alexandrov spaces (e.g. Lemma~\ref{PolyGrad}). 

At least formally, the above $H$ generates a Lagrangian dynamics given by probability measure-valued curves. As in the well-known finite dimensional setting, with reasonable convexity assumptions on the $\sfH$, Lagrangian actions can be introduced as a dual problem and evolution of a corresponding global action minimization problem be studied here. Such minimal action can be described by a Cauchy problem 
\begin{align*}
 \partial_t u(t, \sigma) = H u(t, \cdot)(\sigma).
\end{align*}
Informally, the $u(t)$s are canonical transformations of a Hamiltonian dynamics in $\mathcal P(\R^d)$ generated by the $H$. We will consider such $u$ through an abstraction viscosity solution theory for Hamilton-Jacobi equations in $\mathcal P(\R^d)$. 
\footnote{When action critical point instead of global minimizer is considered, the issue of what notion of solution to use becomes much more subtle. In these cases, the viscosity solution may not always be correct for the context. See further discussions in Section~\ref{Sec1:relHDL}.}
Through dynamical programming argument, the solution $u(t)$ is given by an operator nonlinear semi-group which acts on functions on $\mathcal P(\R^d)$. Through an adaptation of the Crandall-Liggett theory to the viscosity context (e.g. Feng and Kurtz~\cite{FK06}), we may simply consider a resolvent version of the operator equation 
\begin{align}\label{HJmf}
 f - \alpha H f =h.
\end{align}
In the above, the $H$ is a first order nonlinear differential operator in Wasserstein space $(\sfX,\sfd)$, $h \in C_b(\sfX)$ and $\alpha >0$ are given,  $f$ is considered as a solution. Section 4 of Ambrosio and Feng~\cite{AF14} explains a critical role played by a metric geometry nature of the space $(\sfX, \sfd)$, for equations of such type. It is known that the class of test functions \eqref{Poly} is not sufficiently large for making sense of the equation. We leave this point for more discussions in later parts of this introduction section.

Let 
\begin{align*}
f_N - \alpha H_N f_N = h_N.
\end{align*}
We assume convergence of the $h_N$ to $h$ in proper senses. Because that $H_N$ converges to the $H$, we expect the solutions $f_n$s also converge to the $f$. Such result will bring convergence of actions on measure-valued curves as $N \to \infty$.  Within a probability context, Feng and Kurtz~\cite{FK06} has generalized the Barles-Perthame techniques to abstract viscosity solution theory in general metric spaces.  Combined with \cite{AF14}, those techniques can be directly applied to study convergence of a large class of canonical transforms (or actions of minimizer type) in a single-scale mean-field limit context. 

\subsubsection{Hydrodynamic limit}\label{hydroIntro}
In this paper, we consider the hydrodynamic limit problem which is more complex than the above, because it involves more than one scales.  Next, we explain how to obtain convergence of Hamiltonian operators $H_N$ in such situation. For a nice class of test functions $f$, we look for a sequence of test functions $f_N$ such that $f_N \to f$ and $H_N f_N \to Hf$ for some operator $H$. That is,  we verify that $H_N \to H$ in a properly defined operator graph convergence sense. We also identify the limit $H$. 

Let $x_i$s be those in \eqref{Hflow}, we introduce some new coordinates
\begin{align} \label{Sec1:mCoord}
\begin{cases}
& \rho(t; dx):= \frac1N \sum_{i=1}^N \delta_{x_i(t)}(dx), \quad  \sigma(t;dq):= \frac1N \sum_{i=1}^N \delta_{\epsilon^{-1}x_i(t)}(dq), \\
& \qquad \mu(t; dx, dq):= \frac1N \sum_{i=1}^N \delta_{x_i(t), \epsilon^{-1}x_i(t)}(dx, dq).
\end{cases}
\end{align}
It follows that $\pi^1_{\#} \mu(t) = \rho(t)$ and $\pi^2_{\#} \mu(t) = \sigma(t)$. See Section~\ref{appmass} for the notation of $\pi^i$ and $\pi^i_\#$ which are commonly used in optimal transport literature.
For test functions of the form \eqref{Poly} with the $\sigma$ replaced by $\rho$ and the $q$ replaced by $x$,  
we have 
\begin{align*}
H_N f(\rho) &:= H_N\big( {\bf x}, \nabla_{N, {\bf x}} f \big)
=   \int_{\R^d} \Big( \sfH\big( \frac{x}{\epsilon}, \nabla \frac{\delta f}{\delta \rho}(x) \big) - U(x) - V*\rho(x) \Big) \rho(dx). 
\end{align*}
Next, we reveal a hidden separation of scale structure by splitting the coordinates ${\bf x}:=  \big({\bf x},  {\bf q(x)}\big):= \big({\bf x}, \epsilon^{-1} {\bf x}\big)$ according to different scales. 
We take $\varphi_i$s to be of the form
\begin{align*}
\varphi(x):= \varphi_\epsilon(x) :=\phi(x; {\epsilon}^{-1}x), \quad \text{ where } \phi:=\phi(x, q) \in C^1(\R^d \times \R^d), \nabla_x \phi, \nabla_q \phi \in C_b,
\end{align*}
and write
\begin{align*}
\mu(dx,dq):= \rho(dx) \delta_{\epsilon^{-1} x}(dq).
\end{align*}
Note that $\langle \varphi, \mu \rangle = \langle \varphi, \rho \rangle$. 
We re-write test functions 
\begin{align}\label{polyf}
 f(\rho):= f({\bf x}) 
 := \psi(\langle \vec{\varphi}, \rho\rangle) = \psi( \langle \vec{\phi}, \mu\rangle)=: f(\mu).
\end{align}
Then
\begin{align}
 \frac{\delta f}{\delta \rho} (x)&=  \sum_k  \partial_k \psi \big(\langle \vec{\varphi}, \rho \rangle\big)\varphi_k(x) 
= \sum_k  \partial_k \psi \big(\langle \vec{\varphi}, \rho \rangle\big)\phi_k\big(x, \frac{x}{\epsilon}\big) \label{Sec1:GDel1} \\
 \nabla \frac{\delta f}{\delta \rho} (x) 
& = \sum_k \partial_k \psi\big(\langle \vec{\varphi}, \rho \rangle \big) 
 \big(\nabla_x \phi_k(x, \frac{x}{\epsilon}) + \frac{1}{\epsilon} \nabla_q \phi_k (x, \frac{x}{\epsilon})\big) \label{Sec1:GDel2} \\
 & =\big(\nabla_x \frac{\delta f}{\delta \mu}\big) \big(x, \frac{x}{\epsilon}\big)+ \frac{1}{\epsilon} \big(\nabla_q \frac{\delta f}{\delta \mu} \big)\big(x, \frac{x}{\epsilon}\big).
 \nonumber
 \end{align}
Consequently
\begin{align}\label{Heps}
H_N f(\rho) &=  \int_{\R^d}\Big\{ \sfH\Big( q, \nabla_x \frac{\delta f}{\delta \mu}(x, q) 
 + \frac{1}{\epsilon} \nabla_q \frac{\delta f}{\delta \mu}(x,q) \Big) \\
  & \qquad \qquad \qquad \qquad - U(x) -V*\rho(x) \Big\} \mu(dq, dx) + o_N(1). \nonumber
 \end{align}

The above calculation suggests that we should take a class of perturbed test functions
\begin{align}\label{pertfg}
f_N(\rho):=  f(\rho) +\epsilon g(\mu),
\end{align}
where
\begin{align}
f(\rho) & := \psi_0(\langle \varphi_1, \rho\rangle, \ldots, \langle \varphi_K, \rho\rangle), 
\quad  \varphi_k:= \varphi_k(x) \in C^1_c(\R^d), k=1,2,\ldots, K; \label{f0def}\\
g(\mu)&:= \psi_1(\langle \phi_1, \mu\rangle, \ldots, \langle \phi_K, \mu\rangle), \quad \phi_k:=\phi_k(x,q) \in C^2_c(\R^d \times \R^d). \label{gdef}
\end{align}
Note that the $g$ can be understood as both a function of the $\mu$-variable, or a function of the $\rho$-variable:
\begin{align*}
g(\mu) &= \psi_1(\langle \varphi_{\epsilon, 1}, \rho\rangle, \ldots, \langle \varphi_{\epsilon,K}, \rho\rangle), \quad \varphi_{\epsilon,k}(x):=\phi_k\big(x, \frac{x}{\epsilon}\big).
\end{align*}
Therefore,
\begin{align*}
H_N f_N(\rho) &=   \int_{\R^{2d}} \Big\{ \sfH \Big( q, \nabla_x \frac{\delta f}{\delta \rho}(x) 
 + \epsilon (\nabla_x \frac{\delta g}{\delta \mu})\big(x,q\big) 
 +  (\nabla_q \frac{\delta g}{\delta \mu})\big(x,q \big) \Big) \\
 & \qquad \qquad \qquad - U(x) - V*\rho(x) \Big\} \mu(dq,dx)  + o_N(1).  
\end{align*}
 In fact, it is sufficient to consider a special sub-class of the $g$s in the following forms
\begin{align*}
 g(\rho) = \langle \varphi_\epsilon, \rho \rangle, 
\end{align*}
where
\begin{align}\label{varphieps}
 \varphi_\epsilon(x) = \phi\big(x, \nabla \frac{\delta f}{\delta \rho}(x); \frac{x}{\epsilon}\big), \quad \phi(x,P;q) \in C^1(\R^d \times \R^d \times \R^d).
\end{align}
Then  
\begin{align*}
& H_N f_N(\rho) \\
 &= \int_{\R^d} \Big\{ \sfH\Big( q, \nabla_x \frac{\delta f}{\delta \rho}(x) 
 +  \epsilon  (\nabla_x \phi)\big(x,\nabla_x \frac{\delta f}{\delta \rho}(x); q\big) \\
 & \qquad \qquad +   \epsilon\big(  D^2_x \frac{\delta f}{\delta \rho}(x)\big) \cdot \nabla_P \phi\big(x, \nabla_x \frac{\delta f}{\delta \rho}(x); q\big)       \\
 & \qquad \qquad  \qquad + \nabla_q  \phi\big(x, \nabla_x \frac{\delta f}{\delta \rho}(x); q\big)  \Big) - U(x) - V*\rho(x) \Big\} \mu(dx,dq) +o_N(1)    \\
&= \int_{\R^d \times \R^d} \Big(\sfH \big(q, P+ \nabla_q \phi(x,P;q)\big) -U(x) - V*\rho(x) \Big)   \delta_{\nabla_x \frac{{\delta f}}{\delta \rho} (x)}(dP)\mu(dx, dq)
 + o_N(1).
\end{align*}
Suppose that we can solve an auxiliary PDE problem  (where the $q$ is considered as a variable and the $P$ as a parameter) 
\begin{align}\label{cell}
\sfH\big(q, P+ \nabla_q \phi(q)\big) = \bar{\sfH}(P).
\end{align}
By solution, we mean the pair $(\bar{\sfH}, \phi) \in \R \times C(\R^d)$, where  $\bar{\sfH}:=\bar{\sfH}(P)$ is a number and  $\phi:=\phi_{P} :=\phi (P;\cdot)$ is a function, both are indexed by 
$P \in \R^d$. Now, if we write
\begin{align}\label{barHxP}
 \bar{\sfH}(x,P; \rho) := \bar{\sfH}(P) - U(x) - V*\rho(x),
\end{align}
then we can conclude that
\begin{align*}
f_N \to f, \quad H_N f_N \to H f
\end{align*}
with
\begin{align}\label{Sec1:EffH}
 H f(\rho) 
 = \int_{\R^{2d} } \bar{\sfH} (x, P; \rho) \bmu_f(dx,dP),  
  \text{ where }
 \bmu_f (dx, dP) :=  \delta_{\nabla_x \frac{\delta f}{\delta \rho}}(dP) \rho(dx)
\end{align}

The auxiliary PDE \eqref{cell} is known as a ``cell equation". It is at center of homogenization theory for Hamilton-Jacobi equations and weak KAM theory for Hamiltonian dynamical systems in finite dimensions. There is extensive literature on these topics. See, for instance, Lions, Papanicolaou and Varadhan~\cite{LPV87}, Fathi~\cites{Fa97a, Fa97b, Fa98a, Fa98b,FathiBook},  E~\cites{E91,E99}, as well as Evans and Gomez~\cites{EG01, EG02a,EG02b}. There are also unpublished works of Ma\~n\'e,  which were carried out further in Contreras-Iturriaga-Paternain-Paternain~\cite{CIPP98}. See Chapter 9 of \cite{FathiBook} for more on Ma\~n\'e's point of view, and related references.  In general, we know that there is no smooth $C^1$ class of $\phi$ satisfying \eqref{cell} at every $q$, but a generalized  solution in the sense of viscosity solution can be found. Moreover, the effective Hamiltonian $\bar{\sfH}$ is always unique and it has several variational representations. See Appendix~\ref{varbarH} for details.
Viscosity solution $\phi$ for \eqref{cell} can be non-unique. However, the largest critical sub-solution satisfying certain inequalities can be characterized explicitly through some dynamical system quantities. See  Proposition 1.3 of  Davini, Fathi, Iturriaga and Zavidovique~\cite{DFIZ16}, which is summarized in Lemma~\ref{Cell:exist}.  For the above hydrodynamic limit problem, we find a way to proceed rigorously without explicitly using the $\phi$. There are situations where higher order hydrodynamic limits are relevant. In those cases, knowing the form of such particular solution becomes important. We do not pursue higher order hydrodynamics in this paper.   

To proceed rigorously, we will discover that the above derivation looks nicer than it really is. The class of $f$ as defined in \eqref{polyf} is inadequate for studying well-posedness for \eqref{HJmf} as a Hamilton-Jacobi equation in space of probability measures. We need to consider a broader class of functions $f$ and re-run the above asymptotic in a roundabout way in order to make such things rigorous. This is because that the space of probability measures has a singular nature -- it is an Alexandrov metric space with non-negative synthetic curvature. We will make use of a set of well developed first order calculus by Ambrosio-Gigli-Savar\'e~\cite{AGS08}, which is related to a theory of analysis in Alexandrov space as described by recent publication Alexander, Kapovitch and Petrunin~\cite{AKP19}. Due to its importance and technical nature, we devote a whole Section~\ref{AlexCal} to recall and refine some existing results which are needed later in the paper.  Ambrosio and Feng~\cite{AF14} uses metric geometry perspective to formulate and prove well-posedness of Hamilton-Jacobi equation in space of probability measures. See Section 4 of that paper for a closely related Hamilton-Jacobi equation from continuum mechanics. 

Consideration of viscosity solution for Hamilton-Jacobi equation in space of probability measures appeared at least as early as in the late 90s, where Feng~\cite{Fe99} derived a specific PDE model, and called for a need for corresponding comparison principle in order to understand probabilistic large deviation theory for Fleming-Viot stochastic processes. However, a successful comparison principle for similar type of equations did not appear until later. In Feng and Katsoulakis~\cite{FK03}, \footnote{See Reference [66] in the first edition of Ambrosio-Gigli-Savar\'e~\cite{AGS08} and [78] in the second edition.} the authors realized that Otto-calculus in modern mass transport theory was sufficient for putting many estimates in the right order for comparison of solutions. The published form of that work did not appear until years later~\cite{FKa09}. In between, more extensive version of the results (including convergence theories) were developed in context of probabilistic large deviation theory by Feng and Kurtz~\cite{FK06} in a book form --- see Chapters 6,7, 9, and 13.3 there.  See also related publications of Feng and Nguyen~\cite{FN12}, Feng and Swiech~\cite{FS13}, Feng, Mikami and Zimmer~\cite{FMZ21} and references therein. 
The type of equations with Hamiltonian operator~\eqref{Sec1:EffH}, however, is a somewhat different story. Even though the abstract convergence theories still apply, the comparison principle cannot be directly obtained from above mentioned method. 
A missing component, in such context, was found by Ambrosio and Feng~\cite{AF14} using metric-geometry inspired techniques.
This  lead to comparison principle for a new class of equations -- see Section 4 of that paper. In coming sections of this paper, we will further develop key observations and techniques made in the above references.
Regarding different attempts of defining viscosity solution suited for equations in space of probability measures, also see Gangbo, Nguyen and Tadurascu~\cite{GNT08}. 
This definition does not lead to uniqueness. Later improvements were made by Gangbo and collaborators. See in particular Gangbo and Swiech~\cite{GS14} for the use of metric space analysis, which appeared in the same time as \cite{AF14}, although the explicit use of tangent cone techniques in Wasserstein space is absent in \cite{GS14}.

\subsubsection{Relations with other approaches of hydrodynamic limit}\label{Sec1:relHDL}
It is important to emphasize that, as the title reflects, we are not treating hydrodynamic limits for Hamiltonian dynamics of all initial values. Implicitly, we only consider those paths which correspond to global action minimizing Lagrangian dynamics. Such restriction has to do with the viscosity solution techniques we will use, and the finite dimensional weak K.A.M.(Kolomogrov-Arnold-Moser)  results that we will invoke. We also need to emphasize that our approach differs from the traditional program in that we focus on convergence of actions. It is useful to investigate other notion(s) of solution for Hamilton-Jacobi equation which is (are) proper for the general hydrodynamic limit problem. There has been a history on alternative notions of solution studied in Hamiltonian dynamical system literature. We do not digress the topic further in this article to pursue that direction. We mention that, the formal verification of multi-scale Hamiltonian convergence that we described above remains valid, no matter we deal with action minimizers or critical points. Once a new notion of solution is developed in general canonical transform context, we expect to repeat the procedures developed here with new techniques to treat those Hamiltonian dynamics involving critical actions that are not minimal ones.  

With the above points in mind, we introduce relevant Hamilton-Jacobi equations in space of probability measures and study their convergence through an enhanced notion of viscosity solution, in main text of this paper.

Regarding the $2$-scale convergence as mentioned above, we explore physical structure of the problem to reduce such seemingly infinite-dimensional/infinite-particle averaging problem to that of only one-particle/finite dimensional problem. As a result, we could invoke well-developed finite dimensional weak K.A.M. averaging theories to replace the classical ergodic theory step. For connecting the PDE averaging theories with trajectory based ergodic type arguments, see Proposition 3.1 and Theorem 4.1 in Evans and Gomes~\cite{EG01} for brief discussions.  From a physical point of view, it is interesting to note that only ``micro-canonical ensemble" given by the effective Hamiltonian $\bar{\sfH}$ is used here, and that this suffices to characterize the limiting problem.  It is perhaps important to point out that this feature has to do with our initial model assumptions on treating weakly interacting models. Hamiltonian particle models with strong interaction are of course more physically interesting. It is a natural next step to consider.

\subsection{Notations, assumptions and main results}\label{Sec1:main}
Throughout, infimum of a function over empty set is considered $+\infty$ and supremum over empty set is $-\infty$. For a generic metric space $(\sfX,\sfd)$, we denote $M(\sfX)$, $B(\sfX)$,$C(\sfX)$, $C_b(\sfX)$ respectively the spaces of measurable, bounded, continuous, bounded continuous functions. By $\USC(\sfX;\bar{\R})$ and $\LSC(\sfX; \bar{\R})$, we mean upper-semicontinuous and lower-semicontinuous functions on $\sfX$ with value in extended real $\bar{\R}:=\R \cup\{\pm \infty\}$. Similarly, we define $M(\sfX; \bar{\R})$, etc.

In \eqref{defXN}, we use $\sfX_N$ to denote the space of $N$-particle empirical probability measures. It can be used to identify with $(\R^d)^N/\sfG_N$.
Such identification is unique up to an isometry.  We define a surjective projection map $\sfp_N :  (\R^d)^N \mapsto \sfX_N$ by
\begin{align*}
\sigma := \sfp_N({\bf x}):= \frac1N \sum_{i=1}^N \delta_{x_i}, \quad \forall {\bf x}:=(x_1,\ldots, x_N) \in (\R^d)^N.
\end{align*}
The $\sfX_N \subset \sfX$, with the Wasserstein order-2 space $(\sfX, \sfd)$ a separable, metrically complete; and geodesic and non-negatively curved Alexandrov space. For introduction to Alexandrov spaces in general, see Burago, Burago and Ivanov~\cite{BBI01}, and Bridson and Haefliger~\cite{BH99}. For first order calculus and analysis on Alexandrov spaces, see Alexander, Kapovitch and Petrunin~\cite{AKP19}. For specific properties and analysis of the space $(\sfX, \sfd)$ as an Alexandrov space of non-negative curvature, see Ambrosio, Gigli and Savar\'e~\cite{AGS08} and the thesis of Gigli~\cite{Gigli04}.  In Section~\ref{AlexCal}, we give a very brief review on selected techniques for analysis in these spaces, which we also use in later sections.
 
Let $\sfH := \sfH(q,p) : \R^d \times \R^d \mapsto \R$. We assume the following structural conditions. Generalities are not pursued in these conditions, to avoid being side tracked by non-essential issues.  
For instance, periodicity assumption in the following Tonelli type condition can be eliminated under some technical conditions along the lines of Ishii and Siconolfi~\cite{IshiiSi20}. 
\begin{condition}\label{PerCND} \
\begin{enumerate}
\item $\sfH:=\sfH(q,p) \in C(\R^d \times \R^d)$;
\item for each $p \in \R^d$ fixed, $q \mapsto \sfH(q,p)$ is periodic in the sense that
\begin{align*}
 \sfH(q+ k, p) = \sfH(q, p), \quad \forall k:=(k_1, \ldots, k_d) \in \Z^d;
\end{align*}  
\item for each $q \in \R^d$, $p \mapsto \sfH(q,p)$ is convex in $\R^d$;
\item $\sfH$ is uniformly coercive in the $p$-variable 
\begin{align*}
  \liminf_{|p| \to +\infty} \inf_{q \in \R^d }   \sfH(q, p)   = +\infty.
\end{align*}
\item for each $q \in \R^d$ fixed, $p \mapsto \sfH(q, p)$ is super-linear:
\begin{align*}
 \liminf_{|p| \to \infty} \frac{\sfH(q,p)}{1+|p|} = +\infty.
\end{align*}
 \end{enumerate}
\end{condition}
By a periodic function $\varphi : \R^d \mapsto \R$, we mean  
\begin{align*}
 \varphi(q+ k) = \varphi(q), \quad \forall k:=(k_1,\ldots, k_d) \in \Z^d. 
\end{align*} 
We use notations $C_\per(\R^d)$, $\USC_{\per}(\R^d)$ and $\LSC_{\per}(\R^d)$ respectively for set of functions which is continuous and periodic, upper semi-continuous and periodic, lower semicontinuous and periodic, etc.
We also denote $\T^d:= \R^d / \Z^d$.

We define  
\begin{align}\label{sfLdef}
\sfL(q, \xi) := \sup_{p \in \R^d} \big( \xi \cdot p - \sfH(q, p) \big), \qquad \forall (q, \xi) \in \R^d \times \R^d.
\end{align} 
Through properties of Legendre transform, we know that such $\sfL$ is convex and super-linear in $\xi$ as well (provided $p \mapsto \sfH$ is), moreover, $\sfL \in \LSC(\R^d \times \R^d)$ and
\begin{align*}
 \inf_{(q, \xi) \in \R^d \times \R^d} \sfL(q,\xi) \geq - \sup_{q \in \R^d} \sfH(q,0) > -\infty.
\end{align*}
Therefore, (e.g. Theorem 12.2 on page 104 of Rockafellar~\cite{Rock70}),
\begin{align}\label{LegeDual}
 \sfH(q,p) = \sup_{\xi \in \R^d} \big( p \xi - \sfL(q, \xi) \big).
\end{align}
The following concept comes from Bangert~\cite{Bang99}:
\begin{definition}\label{Sec1:closedM}
We call a probability measure $\mu \in \mathcal P(\R^{2d})$ {\em closed}, if 
\begin{enumerate}
\item $\int_{\R^d \times \R^d}   |\xi| \mu(dx, dq) <\infty$
\item $\langle \mu, D_q \varphi \cdot \xi \rangle =0$, for every $\varphi \in C^\infty_\per(\R^d)$.
\end{enumerate}
\end{definition}
We define, for each $v \in \R^d$, 
\begin{align}\label{effsfLdef}
 \bar{\sfL}(v) & := \inf \Big\{ \int_{\R^{2d}} \sfL(q, \xi) \mu(dq, d\xi) : \forall \mu \in {\mathcal P}(\R^{2d}) 
 \text{ is closed with }   \int_{\R^{2d}} \xi \mu(dq,d\xi) = v \Big\}. 
\end{align}
Then, under Condition~\ref{PerCND}, by Proposition~\ref{avgHrep} in Appendix~\ref{varbarH}, we have
\begin{align}\label{Sec1:barHdef}
\bar{\sfH}(P) & := \sup_{v \in \R^d} \big\{ v P - \bar{\sfL}(v) \big\} \\
& = \inf_{\varphi \in C^\infty_{\per}(\R^d)} \sup_{q \in \R^d} \sfH(q, P+ \nabla_q \varphi) = \sup_{\varphi \in C^\infty_{\per}(\R^d)} \inf_{q \in \R^d} \sfH(q, P+ \nabla_q \varphi).
\nonumber
\end{align}
The above implies another equivalent way of introducing the pair $(\bar{\sfL},\bar{\sfH})$. We can first define 
 \begin{align*}
\bar{\sfH}(P)= \inf_\phi \sup_q \sfH(q, P+ \nabla_q \phi),
\end{align*}
 then the $\bar{\sfL}$ with expression~\eqref{effsfLdef} follows through Legendre transform.

There are $U$ and $V$ in the model \eqref{Heps}. Without pursuing generality, we impose the following requirements on them.
\begin{condition}\label{U0CND}
$U \in \Lip(\R^d; \R_+)$ has sub-linear growth at infinity. That is, there exists a concave, increasing, sub-linear function $\beta: \R_+ \mapsto \R_+$ such that 
$0 \leq U(x) \leq \beta(|x|)$.   
\end{condition}
\begin{condition}\label{VCND}
$V \in \Lip(\R^d) \cap L^\infty(\R^d)$.  
\end{condition}
For convenience, we also impose almost quadratic growth conditions (in $p$) on $\sfH$.
\begin{condition}\label{Sec1:Tone2}
There exists $c, C>0$ such that
\begin{align}\label{amQforH}
 -c+  C^{-1}   |p|^2 \leq \sfH(q,p) \leq c+ C |p|^2.
\end{align}
\end{condition}
Through Legendre transform, the above condition implies an almost quadratic growth estimate for the Lagrangian as well:  
\begin{align}\label{Sec1:sfLlb}
   -c +\frac{1}{4C}|\xi|^2 \leq \sfL(q, \xi) \leq  c+ \frac{C}{4} |\xi|^2.
\end{align}
Condition~\ref{Sec1:Tone2} already implies the growth conditions $p \mapsto \sfH(q,p)$ in Tonelli type Condition~\ref{PerCND}. We use it here to confine some technical arguments within the framework of $2$-Wasserstein space.  In principle, one can relax it if needed by using multiple Wasserstein metrics with mixed orders. 

We define one-particle level Lagrangian in an environment given by all particles
${\bf q}:= (q_1, \ldots, q_N)$: 
\begin{align*}
 \sfL_{U,V}(q,\xi; {\bf q}) :=  \sfL(q, \xi) + U(\epsilon q) + \frac1N \sum_{j=1}^N V\big(\epsilon(q-q_j)\big), \quad \forall (q, \xi) \in \R^d \times \R^d.
\end{align*}
We also define $N$-particle level Lagrangian for collective dynamic in the $({\bf q, \dot{q}})$ coordinate
\begin{align*}
 \sfL_N({\bf q}, \bxi) := \frac1N \sum_{i=1}^N  \sfL_{U,V}(q_i, \xi_i; {\bf q});
\end{align*}
or equivalently,  in the $({\bf x, \dot{x}})$-coordinate: 
\begin{align*}
 L_N({\bf x, v}) := \frac1N \sum_{i=1}^N  
    \sfL_{U,V}\big( \frac{x_i}{\epsilon}, v_i ; \frac{\bf{x}}{\epsilon}\big).
\end{align*}
It follows that  
\begin{align*}
L_N({\bf x,v}) = \sup \big\{ \langle {\bf P, v} \rangle_N - H_N({\bf x,P}) : {\bf P} \in (\R^d)^N \big\}.
\end{align*}
\subsubsection{Actions and equations induced by finite particle collective Lagrangian dynamics}
For the collective dynamics in time interval $[0,T]$, {\em action of a path} ${\bf z}(\cdot) :=\{ {\bf z}(t) : t \in [0,T] \} \in C\big([0,T]; (\R^d)^N \big)$ is defined as
\begin{align}\label{Sec1:CollAct}
 A_N[{\bf z}(\cdot),T]:= 
 \begin{cases}
 \int_0^T L_N\big( {\bf z}(t), \dot{\bf z}(t)\big) dt, \quad 
   {\bf z}(\cdot)\in AC\big([0,T]; (\R^d)^N \big), \\
 +\infty, \quad \text{otherwise}.
 \end{cases}
\end{align}
{\em Minimal action} with prescribed initial position ${\bf x}_0$ and terminal position ${\bf x}_1$ is written as  
\begin{align*}
 A_N[{\bf x_0, x_1} ;T]:= \inf \Big\{ A_N[{\bf z}(\cdot),T] : {\bf z}(0) ={\bf x}_0, {\bf z}(T) ={\bf x}_1, {\bf z}(\cdot)  \in C\big([0,T]; (\R^d)^N \big) \Big\}.
\end{align*}

Let ${\mathfrak h}_N \in C\big((\R^d)^N)$ be such that 
$\sup_{(\R^d)^N} {\mathfrak h}_N<\infty$, the following quantity
\begin{align*}
S_N(t,{\bf x}):= S_N(t, {\bf x}; {\mathfrak h}_N):= \sup_{{\bf y} \in (\R^d)^N} \big({\mathfrak h}_N({\bf y})  - A_N[{\bf x,y};t]\big).
\end{align*}
is a viscosity solution to Cauchy problem 
\begin{align}\label{Sec1:MCauN}
\begin{cases}
 \partial_t S_N(t, {\bf x})= H_N\big({\bf x}, \nabla_N S_N(t,{\bf x})\big) \\
  S_N(0,{\bf x}) = {\mathfrak h}_N({\bf x}).
\end{cases}
\end{align}
It is also useful to consider another related quantity: for every $\alpha >0$, we write
\begin{align}\label{Sec1:fiDVal}
 {\mathfrak f}_N({\bf x}) & := {\bf R}_{N;\alpha} {\mathfrak h}_N({\bf x}) \\
 & := \sup \Big\{ \int_0^\infty e^{-\frac{s}{\alpha}}
 \Big( \frac{{\mathfrak h}_N}{\alpha} \big({\bf z}(s)\big) - L_N\big({\bf z}(s), \dot{\bf z}(s)\big)\Big)ds :  \nonumber \\
 & \qquad \qquad {\bf z}(0) ={\bf x},   
  {\bf z}(\cdot) \in AC\big((0, \infty);(\R^d)^N\big) 
  \cap C\big([0,\infty);(\R^d)^N\big) \Big\}. \nonumber
\end{align}
We denote class of continuous functions with sub-linear growth at infinity and bounded from above:
\begin{align*}
{\mathscr C}_N:=\big\{ {\mathfrak h} \in C\big((\R^d)^N\big):
 \sup_{(\R^d)^N} {\mathfrak h}<+\infty,
\quad \lim_{|x| \to +\infty} \frac{|{\mathfrak h}(x)|}{1+|x|} =0 \big\}.
\end{align*}
Suppose ${\mathfrak h}_N \in {\mathscr C}_N$, then it is known that the above ${\mathfrak f}_N \in {\mathscr C}_N$ and it is the unique viscosity solution to 
\begin{align}\label{Sec1:fiDHJ}
 {\mathfrak f}_N({\bf x}) - 
  \alpha H_N\big({\bf x}, \nabla_N {\mathfrak f}_N({\bf x})\big)
    = {\mathfrak h}_N({\bf x}).
\end{align}
Moreover, the ${\bf R}_{N;\alpha} : {\mathscr C}_N \mapsto {\mathscr C}_N$ is a nonlinear contractive map and 
\begin{align}\label{Sec1:SNSG}
S_N(t, {\bf x}; {\mathfrak h}_N) =
  \lim_{n \to \infty} {\bf R}_{N;n^{-1}}^{[nt]} {\mathfrak h}_N({\bf x}).
\end{align}

\subsubsection{Actions and equations by effective collective dynamics of infinite particles}
We recall the definition of $\bar{\sfL}$ as in \eqref{effsfLdef}, and its Legendre transform $\bar{\sfH}$ given thereafter. By convexity arguments, we also have
\begin{align}\label{Sec1:barsfL}
\bar{\sfL}(v) = \sup_{P \in \R^d} \big( Pv - \bar{\sfH}(P)\big).
\end{align}

Let $\sfX:={\mathcal P}_2(\R^d)$ be the space of probability measures with finite second moments, and $\sfd$ be the Wasserstein order-2 metric on $\sfX$ (see Chapter 7.1 of Ambrosio, Gigli and Savar\'e~\cite{AGS08}). The $(\sfX, \sfd) \in \CBB(0)$ is an Alexandrov space with a notion of synthetic curvature bounded from below by $0$ -- See Section~\ref{AlexCal} for more.  
With a slight abuse of notation, we also write
\begin{align}\label{Sec1:LUV}
 \bar{\sfL}_{U,V}(x,v ; \rho) : =  \big(\bar{\sfL}(v) + U(x)\big) + (V* \rho)(x), 
  \quad \forall (x, v) \in \R^d \times \R^d, \forall \rho \in \sfX;
\end{align}
and write
\begin{align}\label{Sec1:Lbnu}
 L(\bnu) := \int_{\R^{2d}} \bar{\sfL}_{U,V}(x,v; \pi_\#^1 \bnu)   \bnu(dx, dv), \quad \forall \bnu \in {\mathcal P}(\R^d \times \R^d),
 \end{align}
where $(\pi_\#^1 \bnu)(dx) = \bnu(dx;\R^d)$. See Section~\ref{appmass} for more regarding optimal mass transport theory and notations.

Let probability-measure-valued curve $\rho(\cdot):=\{ \rho(t) : t \in [0,T]\} \in AC([0,T]; \sfX)$. See Chapter 1 of \cite{AGS08} for definition of absolute continuous curves in such metric space setting. Following analysis in Alexandrov space literature, we 
 introduce tangent cones (see beginning part of Section~\ref{AlexCal}). Following Alexander, Kapovitch and Petrunin~\cite{AKP19}, we define $\frac{d}{dt} \rho(t)$ as an element in this cone (see Definition~\ref{veloCurv} in Section~\ref{AlexCal} in this paper). Following Ambrosio, Gigli and Savar\'e~\cite{AGS08}, we explicitly identify the tangent cone of the $\sfX$ as a subset of $\mathcal P_2(\R^d \times \R^d)$ -- see Lemma~\ref{Sec2:TanIden}) and other related material in Section~\ref{appmass}. Then, we have 
\begin{align}\label{Sec1:tDrho}
 \frac{d}{dt} \rho(t) = \bnu(t):= \bnu(t; dx,dv) \in \Tan_{\rho(t)} \subset {\mathcal P}_2(\R^{2d}), \quad \text{a.e.} \quad t \in [0,T].
\end{align}
We define action of the path $\rho(\cdot)$ by 
\begin{align}\label{Sec1:ActMV}
A[\rho(\cdot); T]:= 
\begin{cases}
 \int_0^T L\big(\bnu(s)\big) ds, \quad \text{ when } \rho(\cdot) \in AC([0,T]; \sfX), \\
 +\infty, \quad \text{ otherwise};
\end{cases}
\end{align}
and action with prescribed initial and terminal boundary conditions: 
\begin{align*}
A[\rho_0, \rho_1;T]:= \inf \Big\{ A[\sigma(\cdot);T] : \sigma(\cdot) \in C([0,T];\sfX), \sigma=\rho_0, \sigma(T) =\rho_1 \Big\}.
\end{align*}
For $h \in C(\sfX)$ with $\sup_\sfX h<+\infty$, we write
\begin{align*}
S(t, \rho):= S(t, \rho; h):= \sup_{\gamma \in \sfX} \big( h(\gamma) - A[\rho, \gamma;t] \big),
\end{align*}
then we expect that the $S$ solves the following Cauchy problem in a proper viscosity solution sense
\begin{align}\label{Sec1:MCau}
\begin{cases}
 \partial_t S(t, \rho)= H \big(\rho, \grad_\rho S(t,\rho)\big) \\
U(0,\rho) = h(\rho),
\end{cases}
\end{align}
where the $H$ is some kind of duality to the above given $L$. Indeed, there is problem here and we will come back to this point, in a bit.

It is also useful to consider another quantity which is related to the action in \eqref{Sec1:ActMV}: for every $\alpha >0$, we write
\begin{align}\label{Sec1:RaMV}
 {\bf R}_\alpha h(\rho) 
 & := \sup \Big\{ \int_0^\infty e^{-\frac{s}{\alpha}}
 \Big( \frac{h}{\alpha} \big(\sigma(s)\big) -   L\big(\bnu(s)\big)\Big)ds : 
  \sigma(\cdot) \in AC\big([0, \infty);\sfX\big)  \\
 & \qquad \qquad \qquad \text{ with } \sigma(0) =\rho,   \frac{d}{ds} \sigma(s) = \bnu(s) \in \Tan_{\sigma(s)} \text{ a.e. } s \Big\}. \nonumber
\end{align}
Using well known optimal control arguments, boosted by abstractions to a metric space setup here (see Lemma 8.18 in \cite{FK06} for instance), we have that 
${\bf R}_\alpha : {\mathscr C} \mapsto {\mathscr C}$ is a nonlinear contractive map on some properly defined subset ${\mathscr C} \subset C(\sfX)$. Moreover,
\begin{align}\label{Sec1:SSG}
 S(t, \rho; h) = \lim_{n \to \infty} {\bf R}_{n^{-1}}^{[nt]} h(\rho).
\end{align}
It is also expected that $f:={\bf R}_\alpha h$ is the unique viscosity solution to 
\begin{align}\label{Sec1:MRes}
 f(\rho) - \alpha H f(\rho) = h(\rho),
\end{align}
for some properly defined Hamiltonian operator $H$.
 
It is a non-trivial issue to rigorously define a PDE in singular  Alexandrov spaces. By singular, we mean space with tangent cone at certain points possibly become not a linear space. That is, a vector in the tangent cone may not have an opposite in the same cone. See Lemmas~\ref{oppolem} and \ref{oppolem2} for more. It turns out, compared with using differential of a function (Definition~\ref{diffDef}) in such cases, we will lose information if we use the notion of gradient (Definition~\ref{GradDef}) -- See Lemmas~\ref{ddist} and \ref{gradDpol}, see also Lemmas~\ref{gdist} ,~\ref{sizesD} and \ref{Sec2:AllEq}. The class of test functions which we can develop calculus also needs to be specified. We will choose distance-squared functions as building blocks -- see Section~\ref{Simfun} and in particular, classes of simple functions $\mathcal S^{\pm}$ and $\mathcal S^{\pm, \infty}$ as specified there.
Let $f \in {\mathcal S}^{+, \infty} \cup {\mathcal S}^{-,\infty}$, by Lemma~\ref{Dmdist} and Remark~\ref{Sec2:Dmdext}, Lemma~\ref{WassDSS} and Remark~\ref{dggamma}, 
$d_\rho f$ exists and can be explicitly expressed. 
It is tempting to introduce yet anther single-valued Hamiltonian operator 
\begin{align}\label{Sec1:bHDef}
{\bH} f(\rho):=  \sup \Big\{ (d_{\rho} f)(\bnu) - \bar{\sfL}(\bnu) : \bnu \in \Tan_\rho \Big\},
\quad \forall f \in{\mathcal S}^{+, \infty} \cup {\mathcal S}^{-,\infty},
\end{align}
and formulate the above PDE problems \eqref{Sec1:MCau} and \eqref{Sec1:MRes} using the operator ${\bf H}$. However, it is difficult to justify a rigorous asymptotic analysis about the limit from $H_N$ to $\bf H$. For this reason, we actually use some estimate of the ${\bf H}f$ from above and below by introducing several pairs of Hamiltonian operators, see Section~\ref{Sec7:GradHam}. In particular, 
let $({\mathbb H}_0, {\mathbb H}_1)$ be defined according to  \eqref{Sec6:bbH0}) and \eqref{Sec6:bbH1}; respectively $({\bf H}_0, {\bf H}_1)$ be defined according to \eqref{Sec7:H0} and \eqref{Sec7:H1}. 
By Lemma~\ref{Sec7:Hcom} and display \eqref{Sec7:IneqH1s}, then   
\begin{align*}
{\bf H}f_0 \leq {\bf H}_0 f_0 \leq {\mathbb H}_0 f_0, \quad \forall f_0 \in {\mathcal S}^{+,\infty}, \text{ and }
{\bf H}f_1 \geq {\bf H}_1 f_1 \geq {\mathbb H}_1 f_1, \quad \forall f_1 \in {\mathcal S}^{-,\infty}.
\end{align*}

\subsubsection{Main results}   
This paper consists of mainly two parts. Sections~\ref{AlexCal}, \ref{VisMetProj}, \ref{CnvHJ} develop some calculus and viscosity solution theories in 
general metric spaces. Sections~\ref{finPartH}, \ref{Sec6}, \ref{Sec7}, \ref{Sec8} and \ref{Sec9} apply these theories to the hydrodynamic limit problem presented in this introduction.  To highlight our main goal,  we only summarize next the two scale hydrodynamic limit results obtained in Theorems~\ref{Sec8:MainThm1} and \ref{Sec9:MainThm2}. See Section~\ref{Sec1:Tech}, however, for some comments concerning abstract arguments in the first part of the paper.

We consider a sequence of functions 
\begin{align*}
 (\{ {\mathfrak h}_N \}_{N \in \N}, h ) \subset C\big((\R^d)^N\big) \times \ldots \times C(\sfX),
\end{align*}
and introduce a special class 
\begin{align*}
{\mathcal C} := \big\{ (\{ {\mathfrak h}_N \}_{N \in \N}, h ) \text{ satisfying those requirements in Definition~\ref{Sec8:ClassC} } \big\}.
\end{align*}
In particular, $(\{ {\mathfrak h}_N \}_{N \in \N}, h ) \in {\mathcal C}$ implies convergence of ${\mathfrak h}_N \to h$ in the following sense: for each $\rho \in \sfX$ and  ${\bf x}:= {\bf x}_N:=(x_1^N, \ldots, x_N^N) \in (\R^d)^N$ with 
 \begin{align*}
 \rho_N(dx):= \frac1N \sum_{j=1}^N \delta_{x_j^N}(dx)
\end{align*}
satisfying $\lim_{N \to \infty} \sfd(\rho_N, \rho)=0$, we have
\begin{align*}
 \lim_N {\mathfrak h}_N ({\bf x}_N) = h(\rho).
\end{align*}

We say a bounded from above function $h: \sfX \mapsto \R \cup\{-\infty\}$
has at most sub-linear growth to $-\infty$, if the following holds: 
\begin{align*}
 h(\rho) \geq -\beta \circ \sfd(\rho, \delta_0), \quad \forall \rho \in \sfX,
\end{align*}
for some non-negative sub-linear function $\beta \in C(\R_+)$. 
If $(\{ {\mathfrak h}_N \}_{N \in \N}, h ) \in {\mathcal C}$, then the $h$ has at most sub-linear growth to $-\infty$.

 \begin{theorem}\label{Main}[Limit theorem to Hamilton-Jacobi PDEs]
Let $\alpha >0$ be arbitrary but fixed. Let $h \in C(\sfX)$ with $\sup_\sfX h<+\infty$, and have at most sub-linear growth to $-\infty$. Moreover, we assume that the $h$ has a modulus of continuity with respect to the $\sfd$-metric, on every $\sfd$-metric balls with finite radius; and that the $h$ is $1$-Wasserstein metric $\sfd_{p=1}$-upper semi-continuous in the $2$-Wasserstein metric space $\sfX$ (see Definition~\ref{Sec6:weakUSC}).
Then there is at most one function $f \in C(\sfX)$, with 
$\sup_\sfX f <+\infty$ and with sub-linearly growing to $-\infty$, such that it is both a sub-solution, in the point-wise strong viscosity sense of Definition~\ref{PtVisDef}, to equation 
 \begin{align*}
 (I - \alpha {\mathbb H}_0) f \leq h,
\end{align*}
as well as a super-solution, in the point-wise strong viscosity sense, to equation 
\begin{align*}
(I -\alpha {\mathbb H}_1) f \geq h.
\end{align*}
Moreover, such $f={\bf R}_\alpha h$ as given by variational representation \eqref{Sec1:RaMV}.

Furthermore, let $H_N:=H_N({\bf x,P})$ be defined according to \eqref{defHN} with the $\epsilon:=\epsilon_N \to 0$ as $N \to \infty$. 
Suppose that $(\{ {\mathfrak h}_N \}_{N \in \N}, h ) \in {\mathcal C}$.
Let ${\mathfrak f}_N:={\mathfrak f}_N({\bf x}):={\bf R}_{N,\alpha} {\mathfrak h}_N$ be the value function defined in \eqref{Sec1:fiDVal}, which is also the unique viscosity solution to partial differential equation \eqref{Sec1:fiDHJ}. Then, for each $\rho \in \sfX={\mathcal P}_2(\R^d)$ and  ${\bf x}:= {\bf x}_N:=(x_1^N, \ldots, x_N^N) \in (\R^d)^N$ with 
 \begin{align*}
 \rho_N(dx):= \frac1N \sum_{j=1}^N \delta_{x_j^N}(dx)
\end{align*}
satisfying $\lim_{N \to \infty} \sfd(\rho_N, \rho)=0$,  we have limit
 \begin{align*}
\lim_{N \to \infty} {\mathfrak f}_N({\bf x}_N) = f(\rho),
\end{align*}
with limiting function the one given by $f={\bf R}_\alpha h$ in the above.
\end{theorem}
 
As an important step of proving the above result, we also prove a comparison principle for sub- and super- solutions of respective equations given by the operators ${\bf H}_0$ and ${\bf H}_1$. For precise statement, see Theorem~\ref{Sec7:CMP}.

\subsection{Comments on physical interests of the main results}
A physical significance of the above result is that it implies convergence of actions.
In the hydrodynamic limit scale, $A_N[\cdot]$ defined in \eqref{Sec1:CollAct} characterizes deterministic finite-particle-level action minimizing collective Lagrangian dynamic, and $A[\cdot]$ in \eqref{Sec1:ActMV} describes a continuum-level effective Lagrangian dynamic defined on probability-measure-valued curves.   From convergence results about 
${\mathfrak f}_N={\bf R}_{N;\alpha} {\mathfrak h}_N \to f := {\bf R}_\alpha h$ (whenever ${\mathfrak h}_N \to h$), we expect convergence of solution semigroups $S_N(\cdot) \to S(\cdot)$ defined in \eqref{Sec1:SNSG} and \eqref{Sec1:SSG} (See Remark~\ref{Sec9:SGRmk}). By dynamical programming principle, and by the arbitrariness of the $h$ and ${\mathfrak h}_N$s, we conclude $A_N[\cdot] \to A[\cdot]$.
Since justification for the above arguments are more or less standard, given that the paper is already long, we do not provide details of these proof, but merely state informally the following expected result:
\begin{align*}
 A_N[\cdot] \to A[\cdot] \text{ in the sense of  $\Gamma$-convergence.}
\end{align*} 

Let $\rho(\cdot)$ be an $A[\cdot]$-action minimizing path. Following \eqref{Sec1:tDrho}, we write $\frac{d}{dt} \rho(t) = \bnu(t) \in \Tan_{\rho(t)}$. 
We define bulk velocity field   
\begin{align}\label{Sec1:bVelo}
u(t,x):= \int_{\mathbb R} v \bnu(t;dv|x).
\end{align}
Assuming $\nabla \bar{\sfL}$ is well defined, we also introduce an enhanced phase space measure
\begin{align}\label{Sec1:Mm}
\bmu(t;dx,dv,dP):=    \delta_{\nabla \bar{\sfL}(v)}(dP)  \bnu(t;dx,dv) ;
\end{align}
and momentum measure (which is a $\R^d$-valued signed-measure)
\begin{align*}
\bm(t,dx):=   \int_{(v,P) \in \R^d \times \R^d} P    \bmu(t; dx,dv,dP);
\end{align*}
and momentum-flux measure (which is a $d \times d$ matrix-valued signed measure):
\begin{align}\label{Sec1:MFm}
\bM(t,dx):=   \int_{(v,P) \in \R^d \times \R^d} (v \otimes P)  \bmu(t;dx,dv,dP) ,
\end{align}
where the $v \otimes P:= \big( v_i P_j\big)_{i,j=1,\ldots,d}$ means a square matrix. 
Following the perturbative computations in Section 3.2 of Feng and Nguyen~\cite{FN12}, at least formally, minimizer of the action $A[\cdot]$ satisfy  hyperbolic system of partial differential equations:  
\begin{align}\label{Sec1:Euler}
\begin{cases}
\partial_t \rho + \div_x (\rho u) =0, \\
\partial_t \bm + \div_x \bM  = \rho \nabla_x (U + 2 V*\rho).
\end{cases}
\end{align}
In the above, by $\div_x \bM$ we mean a vector whose $i$-th component is $\sum_{j=1}^d \partial_j M_{ij}$ where the $\bM=(M_{ij})_{i,j=1,\ldots, d}$.
Of course, $P \mapsto \bar{\sfH}(P)$ may generally not be differentiable at some points, so is the function $v \mapsto \bar{\sfL}(v)$. Such situation corresponds to  phase transition.

In this paper, we only considered {\em globally action minimizing} dynamics. This is because of the use of viscosity solution theory. A challenging task for the future is to generalize the arguments here to those Hamiltonian dynamics which are not global action minimizing, but rather just critical points of the action functional. This requires a new notion of solution for Hamilton-Jacobi equation. Even in the context of finite dimension,  this is an important but under-developed field at current time. Once such a PDE theory is ready, the principal ideas of this paper (namely, multi-scale Hamiltonian convergence implies action convergence, hence corresponding dynamical trajectories), shall still apply. 

In this paper, we also only considered weakly interacting particles. When we take multi-scale strong interacting particles into consideration, non-trivial pressure term in \eqref{Sec1:Euler} will appear. More importantly, we expect the microscopic mechanical energy will be partitioned into two parts -- a slowly varying part which remains to be energy of mechanical nature (described by particle density and pressure), and another highly oscillating part of the energy which will become {\em heat}. This disorganized form of energy is expected to naturally introduce the notions of entropy, temperature, and other forms of free energies etc, into such derivation. Challenging works are still needed for clarifying our understanding on such a picture. For instance, how mathematically rigorous arguments such as weak KAM type averaging on more complicated Hamiltonian operators can be used to justify formal physical arguments historically made using micro-canonical, canonical and grand-canonical ensembles.  We hope the framework proposed in this article provides a testing ground for pursuit of these very interesting directions in the future.   

\subsubsection{A toy example of one dimensional ideal gas}
In the context of our main result, we take the special case of 
\begin{align*}
d=1, \quad {\mathsf H}(q,p) =\frac12 |p|^2 -\sfU(q), \text{ with } \min_{q \in [0,1]} \sfU(q) =0; 
\text{ and } U=0, V=0.
\end{align*}
Let $c_\sfU:= \int_{q \in [0,1]} \sqrt{\sfU}dq$. When $|P|>c_\sfU$, we denote  
$\lambda:=\lambda(P)>0$ the unique solution to 
\begin{align*}
 |P| = \int_{q \in [0,1]}  \sqrt{2(\lambda+ U(q))} dq. 
\end{align*}
Lions, Papanicoulou and Varadhan~\cite{LPV87} identified that
\begin{align*}
\bar{\mathsf H}(P) = 
\begin{cases}
0 \quad & \text{ when } |P| \leq c_{\sfU}, \\
\lambda(P) \quad & \text{ when } |P| > c_{\sfU}.
\end{cases} 
\end{align*}
In particular, if we further simplify by taking $\sfU \equiv0$, then $\bar{\mathsf H}(P) =\frac12 |P|^2$. 
Suppose that we introduce conditional variance of the $\bnu(t;dv|x)$ as defining bulk temperature field
\begin{align*} 
T(t,x) : =  \int_{\R} |v-u(t,x)|^2 {\boldsymbol \nu}(t;dv|x).
\end{align*}
Then $\int_{\R} |v|^2 \bnu(t; dv|x) =  T(x) + u^2(x)$ and \eqref{Sec1:Euler} becomes
\begin{align*}
\begin{cases}
& \partial_t \rho + \partial_x(\rho u) =0, \\
& \partial_t (\rho u) + \partial_x (\rho u^2) + \partial_x p   = 0, 
\text{ where } p:= \rho T  
\end{cases}
\end{align*}
 The pressure-temperature-specific-volume relation $p =\rho T$ verifies the ideal gas law.

Again, we remind readers that the above construction allows us to only infer properties about action{\em-minimizing} path $\rho(\cdot)$.

 \subsubsection{Relation to micro-canonical ensemble in statistical mechanics}
Our paper here is built upon PDE analysis for Hamiltonian asymptotic. There is no explicitly use of ergodic arguments on Lagrangian paths. However, some form of ergodicity is hidden in the background through cell problem \eqref{cell}. We illustrate this next from the point of view of one-particle dynamic in the non-interacting particle model~\eqref{HNODE}. For weakly interacting Hamiltonians given by \eqref{defHN}, by particle permutation symmetry, the one particle argument can be transferred to the infinite particle setting through mass transport techniques.  This gives some heuristics on the form of effective Hamiltonian \eqref{Sec1:EffH} through informal classical arguments concerning statistical mechanics for particles. For simplicity, we assume that the $\sfH \in C^2$. We also write 
$\sfH^P:=\sfH^P(q,p):= \sfH(q,P+p)$.

Let $x_\epsilon(t):= x_{i,\epsilon}(t), P_\epsilon(t):=P_{i,\epsilon}(t)$ be those in \eqref{Hflow}, where the $i$ can be any. 
Sections 3.2 and 3.3 in~\cite{EG01} reveals that, at least along subsequences 
$\epsilon:=\epsilon_k$ with $k \to \infty$ if necessary, the following weak convergence (in $t$) occurs to a limiting measure-valued process 
$\{\sfm(t) \}_{t \geq 0}$:
\begin{align*}
\Phi\big(q(\frac{t}{\epsilon}), p(\frac{t}{\epsilon})\big) \Rightarrow^{w} 
\bar{\Phi}_t :=\langle \Phi, \sfm(t)\rangle, \quad \forall \Phi \in C_c(\R^{2d}).
\end{align*}
Morover,
\begin{align*}
 x_\epsilon(t) \to Q_t, \quad  P_\epsilon(t) \Rightarrow^{w} P_t 
\end{align*}
with
\begin{align*}
 \begin{cases}
  \dot{Q} &\in \partial_P \bar{\sfH}(P), \\
  \dot{P} &= 0.
\end{cases}
\end{align*}
where the $\partial_P \bar{\sfH}$ means sub-gradient for 
the convex (possibly non-smooth) $\bar{\sfH}$.  

 For a given Hamiltonian $\sfH$, the weak K.A.M. theory offers a number of concepts that characterize large time asymptotic sets, in different senses, for global action minimizing dynamics. See Fathi~\cite{FathiBook} for instance. In Appendix~\ref{App:wKAMCnp}, we give a highly condensed summary on Mather Measure ${\mathscr M}_\sfH$, projected Mather Measure ${\mathcal M}_\sfH$, 
Mather set $M_\sfH$, projected Mather set ${\mathbb M}_\sfH$, and projected Aubry set ${\mathbb A}_\sfH$, for definitions as well as their relations.

The measures $\sfm(t) \in {\mathscr M}_{\sfH^P}$ have interesting structures. To clarify, we follow Sections 4.1 and 6 of Evans and Gomez~\cite{EG01} by invoking the Lipschitz graph theorem which originally appeared in Mather~\cite{Mather91}. Such theorem has been further developed by Fathi and Siconolfi~\cite{FS04}, \cite{FS05}, Fathi~\cite{FathiBook}, among others. Let $\phi$ be any viscosity solution to the cell equation
\eqref{cell}. 
Then it is differentiable in the classical sense for every point in the projected Mather set $q \in {\mathbb M}_{\sfH^P}$ (See Appendix~\ref{App:wKAMCnp} for definition).  Moreover, there exists finite constant $C>0$ such that
 \begin{align*}
|\nabla_q \phi - \nabla_{q^\prime} \phi| \leq C|q-q^\prime|, 
 \quad \forall q, q^\prime \in {\mathbb M}_{\sfH^P}.   
\end{align*}
In fact, the above also holds for $q, q^\prime \in {\mathbb A}_{\sfH^P}$ 
the projected Aubry set~\cite{FS04}. 
The projected Aubry set is a larger set than the projected Mather set
-- Lemma~\ref{App:wKAMSets}.
The Lipschitz graph theorem implies that,  for a.e. each $t$ fixed,  
$\sfm^P:=\sfm^P(t)$ has the structure
\begin{align*}
 \sfm^P(dq, d p):= \delta_{\nabla_q \phi}(dp) \sigma(dq) 
  =   \delta_{\{ (q,p) : \text{ s.t.} \sfH^P(q, p) = \bar{\sfH} \}}   \sfm^P(dq,dp),
\end{align*}
for some $\sigma \in {\mathcal M}_{\sfH^P}$ (the set of projected Mather measures).
Note that, for all of the above quantities, the $\bar{\sfH}(P)$ is unique 
(by a comparison principle type argument) and
$ \bar{\sfH} = \langle \sfH, m(t) \rangle$.
However, $\phi$ can be non-unique, and the $\sigma$ can be non-unique.  

On the surface, the above arguments seem to give us a type of ergodic result for the Hamiltonian dynamics \eqref{HNODE}, at least along subsequences $\epsilon:=\epsilon_k \to 0$. Such view point is correct, however, only for those trajectories which are global (in time) action minimizers that satisfy \eqref{HNODE}. Nothing is said for all trajectories with arbitrary initial position-velocity vectors.  Conceptually, the $\sfm^P(t)$ is a kind of micro-canonical measure used in classical statistical mechanics arguments. However, there are subtleties that can destroy some classical informal arguments used in hydrodynamic limit derivations. We explain this point next. 
 
Typical physics textbooks define the micro-canonical ensemble as a uniform measure on the energy shell. First, we introduce a family of phase-volume measures indexed by energy levels $E$: 
with $\chi(r):= {\bf 1}_{[0,+\infty)}(r)$, we define
\begin{align*}
 \Gamma(E; \varphi) & := \int_{\sfH(q,p) \leq E} \varphi(q,p) dq dp \\
& = \int_{(q,p) \in \R^d \times \R^d} \chi (E - \sfH(q,p)) \varphi(q,p) dq dp,
\qquad \forall \varphi \in C_b(\R^d \times \R^d).
 \end{align*} 
 Second, denoting
$\mathcal H^s$ the $s$-dimensional Hausdorff measure, the micro-canonical measure for Hamiltonian $\sfH$ at energy level $E$ is introduced as  
\begin{align*}
\sfm_{m.c.}^{\sfH}(E;\varphi):= 
 \frac{\partial_E \Gamma(E;\varphi) }{\partial_E \Gamma(E;1)}
 = \frac{\int_{\sfH (q,p)=E} \frac{\varphi(q,p)}{|\nabla H|}d {\mathcal H}^{2d-1}}{\int_{\sfH (q,p)=E} \frac{1}{|\nabla H|}d {\mathcal H}^{2d-1}},
 \qquad \forall \varphi \in C_c(\R^d \times \R^d),
\end{align*}
where the last identity follows from the co-area formula (e.g. Proposition 3 on pages 118-119 of Evans and Gariepy~\cite{EG92}).
In general, minimizing Mather measures (and the projected Mather measure $\sigma$) are not unique.
There are also examples where projected Mather measures have non-smooth singular support. Consequently, in general,  
\begin{align*}
 \sfm^P(dq, dp) \neq \sfm_{m.c.}^{\sfH^P}(\bar{\sfH}; dq,dp).
\end{align*}

\subsection{Comments on technical developments}\label{Sec1:Tech}
 In the process of establishing a Hamilton-Jacobi theory for hydrodynamic limit in this paper, we also develop abstract mathematical techniques which can be of independent interests on their own. They include
\begin{enumerate}
\item projections of Hamilton-Jacobi equations given by metric foliation structures/submetry maps (see results with varying levels of generalities in Sections~\ref{projHJ} and~\ref{HJproj} and~\ref{PHJPert});
\item generalized Barles-Perthame half-relaxed limit arguments~\cites{BP87,BP88} for Hamilton-Jacobi equations in metric spaces (Theorem~\ref{Sec4:BPThm} and Lemma~\ref{Sec4:BPAlt} in Section~\ref{CnvHJ});
\item  reduction techniques for averaging of infinite-particle-Hamiltonian to that of single-particle-level (Sections~\ref{finPartH},\ref{Sec6}).  This compliments a method formally introduced in Section 4 in Feng, Mikami and Zimmer~\cite{FMZ21}, which was illustrated in a stochastic model context. 
\item comparison principles for a new class of Hamilton-Jacobi equations in space of probability measures by critical use of Alexandrov space tangent cone structures (see Sections~\ref{Sec7}, \ref{Sec8} and \ref{Sec9}), as well as clarification on relations among different definitions of Hamiltonian operators in such context. 
\end{enumerate}
The first two items on the above list are developed in a general metric space context, which is free of curvature assumptions. The notion of viscosity solution for Hamilton-Jacobi theory is a derivative-free one. In the parts of this paper involving Hamilton-Jacobi equation in space of probability measure, we try to explain everything through an Alexandrov-metric-space perspective. We hope this clarified many issues. However, in a few places, we have to go back relying upon techniques specific to optimal transport theory. We are uncertain if the key properties used, in these places,  can still be extracted into abstract arguments with metric-geometry nature. One such instance is the viscosity regularization-extension techniques in Section~\ref{Sec8} (Lemmas~\ref{Sec8:bbH0bfH0} and \ref{Sec8:bbH1bfH1}), which relied upon optimal transport Lemmas~\ref{Sec2:fdtouch} and \ref{Sec2:gdtouch}. 

 \newpage

\section{First order calculus in Wasserstein space as an Alexandrov metric space}\label{AlexCal}
The discussions in previous section highlight a need in understanding first order calculus on functions defined on Wasserstein space of probability measures.  Therefore, before engaging in an averaging theory for Hamilton-Jacobi equations in such space, we take a detour to recall and improve some results in such direction. In fact, the Wasserstein space is a special case of Alexandrov metric spaces with a notion of synthetic curvature bounded from below by zero. Some of our results in this paper are best presented as properties of Alexandrov spaces to reveal their true natures. In addition, we couldn't find some needed technical tools from existing literature to rigorously realize our above outlined program. Hence, we develop them here, for instance Lemmas~\ref{starmaxagl},\ref{sizesD},\ref{bnustar}, \ref{WassGf}, \ref{Sec2:fdtouch} etc.
Consequently, in the following, we start with Alexandrov metric spaces first, then we focus on Wasserstein spaces. For expositions on Alexandrov spaces, we follow presentations of Bridson and Haefliger~\cite{BH99}, Burago, Burago and Ivanov~\cite{BBI01},   Petrunin~\cite{Pet07}, Ambrosio, Gigli and Savar\'e~\cite{AGS08}, Alexander, Kapovitch and Petrunin~\cite{AKP19}.
For Wasserstein spaces, we follow Ambrosio, Gigli and Savar\'e~\cite{AGS08} and Gigli~\cite{Gigli04}, Villani~\cites{Villani03,Villani09}. 

Following~\cite{AKP19}, we denote $\CBB(\kappa)$ the collection of Alexandrov metric spaces with curvature bounded from below by $\kappa \in \R$.  
For each $x \in \sfX$, we introduce tangent cone $\Tan_x:=\Tan_x \sfX$ as a Euclidean cone over the space of directions (defined below, e.g. \cite{AKP19}), endowed with a cone metric $\sfd_x$. We recall the following basic definitions, concepts and properties. Given $x, y \in \sfX$, a constant speed connecting geodesic is a parameterized path $u : [0, T] \subset \R \mapsto \sfX$ such that $\sfd(u(t), u(s)) =\frac{|t-s|}{T} \sfd(x,y)$. 
In particular, when we take $T:= \sfd(x,y)$, the metric derivative of this curve  (e.g. Chapter 1, \cite{AGS08})  is one, we call it unit speed geodesic. 
Let $G_x$ be the set of all constant speed geodesics starting from $x$.
For each $u,v \in G_x$, we define the following notion of angle (mimicking the cosin law of Euclidean space) 
\begin{align*}
\cos \measuredangle_x (u,v) := \liminf_{s, t \to 0^+} \frac{\sfd^2(u(t), x) 
 + \sfd^2(v(s),x) - \sfd^2(u(t), v(s))}{2 \sfd(u(t), x) \sfd(v(s),x)}.
\end{align*}
Then $\measuredangle_x$ is a pseudo-metric on the space $G_x$. We define an equivalent relation that $u \sim v$ if $\measuredangle_x(u,v) =0$.
Let
\begin{align}\label{Sec2:UpArr}
\Uparrow_x^y := \{ \text{unit speed geodesics from $x$ to $y$} \}.
\end{align}   
The space of  geodesic directions at $x$ is defined by a quotient space 
\begin{align*}
 \Sigma^\prime_x := \cup_{\substack{y \in \sfX\\ y \neq x}} \Uparrow_x^y   \ /\sim. 
\end{align*}
We define space of directions $(\Sigma_x, \measuredangle_x)$ as completion of the
 $(\Sigma^\prime_x, \measuredangle_x)$. We also define the tangent cone $(\Tan_x, \sfd_x)$ as Euclidean cone of the space of directions 
 $(\Sigma_x, \measuredangle_x)$ (e.g. Definition 5.6 on page 59 of \cite{BH99}). We denote apex of the cone by $o:=o_x$ and write $| u |_x := \sfd_x( u, o_x)$. 
Then, for every $u,v \in \Tan_x$,   
\begin{align}\label{cosinelaw}
\sfd_x^2(u,v)  & = |u|_x^2 + |v|_x^2 - 2 |u|_x|v|_x \cos \measuredangle_x(u,v).
\end{align}
 A scalar product on $\Tan_x$ can be introduced by setting
\begin{align}\label{scalarAlex}
 \langle u, v \rangle_x := \frac12 \big( |u|_x^2 +|v|_x^2 - \sfd_x^2(u,v) \big) 
  =|u|_x |v|_x \cos \measuredangle_x(u,v), \quad \forall u, v \in \Tan_x.
\end{align}
For $u, v$ which are constant speed geodesics, the above defining relations also give us
\begin{align}\label{defConeD}
\sfd_x^2(u,v)  = \Big( \lim_{t \to 0^+} \frac{\sfd(u(t),v(t))}{t}\Big)^2.   
\end{align}
In general, $(\Tan_x, \sfd_x)$ may not even be a length metric space (e.g. Halbeisen~\cite{Hal00}) even if the $\sfX$ is. However, Corollary 5.11 on page 62 of \cite{BH99} gives a characterization of Euclidean cone to be geodesic under a geodesic space assumption on the space of directions generating the cone.  

In the above construction, we took completion of the $(\Sigma^\prime, \measuredangle_x)$ to arrive at the space of directions $(\Sigma, \measuredangle_x)$, then took $(\Tan_x, \sfd_x)$ as the Euclidean cone of $(\Sigma, \measuredangle_x)$. If, instead of the above, we take Euclidean cone of the $(\Sigma^\prime, \measuredangle_x)$ as $(\Tan^\prime_x, \sfd_x)$, then metric completion of the $\Tan^\prime_x$, we arrive at the same tangent cone $(\Tan_x, \sfd_x)$. The set $\Tan^\prime_x$ can be identified with geodesics starting at $x$ with arbitrary speed modulo equivalent class given by relation
\begin{align*}
 u \sim v \quad \text{ if and only if } \quad \sfd(u(t), v(t)) = o(t), \text{ as $t \to 0^+$}.
\end{align*}
\cite{AKP19} calls such $\Tan^\prime_x$ {\em space of geodesic tangent vectors} at $x$.

We again recall the definition of $\Uparrow_x^y$, the set of unit speed geodesics connecting $x$ and $y$, as given in \eqref{Sec2:UpArr}. To emphasize explicit parameter dependence of an element $\text{geod}_{[x,y]} \in \Uparrow_x^y$, we write
$\text{geod}_{[x,y]}(t)$ for $t \in [0, \sfd(x,y)]$.
Geodesics do not split in Alexandrov space with curvature bounded from below (e.g. Section 8.37 on page 81 of \cite{AKP19}). For each $u \in \Uparrow_x^y \subset \Sigma_x^\prime \subset \Sigma_x$ with $y \neq x$, we may re-parametrize the curve so that it becomes an arbitrary positive constant speed curve. For $t> 0$, we denote such re-parametrized curve $t \cdot u $ such that $|t \cdot u|_x:= t$. We denote $t \cdot \Uparrow_x^y \subset \Sigma_x^\prime$ the set of such re-parametrized curves. If a metric space $\sfX$ has the property that $\Uparrow_x^y$ is non-empty for every $x,y \in \sfX$, then such $\sfX$ is called a geodesic space.  
Within context of this section,  to simplify, {\it we assume without further mentioning that the space $(\sfX, \sfd) \in {\rm CBB}(\kappa)$ is geodesic, and also that it is a complete metric space.} A number of different definitions on Alexandrov spaces relying on  properties involving angles, triangles, short maps, concavity/convexity, etc etc, become equivalent under such assumption (see Chapter 8 of \cite{AKP19}). When multiple spaces are involved, we may introduce subindex on the metric 
$\sfd:=\sfd_\sfX$ to emphasize dependency on the space $\sfX$~\footnote{Such notation $\sfd_\sfX$ should be distinguished from the $\sfd_x$ -- the latter means a metric on the tangent cone at point $x$.}.

\subsection{First order calculus in Alexandrov metric space}\label{appalex}

For a function $f: \sfX \mapsto \bar{\R}$, we denote its domain
\begin{align*}
\Dom[f] := \{ x \in \sfX : | f(x)| <\infty \}.
\end{align*}

There are versions of semi-convexity -concavity relative to a curvature bound  (above or below) $\kappa \in \R$ that one can introduce -- see Definition 3.17 in \cite{AKP19}. However, for simplicity, we use only the following version. Since all concrete examples that we care about at in the $\CBB(\kappa =0)$ case, such simplification is good enough for our purpose.
\begin{definition}[Semi-concavity / convexity]\label{Sec2:SCCDef}
A lower semi-continuous function $f: \sfX \mapsto \R\cup\{-\infty\}$ is $\lambda$-concave for some $\lambda \in \R$ if: 
for {\em every } unit speed geodesic $z:= \{z(t)\}_{t} \subset \Dom[f]$, 
\begin{align*}
  t \mapsto  f(z(t)) - \lambda\frac{t^2}{2}
\end{align*}
is concave. We also call $f$ to be $\lambda$-convex if $-f$ is $\lambda$-concave.

$0$-convex (resp. $0$-concave) functions are called convex (resp. concave) functions respectively.

If for every $x \in \Dom[f]$, there exists a neighborhood and a $\lambda \in \R$ such that the restriction of $f$ to this neighborhood is $\lambda$-concave (resp. $\lambda$-convex), then $f$ is called semi-concave (resp. semi-convex). We denote $\SCC(\sfX;\bar{\R})$ the collections of semi-concave functions over $\sfX$.
\end{definition}

\subsubsection{Notions of differentials}
This subsection follows \cite{AKP19} and 
\begin{definition}[Velocity of curve]\label{veloCurv}
Let $x(\cdot): [0, \infty) \mapsto \sfX$ with $x(0) =x_0$. We say that $v \in \Tan_{x_0}$ is the right derivative of the curve $x(t)$ at $t=0$, denoted
\begin{align*}
x^{+}(0):=\frac{d}{dt}\Big|_{t=0+} x(t) = v,
\end{align*} 
if the following holds:
For some (and therefore any) sequence of geodesic tangent vectors $v_n \in \Tan^\prime_{x_0}$ with corresponding geodesics $\{ x_n(\cdot) \}_t$, such that $v_n \to v \in (\Tan_{x_0}, \sfd_{x_0})$, we have
\begin{align*}
\limsup_{n \to \infty} \limsup_{\epsilon \to 0}  \frac{1}{\epsilon} \sfd(x(\epsilon), x_n(\epsilon)) =0.
\end{align*}
\end{definition}

\begin{definition}[Differential]\label{diffDef}
Suppose that $f: \sfX \mapsto \bar{\R}$. Let $x_0 \in \Dom[f]$.  
We define a map $d_{x_0} f : \Tan_{x_0} \mapsto \R$ as differential of $f$ at $x_0$, provided this map satisfies the following: 
for every $v \in \Tan_{x_0}$ and every curve $x(\cdot): [0, \infty) \mapsto \sfX$ with $x(0) =x_0$ and $\frac{d}{dt}\big|_{t=0+} x(0) =v \in \Tan_{x_0}$, we have
\begin{align*}
 \frac{d}{dt}\Big|_{t=0+} f \big(x(t)\big) = \big(d_{x_0} f\big) (v).
\end{align*}
\end{definition}

\begin{lemma}[Proposition 6.16 of \cite{AKP19}]
Suppose that $f: \sfX \mapsto \bar{\R}$ is locally Lipschitz and semi-concave, then $d_x f$ is uniquely defined for each $x \in \Dom[f]$.
\end{lemma}

The following Lemma~\ref{ddist} is a metric space analysis version of the first variation formula. See for instance, Section 8.42 on page 84 of \cite{AKP19}  or Chapter 4.5 of \cite{BBI01}. The version in \cite{AKP19} is the most general and clean. However, to state results in such way requires introducing a concept called ``ultra-power" $\sfX^\omega$ of the metric space $\sfX$, which can be quite involved. We extract a useful property in the proof and formulate it as a condition. Such condition always holds when the $\sfX$ is locally compact (see proof of Corollary 4.5.7 in \cite{BBI01}). In general, it follows if a multiple weak-strong topology argument work, which is indeed the case in Wasserstein order-2 metric space example. 
\begin{condition}\label{CNDG}
For every $x_n, x, y_n, y \in \sfX$ with $\lim_{n \to \infty} \sfd(x_n, x)+ \sfd(y_n, y)=0$, and every $u_n \in \Uparrow_{x_n }^{y_n}$, there exists a subsequence $n(k)$ such that $u_{n(k)}:=\{u_{n(k)}(t) \}_{t}$ as a unit speed parametrized curve converge uniformly in $t$ (as $k \to \infty$) to another unit speed parametrized curve $u_0:= \{ u_0(t) \}_t$ with $u_0 \in \Uparrow_x^y$.
\end{condition}
\begin{example}
In Section~\ref{appmass}, we will consider order-2 Wasserstein space $(\sfX:=\mathcal P_2(\R^d), \sfd)$. This is a metrically complete geodesic $\CBB(0)$ space. It is non-locally compact but Condition~\ref{CNDG} is still satisfied. Verification of the condition goes as follow. We refer to the next section for definition of notations and basic results on Wasserstein spaces.

Convergence of geodesic end points $\sfd(\rho_n,\rho)+ \sfd(\gamma_n, \gamma)=0$ implies relative compactness in order-2 Wasserstein metric topology in $\mathcal P_2(\R^{2d})$  of any sequence of optimal plans $\bm_n \in \Gamma^{\opt}(\rho_n, \gamma_n)$. By an explicit geodesic characterization result using maps from the $\bm_n$s  -- See Theorem 7.2.2 of \cite{AGS08}, we conclude convergence along subsequence of geodesic curves as required by Condition~\ref{CNDG}.
\end{example}

\begin{lemma}\label{ddist}
Let $(\sfX, \sfd) \in \CBB(\kappa)$ for some $\kappa \in \R$.
We also assume that it is a geodesic metric space and metrically complete. Then
for every $y \neq x$, 
\begin{align*}
 \big(d_x \sfdist_y\big)(v) \leq \inf_{\xi \in \Uparrow_{x}^y} - \langle \xi ,v \rangle_x, \quad \forall v \in \Tan_x.
\end{align*}
If furthermore, Condition~\ref{CNDG} holds, then 
\begin{align*}
 \big(d_x \sfdist_y\big)(v) = \inf_{\xi \in \Uparrow_{x}^y} - \langle \xi ,v \rangle_x, \quad \forall v \in \Tan_x.
\end{align*}
\end{lemma}
\begin{proof}
With some notational changes, the proof in Theorem 4.5.6. in \cite{BBI01} can be adapted here.
\end{proof}

In Hilbert space situation, if $f$ is semi-concave, then $x \mapsto \nabla_x f$ is an semi-accretive operator. This brings up a host of related variational inequalities.  
Next, we state a result of this kind in $\CBB$ space situation. For simplicity and direct relevance to this article, we assume $(\sfX, \sfd) \in \CBB(0)$. However, general result also hold for any $\kappa \in \R$ by using special functions.  See Section E of Chapter 13 in \cite{AKP19} for details. 
\begin{lemma}\label{dxfCCV}
Suppose $(\sfX, \sfd) \in \CBB(0)$.
Let  $x, y \in \sfX$ with $x \neq y$. Suppose that $f: \sfX \mapsto \R$ is locally Lipschitz and $\lambda$-concave, 
and that any geodesic segment between $x, y \in \sfX$ belongs to $\Dom[f]$.
Then
\begin{enumerate}
\item
\begin{align}\label{dflest}
(d_x f)\big(v \big) \geq \frac{ f(y) - f(x) - \frac{\lambda}{2} \sfd^2(x,y)}{ \sfd(x,y) }, \quad \forall v \in \Uparrow_x^y
\end{align}
\item  
\begin{align*}
 (d_x f) (\xi) + (d_y f)(\eta) \geq    -  \lambda \sfd(x,y), \quad \forall \xi \in \Uparrow_x^y, \eta \in \Uparrow_y^x.
\end{align*}
In particular, 
\begin{align*}
 \big(d_x \frac{\sfdist_y^2}{2} \big) (u) + \big(d_y \frac{\sfdist_x^2}{2} \big)(v) 
  \geq  -  \sfd^2(x,y), 
 \quad \forall u \in \Uparrow_x^y, v \in \Uparrow_y^x.
\end{align*}
\end{enumerate}
 \end{lemma}
\begin{proof}
The proof of Lemma 13.24 in \cite{AKP19} works for the differential inequalities here as well, although the statement of that lemma was a weaker one involving gradient estimate only. 
\end{proof}

\begin{definition}\label{slopDef}
We also define several versions of {\em local Lipschitz constant} which measure respectively upward- downward- and overall- slopes of a function:
\begin{align*}
|D^+_x f|  &:= \limsup_{\substack{y \to x \\ y \neq x}} \frac{\big(f(y) - f(x)\big)\vee 0}{\sfd(x,y)}, \\
 |D^-_x f| & :=\limsup_{\substack{y \to x \\ y \neq x}} \frac{\big(f(x)-f(y)\big) \vee 0}{\sfd(x,y)}, \\
|D_x f| & := \limsup_{\substack{y \to x \\ y \neq x}} \frac{|f(y) - f(x)|}{\sfd(x,y)} = |D^+_x f| \vee |D^-_x f|.
\end{align*}
\end{definition}

\begin{lemma}\label{dDpm}
Let $(\sfX, \sfd) \in \CBB(0)$. Suppose that $f: \sfX \mapsto \bar{\R}$ is locally Lipschitz and semi-concave, then 
\begin{enumerate}
\item $d_x f$ is uniquely defined for each $x \in \Dom[f]$;
\item the map $v \mapsto (d_x f)(v): \Tan_x \mapsto \R$ is Lipschitz;
\item Lipschitz constant of the map $d_x f$ is no bigger than $|D_x f|$:
\begin{align*}
  \sup_{\substack{u \neq v\\ u,v \in \Tan_x}} \frac{| (d_x f)(u) - (d_x f)(v) |}{\sfd_x(u,v)}  \leq |D_x f|;
\end{align*}
\item  the following hold
\begin{align*}
\sup_{\substack{\xi \in \Tan_x, \\ |\xi|_x =1}} \Big( \big(d_x f\big)(\xi) \Big) \vee 0 & = |D_x^+ f|, \\
\sup_{\substack{\xi \in \Tan_x, \\ |\xi|_x =1}} \Big(\big(d_x (-f)\big)(\xi)\Big) \vee 0&= |D_x^- f|, \\
\sup_{\substack{\xi \in \Tan_x, \\ |\xi|_x =1}}|d_x f(\xi)| & = |D_x f|.
\end{align*}
\end{enumerate}
\end{lemma}
\begin{proof} 
For the first three parts, see Proposition 6.16 of \cite{AKP19}. We only prove the fourth property. 
Take $y \in \sfX$ and let $\xi \in \Uparrow_x^y$.  First,
\begin{align*}
(d_x f)( \xi)  & = \lim_{t \to 0^+} \frac{f(\text{geod}_{[x,y]}(t)) - f(x)}{t}  \\
& \leq \limsup_{z \to x} \frac{(f(z) - f(x)) \vee 0 }{  \sfd(x,z)} = |D_x^+ f|.
\end{align*}
Hence $\sup_{\xi} (d_x f)(\xi) \leq |D_x^+ f|$.
Second, by \eqref{dflest}, 
\begin{align*}
 R_f(x,y):=    \frac{f(y) - f(x) - \frac{\lambda}{2} \sfd^2(x,y)}{ \sfd(x,y)} \leq \big(d_x f\big)(\xi) \leq  \sup_{\xi} (d_x f)(\xi).
\end{align*}
Hence, when $|D_x^+f|>0$,
\begin{align*}
\lim_{\epsilon \to 0^+} \sup_{y: \sfd(y,x)<\epsilon} R_f(x,y) 
  \geq \lim_{\epsilon \to 0^+} \sup_{\{y: \sfd(y,x)<\epsilon\} \cap \{ y : f(y) \geq f(x)\}} R_f(x,y) = |D^+_x f|.
\end{align*}
However, when $|D_x^+ f|=0$, we cannot show that the set $\{y: \sfd(y,x)<\epsilon\} \cap \{ y : f(y) \geq f(x)\}$ is non-empty for some $\epsilon>0$ sufficiently small.
But, we can still conclude the trivial inequality $|D_x^+ f| =0\leq  \sup_{\xi} (d_x f)(\xi) \vee 0$.    
\end{proof}
\begin{remark}
We note here that, a notion of gradient $\nabla_x f$ can be defined for semi-concave function $f$ in $\CBB$ spaces, see Definition~\ref{GradDef} next. In particular (Lemma~\ref{UniGrad}), when $|\nabla_x f|_x>0$, it holds that
\begin{align*}
 \sup_{\substack{\xi \in \Tan_x, \\ |\xi|_x =1}}  \big(d_x f\big)(\xi)  = |\nabla_x f|_x .
\end{align*}
\end{remark}

\begin{example}\label{Ddstra}
Let $(\sfX, \sfd)$ be a general complete length metric space (no curvature bound assumption needed). 
Lemma 2.1 in Ambrosio and Feng~\cite{AF14} shows the following
\begin{align*}
|D^+_x \sfdist_y| \leq 1, \quad |D_x^- \sfdist_y| =1, \quad \forall x \neq y.
\end{align*}
In general, it can happen that $|D^+_x \sfdist_y|<1$. 
In fact, if $|D^+_x \sfdist_y|=1$, then $x$ is called $y$-{\em straight} by Definition 8.10 of Alexander, Kapovitch and Petrunin~\cite{AKP19}, and denoted $x \in \Str[y]$. With additional assumption $\sfX \in \CBB$, it is proved (e.g. Theorem 8.11 in \cite{AKP19}) that the set $\Str[y]$ is a dense $G_\delta$ set for every $y \in \sfX$. Moreover, for every $x \in \Str[y]$, there is a unique constant speed connecting geodesic between $x$ and $y$. 
\end{example}

\subsubsection{Gradient calculus}
\begin{definition}[Gradient]\label{GradDef}
Suppose that $f: \sfX \mapsto \bar{\R}$ is locally Lipschitz and semi-concave. Let $x \in \Dom[f]$.
We define the gradient of $f$ at $x$ as an element $\nabla_x f \in \Tan_x$, such that 
\begin{enumerate}
 \item $\big(d_x f\big)(v) \leq \langle \nabla_x f, v \rangle_x$, for every $v \in \Tan_x$;
 \item $\big( d_x f \big)  (\nabla_x f) = | \nabla_x f |_x^2$.
\end{enumerate}
\end{definition}

\begin{lemma}[Direction of steppest ascend]\label{UniGrad}
Let $f: \sfX \mapsto \R$ be locally Lipschitz and semi-concave. 
Then, for every $x \in \sfX$,  there exists a unique gradient $\nabla_x f \in \Tan_x$. Moreover, when $|\nabla_x f|_x >0$,
\begin{enumerate}
\item  there exists a unique maximizer $\xi^*$ of
\begin{align*}
s := \sup\big\{ (d_x f)(\xi) : \xi \in \Tan_x \text{ with }  |\xi|_x =1\big\},
\end{align*}
which is given by 
 \begin{align*}
 \frac{\nabla_x f}{|\nabla_x f|_x} \in \Tan_x. 
\end{align*}
\item    
\begin{align*}
 |\nabla_x f|_x = \sup\big\{ (d_x f)(\xi) : \xi \in \Tan_x \text{ with }  |\xi|_x =1\big\}.
\end{align*}
\end{enumerate}
\end{lemma}
\begin{proof}
See Section 13.E and Lemma 13.20 of \cite{AKP19}.
\end{proof} 

\begin{lemma}[Monotonicity of gradient on semi-concave functions]\label{GradMono}
Let $U: \sfX \mapsto \R$ be semi-concave and locally Lipschitz, $x, y \in \sfX$, and $v \in \Uparrow_x^y$ . Then $\langle v , \nabla_x U \rangle_x \geq  (d_x U)(v) $ and  
\begin{align}\label{metricHWI}
 \langle v, \nabla_x U \rangle_x + \langle u, \nabla_y U \rangle_y
  \geq - \lambda \sfd(x,y),   \qquad   
  \forall x, y \in \sfX, v \in \Uparrow_x^y, u \in \Uparrow_y^x.  
\end{align}  
\end{lemma}

In Alexandrov spaces, the tangent cone can be singular. In particular, suppose $u \in \Tan_x$,  there maybe no notion of opposite direction of $u$ in the tangent cone. The concept of polar helps to clarify such situations. 
\begin{definition}[Polar vector]\label{polar}
Two elements $u,v \in \Tan_x$ are called {\em polar } if
\begin{align*}
 \langle u, w\rangle_x + \langle v, w\rangle_x \geq 0, \quad \forall w \in \Tan_x.
\end{align*}
More generally, $u \in \Tan_x$ is called {\em polar to a set} $V \subset \Tan_x$ if
\begin{align*}
 \langle u, w\rangle_x + \sup_{v \in V} \langle v, w\rangle_x \geq 0, \quad \forall w \in \Tan_x.
\end{align*}
We denote the collection of such $u$ as $V^\circ$.
\end{definition}

The explicit expression in the first variation Lemma~\ref{ddist} gives us the following.
\begin{lemma}
Let $x \neq y$, then 
\begin{align*}
 (d_x \sfdist_y) (u) + (d_x \sfdist_y)(v) \leq 0, 
\end{align*}
for every $u, v \in \Tan_x$ which are polar with respect to each other.
\end{lemma}

Given any element in a tangent cone, we want to find another element in the tangent cone that makes angle between the two elements as wide as possible. Moreover, we want size of the new element to measure the size projected along direction of the original element. These considerations motivate the following concept. If the tangent cone is Hilbertian, we expect the new element to coincide with notion of opposite to the original element. 
\begin{definition}[Extremal polar vector $u^*$]\label{extpolar}
Given a $u \in \Tan_x$, an {\em extremal polar vector } is defined to be an element $u^* \in \Tan_x$ which is polar to $u$ and additionally satisfies
\begin{align*}
\langle u^*, u^* \rangle_x + \langle u^*, u \rangle_x =0.
\end{align*} 
\end{definition}
\begin{lemma}\label{uniqpolar}
For each $u \in \Tan_x$, there exists a unique extremal polar $u^* \in \Tan_x$. In particular, $|u^*|_x \leq |u|_x$.
\end{lemma}
\begin{proof}
The result follow as a special case of the anti-sum lemma in Section F in Chapter 13 of \cite{AKP19}.
\end{proof}
It follows from the above that, if $\Tan_x$ and $u \in \Tan_x$ are such that $\measuredangle_x(u, v) \leq \pi/2$ for every $v \in \Tan_x$, then the unique $u^* = o_x$.

Note that $x \mapsto \sfdist^2_y(x)$ is a 2-concave function in $\CBB(0)$ space (see Corollary 8.24 in \cite{AKP19} for this, as well as Section D in Chapter 8 of that book for general case of $\CBB(\kappa)$ with $\kappa \in \R$ ).
\begin{lemma}\label{gradDpol} Let $x \neq y$ and Condition~\ref{CNDG} holds. Then
\begin{enumerate}
\item  $\nabla_x \sfdist_y \in (\Uparrow_x^y)^\circ$. 
\item if we additionally assume that the set $\Uparrow_x^y:= 
\{ \uparrow_x^y \}$ consists of a singleton \footnote{By Theorem 8.11 of \cite{AKP19}, this assumption holds if $x \in \Str[y]$. Namely, $|D^+_x \sfdist_y|=1$.}, then
$\nabla_x \sfdist_y = (\uparrow_x^y )^*$ is the extremal polar to the $\uparrow_x^y$.
\end{enumerate}
\end{lemma}
\begin{proof}
$\nabla_x \sfdist_y\in (\Uparrow_x^y)^{\circ}$ because of 
Lemma~\ref{ddist}:
\begin{align*}
 \langle \nabla_x \sfdist_y, \eta\rangle_x \geq (d_x \sfdist_y)(\eta) 
 = -\inf_{\xi \in \Uparrow_x^y} \langle \xi, \eta \rangle_x.
\end{align*}

If $\Uparrow_x^y = \{ \uparrow_x^y\}$ is a singleton, then $\nabla_x \sfdist_y$ is polar to the $\uparrow_x^y$. Moreover, from the second defining property of gradient, it follows that
\begin{align*}
 |\nabla_x \sfdist_y|^2 = - \langle \nabla_x \sfdist_y, \uparrow_x^y \rangle_x.
\end{align*}
Hence it is the extremal polar.
\end{proof}

The following property made it clear that the ``extremal" in the definition of $u^*$ means ``maximal angle" that $u^*$ can open with respect to $u$, within the tangent cone $\Tan_x$.
\begin{lemma}\label{starmaxagl}
Let $u \in \Tan_x$ be such that $|u^*|_x > 0$. Then
\begin{align*}
\sup_{\substack{\xi \in \Tan_x \\ |\xi|_x=1}}  -\langle u, \xi \rangle_x   = - \langle u, 
 \frac{u^*}{|u^*|_x}\rangle_x = |u^*|_x.
\end{align*}
That is, for every $u \in \Tan_x$ with $|u|_x=1$ and $|u^*|_x>0$,  
\begin{align*}
\sup_{\substack{\xi \in \Tan_x \\ |\xi|_x=1}}   \measuredangle_x (u, \xi)  
= \measuredangle_x (u, \frac{u^*}{|u^*|_x}) =\arccos\big( - |u^*|_x \big).
\end{align*}
\end{lemma}
\begin{proof}
First of all, by definition of polarity,
\begin{align*}
 \sup_{\xi \in \Tan_x, |\xi|_x =1} - \langle u, \xi \rangle_x
     \leq \sup_{\xi \in \Tan_x, |\xi|_x =1} \langle u^*, \xi \rangle_x \leq |u^*|_x |\xi|_x \cos \measuredangle_x (u^*, \xi) \leq | u^*|_x.
\end{align*}
Second, by definition of extremal polar vector,
\begin{align*}
\max_{\xi \in \Tan_x, |\xi|_x=1} - \langle u, \xi \rangle_x \geq - \langle u, \frac{u^*}{|u^*|}\rangle_x = |u^*|_x.
\end{align*}
\end{proof}

In general, $\measuredangle_x (u, u^*) \neq \pi$ unless they becomes {\em opposite}.  
\begin{definition}[Opposite]
We say that $u, v \in \Tan_x$ are {\em opposite } to each other (symbolically written $u + v =0$) in either of the following situations
\begin{enumerate}
\item $|u|_x=|v|_x=0$;
\item $\measuredangle_x (u,v) =\pi$ and $|u|_x =|v|_x$.
\end{enumerate}
\end{definition}
By definition, $\langle u^*, u \rangle_x =-|u^*|_x^2$ always holds. If the $u, u^*$ are opposite, we also get the extra property 
\begin{align*}
 \langle u^*, u\rangle_x = -|u|_x^2, \quad  |u|_x=|u^*|_x.
\end{align*}

\begin{lemma}[Proposition 13.37 of \cite{AKP19}]\label{oppolem}
For $u,v \in \Tan_x$ to be opposite is equivalent to 
  $\langle u, w\rangle_x + \langle v, w\rangle_x =0$ for every $w \in \Tan_x$
\end{lemma}

\begin{lemma} \label{oppolem2}
Let $u \in \Tan_x$. Then the following are equivalent
\begin{enumerate}
\item the $u^*$ and $u$ are opposite.
\item  $(u^*)^* = u$. 
\item $|u^*|_x=|u|_x$.
\end{enumerate}
\end{lemma}
\begin{proof} 
Suppose that $u^*$ and $u$ are opposite, by Lemma~\ref{oppolem}, taking $w=u$,
\begin{align}\label{u2star}
| u |_x^2+  \langle u^*, u \rangle_x =0.
\end{align}
Hence $(u^*)^*=u$ (in view of the uniqueness result in Lemma~\ref{uniqpolar}).
Suppose that $(u^*)^*=u$, then \eqref{u2star} holds. Combined with the defining relation of $u^*$, we have
\begin{align}\label{eqnormuu}
 |u|_x^2 = -  \langle u^*, u \rangle_x =|u^*|_x^2.
\end{align}
Suppose that $|u|_x=|u^*|_x$. By definition of $u^*$, \eqref{eqnormuu} holds.  By the cosin law, 
\begin{align*}
\cos \measuredangle_x (u, u^*)= -1.
\end{align*} 
Therefore $u^*$ and $u$ are opposite. 
\end{proof}

In the following, we develop a set of new results illustrating relation among the earlier introduced notions of differential, sub- super- gradient as well as gradient of a Lipschitz semi-concave function. In particular, their relations with polar of certain geodesics when the function becomes a distance.   
  
\begin{definition}[Fr\'echet super- and sub-gradients]\label{Sec2:Frech}
Let $f : \sfX \mapsto \bar{\R}$ with $f(x) \in \R$. We respectively denote super-, sub-gradients of $f$ at $x$ by 
$\bpartial^+_x f, \bpartial^-_x f$. These are subsets of $\Tan_x$ satisfying the following properties. 

We say $u \in \bpartial^+_x f \subset \Tan_x$, if there exists a modulus of continuity $\omega_u$ such that 
\begin{align*}
f(y) - f(x) \leq \sup_{v \in  \sfd(x,y) \cdot \Uparrow_x^y } \langle u, v \rangle_x + \sfd(x,y) \omega_u( \sfd(x,y)), \quad \forall y \in \sfX.
\end{align*}
If $u \in \Tan_x$ satisfies that, there exists modulus of continuity  $\omega_u$ with
\begin{align*}
f(y) - f(x) \leq \inf_{v \in \sfd(x,y) \cdot \Uparrow_x^y} \langle u, v \rangle_x + \sfd(x,y) \omega_u( \sfd(x,y)), \quad \forall y \in \sfX,
\end{align*}
then we say $u$ belongs to a strong super-gradient  $\bpartial^{s,+}_x f$. It follows then $\bpartial^{s,+}_x f \subset \bpartial^{s}_x f$.
 
Analogously, we say $u \in \bpartial^-_x f \subset \Tan_x$, if there exists a modulus of continuity $\omega_u$ such that 
\begin{align*}
f(y) - f(x) \geq \inf_{v \in\sfd(x,y) \cdot \Uparrow_x^y} \langle u, v \rangle_x - \sfd(x,y) \omega_u (\sfd(x,y)), \quad \forall y \in \sfX.
\end{align*}
and   $u \in \bpartial^{s,-}_x f \subset \Tan_x$, if there exists a modulus of continuity $\omega_u$ with 
\begin{align*}
f(y) - f(x) \geq  \sup_{v \in \sfd(x,y) \cdot \Uparrow_x^y} \langle u, v \rangle_x - \sfd(x,y) \omega_u (\sfd(x,y)), \quad \forall y \in \sfX.
\end{align*}
It follows that $\bpartial^{s,-}_x f \subset \bpartial^{-}_x f$.
\end{definition}

\begin{lemma}\label{GsubG}
Let $f \in \Lip_{\loc}( \sfX) \cap \SCC(\sfX)$, then
\begin{align*}
 \nabla_x f \in \bpartial_x^{s,+} f.
\end{align*}
\end{lemma}
\begin{proof}
This follows from \eqref{dflest} and the first defining property of gradient of a semi-concave function.
\end{proof}

\begin{lemma}\label{SlopDiff}
Let $f \in \SCC(\sfX) \cap \Lip_{\loc}(\sfX)$ and $x \in \sfX$.
Then
\begin{align*}
| D_x^+ f| \leq \inf \big\{ |u|_x :  u \in \bpartial_x^+ f \big\}.
\end{align*}
Note that, by convention, $\inf$ over empty set is $+\infty$.
\end{lemma}
\begin{proof}
Let $u \in \bpartial_x^{+} f$. Then there exists a modulus of continuity $\omega$ such that
\begin{align*}
 \frac{f(z) - f(x)}{\sfd(x, z)}   \leq \sup_{\eta \in \Uparrow_x^z } \langle  u, \eta \rangle_x    + \omega(\sfd(x,z)) \leq |u|_x + \omega(\sfd(x,z))
  \qquad  \forall  z \in \sfX.
\end{align*}
Therefore the conclusion follows. 
\end{proof}

\begin{lemma}\label{sizesD}
Let $x,y \in \sfX$, $x \neq y$ and Condition~\ref{CNDG} holds. 
Then  
\begin{enumerate}
\item  
\begin{align}\label{polsubD}
\{ u^* : u \in \Uparrow_x^y\} \subset \big(\Uparrow_x^y\big)^\circ = \big(\bpartial^{s,+}_x \sfdist_y \big).
\end{align}
\item  
\begin{align}\label{allDnsame}
 |\nabla_x \sfdist_y|_x& =|D^+_x \sfdist_y|   =\sup_{\substack{|v|_x=1, \\ v \in \Tan_x}}  \big( (d_x \sfdist_y)(v) \big)\vee 0 =
 \inf_{w \in \bpartial^+_x \sfdist_y} |w|_x  \\
 & \qquad \quad = \inf_{w \in (\Uparrow_x^y)^\circ} |w|_x
= \inf_{w \in \bpartial^{s,+}_x \sfdist_y}  |w|_x   \leq \inf_{u \in \Uparrow_x^y}  |u^*|_x .  \nonumber
\end{align}
\item  $\nabla_x \sfdist_y$ is a minimal element in $\bpartial_x^+ \sfdist_y$ (as well as a minimal element in  $\bpartial_x^{s,+} \sfdist_y$) in the sense that
\begin{align*}
 \nabla_x \sfdist_y \in \bpartial_x^+ \sfdist_y \quad \text{ and } \quad
|\nabla_x \sfdist_y|_x= \inf_{w \in \bpartial^{s,+}_x \sfdist_y}  |w|_x.
\end{align*}
\end{enumerate}
\end{lemma}
\begin{proof}
It follows from definitions that $\{ u^* : u \in \Uparrow_x^y\} \subset \big(\Uparrow_x^y\big)^\circ$. 

Denote $f:=\sfdist_y$, then $f$ is locally concave when bounded away from the point $y$. Since the $f$ is at most linearly growing with respect to the metric $\sfdist_y$, by selecting a large enough $\lambda \in \R_+$ which may depend on $y$, through the results in \eqref{dflest} and Lemma~\ref{ddist},  we have that 
\begin{align*}
 f(z) - f(x) & \leq \big( \inf_{u \in \Uparrow_x^y} -\langle  u , \xi \rangle_x \big) \sfd(x,z) + \frac{\lambda}{2} \sfd^2(x, z), \\
  & \leq   \big( \inf_{v \in (\Uparrow_x^y)^\circ} \langle  v , \xi \rangle_x \big) \sfd(x,z) + \frac{\lambda}{2} \sfd^2(x, z), 
  \quad  \forall \xi \in \Uparrow_x^z, z \in \sfX.
\end{align*}
Hence $(\Uparrow_x^y)^\circ \subset \bpartial^{s,+}_x f$. 
Next, let $w \in \bpartial^{s,+}_x f$, then in view of Lemma~\ref{ddist} and the defining inequalities of super-gradients of $f$ at $x$,
\begin{align}\label{ddistwxi}
\inf_{u \in \Uparrow_x^y} \big( -\langle u, \xi \rangle_x \big) = (d_x f)(\xi)  \leq \langle w, \xi \rangle_x , \quad \forall \xi \in \Uparrow_x^z, \forall z \in \sfX.
\end{align}
That is, $w \in \big(\Uparrow_x^y\big)^\circ$. Hence  
$\bpartial^{s,+}_x f\subset \big(\Uparrow_x^y\big)^\circ$. 
Therefore, \eqref{polsubD} is verified.

Next,  we prove \eqref{allDnsame}. First,  
\begin{align*}
  \inf_{w \in \bpartial^{s,+}_x f}  |w|_x
\leq |\nabla_x f|_x   \leq 
 \sup_{\substack{|v|=1,\\ v \in \Tan_x}}  (d_x f)(v) \vee 0 = |D^+_x f|.
\end{align*}
In the above, the first inequality follows because of $\nabla_x f \in \bpartial^{s,+}_x f$ (Lemma~\ref{GsubG}), the last equality holds because of Lemma~\ref{dDpm}. To verify the second inequality, we only need to verify the non-trivial case when $|\nabla_x f|_x>0$. 
From Lemma~\ref{UniGrad},
\begin{align*}
 |\nabla_x f|_x = (d_x f)\big(\frac{\nabla_x f}{|\nabla_x f|_x}\big) 
 =\sup_{\substack{|v|=1,\\ v \in \Tan_x}}  (d_x f)(v).
 \end{align*}
Second, in view of Lemma~\ref{SlopDiff} and \eqref{polsubD}, we have that
\begin{align*}
 |D^+_x f|   \leq
 \inf_{w \in \bpartial^+_x f} |w|_x 
\leq  \inf_{w \in \bpartial^{s,+}_x f}  |w|_x  
= \inf_{w \in (\Uparrow_x^y)^\circ} |w|_x
\leq \inf_{u \in \Uparrow_x^y}  |u^*|_x .   
\end{align*}
Combine the above first and second points, we conclude that \eqref{allDnsame} holds.

Finally, in view of  Lemma~\ref{GsubG} and the identities in \eqref{allDnsame}, $\nabla_x f \in \bpartial_x^{s,+} f \subset \bpartial_x^+ f$ exists as a minimal element in $\bpartial_x^+ f$ (also in  $\bpartial_x^{s,+} f$).
\end{proof}

\begin{lemma}\label{Sec2:AllEq}
If the following minimax equality holds
\begin{align}\label{minmaxCND}
\sup_{\substack{|v| =1\\ v \in \Tan_x}} \inf_{u \in \Uparrow_{x}^y} - \langle u ,v \rangle_x 
=  \inf_{u \in \Uparrow_{x}^y} \sup_{\substack{|v| =1\\ v \in \Tan_x}}  - \langle u ,v \rangle_x.
\end{align}
Then all quantities in \eqref{allDnsame} are equal.
 \end{lemma}
\begin{proof}
To show that \eqref{minmaxCND} implies all quantities are equal in \eqref{allDnsame},
we notice 
\begin{align*}
\inf_{u \in \Uparrow_x^y} |u^*|_x = \inf_{u \in \Uparrow_x^y} 
 \sup_{\substack{|v|=1\\ v \in \Tan_x}} - \langle u, v \rangle_x
= \sup_{\substack{|v|=1\\ v \in \Tan_x}} \inf_{u \in \Uparrow_x^y} 
 - \langle u, v \rangle_x
=\sup_{\substack{|v|=1\\ v \in \Tan_x}}  (d_x \sfdist_y)(v).
\end{align*}
\end{proof}

The inclusion relation \eqref{polsubD} makes us wonder if the minimal element in $\bpartial^{s,+}_x \sfdist_y$, which gives the gradient, can also be selected from the subset $\{ u^* : u \in \Uparrow_x^y\}$.  We have the following result.

\begin{lemma}\label{gdist}
Let $x,y \in \sfX$, $x \neq y$ and Condition~\ref{CNDG} hold.  Then there exists a minimizer $u_0 \in \Uparrow_x^y$ solving the minimization problem 
\begin{align*}
|u_0^*|_x=  \inf_{u \in \Uparrow_x^y} |u^*|_x. 
\end{align*}
Assume that all the quantities in \eqref{allDnsame} are equal.
Suppose furthermore that minimal element in either $\bpartial_x^{s,+} \sfdist_y$ or $\bpartial_x^+ \sfdist_y$ is unique, then
\begin{align}\label{gradDpol2}
 \nabla_x \sfdist_y =  u_0^*.
 \end{align}
 \end{lemma}
\begin{proof}
We take a sequence $u_n \in \Uparrow_x^y$ such that 
\begin{align*}
 \lim_{n \to \infty}    |u_n^*|_x =   \inf_{u \in \Uparrow_{x}^y} |u^*|_x. 
\end{align*}
By Condition~\ref{CNDG}, there exists another unit speed curve $u_0 \in \Uparrow_x^y$ such that, up to selection of a subsequence (still denoted using the $n$), the $u_n$s as curves converges uniformly in time along subsequences.
For $\CBB$ spaces,  angle between hinges is lower semicontinuous with respect to convergence of hinges (e.g. Section 8.40 on page 82 of \cite{AKP19} or Theorem 4.3.11 of \cite{BBI01}). Consequently,
\begin{align*}
\liminf_{n \to \infty} - \langle u_n, v \rangle_x = \liminf_{n \to \infty} -\cos \measuredangle (u_n,v) \geq -\cos \measuredangle (u_0,v)
= -\langle u_0, v \rangle_x, \quad \forall v \in \Sigma^\prime_x.
\end{align*}
Therefore, by Lemma~\ref{starmaxagl},
\begin{align*}
\liminf_{n \to \infty} |u_n^*|_x = \liminf_{n \to \infty} \sup_{\substack{|v| =1\\ v \in \Sigma^\prime_x}}   - \langle u_n, v \rangle_x
   \geq \sup_{\substack{|v| =1\\ v \in \Sigma^\prime_x}}  -\langle u_0, v \rangle_x =|u_0^*|_x .
\end{align*}
Hence  $\inf_{u \in \Uparrow_{x}^y} |u^*|_x = |u_0^*|$.

By \eqref{polsubD}, the $u_0^* \in \bpartial_x^{s,+} \sfdist_y \subset \bpartial_x^+ \sfdist_y$.
If all quantities in \eqref{allDnsame} are equal, then the $u_0^*$ is a minimal element in $\bpartial_x^{s,+} \sfdist_y$ 
as well as a minimal element in $\bpartial_x^+ \sfdist_y$. 
Since the $\nabla_x \sfdist_y$ is also a minimal element in $\bpartial^+_x \sfdist_y$ and in $\bpartial^{s,+}_x \sfdist_y$ (Lemma~\ref{sizesD}), 
the uniqueness assumption on minimal element implies \eqref{gradDpol}.
\end{proof}

\begin{remark}\label{Sec2:GradDiff}  
Combining all the above results, we discover that differential $d_x \sfdist_y$ could contain strictly more information than gradient $\nabla_x \sfdist_y$, when the $x$ becomes a singular point in the sense that $|D^+_x \sfdist_y| \neq 1$. Because of this, it makes sense for us to work mostly differentials when formulating first order Hamilton-Jacobi partial differential equations in terms of Hamiltonian operators $Hf(x) := H(x, d_x f)$. 
\end{remark}
\subsubsection{Simple functions}\label{Simfun}
Following Petrunin~\cite{Pet07}, we consider some classes of simple smooth test functions in $(\sfX, \sfd)$ and their differential properties.  

Let
\begin{align}\label{defPsiK}
{\bf \Psi}_K & := \big\{ \psi \in C^2(\R^K):  \psi \geq 0, \psi \text{ semi-concave},
  \partial_k \psi > 0, \forall k =1,2,\ldots,K \big\}.
\end{align}
We write
\begin{align} 
{\mathcal S}^+& := \big\{f:= f(x):= \psi\big(\sfdist^2_{y_1}(x), \ldots, \sfdist^2_{y_K}(x)\big), 
\quad \forall y_k \in \sfX,  \psi \in {\bf \Psi}_K, K \in \N \big\}, \label{SS+} \\
{\mathcal S}^-& := \big\{ g:=g(y):= - \psi\big( \sfdist^2_{x_1}(y), \ldots, \sfdist^2_{x_K}(y)\big), 
\quad \forall x_k \in \sfX,  \psi \in {\bf \Psi}_K, K \in \N\}.  \label{SS-}
\end{align}
If situation warrants, we may also write ${\mathcal S}^+_\sfX$ and ${\mathcal S}^-_\sfX$ to emphasize the underlying space $\sfX$.

\begin{lemma}\label{S+SCC}
Every function in ${\mathcal S}^+$ is locally semi-concave in $(\sfX, \sfd)$. Respectively, every function in ${\mathcal S}^-$ is locally semi-convex.
\end{lemma}
\begin{proof}
We only verify the claim regarding $f \in {\mathcal S}^+$, the other can be similarly proved.

Let $x, z \in \sfX$ and $x(t)$ be a unit speed connecting geodesic with $x(0)=x$ and $x(1)=z$. Since $\sfX \in \CBB$,   $t \mapsto \sfdist_y^2(x(t)) $ is semi-concave. That is, $D^2_t \sfdist_y^2(x(t)) \leq \kappa$ for $t$ almost everywhere. 
We note that for a $\lambda$-concave $\psi$,  
\begin{align*}
 \sum_{i,j=1}^K (\partial_{ij} \psi) \eta_i \eta_j  \leq \lambda \sum_{i=1}^K |\eta_i|^2, \quad \forall (\eta_1,\ldots, \eta_K) \in \R^K.
\end{align*}
Through regularization and approximation, therefore the following holds for $t$ almost everywhere,
\begin{align*}
D^2_t f(x(t)) & = \sum_{k=1}^K \partial_k \psi\big(\sfdist^2_{y_1}(x(t)), \ldots, \sfdist^2_{y_K}(x(t))\big) D^2_t \sfdist^2_{y_k}(x(t)) \\
& \qquad \qquad 
  + \sum_{i,j=1}^K \partial_{ij} \psi\big(\sfdist^2_{y_1}(x(t)), \ldots, \sfdist^2_{y_K}(x(t))\big) D_t \sfdist^2_{y_i}(x(t)) D_t \sfdist^2_{y_j}(x(t)) \\
   & \leq \kappa \sum_{k=1}^K  \partial_k \psi\big(\sfdist^2_{y_1}(x(t)), \ldots, \sfdist^2_{y_K}(x(t))\big)
     + 4 \lambda \sum_{k=1}^K |D_t \sfdist_{y_k}^2(x(t))|^2.
\end{align*}
Hence $f$ is locally semi-concave.
We note that the $D_t \sfdist_y(x(t))$ has an explicit expression given by the first variation formula.
\end{proof}

From Lemma~\ref{ddist}, we have the following.
\begin{lemma}\label{Dmdist}
Let $f \in {\mathcal S}^+$, then
\begin{align*}
( d_x f)(v) &= 2 \inf_{ \substack{\xi_k \in \Uparrow_x^{y_k} \\ k=1,\ldots,K}}  
\sum_{k=1}^K  \partial_k \psi\big(\sfdist^2_{y_1}(x), \ldots, \sfdist^2_{y_K}(x)\big) \sfdist_{y_k}(x) 
 \big( - \langle \xi_k, v \rangle_x\big), \quad \forall v \in \Tan_x.
\end{align*}
Let $g \in {\mathcal S}^-$, then 
\begin{align*}
( d_y g)(u) & = -\big(d_y (-g)\big)(u)\\
& =2 \sup_{ \substack{\eta_k \in \Uparrow_y^{x_k} \\ k=1,\ldots,K}}  
\sum_{k=1}^K  \partial_k \psi\big(\sfdist^2_{x_1}(y), \ldots, \sfdist^2_{x_K}(y)\big) \sfdist_{x_k}(y) 
 \big(  \langle \eta_k, u \rangle_y\big), \quad \forall u \in \Tan_y.
\end{align*}
\end{lemma}

We also introduce two slightly larger classes of test functions than the ${\mathcal S^+}$ and ${\mathcal S}^-$. Let $\R^\infty$ be countable infinite product space of $\R$ with the usual product topology. For $r:= (r_1, \ldots, r_k,\ldots) \in \R^\infty$, we denote the usual sequence space norm $|r|_{l^p}$ for $p \in [1,+\infty]$, in particular, $|r|_{l^\infty}:= \sup_{k \in \N} |r_k|$.
For a function $\psi : \R^\infty \mapsto \R_+$ with $\partial_k \psi \in C(\R^\infty)$, we denote $\nabla \psi := (\partial_1 \psi, \ldots, \partial_k \psi, \ldots)$.  Let 
\begin{align} \label{Sec2:defPsi}
{\bf \Psi} &:= \Big\{ \psi :  \psi, \partial_k \psi, \partial_{ij} \psi \in C(\R^\infty), \\
 & \qquad \qquad \psi_K:=\psi_K(r_1,\ldots, r_K):=\psi(r_1,\ldots, r_K, 0,0,\ldots) \in \Psi_K, \nonumber \\
 & \qquad \qquad \quad \text{ and } \sup_{r \in \R^\infty, |r|_{l^\infty} \leq C}  
 \big(|\nabla \psi(r)|_{l^1}+|\nabla \psi(r)|_{l^2}\big) <+\infty, \forall C \in \R_+,
  \nonumber \\
 & \qquad \qquad \qquad \text{ and } \lim_{k \to \infty} \psi(r_1,  \ldots, r_k, 0, 0, \ldots) =  \psi(r_1,  \ldots, r_k, r_{k+1}, \ldots) \Big\}, \nonumber
\end{align}
we write
\begin{align} 
{\mathcal S}^{+,\infty}
& := \Big\{f:= f(x):= \psi\big(\sfdist^2_{y_1}(x), \ldots, \sfdist^2_{y_k}(x), \ldots \big), \label{extSS+}  \\
& \qquad \qquad \qquad 
\forall y_k \in \sfX \text{ with } \sup_k \sfdist_{y_1}(y_k)<\infty, \psi \in {\bf \Psi} \Big\}, \nonumber \\
{\mathcal S}^{-,\infty}& := \Big\{ g:=g(y):= - \psi\big( \sfdist^2_{x_1}(y), \ldots, \sfdist^2_{x_k}(y), \ldots \big), \label{extSS-} \\
&\qquad \qquad \qquad 
\forall x_k \in \sfX \text{ with } \sup_k \sfdist_{x_1}(x_k) <\infty,  \psi \in {\bf \Psi} \Big\}. 
\nonumber 
\end{align}

\begin{remark}\label{Sec2:Dmdext}
The conclusions of Lemma~\ref{Dmdist} can be extended for $f \in {\mathcal S}^{+,\infty} $ and $g \in {\mathcal S}^{-,\infty}$ as well.
\end{remark}

\subsection{First order calculus in Wasserstein space of probability measures}\label{appmass}
Next, we study a situation where $\sfX:= \mathcal P_2(\R^d)$ is the space of probability measures over $\R^d$ with finite second moments, and $\sfd$ is the Wasserstein order-2 metric (e.g. Chapter 7.1 of \cite{AGS08}). We call such metric space order-2 Wasserstein space.  By Theorem 7.3.2 in Ambrosio-Gigli-Savar\'e~\cite{AGS08},
the space $({\mathcal P}_2(\R^d), \sfd) \in \CBB(0)$, and it is a geodesic and complete separable metric space.  
 Moreover, tangent cone $\Tan_\rho$ of this space can be identified explicitly using probability-measure-coupling techniques. 
This also leads to more probability-measure based representation of differentials and gradients for simple smooth test functions given in \eqref{SS+} and \eqref{SS-}.

We introduce a few additional mass transport notations in the following: Let $\pi^{i_1,i_2,\ldots,i_l} : (\R^d)^K \mapsto (\R^d)^l$ be a projection 
\begin{align*}
 \pi^{i_1,i_2,\ldots, i_l} (x_1, \ldots, x_K) = (x_{i_1}, x_{i_2}, \ldots x_{i_l}), \quad \forall x_i \in \R^d.
\end{align*}
For $\rho, \gamma \in {\mathcal P}_2(\R^d)$ and $\bm \in {\mathcal P}_2(\R^{2d})$, 
\begin{enumerate}
\item $\Gamma(\rho,\gamma):= \big\{ \bmu \in {\mathcal P}(\R^{2d}) :  
    \pi^{1}_\# \bmu = \rho, \pi^{2}_{\#} \bmu =\gamma \big\}$;
\item $\Gamma^{\opt}(\rho,\gamma) := \big\{ \bmu \in {\mathcal P}(\R^{2d}) :  
    \sfd^2(\rho,\gamma)=\int_{\R^{2d}} |x-y|^2 \bmu(dx, dy) \big\}$; 
    \item $\Gamma^{\opt} (\bm, \gamma):=\big\{ \bmu \in {\mathcal P}(\R^{3d}) :  
    \pi^{1,2}_\# \bmu = \bm, \pi^{1,3}_{\#} \bmu \in \Gamma^{\opt}(\pi^1_\# \bm, \gamma) \big\}$;
\end{enumerate}
Throughout this section, we write a regular conditional probability decomposition (also known as slicing measure decomposition) for 
 any $\bpi \in \mathcal P(\R^d \times \R^d)$ with $\pi^1_\#\bpi =\rho$ as
 \begin{align*}
  \bpi(dx,dy):= \bpi(dy|x) \rho(dx).
\end{align*}
In the mass transport context, we write $\grad_\rho f$ for gradient of a semi-concave function $f$ instead of $\nabla_\rho f$ 
because the notation $\nabla_x \varphi$ is reserved for gradient of a function on Euclidean space $\varphi: \R^d \mapsto \R$.


 \subsubsection{Tangent cone identification}
In this section, following original arguments in \cites{Gigli04,AGS08}, we identify the tangent cones $\Tan_\rho$ as defined in general abstract sense at the beginning of  Section~\ref{AlexCal}. 
 
We define
\begin{align*}
G(\rho)& := \big\{ {\bnu}:= {\bnu}(dx, d \xi) \in {\mathcal P}_2(\R^{2d})_\rho : \pi^1_\# {\bnu} = \rho, (\pi^1, \pi^1+ \epsilon \pi^2)_{\#} {\bnu}
  \in \Gamma^{\opt}(\rho, \gamma) , \\
  & \qquad  \text{for some }  \gamma \in \sfX, \epsilon >0\big\}.
\end{align*}
For each ${\bnu}_i \in {\mathcal P}_2(\R^{2d})$ with $\pi^1_{\#} {\bnu}_i = \rho$, $i=1,2$, we define a metric
\begin{align}\label{Sec2:drho}
 \sfd_\rho({\bnu}_1, {\bnu}_2) &:= \inf \big\{ \int_{\R^{3d}} |\xi -\eta|^2 \bM(dx; d \xi, d \eta) : \bM \in {\mathcal P}_2(\R^{3d})  \\
  & \quad \qquad \qquad \pi^{1,2}_\# \bM = \nu_1, \pi^{1,3}_\# \bM = \nu_2 \big\},
   \nonumber
\end{align}
and a scalar product
\begin{align}\label{Sec2:ScaDef}
\langle {\bnu}_1, {\bnu}_2 \rangle_\rho &:= \max \big\{ \int_{\R^d \times \R^d \times \R^d} (\xi \cdot \eta) \bM(dx; d \xi , d \eta) :
  \bM \in {\mathcal P}_2(\R^{2d}), \\
 & \qquad \qquad  \qquad \pi^{1,2}_\# \bM = {\bnu}_1, \pi^{1,2}_\# \bM = {\bnu}_2  \big\}. \nonumber
\end{align}
In particular, when ${\bnu}_1 = {\bnu}_2=: {\bnu}$, the above maximum is attained at 
\begin{align*}
 \bM(dx; d \xi, d \eta) :=\delta_{\xi}(d \eta) \nu(dx, d \xi);
\end{align*}
and
\begin{align*}
\Vert {\bnu}\Vert_{\rho}^2:= \langle {\bnu}, {\bnu}\rangle_\rho = \int_{\R^{2d}} |\xi|^2 \nu(dx, d\xi). 
\end{align*}
We now define
\begin{align*}
 \Tan_\rho:= \overline{G(\rho)}^{\sfd_\rho(\cdot, \cdot)}, \quad \Tan:= \bigsqcup_{\rho \in \sfX} \Tan_\rho.
\end{align*}

\begin{lemma}[Proposition 12.4.2 of \cite{AGS08}, Theorem 4.12 of~\cite{Gigli04}]
\label{Sec2:TanIden}
The tangent space $(\Tan_\rho, \sfd_\rho)$ defined above coincides (up to isometry) with the tangent cone introduced in abstract Alexandrov metric space setting at the beginning of  Section~\ref{AlexCal}, where the $(\sfX, \sfd)$ is viewed as a geodesic $\CBB(0)$ space with complete metric $\sfd$. 
 \end{lemma}

\subsubsection{Identification of tangent cones, polar and extremal polar vectors}
For each $\rho \in \sfX$, we denote
\begin{align*}
 {\mathcal P}_2(\R^{2d})_\rho&:= \big\{ {\bnu}:={\bnu}(dx, d\xi) \in \mathcal P_2(\R^{2d}) : \pi^1_\# {\bnu} = \rho \big\}, \\
 \sfD^2_\rho(\bnu_1, \bnu_2)&:= \int_{\R^d} \sfd^2\big({\bnu}_1(\cdot|x), {\bnu}_2(\cdot|x)\big)  \rho(dx), 
   \quad \forall \bnu_i \in {\mathcal P}_2(\R^{2d})_\rho,
\end{align*}
where the $\bnu_i(\cdot|x)$s are the disintegrations of the $\bnu_i$s with respect to $\rho(dx)$ -- i.e. $\bnu_i(dx,d \xi) = {\bnu}_i(d \xi |x) \rho(dx)$. It follows then 
\begin{align*}
\Tan_\rho \subset {\mathcal P}_2(\R^{2d})_\rho \subset {\mathcal P}_2(\R^{2d}).
\end{align*}
\begin{lemma}[Proposition 12.4.6 of \cite{AGS08} or Proposition 4.2 in \cite{Gigli04}]
We have 
\begin{align*}
\sfd_\rho( \bnu, \bmu) = \sfD_\rho(\bnu, \bmu), \quad \forall \bnu, \bmu \in  {\mathcal P}_2(\R^{2d})_\rho.
\end{align*}
\end{lemma}

\begin{lemma}[Theorem 4.5. in \cite{Gigli04}]
For each $\rho \in {\mathcal P}_2(\R^d)$, the $({\mathcal P}_2(\R^{2d})_\rho, \sfD_\rho )$ is a complete metric space.
\end{lemma}

 \begin{definition}\label{LinPTan}
 For any $t \in \R$ and $\bmu, {\boldsymbol \nu} \in {\mathcal P}_2(\R^{2d})_\rho$, we define 
 \begin{align*}
 t \cdot \bmu & := ( \pi^1, t \pi^2)_\# \bmu, \\
 \bmu \oplus {\boldsymbol \nu} &:= \big\{ (\pi^1, \pi^2 + \pi^3)_\# \bM : \exists \bM \in {\mathcal P}_2(\R^{3d}),
  \pi^{1,2}_\# \bM =\bmu, \pi^{1,3}_\# \bM = {\boldsymbol \nu} \big\} .
\end{align*}
\end{definition}

\begin{lemma}[Proposition 4.25 of \cite{Gigli04}]\label{addmeas}
For $\bmu, {\boldsymbol \nu} \in \Tan_\rho$, we have $\bmu \oplus {\boldsymbol \nu} \subset \Tan_\rho$.
\end{lemma}
 \begin{lemma}[Proposition 4.27 of \cite{Gigli04}]\label{omax}
 For every $\alpha \geq 0$ and $\bmu, {\boldsymbol \nu}, {\boldsymbol \nu}_i \in {\mathcal P}_2(\R^{2d})_\rho$ with $i=1,2,3$, we have
 \begin{align*}
 & \langle \alpha \cdot {\boldsymbol \nu}, \bmu\rangle_\rho = \langle {\boldsymbol \nu}, \alpha \cdot \bmu\rangle_\rho = \alpha \langle {\boldsymbol \nu}, \bmu\rangle_\rho, 
  \\
& \langle {\boldsymbol \nu}_1, {\boldsymbol \nu}_3 \rangle_\rho + \langle {\boldsymbol \nu}_2, {\boldsymbol \nu}_3 \rangle_\rho =\max  \langle {\boldsymbol \nu}_1 \oplus {\boldsymbol \nu}_2, {\boldsymbol \nu}_3 \rangle_\rho,
\end{align*}
where the max is over the set ${\boldsymbol \nu}_1 \oplus {\boldsymbol \nu}_2$.
 \end{lemma}

\begin{lemma}\label{negnu}
Let ${\boldsymbol \nu} \in \Tan_\rho$. Then 
\begin{enumerate}
\item $(-1)\cdot {\boldsymbol \nu} \in \Tan_\rho$ is polar to ${\boldsymbol \nu}$ (Definition~\ref{polar}).
\item $\Vert (-1)\cdot {\boldsymbol \nu} \Vert_{\rho} = \Vert {\boldsymbol \nu} \Vert_\rho$.
\end{enumerate}
\end{lemma}
\begin{proof}
The fact that $(-1) \cdot {\boldsymbol \nu} \in \Tan_\rho$ was proved in Proposition 4.29 of \cite{Gigli04}. Take $\bmu \in \Tan_\rho$, then
\begin{align*}
&   \langle {\bnu}, \bmu\rangle_\rho + \langle (-1)\cdot {\bnu}, \bmu\rangle_\rho  \\
  &\geq \max_{\substack{\bM \in {\mathcal P}_2(\R^{4d}), \\ 
                         \pi^{1,2}_\# \bM = \bmu, \pi^{1,3}_\# \bM = {\bnu}\\
                         \pi^{1,4}_\# \bM = (-1) \cdot {\bnu}}} 
  \big\{ \int_{\R^{3d}} (\xi_1 \cdot \xi_2 + \xi_1 \cdot \xi_3) \bM(dx; d \xi_1, d \xi_2, d \xi_3) 
  \big\} \\
  & = \max_{\substack{\bM \in {\mathcal P}_2(\R^{3d}), \\ 
                         \pi^{1,2}_\# \bM = \bmu, \pi^{1,3}_\# \bM = {\bnu}}} 
  \big\{ \int_{\R^{3d}} \big(\xi_1 \cdot \xi_2 + \xi_1 \cdot (-\xi_2)\big) \bM(dx; d \xi_1, d \xi_2) \big\} \\
 & =0.
\end{align*}
In the above, the inequality follows from \eqref{Sec2:ScaDef}, the equality follows from the definition $(-1) \cdot \bnu$.

Therefore, $(-1) \cdot {\bnu} $ and ${\bnu}$ are polar. 

$\Vert (-1) \cdot {\bnu} \Vert_{\rho} = \Vert {\bnu} \Vert_{\rho}$ follows from definition. 
\end{proof}

In general, $({\bnu})^* \neq (-1) \cdot {\bnu}$ in the sense of Definition~\ref{extpolar}. Otherwise,  with the property $\Vert (-1) \cdot {\bnu} \Vert_\rho = \Vert {\bnu} \Vert_\rho$, 
the two tangent elements become opposite (Lemma~\ref{oppolem2}).
However, as Remark 4.28 in Gigli~\cite{Gigli04} pointed out,
\begin{align*}
- \langle (-1)\cdot {\bnu}, \bmu \rangle_\rho = \min_{\substack{\bM \in {\mathcal P}_2(\R^{3d}), \\ 
                         \pi^{1,2}_\# \bM = {\bnu}, \pi^{1,3}_\# \bM = \bmu}} \int_{\R^{3d}} (\xi_1 \cdot \xi_2) \bM(dx; d \xi_1, d \xi_2), \\        
  \langle {\bnu}, \bmu \rangle_\rho = \max_{\substack{\bM \in {\mathcal P}_2(\R^{3d}), \\ 
                         \pi^{1,2}_\# \bM = {\bnu}, \pi^{1,3}_\# \bM = \bmu}} \int_{\R^{3d}} (\xi_1 \cdot \xi_2) \bM(dx; d \xi_1, d \xi_2).                                
\end{align*}
The above two quantities are not the same for a generic $\bmu \in \Tan_\rho$( Lemma~\ref{oppolem}). 
 
 \begin{lemma}[Proposition 4.25 of \cite{Gigli04}]\label{Gig4-25}
 For every ${\boldsymbol \nu}, \bmu \in \Tan_\rho$, ${\boldsymbol \nu} \oplus \bmu \in \Tan_\rho$ and $t \cdot {\boldsymbol \nu} \in \Tan_\rho$ for all $t \in \R_+$.
 \end{lemma}

\begin{lemma}\label{bnustar}
Let ${\boldsymbol \nu} \in \Tan_\rho$. Then its extremal polar is given by 
\begin{align*}
({\boldsymbol \nu})^*(dx, dv) = \delta_{-v(x)}(dv) \rho(dx), \text{ where } v(x) := \int_{\R^d} w {\boldsymbol \nu}(dw|x).
\end{align*}
\end{lemma}
\begin{proof}
It follows from Theorem 12.4.4 in \cite{AGS08} that ${\boldsymbol \nu}^* \in \Tan_\rho$. Therefore, 
conclusion of the lemma follows from two more observations. 
One, we always have
\begin{align*}
 \langle {\boldsymbol \nu}, {\boldsymbol \nu}^*\rangle_\rho + \Vert {\boldsymbol \nu}^* \Vert_\rho^2
 & =   \int_{\R^{2d}} w  \big(-v(x)\big)    {\boldsymbol \nu}(dw|x) \rho(dx) 
  +  \int_{\R^d} |v(x)|^2 \rho(dx)  \\
 & =  \int_{\R^d}\big( - |v(x)|^2 + |v(x)|^2 \big) \rho(dx) =0.
\end{align*}
Two, let $\bmu \in \Tan_\rho$, we define a 3-variable probability measure using conditional independence:
\begin{align*}
\bM(dx; dv, du) :=  {\boldsymbol \nu}(dv| x) \bmu(du|x) \rho(dx) = {\boldsymbol \nu}(dv|x) \bmu(dx,du) \in \mathcal P_2(\R^{3d}).
\end{align*}
 Then $\pi^{1,2}_\# \bM = {\boldsymbol \nu}$ and $\pi^{1,3}_\# \bM =\bmu$, and
\begin{align*}
 \langle{\boldsymbol \nu}, \bmu\rangle_\rho + \langle {\boldsymbol \nu}^*, \bmu \rangle_\rho
 & \geq   \int_{\R^{3d}}  v u   \bM(dx,dv,du) - \int_{\R^{2d}} v(x) u \bmu(dx, du) \\
 & = \int_{\R^{2d}}  \big(\int_{\R^d} v {\boldsymbol \nu}(dv|x)\big)  u   \bmu(dx,du) - \int_{\R^{2d}} v(x) u \bmu(dx, du) =0.
 \end{align*}
 \end{proof}

\subsubsection{Differentials and gradients}
\begin{lemma}\label{PolyGrad}
Let $f$ be as in \eqref{Poly}:
\begin{align*}
f(\sigma):= \psi(\langle \varphi_1, \sigma \rangle, \ldots, \langle \varphi_K, \sigma \rangle), \quad \forall \varphi_i \in C^1(\R^d), \psi \in C^2(\R^K).
\end{align*}
Then
\begin{align*}
 \grad_\rho f  = {\boldsymbol \nu}(dx, du) &:= \delta_{\nabla\frac{\delta f}{\delta \rho}(x)}(du) \rho(dx), 
\end{align*}
where the
\begin{align*}
\frac{\delta f}{\delta \rho} & :=\sum_{k=1}^K \partial_k \psi( \langle \varphi_1, \rho \rangle, \ldots, \langle \varphi_K, \rho \rangle )  \varphi_k.
\end{align*}
\end{lemma}
\begin{proof} We note that $f$ is semi-concave in the Wasserstein space $(\sfX, \sfd)$. 
Direct calculation gives
\begin{align*}
(d_\rho f)(\bmu) =   \int_{\R^{2d}} \big( v \cdot \nabla_x\frac{\delta f}{\delta \rho} \big) \bmu(dx, dv), \quad \forall \bmu \in \Tan_\rho.
\end{align*}
Consequently,
\begin{align*}
\langle {\boldsymbol \nu}, \bmu\rangle_\rho  = \int_{\R^{2d}}  \big( v \cdot   \nabla_x\frac{\delta f}{\delta \rho}(x) \big)\bmu(dx, dv) 
 = \big( d_\rho f\big) (\bmu), \quad \forall \bmu \in \Tan_\rho.
\end{align*}
In particular, 
\begin{align*}
\Vert {\boldsymbol \nu} \Vert^2_\rho = \int_{\R^{2d}} \Big| \nabla_x\frac{\delta f}{\delta \rho} \Big|^2 \rho(dx) = \big( d_\rho f\big) ({\boldsymbol \nu}).
\end{align*}
\end{proof}

Next, following Chapter 12.4 of Ambrosio, Gigli and Savar\'e~\cite{AGS08}, we define a concept of exponential map at least on $G(\rho)$, 
 a dense subset of the tangent cone. In a similar way, we also define a notion of (right) inverse exponential maps.
 \begin{definition} [Exponential, inverse exponential maps]\label{Emap}
 We define
\begin{align*}
\exp_\rho(\bnu):= \big\{ (\pi^1 + \pi^2)_\# \bnu \big\},   \qquad \forall \bnu \in G(\rho), 
\end{align*}
and for every $\gamma \in \sfX$,
\begin{align*}
\exp_\rho^{-1} (\gamma):=\big\{ \bnu \in G(\rho): \gamma= \exp_\rho(\bnu)  \big\} 
 = \big\{ \bnu \in G(\rho) : (\pi^1, \pi^1 + \pi^2)_\# \bnu \in 
  \Gamma^{\opt}(\rho,\gamma) \big\} .
\end{align*}
\end{definition}
 
 From abstract results in Lemma~\ref{UniGrad}, we see that $\grad_\rho \sfdist_\gamma^2 \in \Tan_\rho$ exists uniquely. Next, we find a probability-measure representation of such quantity.  We also elaborate on an explicit selection criteria of it from the sets of super-, sub-differentials.
 
%
%
%
%
%

By Theorem 10.2.2 of \cite{AGS08}, we have
\begin{lemma}\label{supsubD}
\begin{align*}
 \big\{ (\pi^1, (\pi^2-\pi^1))_\# \bmu : \bmu \in \Gamma^{\opt}(\rho, \gamma)  \big\} 
 & \subset \bpartial^{s,+}_\rho \big( \frac12 \sfdist_\gamma^2 \big), \\
 \big\{ (\pi^2, (\pi^1-\pi^2))_\# \bmu : \bmu \in \Gamma^{\opt}(\rho, \gamma)  \big\}
 &  \subset \bpartial^{s,-}_\gamma \big( -\frac12 \sfdist_\rho^2 \big).
\end{align*}
\end{lemma}
Next, noting $\sfd(\rho,\gamma) \cdot \Uparrow_\rho^\gamma = \exp_\rho^{-1}(\gamma)$, we consider a minimization problem
\begin{align*}
s := \inf \Big\{ \Vert \bmu^* \Vert_\rho : \bmu \in \exp_\rho^{-1}(\gamma)  \Big\}.
\end{align*}
For every $\bpi \in \Gamma^{\opt}(\rho,\gamma)$, making a change of variable
 \begin{align}\label{bpi2bmu}
\bmu:= \big( \pi^1, (\pi^2-\pi^1)\big)_\# \bpi \in  \exp_\rho^{-1}(\gamma)  \subset  \Tan_\rho.
\end{align}
By Lemma~\ref{bnustar},
\begin{align*}
 \bmu^*= \delta_{-v(x)}(dv)\rho(dx), \quad  v(x) := \int_{\R^d} v \bmu(dv|x)= \Big(\int_{\R^d} y \bpi(dy|x)\Big) -x.
\end{align*}
Consequently, 
 \begin{align*}
s^2 & = \min \Big\{ \int_{\R^d} |\int_{\R^d} y \bpi(dy|x) -x|^2 \rho(dx) : \bpi \in \Gamma^{\opt}(\rho,\gamma) \Big\} \\
 & = \min \Big\{ \int_{\R^d} |v|^2 d\rho : (\rho, v) \text{ given as above} \Big\}.
  \end{align*}

Following Definition~\ref{GradDef}, gradient of every semi-concave locally Lipschitz function $f$ is well defined. In such mass transport situation, we use $\grad_\rho f$ for such notation, to distinguish the $\nabla$ notation which could appear as in $\nabla_x \rho$. For semi-convex locally Lipschitz function $g$, we define $\grad g = - \grad (-g)$ as the $-g$ is semi-concave. Therefore, the notion $\grad f$ is well defined for $f \in {\mathcal S}^{+,\infty} \cup {\mathcal S}^{-,\infty}$.

 \begin{lemma}[Identification of Gradient]\label{GradS} 
It holds that 
\begin{align}\label{sDeqmass}
s &=  \Big|D^+_\rho \frac12 \sfdist_\gamma^2\Big| = \Vert \grad_\rho \frac12 \sfdist_\gamma^2 \Vert_\rho  \\
&= \inf\Big\{\Vert \bmu\Vert_\rho : \bmu \in \bpartial^{+}_\rho \frac12 \sfdist_\gamma^2  \Big\}
= \inf \Big\{ \Vert \bmu\Vert_\rho : \bmu \in \bpartial^{s,+}_\rho \frac12 \sfdist_\gamma^2 \Big\}. \nonumber
\end{align}
Moreover, there exists $\bpi_0 \in \Gamma^{\opt}(\rho,\gamma)$ which is the unique minimizer  
 \begin{align}\label{Sec2:pimini}
  &  \int_{\R^d} |\int_{\R^d} y \bpi_0(dy|x) -x|^2 \rho(dx) \\
  & = 
   \min \big\{  \int_{\R^d} |\int_{\R^d} y \bpi(dy|x) -x|^2 \rho(dx) : \bpi \in \Gamma^{\opt}(\rho,\gamma) \big\}, \nonumber
  \end{align}
and the above quantity equal to $s^2$.
From $\bpi_0$ we define $\bmu_0$ according to \eqref{bpi2bmu}, then 
  \begin{align*}
 \grad_\rho \big(\frac12 \sfdist_\gamma^2\big) = (\bmu_0)^*
  = \delta_{-v_0(x)}(dv)\rho(dx), \quad 
   v_0(x) := \int_{\R^d} v \bmu_0(dv|x).
\end{align*}
\end{lemma}
\begin{proof}
By Theorem 10.4.12 of \cite{AGS08} and Lemma~\ref{bnustar} for identification of $\bnu^*$ below, we have that
\begin{align*}
 \Big|D^+_\rho  \frac12 \sfdist_\gamma^2\Big|^2 =s^2 = \inf_{\bnu \in \Uparrow_\rho^\gamma} |\bnu^*|_\rho^2,
\end{align*}
and that the minimizer $\bpi_0$ in \eqref{Sec2:pimini} is unique. Therefore, identity \eqref{sDeqmass} holds because of Lemma~\ref{sizesD} \footnote{Note that, in the notations of that lemma, 
$\inf_{u \in \Uparrow_x^y} |u^*| = |D^+_x \sfdist_y|$, implying equalities for all the quantities.}.

By Theorem 10.3.11 of \cite{AGS08}, minimal selection of 
$\bpartial_\rho^+ \frac12 \sfdist_\gamma^2$  is unique 
\footnote{Note that the theorem in \cite{AGS08} is applied to negative distance squared function, which is stated for sub-differentials. We converted the results to super-differentials by getting rid of the negative sign.}.
Hence, by Lemma~\ref{gdist},
\begin{align*}
 \grad_\rho  \frac12 \sfdist_\gamma^2  = \bmu_0^*.
\end{align*}
\end{proof}
 
 \subsubsection{Simple smooth test functions}\label{SSmass}
We recall the definition of several classes of simple smooth test functions $\mathcal S^{\pm}$ in \eqref{SS+}, \eqref{SS-} and ${\mathcal S}^{\pm,\infty}$ in \eqref{extSS+}, \eqref{extSS-}. We adapted them to the Wasserstein situation, in particular,
\begin{align}
{\mathcal S}^+ \ni f &:=f(\rho) := \psi\big(\sfdist_{\gamma_1}^2(\rho), \ldots, \sfdist_{\gamma_K}^2(\rho)\big) , \label{mtranf}\\
{\mathcal S}^- \ni g &:= g(\gamma) := -\psi\big(\sfdist_{\rho_1}^2(\gamma), \ldots, \sfdist_{\rho_K}^2(\gamma)\big). \label{mtrang}
\end{align}
By the  semi-concavity and local Lipschitzness of $f$ (Lemma~\ref{S+SCC}), $d_\rho f$ is well defined, and $\mu =\grad_\rho f$ exists uniquely in $\Tan_\rho$ (Lemma~\ref{UniGrad}). We identify these quantities explicitly next.

\begin{lemma}\label{WassDSS}
Let $f \in {\mathcal S}^+$ as in \eqref{mtranf}, or more generally $f \in {\mathcal S}^{+,\infty}$ which can be considered notation-wise as $K:=+\infty$. We denote
\begin{align}\label{alphaDef}
\alpha_k:= \alpha_k(\rho; \gamma_1, \ldots, \gamma_K):= \partial_k  \psi\big(\sfdist^2_{\gamma_1}(\rho), \ldots, \sfdist^2_{ \gamma_K}(\rho)\big) \geq 0.
\end{align}
Then (allowing the case $K=+\infty$)
\begin{align*}
 \big(d_\rho f\big)(\bnu) 
 = 2 \inf_{\substack{\bM \in \mathcal P_2(\R^{(K+2)d})\\ 
   \pi^{1, k+1}_\# \bM \in \Gamma^{\opt}(\rho,\gamma_k),  k=1,\ldots,K \\
   \pi^{1,K+2}_\#\bM = \bnu} }
  \int_{\R^{(K+2)d}} \big(\sum_{k=1}^K \alpha_k (x-y_k) \cdot v\big) \bM(dx, dy_1, \ldots, d y_K; dv)
\end{align*}
holds for every $\bnu \in \Tan_\rho$.  In particular,
\begin{align}\label{sec2:ddsqr}
\big(d_\rho \sfdist_\gamma^2\big) (\bnu)
 =   2 \inf_{ \bmu \in \exp_\rho^{-1}(\gamma) } 
 \big( -\langle \bmu, \bnu \rangle_\rho \big), \quad \forall \bnu \in G(\rho).
\end{align}
\end{lemma}
\begin{proof}
We look at the case of $f \in {\mathcal S}^+$ first.
Apply Lemma~\ref{Dmdist} to the special case of Wasserstein space, and take into account of results in Lemmas~\ref{addmeas} and \ref{omax}, then 
\begin{align*}
 \big(d_\rho f\big)({\boldsymbol \nu}) & = 2  \inf_{\substack{\bmu_k \in \exp_\rho^{-1}(\gamma_k) \\ k=1,2,\ldots,K}}
  \sum_{k=1}^K - \alpha_k  \langle \bmu_k,  \bnu \rangle_{\rho},   \\
 &  = 2  \inf_{\substack{\bmu_k \in \exp_\rho^{-1}(\gamma_k) \\ k=1,2,\ldots,K}}
  - \sup \langle  \oplus_{k=1}^K \alpha_k \cdot \bmu_k,  {\bnu} \rangle_{\rho},  , \\
  &= 2  \inf_{\substack{\bmu_k \in \exp_\rho^{-1}(\gamma_k) \\ k=1,2,\ldots,K \\
             \bmu \in \oplus_{k=1}^K \alpha_k \cdot \bmu_k}} 
              - \langle \bmu,  \bnu \rangle_{\rho}, \qquad \forall {\bnu} \in \Tan_\rho.
\end{align*}
This gives the conclusion.

The general case of $f \in {\mathcal S}^{+,\infty}$ (by taking $K=\infty$ in the above expressions) follows from Remark~\ref{Sec2:Dmdext}.
\end{proof}

To help with presentation, we also introduce notation for a particular type of optimal multi-plans.
 \begin{definition}  \label{MultOpt}
Let  $\rho, \gamma_k\in \mathcal P_2(\R^d)$, $k=1,2,\ldots,K$ where the $K \in \N \cup\{+\infty\}$. We denote
\begin{align*}
 \Gamma^{\opt} (\rho; \gamma_1, \ldots, \gamma_K) 
& :=  \Big\{ \bM := \bM(dx; dy_1,\ldots, dy_K) \in \mathcal P_2\big(\R^{(1+K)d} \big) \text{ such that }\\
 &   \qquad \qquad \pi^{1, 1+k}_\# \bM \in \Gamma^{\opt}(\rho, \gamma_k), \quad k=1,\ldots,K \Big\}.
\end{align*}
 \end{definition}  
For the $f\in \mathcal S^+$, and $\bM  \in \Gamma^{\opt}(\rho; \gamma_1, \ldots, \gamma_K)$, 
we write
\begin{align}\label{bnu0Def}
 \bnu^{\bM}_f(dx, dP) :=   \int_{(y_1,\ldots, y_K) \in \R^{Kd}} 
  \delta_{\sum_{k=1}^K  2 \alpha_k  (x-y_k)} (dP)\bM(dx, dy_1,\ldots, dy_K),
 \end{align}
 where the $\alpha_k$s are defined according to \eqref{alphaDef}.
 In the case of $f \in \mathcal S^{+,\infty}$ ($K=+\infty)$, at least when $\sup_k \sfd(\gamma_k, \rho) < \infty$, by those uniform summability requirements on $\nabla \psi(r) \in l^1 \cap l^\infty$ {\color{red}check here again} in the definition of test functions ${\bf \Psi}$ in \eqref{Sec2:defPsi}, $\bnu^{\bM}_f \in \mathcal P_2(\R^d \times \R^d)$ is well defined. 
 
 \begin{lemma}\label{Sec2:nupar}
 For $f \in {\mathcal S}^+$, $ \bnu^{\bM}_f \in \bpartial_\gamma^{s,+} f$.
 \end{lemma}

\begin{remark}\label{dggamma}
Using $d(-f) = - d f$, we have that for every $g \in {\mathcal S}^-$ as in \eqref{mtrang} 
and for every $\bnu \in \Tan_\gamma \sfX$,
\begin{align*}
 \big(d_\gamma g\big)(\bnu) = 2  \sup_{\substack{\bM \in \mathcal P_2(\R^{(K+2)d})\\ 
   \pi^{1, k+1}_\# \bM \in \Gamma^{\opt}(\gamma,\rho_k),  k=1,\ldots,K \\
   \pi^{1,K+2}_\#\bM = \bnu} }
  \int_{\R^{(K+2)d}} \big(\sum_{k=1}^K \beta_k (x_k-y) \cdot v\big) \bM(dy, dx_1, \ldots, d x_K; dv),
\end{align*}
where the 
\begin{align}\label{betaDef}
\beta_k := \beta_k(\gamma; \rho_1, \ldots, \rho_K):= \partial_k  \psi\big(\sfdist^2_{\rho_1}(\gamma), \ldots, \sfdist^2_{ \rho_K}(\gamma)\big) \geq 0.
\end{align}
In the same way, the above expression also holds (by setting $K=+\infty$) for 
$g \in {\mathcal S}^{-,\infty}$.

For the $g \in {\mathcal S}^-$ as in \eqref{mtranf} and 
$\bM  \in \Gamma^{\opt}(\gamma;\rho_1, \ldots, \rho_K)$,
we write
\begin{align}\label{Sec2:bnu1Def}
 \bnu^{\bM}_g(dy, dP) :=   \int_{(x_1,\ldots, x_K) \in \R^{Kd}} 
  \delta_{\sum_{k=1}^K  \beta_k 2 (x_k-y)} (dP)\bM(dy, dx_1,\ldots, dx_K),
 \end{align}
 where the $\beta_k$s are defined according to \eqref{betaDef}.
It follows then $ \bnu^{\bM}_g \in \bpartial_\gamma^{s,-} g$.  
Analogous expression and result also hold for $ g \in {\mathcal S}^{-,\infty}$ by setting $K=+\infty$.
\end{remark}

We also have the following.\footnote{We could use Lemma~\ref{GradS} to give a proof, but we choose to give a different one based upon direct calculations.}
\begin{lemma}\label{WassGf}
In the context of Lemma~\ref{WassDSS}, let  $\bpi_{0,k} \in \Gamma^{\opt}(\rho,\gamma_k)$ 
denote the unique minimizer (Lemma~\ref{GradS}) of 
\begin{align}\label{uniMini}
 s_k^2 & := \int_{\R^d} |\int_{\R^d} y \bpi_{0,k}(dy|x) - x|^2 \rho(dx) \\
 & = 
 \inf \Big\{ \int_{\R^d} |\int_{\R^d} y \bpi(dy|x) - x|^2 \rho(dx) : 
    \bpi \in \Gamma^{\opt}(\rho,\gamma_k)\Big\}. \nonumber
\end{align}
We define  
\begin{align*}
 v_k(x) :=\int_{\R^d} (y -x) \bpi_{0,k}(dy|x), \quad k=1,2,\ldots.
\end{align*}
where the $\alpha_k$s are defined as in \eqref{alphaDef}.
Then 
\begin{align}\label{Sec2:gradf}
 \grad_\rho f =  \delta_{- \sum_{k=1}^K 2 \alpha_k v_k(x) }(du) \rho(dx), \quad
 \forall f \in {\mathcal S}^+.
\end{align}
\end{lemma}
\begin{proof} For convenience,  we denote 
\begin{align}\label{defnuu0}
{\boldsymbol \nu}_0(dx,du)  := \delta_{u_f(x)}(du) \rho(dx), \quad 
u_f(x):=- \sum_{k=1}^K 2 \alpha_k v_k(x).
\end{align}
Then
\begin{align*}
 \Vert {\bnu}_0\Vert_\rho \leq \sum_{k=1}^K 2 \alpha_k \sqrt{\int_{\R^d} |v_k|^2 d \rho} <\infty.
\end{align*}

By Lemma~\ref{WassDSS}, for every ${\bnu} \in \Tan_\rho$,
\begin{align*}
  \quad (d_\rho f) ({\boldsymbol \nu}) 
  & =  2 \inf_{\substack{\bM \in \mathcal P_2(\R^{(K+2)d})\\
 \pi^{1, k+1}_\# \bM \in \Gamma^{\opt}(\rho,\gamma_k), \\ k=1,\ldots,K \\ \pi^{1,K+2}_\# \bM ={\bnu}}}
  \int_{\R^{(K+1)d}} \Big(\big(\sum_{k=1}^K \alpha_k (x-y_k)\big) \cdot v \Big) d \bM  \\
 & \leq 2 \int_{\R^d} \Big( v\cdot \big(\sum_{k=1}^K \alpha_k \int_{\R^d} (x-y_k) \bpi_{0,k}(dy_k|x)\big)  \Big) {\bnu}(dx,dv) \\
& =  \langle {\bnu}_0, {\bnu}\rangle_\rho.
\end{align*}

Next, we show that $ (d_\rho f) ({\bnu}_0) = \Vert {\bnu}_0 \Vert_\rho^2$.  

Let smooth vector field $\xi:= \xi(x) \in C_c^1(\R^d; \R^d)$.
 We consider the following continuity equation in a weak (Schwartz distributional) solution sense:
\begin{align*}
 \partial_t \sigma + \div_x (\sigma \xi) & =0,  \quad \forall (t,x) \in (0,1) \times \R^d, \\
  \sigma(0) & = \rho.
\end{align*} 
By classical PDE theory, there exists a unique solution $\sigma:=\sigma(t):=\sigma(t,dx)$ where the 
curve $[0,T] \ni t \mapsto \sigma(t) \in {\mathcal P}_2(\R^d) $ is continuous in $t \in [0,T]$ for any $T>0$. 
 Moreover, by a special property of the Wasserstein space (Proposition 7.3.6 of \cite{AGS08}),
\begin{align*}
\frac{d}{dt}\Big|_{t=0+}\sfdist_{\gamma_k}^2(\sigma(t)) = (-2) \int_{\R^{2d}} (y-x)\xi(x) \bpi(dx, dy), 
\quad \forall \bpi \in \Gamma^{\opt}(\rho, \gamma_k).
\end{align*}
Writing ${\bnu}^\xi(dx, dv):= \delta_{\xi(x)}(dv) \rho(dx)$, then we arrive in particular  
\begin{align*}
 (d_\rho \sfdist_{\gamma_k}^2)({\bnu}^\xi) = (- 2) \int_{\R^d} \xi(x) \cdot \big(\int_{\R^d} (y-x) \bpi_{0,k}(dy|x)\big)   \rho(dx) 
 = (-2) \int_{R^d} \big(\xi(x) \cdot v_k(x)\big)\rho(dx).
\end{align*}
Using Lipschitz continuity of ${\bnu} \mapsto (d_\rho \sfdist_{\gamma_k})({\bnu})$ in $\Tan_\rho$ (see Lemma~\ref{dDpm}), we 
approximate the ${\bnu}_0$ in \eqref{defnuu0} by those ${\bnu}^\xi$s (equivalently, approximate the $u_f$ by $\xi$s), giving
\begin{align*}
 (d_\rho \sfdist_{\gamma_k}^2)({\bnu}_0)  =(-2) \int_{R^d} \big(u_f(x) \cdot v_k(x)\big)  \rho(dx),\quad \text{ where }  \quad 
 u_f(x) := - \sum_{k=1}^K 2 \alpha_k v_k(x).
\end{align*}
Consequently,
\begin{align*}
    (d_\rho f) ({\bnu}_0) & =  \sum_{k=1}^K  \alpha_k \big(d_\rho \sfdist_{\gamma_k}^2\big) ({\bnu}_0) = 
  \int_{\R^d}  \Big(u_f \cdot (\sum_{k=1}^K  \big(-2\alpha_k) v_k\big) \Big)d \rho = \int_{\R^d} |u_f|^2 d \rho = \Vert {\bnu}_0 \Vert_\rho^2.
\end{align*}
\end{proof}

\subsubsection{A special linear subspace of the tangent cone $\Tan_\rho$}
Following Ambrosio, Gigli and Savar\'e~\cite{AGS08} (see also Appendix D.5 of Feng and Kurtz~\cite{FK06}), we define
\begin{align}\label{defLnab}
 L_{\nabla,\rho}^2 := L^2_{\nabla,\rho}(\R^d;\R^d)
 := \overline{\{ \nabla_x \varphi : \varphi \in C^\infty_c(\R^d) \}}^{L^2_\rho}
 \subset L^2_\rho:=L^2_{\rho}(\R^d;\R^d), \quad \forall \rho \in \mathcal P_2(\R^d).
\end{align}
Note that the $L^2_{\nabla, \rho}$ should be thought of as equivalent class of functions especially when the $\rho$ becomes singular (e.g. without full support on $\R^d$ or becomes non-diffusive on $\R^d$ etc etc...)
This $L^2_{\nabla, \rho}$ is a special linear subspace of the cone structure $\Tan_\rho$. 

\begin{lemma}\label{defPiproj}
Let $\rho \in \mathcal P_2(\R^d)$. For every $v \in L^2_\rho(\R^d; \R^d)$, there is a unique $\Pi_\rho(v) \in L_{\nabla, \rho}^2$ such that
  \begin{align*}
 \int_{\R^d} \big(v-\Pi_\rho(v)\big) \big(\nabla \varphi \big) d \rho =0, \quad \forall \varphi \in C_c^\infty(\R^d).
\end{align*}
In addition, the map of $\Pi_\rho: v \to \Pi_\rho(v)$ is a linear projection operator in $L_\rho^2(\R^d;\R^d)$.
 \end{lemma}
 \begin{proof}
 This is a restatement of Lemma 8.4.2 in \cite{AGS08}. See also Lemma D.49 of \cite{FK06}.
 \end{proof}
 
 The above result implies a probabilistic representation of the $\Pi_\rho$ as follow. Let $X$ be a random variable defined in some ambient probability space $(\Omega, {\mathcal F}, P)$ such that its probability distribution is equal to $\rho$. That is, $P(X \in dx) = \rho(dx)$. We have $E[|X|^2]<\infty$. Let sub-sigma field 
 \begin{align*}
{\mathcal G} :=\sigma\big\{ \nabla \varphi(X): \forall \varphi \in C^\infty_c(\R^d) \big\}
 \subset {\mathcal F}.
\end{align*}
Then the following conditional expectation representation holds
\begin{align}\label{Sec2:BPrRep}
 \Pi_\rho\big(v\big)(X) = E[ v(X) | {\mathcal G}], \quad \text{a.s.}
\end{align}
More generally, we introduce the next concept.
 \begin{definition}[Barrycentric projection]\label{BarPDef}
 For each $\bmu := \bmu(dx,dv) \in {\mathcal P}_1(\R^{2d})$ admitting the disintegration $\bmu(dx,dv)= \bmu(dv|x) \rho(dx)$, we denote its Barrycentric projection (with respect to the first marginal $\rho:= \pi^1_\# \bmu$) as
 \begin{align*}
\big({\mathcal B} \bmu\big) (x):=  \int_{\R^d} v \bmu(dv|x), \quad \text{ for } \rho \text{ \rm- a.e. in } x \in \R^d.
\end{align*} 
In the above, the disintegration is chosen so that $\big({\mathcal B} \bmu\big)$ is measurable in $x$.
 \end{definition}
 
 It follows from Lemma~\ref{WassGf}, then
 \begin{align*}
\big( {\mathcal B} \grad_\rho f\big)(x) = u_f(x), \quad \forall f \in {\mathcal S}^+, \text{ with the }  u_f \text{ given by } \eqref{defnuu0}.
\end{align*}

\begin{lemma}\label{Bproj}
For every $\rho \in \mathcal P_2(\R^d)$, ${\mathcal B}(\Tan_\rho) = L^2_{\nabla,\rho}$. 
That is, the $L^2_{\nabla, \rho}$ is image of $\Tan_\rho$ through barycentric projection.
In particular, when $\rho(dx) =\rho(x) dx$ admits a Lebesgue density, the $L^2_{\nabla, \rho}$ and $\Tan_\rho$ are identical up to an isometry.  
\end{lemma} 
\begin{proof}
This is Theorem 12.4.4 of \cite{AGS08}.
\end{proof}

\begin{lemma}\label{JensenLike}
Suppose that $\sfL :\R^d \times \R^d \mapsto \R_+$ and that
$v \mapsto \sfL(x, v)$ is lower-semicontinuous and convex for each $x \in \R^d$ fixed.  
Let $\rho \in \mathcal P_2(\R^d)$ and $v:=v(x) \in L^2_\rho(\R^d;\R^d)$. Then
\begin{align*}
 \int_{\R^d \times \R^d} \sfL\Big(x, \Pi_\rho\big(v\big)(x)\Big) \rho(dx) 
  \leq \int_{\R^d \times \R^d} \sfL\big(x, v(x)\big) \rho(dx). 
\end{align*} 
\end{lemma}
\begin{proof} 
Using the probabilistic representation \eqref{Sec2:BPrRep}, this lemma follows from Jensen's inequality. 
\end{proof}

\begin{lemma}\label{WassNisio}
Let $\bar{\sfH}$ be defined as in \eqref{barHxP}, and $L$ by \eqref{Sec1:LUV} and \eqref{Sec1:Lbnu}. 
For every $\bmu:= \bmu(dx,dP)=: \bmu(dP|x) \rho(dx) \in \Tan_\rho \subset \mathcal P_2(\R^{2d})$, we have
\begin{align*}
 \int_{\R^d}  \bar{\sfH}\Big(x, \big({\mathcal B} \bmu\big)(x); \rho\Big) \rho(dx) 
 &  \leq \sup  \big\{  \langle \bmu, {\boldsymbol \nu} \rangle_\rho  - L({\boldsymbol \nu}) : {\boldsymbol \nu} \in \Tan_\rho \big\} \\
 & \leq  \int_{\R^{2d}}  \bar{\sfH}(x, P;\rho) \bmu(dx; d P).
\end{align*}
In particular, if $\bmu(dP|x):= \delta_{u(x)}(dP)$ for some measurable function $u$, then all the above inequalities become equalities.
 \end{lemma}
\begin{proof}
We denote, for every ${\boldsymbol \nu}, \bmu \in \Tan_\rho$,  
\begin{align*}
I[\bmu, {\boldsymbol \nu}]&:= \sup  \big\{  \int_{\R^{3d}} 
 \big( v \cdot P   - \bar{\sfL}_{U,V} \big) \bM(dx; dP, d v)
  : \bM \in {\mathcal P}_2(\R^{3d}), \pi^{1,2}_\# \bM = \bmu, \pi^{1,3}_\# \bM = {\boldsymbol \nu} \big\}, \\
 I[\bmu]&:= \sup_{{\boldsymbol \nu} \in \Tan_\rho} I[\bmu, {\boldsymbol \nu}] 
 =  \sup_{{\boldsymbol \nu} \in \Tan_\rho}  \big (  \langle \bmu, {\boldsymbol \nu} \rangle_\rho  - L({\boldsymbol \nu}) \big)  .
\end{align*}
Then 
\begin{align*}
I[\bmu,{\boldsymbol \nu}] & \leq \sup\big\{  \int_{\R^{3d}} \bar{\sfH}(x, P;\rho) \bM(dx; dP, d v)
 : \bM \in {\mathcal P}_2(\R^{3d}), \pi^{1,2}_\# \bM = \bmu, \pi^{1,3}_\# \bM = {\boldsymbol \nu} \big\} \\
& =  \int_{\R^{2d}}  \bar{\sfH}(x, P;\rho) \bmu(dx; d P).
\end{align*}
To finish the proof, we only need to show that 
\begin{align}\label{Sec7:IGEB}
I[\bmu] \geq \int_{\R^d} \sfH\big(x, ({\mathcal B}\bmu)(x);\rho\big) \rho(dx).
\end{align}

For every $\varphi \in C_c^\infty(\R^d)$, we take
${\boldsymbol \nu}_\varphi(dx, dv) := \delta_{\nabla_x \varphi}(dv) \rho(dx)$.
Then by Chapter 4.4 of Gigli~\cite{Gigli04} or Theorem 12.4.4 of \cite{AGS08}, ${\boldsymbol \nu}_\varphi \in \Tan_\rho$. 
Consequently,  we have  
\begin{align*}
I[\bmu] & \geq \sup_{\varphi \in C^\infty_c(\R^d)} I[\bmu, {\boldsymbol \nu}_\varphi] = \sup_{\varphi \in C^\infty_c(\R^d)}
  \big\{ \int_{\R^{2d}} \Big(P \cdot \nabla_x  \varphi 
   - \bar{\sfL}_{U,V}\big(x,  \nabla \varphi(x)\big)\Big) \bmu(dx, dP) \big\}  \\
  & =  \sup_{\varphi \in C^\infty_c(\R^d)}
  \big\{ \int_{\R^{2d}} \Big( ({\mathcal B}\bmu)(x) \cdot \nabla_x  \varphi 
   - \bar{\sfL}_{U,V}\big(x,  \nabla \varphi(x)\big)\Big) \rho(dx) \big\}.
\end{align*}
On the other hand, for each $\epsilon>0$, there exits measurable function $u_\epsilon:=u_\epsilon(x)$ such that 
\begin{align*}
\bar{\sfH}\big(x, ({\mathcal B}\bmu)(x);\rho \big) \leq \epsilon +  ({\mathcal B}\bmu)(x) u_\epsilon(x) - \bar{\sfL}_{U,V}\big(x, u_\epsilon(x);\rho\big).
\end{align*}
We recall the space $L^2_{\nabla, \rho}$ defined in \eqref{defLnab}, and the projection operator $\Pi_\rho$ in Lemma~\ref{defPiproj}. By Lemma~\ref{Bproj}, ${\mathcal B}\bmu \in L^2_{\nabla, \rho}$. In view of Lemma~\ref{JensenLike},
\begin{align*}
&  \int_{\R^d} \Big(({\mathcal B}\bmu)(x) u_\epsilon(x) - \bar{\sfL}_{U,V}\big(x, u_\epsilon(x);\rho\big)\Big) \rho(dx) \\
& \leq \int_{\R^d} \Big(({\mathcal B}\bmu)(x) \Pi_\rho(u_\epsilon)(x) - \bar{\sfL}_{U,V}\big(x, \Pi_\rho(u_\epsilon)(x);\rho\big)\Big) \rho(dx).
\end{align*}
Combine the last three inequalities, we arrive at \eqref{Sec7:IGEB}.
\end{proof}

\subsubsection{Projection of ${\mathcal P_2}(\R^{2d})_\rho$ onto tangent cone $\Tan_\rho$}\label{projTanS}
We take some results from Section 6 of Gigli's thesis~\cite{Gigli04}. 
\begin{lemma}
Let $\bm \in {\mathcal P}_2(\R^{2d})_\rho$. Then there exists a unique $\bmu_\rho \in \Tan_\rho \sfX$ which is a minimizer in the following sense
\begin{align*}
   \sfD_\rho(\bm, \bmu_\rho)  = \inf_{\bmu \in \Tan_\rho \sfX} \sfD_\rho(\bm, \bmu).
\end{align*}
\end{lemma}
\begin{proof}
This is Propositions 4.30 in \cite{Gigli04}. 
\end{proof}
\begin{definition}[Projection onto the tangent space]\label{scrPDef}
We call the above $\bmu_\rho$ the projection of $\bm$ onto the tangent space $\Tan_\rho \sfX$, and denote it as ${\mathscr P}_\rho \bm :=\bmu_\rho$. 
\end{definition}

\begin{lemma}\label{projTan}
For each $\bm := \bm(dx; d\xi)\in {\mathcal P}_2(\R^{2d})$, 
there is a unique   $\bsigma := \bsigma(dx; d\xi, dv) \in \mathcal P_2(\R^{3d})$ such that
\begin{align}\label{sec2:projDist}
 \sfD_\rho (\bm, {\mathscr P}_\rho \bm) =\Big( \int_{(x, \xi, v) \in \R^{3d}} |\xi - v|^2 \bsigma(dx; d\xi, dv) \Big)^{1/2}.
\end{align} 
Moreover, such $\bsigma$ is given by a map in the following sense: there exists a Borel map $p:= p_\rho(x, \xi) : \R^d \times \R^d \mapsto \R^d$ such that
\begin{align*}
  \bsigma(dx; d\xi, dv) = \delta_{p(x,\xi)}(dv) \bm(dx, d \xi). 
\end{align*}
\end{lemma}
\begin{proof}
This is an adaptation of Proposition 4.32 in \cite{Gigli04}. 
\end{proof}

\begin{lemma}\label{p4coupling}
Given $\rho:=\rho(dx) \in \sfX:= \mathcal P_2(\R^d)$, $\bmu:=\bmu(dx, dP) \in \Tan_\rho \sfX$ and $\bm:= \bm(dx; d\xi) \in {\mathcal P}_2(\R^d)_\rho$. Let $\bsigma:=\bsigma(dx; d\xi, dv) \in {\mathcal P}_2(\R^{3d})$ be the lifted probability measure of $\bm$ uniquely defined in Lemma~\ref{projTan}. Then for any  $\btau:=\btau(dx; d\xi, dv, dP) \in {\mathcal P}_2(\R^{4d})$ (which is a further lift of the $\bsigma$) satisfying 
\begin{enumerate}
\item $\pi^{1,2,3}_\# \btau = \bsigma$ ,
\item $\pi^{1,4}_\# \btau = \bmu  \in \Tan_\rho \sfX$, 
\end{enumerate}
it holds that 
\begin{align}\label{xivID}
\int_{\R^{4d}} (\xi \cdot P)\btau(dx; d\xi, d v, d P) = \int_{\R^{4d}} ( v \cdot P) \btau(dx; d \xi, dv, dP).
\end{align} 
\end{lemma}
\begin{proof}
This is an adaptation of Proposition 4.33 in \cite{Gigli04}. 
\end{proof}

\subsubsection{Inferring first order derivative of simple functions through touching by distance-squared functions}
We strengthen a result stated in Lemma D.55 on page 401 of Feng and Kurtz~\cite{FK06}.   
First, Feng and Kurtz used a special notion of gradient for semi-continuous functions~\footnote{Such special notion can be generally different than the notion of gradient defined using Alexandrov space analysis techniques.} in the Wasserstein space setting. Second, the result showed that such gradient can be identified using another ``smoother" function touching from the semi-continuous function from either above or below. 
In earlier part of this section, we recalled and developed certain aspects of first order calculus in Alexandrov space. In particular, we concluded that differentials can give more information than gradients in Remark~\ref{Sec2:GradDiff}.
Next, we generalize Lemma D.55 of \cite{FK06} using the language of Alexandrov space differentials. Differentials are determined by their actions on geodesic directions. 
These directions are generated by geodesic curves connecting two given points. However, such geodesics are generally non-unique (i.e. multiple optimal plans may exist for the Kantorovich formulation of optimal transport problems). Therefore, within the scope of applications of this paper, we need to recognize a subtle distinction between  ``along {\em some} geodesic direction" vs ``along {\em every} geodesic direction"  generated by ``straight" path connecting two points \footnote{Such subtleties were already noted in numerous statements and formulation of concepts and theorems in \cite{AGS08}.}.
 The following results offer a key step for relating these two statements -- see the proof of Lemma~\ref{Sec8:bbH0bfH0}. 
 
\begin{lemma}\label{Sec2:fdtouch}
Let  $\rho_\epsilon, \gamma_\epsilon \in {\mathcal P}_2(\R^d)$, $\epsilon>0$, and $f := f_{\gamma_1, \ldots, \gamma_k,\ldots}  \in {\mathcal S}^{+,\infty}$ be a simple function as defined in \eqref{mtranf} with $K=+\infty$ (see also the general notations in \eqref{extSS+}). Suppose that
$\rho_\epsilon$ is a local maximizer of 
\begin{align}\label{fmax1}
 \rho \mapsto  -\frac{\sfd^2(\rho,\gamma_\epsilon)}{2\epsilon} - f(\rho).
\end{align}
and that $\bmu:=\bmu(dx, dy) \in \Gamma^{\opt}(\rho_\epsilon,\gamma_\epsilon)$. Then for every $\bN:= \bN(dx, dy; dy_1, \ldots, dy_K) \in {\mathcal P}_2(\R^{(K+2)d})$ which is a lift 
\footnote{The collection of such $\bN$s is non-empty, since we can at least use conditionally independent random variables to construct such coupling.}
of the $\bmu$ in the sense that 
\begin{align}\label{sec2:pikDef}
 \pi^{1, 2}_\#\bN =  \bmu, \text{ and } 
   \pi^{1, k+2}_\# \bN \in \Gamma^{\opt}(\rho_\epsilon,\gamma_k),
    \quad \text{ for } k=1,\ldots
\end{align}
we have
\begin{align*}
 \frac{x-y}{\epsilon} = 2 \sum_k \alpha_k(\rho_\epsilon;\gamma_1, \ldots, \gamma_K) \big(y_k-x\big), \quad \bN \text{\rm - almost everywhere},
\end{align*}
with the $\alpha_k:=\alpha_k(\rho_\epsilon; \gamma_1, \ldots, \gamma_K)$s defined as in \eqref{alphaDef}.
\end{lemma}
\begin{remark}\label{Sec2:fdRMK}
To streamline main arguments in the proof, we present the case as if $K$ is finite. 
There is no essential changes in the case when $K=\infty$ (i.e. countable $\alpha_k$s), as long as we have an extra property
\begin{align*}
\sup_{\rho} \sum_{k=1}^\infty \alpha_k^2(\rho; \gamma_1, \ldots) < \infty,
\end{align*}
where the $\sup_\rho$ above is over metric-balls of arbitrary but finite radius in $\sfX:={\mathcal P}_2(\R^d)$.
\end{remark} 
 
\begin{proof}
The maximizing property in \eqref{fmax1} implies that
\begin{align}\label{sec2:dfdD}
\big( d_{\rho_\epsilon} f\big)(\bnu) \geq  \Big(d_{\rho_\epsilon} \big(- \frac{\sfdist^2_{\gamma_\epsilon}}{2\epsilon}\big)\Big) (\bnu), \quad \forall \bnu \in \Tan_{\rho_\epsilon} \sfX.
\end{align}
We claim that this further implies that, for each $\bN$ satisfying \eqref{sec2:pikDef},  the following must holds:
\begin{align}\label{sec2:xyxi}
 & \int_{\R^{(K+2)d}} \big(2 \sum_{k=1}^K \alpha_k (x-y_k) \cdot \xi(x,y,y_1, \ldots, y_K) \big) \bN(dx, dy, dy_1, \ldots, d y_K) \\
& \qquad \geq
 \int_{\R^{(K+2)d}} \Big(\frac{(y-x)}{\epsilon} \cdot \xi(x,y,y_1, \ldots, y_K) \Big) \bN(dx,dy, dy_1, \ldots, dy_K), \nonumber \\
 & \qquad \qquad \qquad 
  \quad \forall \xi :=\xi(x,y, y_1,\ldots, y_K)\in C^\infty_c(\R^{(2+K)d};\R^d). \nonumber
\end{align}
Next, by arbitrariness of the $\xi$, the above inequality holds with $\xi$ replaced by $-\xi$ as well. Hence the inequality is indeed an equality, giving
\begin{align*}
 \int_{\R^{(K+2)d}} \Big(\big(2 \sum_{k=1}^K \alpha_k (x-y_k) - \frac{y-x}{\epsilon} \big)\cdot \xi(x,y,y_1, \ldots, y_K) \Big) \bN(dx, dy, dy_1, \ldots, d y_K)=0.
\end{align*}
Consequently 
\begin{align*}
 \int_{\R^{(K+2)d}} \Big|\big(2 \sum_{k=1}^K \alpha_k (x-y_k) - \frac{y-x}{\epsilon} \big|^2\bN(dx, dy, dy_1, \ldots, d y_K)   =0,
\end{align*}
giving conclusion of the lemma.

Next, we verify the claim \eqref{sec2:xyxi} in six steps. 

First, we define a lift of the $\bN$ by attaching two more variables $\xi$ and $P$:
\begin{align*}
& \widehat{\bN} (dx, dy, dy_1, \ldots, d y_K; d \xi, dP) \\
&  \qquad :=  \delta_{\xi(x,y, y_1, \ldots, y_K)}(d\xi)   
     \delta_{\sum_{k=1}^K \alpha_k (x-y_k)  }(d P)   
 \bN(dx,dy, dy_1, \ldots, dy_K) \in \mathcal P_2(\R^{(4+K)d}).
\end{align*}
Its projection into the $(x,\xi)$-variable gives
\begin{align*}
\bm:= \bm(dx, d\xi):= (\pi^{1, 3+K}_\# \widehat{\bN} )(dx, d\xi) 
  \in {\mathcal P}_2(\R^{2d})_{\rho_\epsilon} \subset {\mathcal P}_2(\R^{2d}).
\end{align*} 
In general  $\bm \not \in \Tan_{\rho_\epsilon} \sfX$, but we may consider its projection to the tangent cone as given by Definition~\ref{scrPDef}:
\begin{align}\label{nuprojm}
\bnu := \bnu_{\rho_\epsilon}(dx, dv):= {\mathscr P}_{\rho_\epsilon} \bm \in \Tan_{\rho_\epsilon} \sfX.
\end{align}
We also introduce notations for other projected marginal measures
\begin{align}\label{Sec2:projbmu}
 \bmu:=\bmu(dx, d P):=  \pi^{1, K+4}_\# \widehat{\bN},
\end{align}
and
\begin{align*}
 \bmu_k(dx,dP):= \int_{ y_1,\ldots, y_K}\delta_{x-y_k}(dP) \bN(dx, dy_1,\ldots, dy_K). 
\end{align*}
Since $\pi^{1,k+2}_\# \bN \in \Gamma^{\opt}(\rho_\epsilon, \gamma_k)$ 
 (see \eqref{sec2:pikDef}), by the first part of Lemma~\ref{negnu}, 
we have
\begin{align*}
 \bmu_k = (-1) \cdot \int_{ y_1,\ldots, y_K}\delta_{y_k-x}(dP) \bN(dx, dy_1,\ldots, dy_K) \in \Tan_{\rho_\epsilon} \sfX.
\end{align*}
Consequently,  by Lemmas~\ref{addmeas}, \ref{negnu} and \ref{Gig4-25},
\begin{align*}
\bmu \in \oplus_k \big( \alpha_k \cdot \bmu_k\big) \subset \Tan_{\rho_\epsilon}\sfX.
\end{align*}
 
Second, we identify structure of the $\bnu$ in \eqref{nuprojm} more explicitly. By Lemma~\ref{projTan}, there exists a Borel map $p:=p_{\rho_\epsilon}(x, \xi) :\R^d \times \R^d \mapsto \R^d$ such that 
\begin{align*}
 \bsigma(dx; d \xi, d v) := \delta_{p(x,\xi)}(dv) \bm(dx, d \xi).
\end{align*}
We also have $\pi^{1,3}_\# \bsigma = {\mathscr P}_{\rho_\epsilon}\bm = \bnu$.
 With such $p$, we further lift the $\widehat{\bN}$ by attaching one more variable $v$:
\begin{align*}
 \widehat{\widehat{\bN}}(dx,dy, dy_1, \ldots, dy_K; d \xi, dP, dv)
 & :=  \delta_{p(x,\xi)}(dv)  \widehat{\bN}(dx,dy, dy_1, \ldots, dy_K; d \xi, dP) \\
 & \qquad \in {\mathcal P}_2(\R^{(5+K)d}).
\end{align*}
Let 
\begin{align*}
\bM_\epsilon := \bM_\epsilon(dx, dy_1, \ldots, dy_K, dv) :=\pi^{1,3,\ldots, K+2, K+5}_\#  \widehat{\widehat{\bN}}. 
\end{align*}
Then
\begin{align}\label{Mepspj1}
 \pi^{1,k+1}_\# \bM_\epsilon = \pi^{1,k+2}_\#  \widehat{\widehat{\bN}} \in \Gamma^{\opt}(\rho_\epsilon, \gamma_k), \quad k=1,2,\ldots, K,
\end{align}
and
\begin{align}\label{Mepspj2}
 \pi^{1,K+2}_\# \bM_\epsilon = \pi^{1, K+5}_\#  \widehat{\widehat{\bN}} = \pi^{1,3}_\# \bsigma = \bnu.
\end{align}
Relationship among the above various marginal probability measures is rather involved. However, the intuition is simple: We are merely introducing more and more random variables living in a same probability space whose joint-distributions is compatible with the various marginal measures. We do this so that integrations with respect to these measures can be viewed as expectations in a fixed ambient probability space. Such probabilistic coupling techniques can be graphically represented using Figure~\ref{sec2:figRVrep}. Readers are invited to re-write our proof using expectations of random variables in the lifted probability space $(\Omega, {\mathcal F}, \P)$, and verify that the measure-theoretic arguments here have an alternative presentation using submetry projection arguments from the lifted probabilistic formulation. 
  
\begin{figure}[h]\label{sec2:figRVrep}
\centering

\begin{align*}
\begin{tikzpicture}
\draw[black, thick] (-4,0) --(6,0);
\draw (-3,4) -- (-3,0);
\filldraw (-3,2) circle (1pt) node [anchor= east] {$X_\epsilon$} ;
\draw (-2,4) -- (-2,0);
\filldraw (-2,2) circle (1pt) node [anchor = east] {$Y_\epsilon$}; 
\draw (-1,4) -- (-1,0);
\filldraw (-1,2) circle (1pt) node [anchor=east] {$Y_1$};
\draw (0,4) -- (0,0);
\filldraw (0,2) circle (1pt) node [anchor=east] {$\dots$};
\draw (1,4) -- (1,0);
\filldraw (1,2) circle (1pt) node [anchor=east] {$Y_K$};
\draw (2,4) -- (2,0);
\filldraw (2,3) circle (1pt) node [anchor=west] {$\Xi$};
\draw (3,4) -- (3,0);
\filldraw (3,3.7) circle (1pt) node [anchor=west] {$P$};
 \draw (4,4) -- (4,0);
\filldraw (4,1.3) circle (1pt) node [anchor=west] {$V$};
\node (a)  at (-3,-0.3) {$\rho_\epsilon$};
\node (b) at  (-2, -0.3) {$\gamma_\epsilon$};
\node (c)  at (-1,-0.3) {$\gamma_1$};
\node (d)  at (0,-0.3) {\ldots};
\node (e)  at (1,-0.3) {$\gamma_K$};
\node (f)  at (2, -0.3) {$\pi^2_\# \bm $};
\node (g)  at (3, -0.3) {$\pi^2_\# \bmu$};
\node (h)  at (4, -0.3) {$\pi^3_\# \bsigma$};
\node (Y) at (-5, 4.0) {$ L^2(\Omega, {\mathcal F}, {\mathbb P})$};
\node (X) at (-5, 0.0) {${\mathcal P}_2(\R^d)$};
\draw [->] (Y) -- (X)  node[midway, left] {$\sfp$};
\end{tikzpicture}
\end{align*}

\caption{A graphical representation of marginal probability measures as projection, through a submetry map $\sfp$, of square integrable random variables, where the probability space $(\Omega, {\mathcal F}, {\mathbb P}):= ([0,1], {\mathcal B}_{[0,1]}, \text{Leb})$ and random variables
$\Xi:=\xi(X_\epsilon,Y_\epsilon, Y_1,\ldots, Y_K)$, $P:=\sum_k \alpha_k (X_\epsilon -Y_k)$ and $V:= p(X_\epsilon, \Xi)$. The measure $\widehat{\widehat{\bN}}$ is the joint distribution of $(X_\epsilon, Y_\epsilon, Y_1,\ldots, Y_K, \Xi,P, V)$, which can be viewed as a ``section" in the graph.} 
\end{figure}
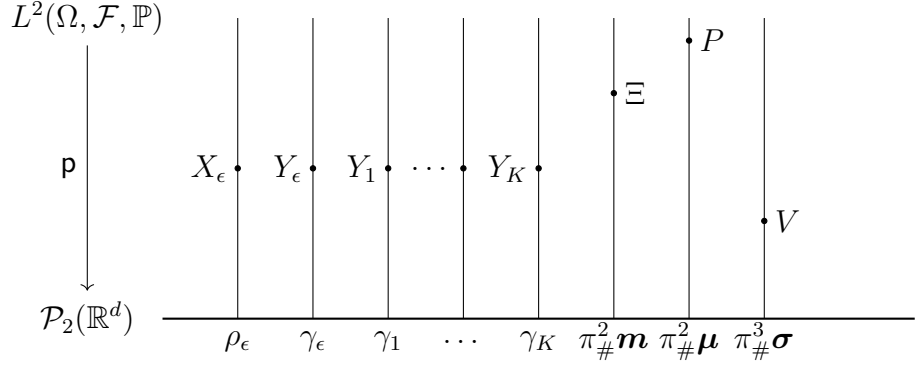

Third, denoting
\begin{align*}
\btau := \btau(dx; d \xi, d v, dP):= \pi^{1,3+K, 5+K,4+K}_\#  \widehat{\widehat{\bN}}, 
\end{align*}
then 
\begin{align*}
  \pi^{1,2,3}_\# \btau = \bsigma, \quad \pi^{1,4}_\# \btau = \bmu   \in \Tan_{\rho_\epsilon} \sfX.
\end{align*}
where the $\bmu$ is defined in \eqref{Sec2:projbmu}. 
By Lemma~\ref{p4coupling},  the following holds
\begin{align*}
\int_{\R^{4d}} (\xi \cdot P)\btau(dx; d\xi, d v, d P) = \int_{\R^{4d}} ( v \cdot P) \btau(dx; d \xi, dv, dP).
\end{align*}
Next, we note that, on one hand,
\begin{align*}
\int_{\R^{4d}} (P \cdot v ) \btau(dx; d \xi, dv, dP) &= \int_{\R^{(5+K)d}} (P \cdot v )  \widehat{\widehat{\bN}}(dx, dy, dy_1, \ldots, dy_K; d \xi, dP, dv)  \\
& = \int_{\R^{(2+K)d}} \big( 2 \sum_{k=1}^K \alpha_k (x-y_k)  \cdot v \big) \bM_\epsilon(dx, dy_1, \ldots, dy_K, dv).
\end{align*}
On the other hand,
\begin{align*}
& \int_{\R^{4d}} (P \cdot \xi ) \btau(dx; d \xi, dv, dP) \\
&= \int_{\R^{(5+K)d}} (P \cdot \xi )  \widehat{\widehat{\bN}}(dx, dy, dy_1, \ldots, dy_K; d \xi, dP, dv)  \\
& = \int_{\R^{(2+K)d}} \big( 2 \sum_{k=1}^K \alpha_k (x -y_k)  \cdot \xi(x,y,x_1, \ldots, x_K) \big) \bN(dx,dy, dy_1, \ldots, dy_K).
\end{align*}
Consequently
\begin{align}\label{veqnxi}
 & \int_{\R^{(2+K)d}} \big( 2 \sum_{k=1}^K \alpha_k (x-y_k)  \cdot v \big) \bM_\epsilon(dx, dx_1, \ldots, dx_K, dv) \\
& =  \int_{\R^{(2+K)d}} \big( 2 \sum_{k=1}^K \alpha_k (x-y_k)  \cdot \xi(x,y,y_1, \ldots, y_K) \big) \bN(dx,dy, dy_1, \ldots, dy_K). \nonumber
\end{align}

In fourth step, we have
\begin{align*}
\big(d_{\rho_\epsilon} f \big)(\bnu)  
& =   \inf_{\substack{\bM \in \mathcal P_2(\R^{(K+2)d})\\ 
   \pi^{1, k+1}_\# \bM \in \Gamma^{\opt}(\rho_\epsilon,\gamma_k),  k=1,\ldots,K \\
   \pi^{1,K+2}_\#\bM = \bnu} }
  \int_{\R^{(K+2)d}} 2 \big(\sum_{k=1}^K \alpha_k (x-y_k) \cdot v\big) \bM(dx, dy_1, \ldots, d y_K; dv),
 \\
 & \leq \int_{\R^{(K+2)d}} 2\big(\sum_{k=1}^K \alpha_k (x-y_k) \cdot v  \big) \bM_\epsilon(dx, dy_1, \ldots, d y_K; dv) \\
 & =  \int_{\R^{(2+K)d}}  \big( 2 \sum_{k=1}^K \alpha_k (x -y_k)  \cdot \xi(x,y,y_1, \ldots,y_K) \big) \bN(dx,dy, dy_1, \ldots, dy_K).
\end{align*}
In the above, the first equality follows from by Lemma~\ref{WassDSS}, the first inequality follows from \eqref{Mepspj1} and \eqref{Mepspj2}, and the last equality follows from \eqref{veqnxi}.

In the fifth step, invoking \eqref{sec2:ddsqr} (or Remark~\ref{dggamma}) and proceed with similar derivations as the above four steps, we conclude the following. For the $\bnu$ constructed in \eqref{nuprojm},  in view of \eqref{sec2:ddsqr} 
and \eqref{Sec2:ScaDef}, we have  
\begin{align*}
\Big(d_{\rho_\epsilon} \big(- \frac{\sfdist_{\gamma_\epsilon}^2}{2\epsilon}\big)\Big) (\bnu)
&  =  \sup_{\substack{\pi^{1,2}_\#\bM \in \Gamma^{\opt}(\rho_\epsilon, \gamma_\epsilon)\\ \pi^{1,3}_\# \bM = \bnu }} \int_{\R^{3d}} \Big(\frac{(y-x)}{\epsilon} \cdot v \Big) \bM(dx,dy;dv) \\
& \geq \int_{\R^{(K+2)d}} \Big(\frac{(y-x)}{\epsilon} \cdot \xi(x,y, y_1, \ldots, y_K) \Big) \bN(dx,dy, dy_1, \ldots, d y_K).
\end{align*}

In the sixth step, we combine estimates in the fourth and fifth steps with \eqref{sec2:dfdD} to conclude \eqref{sec2:xyxi}.
 \end{proof}

In a similar way, we can prove the following.
\begin{lemma}\label{Sec2:gdtouch}
Let $\rho_\epsilon, \gamma_\epsilon \in \mathcal P_2(\R^d)$, $\epsilon >0$, and
 $g:= g_{\rho_1, \ldots, \rho_k,\ldots} \in \mathcal S^{-,\infty}$ be the simple function defined in \eqref{mtrang} with $K=+\infty$ (see \eqref{extSS-} for more general situation). Suppose that
  $\gamma_\epsilon$ is a local maximizer of
\begin{align}\label{gmax2}
 \gamma \mapsto g(\gamma) -\frac{\sfd^2(\rho_\epsilon,\gamma)}{2\epsilon},
\end{align}
and that $\bmu:=\bmu(dx, dy) \in \Gamma^{\opt}(\rho_\epsilon,\gamma_\epsilon)$. 
Then for any $\bN:= \bN(dx, dy; dx_1, \ldots, dx_K) \in \mathcal P_2(\R^{(K+2)d})$ which is a lift 
 of the $\bmu$ in the sense that 
\begin{align*}
 \pi^{1, 2}_\#\bN =  \bmu, \text{ and } \pi^{2, k+2}_\# \bN \in \Gamma^{\opt}(\gamma_\epsilon,\rho_k),  k=1,\ldots,K,
\end{align*}
we have
\begin{align*}
 \frac{(y-x)}{\epsilon} = 2 \sum_{k=1}^K \beta_k(\gamma_\epsilon;\rho_1, \ldots, \rho_K) \big(x_k-y\big), \quad \bN \text{\rm - almost everywhere}.
\end{align*}
where the $\beta_k:=\beta_k(\gamma_\epsilon; \rho_1, \ldots, \rho_K)$ are defined as in \eqref{betaDef}.
\end{lemma}

 \begin{remark}
In the above two lemmas, we used $\mathcal S^{\pm, \infty}$ instead of just $\mathcal S^{\pm}$. This is because that, using Borwein-Preiss smooth perturbed optimization principle (Lemma~\ref{BorPre}), we are guaranteed to have extremal points $\rho_\epsilon$ in Lemma~\ref{Sec2:fdtouch} and $\gamma_\epsilon$ in Lemma~\ref{Sec2:gdtouch} for a large sub-class of functions in $\mathcal S^{\pm,\infty}$.
 \end{remark}
 
\newpage
 
 \section{Viscosity solution theory in metric spaces, and projection of equations from metrically foliated spaces}\label{VisMetProj}  
We are interested in viscosity solutions in quotient metric spaces (see heuristic discussions in Section~\ref{intro}).
For such purpose, we develop abstract results concerning projection of viscosity solutions through submetry maps. This is done in Sections~\ref{projHJ} and \ref{HJproj}, after we introduce (generalized) notions of viscosity solution for equations in metric spaces in Section~\ref{Sec:VisDef}. 

For a quick introduction on the concept and properties of submetry, see Appendix~\ref{App:SubM}.

\subsection{Definitions of viscosity solution}\label{Sec:VisDef}
Let $(\sfX, \sfd)$ be a metric space,
$\alpha >0$ and $h \in B(\sfX)$. We allow operators be multivalued and identify them with their graphs.
We consider an operator $H \subset M(\sfX; \bar{\R}) \times M(\sfX; \bar{\R})$, 
and sub-solution $\overline{f}$ and super-solution $\underline{f}$ respectively to equations formally written as inequalities
\begin{align}
 \overline{f} - \alpha H \overline{f} & \leq h, \label{sub} \\
 \underline{f} - \alpha H \underline{f} & \geq h. \label{sup}
\end{align}
Next, we introduce versions of viscosity solution in this context. The following is an adaptation of Definition 7.1 of Feng and Kurtz~\cite{FK06}. 
Motivations for such definition came from representation theorems on dissipative operators in function spaces (see Sato~\cite{Sato68} for details and Appendix A.3 of \cite{FK06} for a summary).

\begin{definition}[Sequential definition of viscosity solution]\label{SeqVisDef}
We call $\overline{f}$ a viscosity sub-solution to \eqref{sub} {\it in the sequential sense}, if 
\begin{enumerate}
\item $\overline{f} \in M(\sfX; \bar{\R})$;
\item for every $(f_0,g_0) \in H$ with 
$\sup_\sfX (\overline{f} - f_0)<\infty$, {\em there exists} $\{ x_n \}_{n=1,2,\ldots} \subset \sfX$ satisfying 
\begin{align}\label{seqffbar}
 \lim_{n \to \infty} (\overline{f} - f_0)(x_n) = \sup_\sfX (\overline{f} - f_0)
\end{align}
and
\begin{align}\label{seqgg0}
 \liminf_{n \to \infty} \big( \overline{f} - h - \alpha g_0 \big)(x_n) \leq 0.
\end{align}
\end{enumerate}
In the above definition, if the (2) is replaced by the following (2a), then we call $\overline{f}$ 
a {\em strong} viscosity sub-solution {\it in the sequential sense}:
\begin{enumerate}
\item [(2a)]  for every $(f_0,g_0) \in H$ with 
$\sup_\sfX (\overline{f} - f_0)<\infty$, and {\em for every} $\{ x_n \}_{n=1,2,\ldots} \subset \sfX$ satisfying 
\eqref{seqffbar}, we have \eqref{seqgg0}.
\end{enumerate}
 
We call $\underline{f}$ a viscosity super-solution to \eqref{sup} {\em in the sequential sense}, if 
 \begin{enumerate}
 \item $\underline{f} \in M(\sfX; \bar{\R})$ ;
 \item for every $(f_1, g_1) \in H$ with 
$\sup_\sfX(f_1 - \underline{f})<\infty$, {\em there exists} $\{ x_n \}_{n =1,2,\ldots} \subset \sfX$ satisfying 
\begin{align*}
 \lim_{n \to \infty} (f_1 - \underline{f})(x_n) = \sup_\sfX (f_1 - \underline{f})
\end{align*}
and
\begin{align*}
 \liminf_{n \to \infty} \big( \underline{f} - h - \alpha g_1 \big)(x_n) \geq 0.
\end{align*}
\end{enumerate}
Similarly, we define {\it strong} viscosity super-solution {\it in the sequential sense}.

If a function is both a sub-solution and a super-solution in the sequential sense, it is called a solution in the sequential sense. Similarly, we define strong solution in the sequential sense.
\end{definition}

Given their connections to dissipative and strongly dissipative operators in the Banach space of bounded functions (see Appendix A.3 of \cite{FK06} for references), the above definitions are natural and more convenient to use, when the underlying metric space $\sfX$ is non-locally compact. In particular, we recall that strong (in norm) infinitesimal generator of a possibly nonlinear contraction operator semigroup in Banach space is strongly dissipative. If, however, the $\sfX$ is locally compact, definitions and techniques can simplify. We are therefore led to the following concepts, which are most frequently used in partial differential equation literature.  

\begin{definition}[Pointwise definition of viscosity solution]\label{PtVisDef}
We call $\overline{f}$ a viscosity sub-solution to \eqref{sub} {\it in the point-wise definition sense}, if 
\begin{enumerate}
\item $\overline{f} \in M(\sfX; \bar{\R})$; 
\item for every $(f_0,g_0) \in H$ with 
$\sup_\sfX (\overline{f} - f_0)<\infty$, there {\em exists} $x_0 \in \sfX$ satisfying 
\begin{align}\label{vDefmaxf}
(\overline{f} - f_0)(x_0) = \sup_\sfX (\overline{f} - f_0)
\end{align}
and
\begin{align}\label{vsDefHf}
  \overline{f}(x_0) - h(x_0) \leq  \alpha g_0(x_0).
\end{align}
\end{enumerate}
In the above definition, if the second point is replaced by the following (2A), 
then we call $\overline{f}$ a {\em strong} viscosity sub-solution {\it in the point-wise definition sense}:
\begin{enumerate}
\item [(2A)] for every $(f_0,g_0) \in H$ with 
$\sup_\sfX (\overline{f} - f_0)<\infty$, and {\em for every} $x_0 \in \sfX$ satisfying 
\begin{align*}
(\overline{f} - f_0)(x_0) = \sup_\sfX (\overline{f} - f_0)
\end{align*}
we have
\begin{align*}
  \overline{f}(x_0) - h(x_0) \leq  \alpha g_0(x_0).
\end{align*}
\end{enumerate}

We call $\underline{f}$ a super-solution to \eqref{sup} {\it in the point-wise definition sense}, if
 \begin{enumerate}
 \item $\underline{f} \in M(\sfX; \bar{\R})$ 
 \item for every $(f_1, g_1) \in H$ with 
$\sup_\sfX(f_1 - \underline{f})<\infty$, there exists $x_1 \in \sfX$ satisfying 
\begin{align*}
 (f_1 - \underline{f})(x_1) = \sup_\sfX (f_1 - \underline{f})
\end{align*}
and
\begin{align*}
  \underline{f}(x_1) - h(x_1) \geq \alpha g_1 (x_1).
\end{align*}
\end{enumerate}
{\em Strong} viscosity super-solution {\it in the point-wise definition sense} is defined similarly.

In the above point-wise definitions of viscosity solutions, if a function is both a sub-solution and a super-solution, it is called a solution. Similarly, being both strong sub-solution and strong super-solution defines a strong solution.
\end{definition}
 
\begin{remark}
We will frequently work with upper semicontinuous sub-solution and lower semicontinuous super-solution. When that is the case, we will make these restrictions explicit in respective statements.

Given a function, while we can always find a sequence of points $\{ x_n\}_{n=1,2,\ldots}$ approximating supremum or infimum, there is no guarantee that we can always find a point $x_0$ that attains the extreme. Hence the above notion of point-wise viscosity definition has the risk of being vacuous. Note also that the defining property for strong point-wise solution only needs to hold when extremizing point exists. Therefore, strong sub- (super-) solution in the point-wise sense does not necessarily imply sub- (super-) solution in the point-wise sense.  However, with careful construction of test functions, one can frequently guarantee existence of extremal point(s). One possibility is that we require test functions to have compact finite level sets and proper semicontinuity properties. Note that this means that, in the case of non-locally compact $\sfX$, we are forced to consider test functions which are not continuous but merely semi-continuous. Feng and Katsoulakis~\cite{FKa09}, Feng and Kurtz~\cite{FK06},  Feng and Nguyen~\cite{FN12} give examples of this kind. Lemmas~\ref{seq2pw} and \ref{seq2pw2} next illustrate some basic properties.  There is another approach to guarantee existence of extremal points. It has the benefit of allowing us to use continuous test functions in non-locally compact metric state space settings. Such approach has a longer history. 
Starting from 1985, Crandall and Lions~\cites{CL85,CL86,CL86b,CL90,CL91,CL94,CL94b} published a series of works developing viscosity solution in Banach (mostly Hilbert) spaces.  The first two papers in that series introduced a perturbative method for constructing test functions by adding small perturbation of a distance function (e.g. Ekeland's principle~\cite{Eke79}). In fact, one can also use a smooth version of perturbation by adding combinations of distance-squared functions (e.g. Borwein-Preiss~\cite{BorPre87}). See Feng and Swiech~\cite{FS13}  and Ambrosio and Feng~\cite{AF14} for illustrations.
\end{remark}

\subsubsection{Sequential solution versus point-wise solution}
From the defining relations of viscosity solutions, sub- (resp. super-)solution in the point-wise viscosity sense always implies sub- (resp. super-) viscosity solution property in the sequential sense. Under the following extra conditions, a type of converse also holds. We show this next.
 \begin{condition}\label{LypSub}
For every $(f_0, g_0) \in H$ and every $C \in \R$, the following sub-level is compact:
 \begin{align*}
  \{ x \in \sfX : \big( f_0 - \alpha g_0\big) (x) \leq C\} \subset \subset \sfX.
\end{align*}
 \end{condition}

\begin{lemma}\label{seq2pw}
Suppose that $H \subset \LSC(\sfX; \bar{\R}) \times \USC(\sfX; \bar{\R})$, 
that $h \in C(\sfX)$ and is bounded from above, and that Condition~\ref{LypSub} holds. 

Then for every $\overline{f} \in \USC(\sfX;\bar{\R})$ which is a sequential viscosity sub-solution, it is also a point-wise viscosity sub-solution (with guaranteed existence of the extremal points). 
\end{lemma}
\begin{proof} Suppose that $\overline{f}$ is a sequential sub-solution. Then for every $(f_0, g_0)  \in H$, 
there exists $\{x_n \} \subset \sfX$ such that  
\begin{align*}
 c:= \lim_{n \to \infty} ( \overline{f}-f_0)(x_n) = \sup_\sfX (\overline{f} - f_0), \quad
  \limsup_{n \to \infty} \big( \overline{f} - h - \alpha g_0 \big)(x_n) \leq 0.
\end{align*}
Hence $\limsup_{n \to \infty} \big( f_0 - \alpha g_0 \big) (x_n) \leq \limsup_{n \to \infty} h(x_n) -c$.
By Condition~\ref{LypSub},  there exists $x_0 \in \sfX$ with (relabeling a sequence if necessary) $\lim_{n \to \infty} x_n = x_0$.

In addition, since $\overline{f} \in \USC(\sfX; \bar{\R})$ and $f_0 \in \LSC(\sfX; \bar{\R})$, 
\begin{align*}
 \lim_{n \to \infty} (\overline{f} - f_0)(x_n) = \sup_\sfX (\overline{f} - f_0) =   (\overline{f} - f_0)(x_0).
\end{align*}
This further implies that 
\begin{align*}
\liminf_{n \to \infty} \overline{f} (x_n) = 
\lim_{n \to \infty} (\overline{f} - f_0)(x_n) + \liminf_{n \to \infty} f_0(x_n) 
\geq (\overline{f} - f_0)(x_0) +  f_0(x_0) = \overline{f}(x_0),
\end{align*}
giving $\lim_{n \to \infty} \overline{f} (x_n) =  \overline{f}(x_0)$. 
Similarly, we can derive $\limsup_{n \to \infty}  f_0(x_n) \leq f_0(x_0)$,
 hence $\lim_{n \to \infty}  f_0(x_n) = f_0(x_0)$. Since $g_0 \in \USC(\sfX; \bar{\R})$ and $h\in C(\sfX)$,
\begin{align*}
\overline{f}(x_0) - h(x_0) = \limsup_{n \to \infty} \big(\overline{f}(x_n) - h(x_n)\big) \leq \alpha \limsup_{n \to \infty} g_0(x_n) \leq \alpha g_0(x_0).
\end{align*}
That is, $\overline{f}$ is a sub-solution in the point-wise definition sense.
\end{proof}

In a similar vein, we can prove the following.
 \begin{condition}\label{LypSup}
 For every $(f_1, g_1) \in H$ and every $C \in \R$, the following sup-level is compact
 \begin{align*}
\{ x \in \sfX : \big(f_1 - \alpha g_1\big)(x) \geq C\} \subset \subset \sfX.
\end{align*}
 \end{condition}
 
\begin{lemma}\label{seq2pw2}
Suppose that $H \subset \USC(\sfX; \bar{\R}) \times \LSC(\sfX; \bar{\R})$, 
that $h \in C(\sfX)$ and is bounded from below, and that Condition~\ref{LypSup} holds. 

Then for every $\underline{f} \in \LSC(\sfX; \bar{\R})$ which is a sequential viscosity super-solution, it is also a point-wise viscosity super-solution (with guaranteed existence of extremal points).
\end{lemma}
 
\subsubsection{Point-wise strong solution implies point-wise solution}
The relationship between point-wise strong solution vs. point-wise solution is a bit subtle. Note that, in our definitions, the {\em strongness} of point-wise solution does not require {\em a priori} extremizing point always exist. It merely requires relevant inequalities to hold at those extremizing points, {\em once they exist}. In contrast, the definition of point-wise viscosity solution requires existence of extremizing point, always. In other words, verifying a function is a point-wise viscosity solution requires construction of the extremizing point(s) first.
However, if we assume that domain of the Hamiltonian operator is chosen so that there will always be extremizing point, then strong point-wise sub- (resp. super-) solution imply point-wise sub- (resp. super-) solution. Next, we give a condition so that such assumption can be readily verified. 

We consider the case of sub-solutions. Let $(\sfX, \sfd)$ be a metric space, $H$ is a possibly multi-valued operator with its graph $H \subset \LSC(\sfX) \times \USC(\sfX; \bar{\R})$. 

\begin{condition}\label{cndKlev}
For each $f_0 \in D(H)$,
\begin{enumerate}
\item there exists a sub-linear function $\beta:=\beta_{f_0}: \R \mapsto \R$ 
such that $\bar{f} \leq \beta\circ f_0$;
\item  $f_0$ has compact finite sub-levels. 
\end{enumerate}
\end{condition}

\begin{lemma}\label{spw2pw}
Let $\bar{f} \in \USC(\sfX; \bar{\R})$ be a point-wise strong viscosity sub-solution to \eqref{sub}.  
Suppose that $\bar{f}$ and $H$ satisfies Condition~\ref{cndKlev}.
Then the $\bar{f}$ is also a point-wise viscosity sub-solution.
\end{lemma}

Similarly, one can state a result for the super-solution case. 

\subsubsection{Sequential viscosity solution implies strong point-wise solution, for local Hamiltonian operators}
Next, we prove that point-wise solution can become strong point-wise solution in a wide variety of situations.
Without pursuing generality, we only consider scenarios where the Hamiltonian operator has a special property \eqref{Sec3:gapp} which is natural for local operators such as differential operators.  In fact, we will prove a stronger result that, with such special property, sequential viscosity solution becomes strong point-wise viscosity solution.  

The following is an adaptation of Lemma 3.6 in Feng~\cite{Feng06}.
\begin{lemma}\label{Sec3:locF}
 Let $\overline{f} \in \USC(\sfX; \bar{\R})$, $f_0 \in D(H)$ and $x_0 \in \sfX$ satisfy $(\overline{f} - f_0)(x_0) = \sup_{\sfX} (\overline{f} - f_0)$. We introduce a perturbation of the $f_0$ by
\begin{align*}
f_{\theta}(x) := f_0(x) + \theta \sfd^2(x, x_0), \quad \forall \theta >0.
\end{align*}
 Then the following properties hold:
\begin{enumerate}
\item  $(\overline{f} - f_{\theta})(x_0) > (\overline{f} - f_{\theta})(x)$ for every $x \neq x_0$.
\item  for every $\{ x_{n,\theta} \} \subset \sfX$ satisfying 
\begin{align}\label{Sec3:seqmax}
 \lim_{n \to \infty} (\overline{f} - f_{\theta})(x_{n,\theta}) = \sup_{\sfX} (\overline{f} - f_{\theta}),
\end{align}
we have $\lim_{n \to \infty} \sfd(x_{n,\theta}, x_0) =0$, and 
\begin{align*}
\lim_{n \to \infty} \overline{f}(x_{n,\theta}) = \overline{f}(x_0), \quad
 \lim_{n \to \infty}  f_0(x_{n,\theta}) =  f_0(x_0).
\end{align*}
\end{enumerate}
\end{lemma}
\begin{proof}
We only need to prove the second property. 

By the assumptions,
\begin{align*}
(\overline{f} - f_0)(x_0)   = (\overline{f} - f_\theta)(x_0)  
  & \leq \sup_\sfX( \overline{f} - f_{\theta})
 = \lim_{n \to \infty} (\overline{f} - f_\theta)(x_{n,\theta}) \\
 & = \lim_{n \to \infty} \big( (\overline{f} - f_0)(x_{n,\theta}) - \theta \sfd^2(x_{n,\theta}, x_0)\big) \\
 & \leq \lim_{n \to \infty}   (\overline{f} - f_0)(x_{n,\theta}) \leq \sup_\sfX  (\overline{f} - f_0) 
   =    (\overline{f} - f_0)(x_0).
\end{align*}
Consequently, 
\begin{align*}
(\overline{f} - f_\theta)(x_0)= \sup_\sfX (\overline{f} - f_\theta) = \lim_{n \to \infty} (\overline{f} - f_\theta)(x_{n,\theta}) = \lim_{n \to \infty}(\overline{f} - f_0)(x_{n,\theta}) = (\overline{f} - f_0)(x_0),
\end{align*}
giving
\begin{align*}
\lim_{n \to \infty} \theta \sfd^2(x_{n,\theta}, x_0) = 
  \lim_{n \to \infty} \Big( \big( \overline{f} - f_0\big)(x_{n,\theta})  - \big(\overline{f} - f_\theta\big)(x_{n,\theta}) \Big) =0.
\end{align*}
From the above, we also conclude that $\lim_{n \to \infty}(\overline{f} - f_0)(x_{n,\theta})= (\overline{f} - f_0)(x_0)$. By lower semicontinuity of $f_0$ and upper-semicontinuity of $\overline{f}$, therefore we have
\begin{align*}
\overline{f}(x_0) \geq \limsup_{n \to \infty} \overline{f}(x_{n,\theta}) 
 & = \limsup_{n \to \infty} \big((\overline{f} - f_0)(x_{n,\theta}) + f_0(x_{n,\theta})\big) \\
 & \geq (\overline{f} - f_0)(x_0) + f_0(x_0)= \overline{f}(x_0) .
\end{align*}
The above also holds with $\limsup_n$ replaced by $\liminf_n$. Hence we conclude.
\end{proof}

\begin{lemma}\label{Sec3:seq2spw}
Let $\bar{f} \in \USC(\sfX; \bar{\R})$ be a sequential viscosity sub-solution to \eqref{sub} with $h \in C(\sfX)$.
Suppose that for each $(f_0, g_0) \in H$ and every $\theta \in (0, \theta_0)$ for some $\theta_0>0$, there exists $g_\theta \in H f_\theta$, such that in context of Lemma~\ref{Sec3:locF}, for the $\lim_{n \to \infty} x_{n,\theta} \to x_0$, we have 
\begin{align}\label{Sec3:gapp}
\limsup_{\theta \to 0^+} \limsup_{n \to \infty} g_\theta(x_{n,\theta}) \leq g_0(x_0)
\end{align}
 Then $\bar{f}$ is also a point-wise strong viscosity sub-solution.
 
 If, in addition, we assume that Condition~\ref{LypSub} holds, then for each $f_0 \in D(H)$ at least one maximizer $x_0 \in \sfX$ in 
 $(\overline{f} - f_0)(x_0) = \sup_\sfX(\overline{f} - f_0)$ is guaranteed to exist.  
\end{lemma}
\begin{proof}
Let $f_0 \in D(H)$ and let $x_0 \in \sfX$ be such that $(\overline{f} - f_0)(x_0) = \sup_{\sfX} (\overline{f} - f_0)$.  We introduce $f_\theta$ as in the previous lemma. By the defining property of sequential viscosity sub-solution, there exists $x_{n,\theta} \in \sfX$ such that \eqref{Sec3:seqmax} holds and that
\begin{align*}
\limsup_{n \to \infty}\big(  (\overline{f} - h)  - \alpha g_{\theta} \big)(x_{n,\theta})  \leq 0. 
\end{align*}
By Lemma~\ref{Sec3:locF} and \eqref{Sec3:gapp}, then
\begin{align*}
 \big( \overline{f} - h\big)(x_0) = 
\limsup_{\theta \to 0^+}  \limsup_{n \to \infty} (\overline{f} - h)(x_{n,\theta}) 
\leq \alpha g_0(x_0). 
\end{align*}
We conclude that $\overline{f}$ is a strong point-wise viscosity solution.

With Condition~\ref{LypSub}, the existence of maximizer $x_0$ is a consequence of Lemma~\ref{seq2pw}.
\end{proof}
\begin{remark}\label{Sec3:seq2SUP}
Again, in a similarly way, we have corresponding result for the case of super-solution.
\end{remark}

 \subsubsection{Viscosity extension for test functions which are sup- inf- envelop of simpler test functions }\label{Visc-ext}
The notions of viscosity sub- super- solutions are stable under certain (possibly) non-smooth variational extensions of the Hamiltonians.   Lemmas 7.7 and 13.21 in Feng and Kurtz~\cite{FK06} presented one such type of situation. Next, we consider situation of a related but different type. 
 
Let $\Lambda$ be an index set such that $(f_\lambda, g_\lambda) \in H$ 
for every $\lambda \in \Lambda$. For each $x \in \sfX$ fixed, considering  $\lambda \mapsto f_\lambda(x)$ as a function of the $\lambda$, we define 
{\em set of extremal parameters}:
\begin{align}
{\mathcal E}_{\Lambda}^-[f_\cdot(x)]:=\big\{ \alpha \in \Lambda :  f_\alpha(x) = \inf_{\lambda \in \Lambda} f_\lambda(x)   \big\} \subset \Lambda,   
\label{E-def}\\
{\mathcal E}_{\Lambda}^+[f_\cdot(x)]:=\big\{ \alpha \in \Lambda :  f_\alpha(x) = \sup_{\lambda \in \Lambda} f_\lambda(x)  \big\} \subset \Lambda.
\label{E+def}
\end{align}
Note that these sets can be empty in general.
\begin{lemma}\label{parextHJ}
Suppose that $\overline{f} \in M(\sfX; \bar{\R})$ is a strong point-wise sub-solution to \eqref{sub}. 
We define an extension of $H$ by
\begin{align*}
H_0:=  H \cup \big\{ (f,g) :  f  = \inf_{\lambda \in \Lambda} f_\lambda, 
 g(x):= \inf_{\lambda \in {\mathcal E}_{\Lambda}^-[f_{\cdot}(x)]} g_\lambda(x);  \text{ where the } (f_\lambda, g_\lambda) \in  H  \big\}.  
\end{align*}
Then $\overline{f}$ is also a strong pointwise  sub-solution with the operator $H$ replaced by $H_0$.

Suppose that $\underline{f} \in M(\sfX; \bar{\R})$ is a strong point-wise super-solution to 
\eqref{sup}.  We define
\begin{align*}
H_1  &:= H \cup \big\{ (f,g) : f = \sup_{\lambda \in \Lambda}  f_\lambda , 
  g(x):= \sup_{\lambda \in {\mathcal E}_{\Lambda}^+[f_{\cdot}(x)]} g_\lambda(x)
  \big\}.
\end{align*}
Then $\underline{f}$ is also a strong point-wise super-solution with the operator $H$ 
replaced by $H_1$.
\end{lemma}
In the above, we follow the convention that $\inf$ over empty set is considered $+\infty$, 
and $\sup$ over empty set is $-\infty$.
\begin{proof} We only prove the sub-solution property. The super-solution case is similar.

Let $f :=\inf_{\lambda \in \Lambda} f_\lambda$ and $x_0 \in \sfX$ be such that 
$\sup_\sfX (\overline{f} -f ) = (\overline{f} - f)(x_0)$. Next, we verify that 
\begin{align*}
( \overline{f} - h)(x_0) \leq \alpha 
 \inf_{\lambda \in {\mathcal E}_{\Lambda}^-[f_\cdot(x_0)]} g_\lambda(x_0). 
\end{align*}
With the convention that inf over empty set is $+\infty$, 
we only need to prove the case when 
${\mathcal E}_{\Lambda}^-[f_\cdot(x_0)] \neq \emptyset$.  First, we note that
\begin{align*}
\sup_\sfX  (\overline{f} -f ) = \sup_{\lambda \in \Lambda} \sup_\sfX  (\overline{f} -f_\lambda ) 
   \geq \sup_\sfX  (\overline{f} -f_\lambda ), \quad \forall \lambda \in \sfX.
\end{align*}
Second, for each $\lambda_0 \in {\mathcal E}_{\Lambda}^-[f_{\cdot}(x_0)]$ (that is, 
$f_{\lambda_0}(x_0) = f(x_0)$ holds),  the above implies that
\begin{align*}
 (\overline{f} - f_{\lambda_0})(x_0) = (\overline{f} - f)(x_0) = \sup_\sfX(\overline{f} -f ) 
 \geq \sup_\sfX(\overline{f} - f_{\lambda_0}).
\end{align*} 
By the point-wise strong viscosity sub-solution property,
\begin{align*}
( \overline{f} - h)(x_0) \leq \alpha g_{\lambda_0}(x_0). 
\end{align*}
By arbitrariness of the $\lambda_0$ within the set ${\mathcal E}_{\Lambda}^-[f_{\cdot}(x_0)]$, we conclude. 
\end{proof}

\subsection{Projection by submetry of viscosity solutions}\label{projHJ}
In many situations, we are interested in Hamilton-Jacobi equation in a metric space $(\sfX, \sfd_\sfX)$, where the $\sfX$ is base space from a metrically foliated space $(\sfY, \sfd_\sfY)$ (e.g. Definition~\ref{MFoliate}). We discuss such issue in the next two subsections.
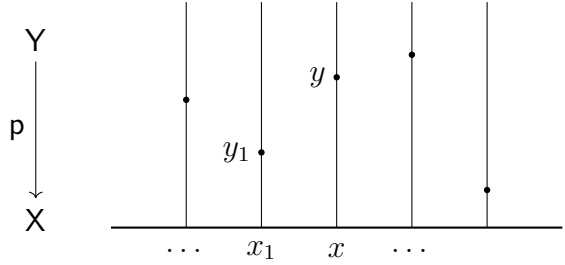
\begin{align*}
\begin{tikzpicture}
\node (Y) at (-4, 2.5) {$\sfY$};
\node (X) at (-4, 0.1) {$\sfX$};
\draw[black, thick] (-3,0) --(3,0);
\draw (-2,3) -- (-2,0);
\filldraw (-2,1.7) circle (1pt);
\draw (-1,3) -- (-1,0);
\filldraw (-1,1.0) circle (1pt) node [anchor=east] {$y_1$};
\draw (0,3) -- (0,0);
\filldraw (0,2) circle (1pt) node [anchor=east] {$y$};
\draw (1,3) -- (1,0);
\filldraw (1,2.3) circle (1pt)   ;
\draw (2,3) -- (2,0);
\filldraw (2,0.5) circle (1pt)    ;
  \node  at  (-2, -0.3) {\ldots};
\node     at (-1,-0.3) {$x_1$};
\node    at (0,-0.3) {$x$};
\node    at (1,-0.3) {\dots};
\draw [->] (Y) -- (X)  node[midway, left] {$\sfp$};
 \end{tikzpicture}
\end{align*}

In such context, the natural projection map $\sfp : \sfY \mapsto \sfX$ is a submetry (see Definition~\ref{subMDef} and Lemma~\ref{subMFol}). Usually, we can write down a Hamilton-Jacobi equation in $\sfY$ which represents physical model defined with finer details. If the physical situation suggests invariance or symmetry along each $\sfp^{-1}(x)$, then we expect a projected equation exists at the reduced level state space $\sfX$. 
Next, we proceed more generally by working with a setup where only approximate versions of invariance or symmetry exist along $\sfp^{-1}(x)$. Thus we are lead to consider perturbed test functions  (e.g. Sections~\ref{PHJPert}).  We also separately discuss the sub- super-solution cases by using possibly different Hamiltonians.
 
Throughout this subsection, let  $( \sfY, \sfd_{\sfY})$ and $(\sfX, \sfd_\sfX)$ be generic metric spaces.
We assume that $\sfp : \sfY \mapsto \sfX$ is a submetry (Definition~\ref{subMDef}). 
We start with a pair of operators in $\sfY$:
\begin{align*}
{\mathcal H}_0 \subset M(\sfY;\bar{\R}) \times M(\sfY; \bar{\R}), \text{ and } 
 {\mathcal H}_1 \subset M(\sfY;\bar{\R}) \times M(\sfY; \bar{\R});
\end{align*}
and consider respectively sub- super-solutions to 
\begin{align}
\bar{\mathfrak f} - \alpha {\mathcal H}_0 \bar{\mathfrak f} & \leq {\mathfrak h}_0, \label{subHY}  \\
\underline{\mathfrak f} - \alpha {\mathcal H}_1 \underline{\mathfrak f} & \geq {\mathfrak h}_1. \label{supHY} 
\end{align}
We are interested in projecting these equations and solutions in $\sfY$ to 
sub- and super-solutions to equations in $\sfX$
\begin{align}
\bar{f} - \alpha H_0 \bar{f} & \leq h_0, \label{H0proj} \\
\underline{f} - \alpha H_1 \underline{f} & \geq h_1, \label{H1proj}
\end{align}
defined with a new set of Hamiltonians 
\begin{align}\label{abprojH0H1}
H_0 \subset M(\sfX;\bar{\R}) \times M(\sfX; \bar{\R}), \text{ and } 
H_1 \subset M(\sfX;\bar{\R}) \times M(\sfX; \bar{\R}).
\end{align}
A natural question is getting sharp estimates motivating definitions of $H_0$ and $H_1$. 

The following notion of projections are useful in our context.
 \begin{definition}[Inf- and sup- projections]\label{LLPdef}
Let $\sfp : \sfY \mapsto \sfX$ be a  submetry and ${\mathfrak f}: \sfY \mapsto \bar{\R}$ . The {\em inf-projection} of the $\mathfrak f$ is a function on $\sfX$ defined by
\begin{align}\label{Linfproj}
 \sfp_{\inf}{\mathfrak f}(x) := \inf\{ {\mathfrak f}(y) : y \in \sfp^{-1}(x) \}.
\end{align}
Similarly, we define the {\em sup}-projection of the $\mathfrak f$ by
\begin{align}\label{Lsupproj}
{\sfp}_{\sup} {\mathfrak f}(x) := \sup \{ {\mathfrak f}(y) : y \in \sfp^{-1}(x) \}.
\end{align}
\end{definition}

Let $\sfp: \sfY \mapsto \sfX$ be a submetry.   For test functions
${\mathfrak f}_0(y):=\beta \sfd_\sfY^2(y,z)$ with $\beta >0, y,z \in \sfY$,
we have 
\begin{align*}
\sfp_{\inf} {\mathfrak f}_0(x) =\beta  \inf_{y \in \sfp^{-1}(x)} \sfd_\sfY^2(y, z)
 = \beta \sfd_\sfY^2(\sfp^{-1}(x), z) = \beta  \sfd_\sfX^2\big(x,\sfp(z)\big),
\end{align*}
where the last equality follows from equi-distant property of metric foliation (see Lemma~\ref{subMFol}). 
Similarly, for ${\mathfrak f}_1(z) := - \beta \sfd_\sfY^2(y,z)$ with $\beta >0, y,z \in \sfY$, we have 
\begin{align*}
\sfp_{\sup} {\mathfrak f}_1(x) = -\beta \sfd_\sfX^2\big(\sfp(y),x\big).
\end{align*}
For general test functions, we have the following regularity result. 
\begin{lemma}
Suppose that ${\mathfrak f} \in \LSC(\sfY; \R \cup \{+\infty\})$ and every finite sub-level set is compact, then $f:=\sfp_{\inf} {\mathfrak f} \in \LSC(\sfX; \R \cup\{+\infty\})$.
Similarly, if ${\mathfrak f} \in \USC(\sfY; \R \cup \{-\infty\})$ with $-{\mathfrak f}$ having compact finite sub-levels, then $f:=\sfp_{\sup} {\mathfrak f} \in \USC(\sfX; \R \cup\{-\infty\})$.  
\end{lemma}
\begin{proof} We prove the lower semicontinuity case. 
Let $x_n, x_0 \in \sfX$ be $\sfd_\sfX(x_n, x_0 ) \to 0$. We only need to show that $\liminf_{n \to \infty} f(x_n) \geq f(x_0)$ in the special case when left hand side is finite. 

By the definition in \eqref{Linfproj}, there exists $y_n \in \sfp^{-1}(x_n)$ with $f(x_n) \geq {\mathfrak f}(y_n) -\frac1n$. By the compact finite sub-level assumption, up to selection of subsequence, there exists a $y_0 \in \sfY$ with $\lim_{k \to \infty}\sfd_\sfY(y_{n(k)}, y_0) = 0$. By the equi-distant property of metric foliation induced by the $\sfp$ (Lemma~\ref{subMFol}), we also have $\lim_n \sfd_\sfY(y_n, \sfp^{-1}(x_0)) =0$.
Since the leaf $\sfp^{-1}(x_0)$ is closed, we have $y_0 \in \sfp^{-1}(x_0)$, giving 
\begin{align*}
\liminf_{n \to \infty} f(x_n) \geq \liminf_{n \to \infty} {\mathfrak f}(y_n)  \geq {\mathfrak f}(y_0) \geq f(x_0).
\end{align*}
\end{proof}

\subsubsection{Projected Hamiltonian operators and viscosity solutions -  Multi-valued operators}
We recall the convention that $\inf$ over an empty set is $+\infty$, and $\sup$ over an empty set is $-\infty$.

In the context of $\sfX= \sfp(\sfY)$ where the $\sfp$ is a submetry,
 we define another type of extremal sets similar to the use of \eqref{E-def} and \eqref{E+def} in Subsection~\ref{Visc-ext}: 
 for every ${\mathfrak f} : \sfY \mapsto \bar{\R}$, $x \in \sfX$ and $\epsilon\geq 0$, we define
 \begin{align}
 E_\epsilon^-[{\mathfrak f};x] &:= 
      \big\{ y \in \sfp^{-1}(x) : {\mathfrak f}(y) -\epsilon 
      \leq \inf_{\sfp^{-1}(x)} {\mathfrak f} \big\}, 
       \label{defEe-}\\
 E_\epsilon^+[{\mathfrak f};x] &:=
           \big\{ y \in \sfp^{-1}(x) : {\mathfrak f}(y) +\epsilon
            \geq \sup_{\sfp^{-1}(x)} {\mathfrak f} \big\}. 
              \label{defEe+}
\end{align}
If $\mathfrak f \in \LSC(\sfY; \bar{\R})$, then $E_\epsilon^-[{\mathfrak f};x]$ is a closed set. 
Similarly, if $\mathfrak f \in \USC(\sfY; \bar{\R})$, then $E_\epsilon^+[{\mathfrak f};x]$ is a closed set.
We define two multi-valued operators in $\sfX$, through their graphs, by
\begin{align}
H_0    &:=\big\{ (f,g) :  f= \sfp_{\inf} {\mathfrak f}, 
g= g_\epsilon(x)= \sup_{y \in E_\epsilon^-[\mathfrak f; x]} {\mathfrak g}(y), 
  \forall ({\mathfrak f},{\mathfrak g}) \in  {\mathcal H}_0, \forall \epsilon>0  \big\}, \label{H0insuDef} \\
H_1  &:= \big\{ (f,g) : f = \sfp_{\sup}{\mathfrak f} , 
g= g_\epsilon(x)= \inf_{y \in E_\epsilon^+[\mathfrak f; x]}{\mathfrak g}(y),
\forall {(\mathfrak f}, {\mathfrak g}) \in {\mathcal H}_1, \forall \epsilon >0 \big\}.\label{H1suinDef}
\end{align}
 
 \begin{lemma}\label{sprojHJ}
Let ${\mathfrak h}_i \in M(\sfY)$ be bounded for $i=0,1$. Suppose that 
\begin{enumerate}
\item $\bar{\mathfrak f} \in M(\sfY; \bar{\R})$ is a sequential viscosity sub-solution to \eqref{subHY},
\item $\bar{\mathfrak f}$ is invariant over $\sfp^{-1}(x)$ for every $x \in \sfX$ fixed. That is, 
\begin{align}
\bar{\mathfrak f}(y) = \text{constant}=: \bar{f}(x), \quad \forall y \in \sfp^{-1}(x). \label{finv0}
\end{align} 
\end{enumerate}
Then $H_0 \subset M(\sfX; \bar{\R}) \times M(\sfX; \bar{\R})$ and
the $\bar{f} \in M(\sfX; \bar{\R})$ in \eqref{finv0} is a sequential viscosity sub-solution to \eqref{H0proj} with the $h_0 := \sfp_{\sup} {\mathfrak h}_0$. 
Moreover, if the $\bar{\mathfrak f} \in \USC(\sfY; \bar{\R})$ and the $\sfp$ is a strong submetry (see Appendix~\ref{App:SubM}), then $\overline{f} \in \USC(\sfX; \bar{\R})$.

Similarly, suppose that 
\begin{enumerate}
\item $\underline{\mathfrak f} \in M(\sfY; \bar{\R})$ is a sequential viscosity super-solution to \eqref{supHY}.
\item $\underline{\mathfrak f}$ is invariant over $\sfp^{-1}(x)$ for every $x \in \sfX$ fixed:
\begin{align*}
 \underline{\mathfrak f}(y) = \text{constant} =: \underline{f}(x), \quad \forall y \in \sfp^{-1}(x).
\end{align*}
\end{enumerate}
Then  $H_1 \subset M(\sfX; \bar{\R}) \times M(\sfX; \bar{\R})$ and 
the $\underline{f} \in M(\sfX; \bar{\R})$ is a sequential viscosity super-solution to \eqref{H1proj} with the $h_1 := \sfp_{\inf} {\mathfrak h}_1$. Moreover, if the $\underline{\mathfrak f} \in \LSC(\sfY; \bar{\R})$ and the $\sfp$ is a strong submetry, then $\underline{f} \in \LSC(\sfX; \bar{\R})$.
\end{lemma}
\begin{proof} 
We only prove the sub-solution case. The super-solution case is similar.
   
Let $\epsilon>0$ and $(f,g) \in H_0$ be such that $f= \sfp_{\inf} {\mathfrak f}$ and $g(x)= g_\epsilon(x) = \sup_{y \in E^-_\epsilon[\mathfrak f;x]} {\mathfrak g}(y)$ with the $({\mathfrak f}, {\mathfrak g}) \in {\mathcal H}_0$. 
Since the $\bar{{\mathfrak f}}$  is a sequential viscosity sub-solution to \eqref{subHY}, 
there exists $y_n:= y_n^{\mathfrak f, \mathfrak g} \in \sfY$ such that 
\begin{align}
 & \lim_{n \to \infty} (\overline{\mathfrak f} - {\mathfrak f})(y_n) 
     = \sup_{\sfY}(\overline{\mathfrak f} - {\mathfrak f}), \label{Sec3:supyf} \\
 & \limsup_{n \to \infty} (\overline{\mathfrak f} - {\mathfrak h}_0 - \alpha {\mathfrak g})(y_n)   \leq 0. \label{Sec3:Subeqn}
\end{align}
We define $x_n:=x_n^{\mathfrak f, \mathfrak g}:= \sfp(y_n) \in \sfX$.  Then 
\begin{align*}
  \sup_{\sfY} (\bar{\mathfrak f} - {\mathfrak f}) &
  = \sup_{x \in \sfX} \sup_{y \in \sfp^{-1}(x)}  \big(\bar{\mathfrak f}(y) - {\mathfrak f}(y)\big) 
  = \sup_{x \in \sfX}  \big(\bar{f}(x)- (\sfp_{\inf}{\mathfrak f})(x)\big) 
   =  \sup_{\sfX}   (\bar{f} - f),
\end{align*}
and 
\begin{align*}
 \sup_{\sfY} (\bar{\mathfrak f} - {\mathfrak f}) 
=\lim_{n \to \infty}  (\bar{\mathfrak f} - {\mathfrak f}) (y_n) 
\leq \lim_{n \to \infty}   (\bar{f} - f) (x_n),
\end{align*}
implying $\lim_{n \to \infty}(\bar{f} - f) (x_n)= \sup_{\sfX}(\bar{f} - f)$.
Moreover, because of \eqref{Sec3:supyf}, for the $\epsilon>0$, there exists $N:=N(\epsilon)$ large enough so that, for $n>N$, we have  
\begin{align*}
 \bar{f}(x_n) - {\mathfrak f}(y_n) = (\bar{\mathfrak f} - {\mathfrak f})(y_n)
 &  \geq \sup_\sfY(\bar{\mathfrak f} - {\mathfrak f})  - \epsilon \\
&  \geq  \sup_{\sfp^{-1}(x_n)}(\bar{\mathfrak f} - {\mathfrak f})  - \epsilon  
   =   \bar{f}(x_n) - \inf_{\sfp^{-1}(x_n)}{\mathfrak f} - \epsilon.
\end{align*}
Therefore $y_n \in E^-_\epsilon[\mathfrak f;x_n]$. Consequently, \eqref{Sec3:Subeqn} gives
\begin{align*}
 \limsup_{n \to \infty} \big(\overline{f}(x_n) - h_0(x_n)  - \alpha  g_\epsilon(x_n)\big) \leq  
 \limsup_{n \to \infty} (\overline{\mathfrak f} - {\mathfrak h}_0 - \alpha  {\mathfrak g})(y_n) \leq 0.
\end{align*}

Now we put the additional assumptions that the $\bar{\mathfrak f} \in \USC(\sfY; \bar{\R})$ and that the $\sfp: \sfY \mapsto \sfX$ is a {\em strong} submetry. Let $x_n , x_0 \in \sfX$ be $\sfd_\sfX(x_n, x_0) \to 0$. Take a $y_0 \in \sfp^{-1}(x_0)$, by the 2-point lifting property (Lemma~\ref{2pts}), there exists $y_n \in \sfY$ with $\sfd_\sfY(y_n, y_0) = \sfd_\sfX(x_n, x_0)$. Hence
\begin{align*}
\limsup_{n \to \infty} \overline{f}(x_n) = \limsup_{n \to \infty} \overline{\mathfrak f}(y_n) \leq \overline{\mathfrak f}(y_0) = \overline{f}(x_0).
\end{align*}
\end{proof}

One can view $H_i$ as projected operators from the ${\mathcal H}_i$s, $i=0,1$. If we strengthen some assumptions in Lemma~\ref{sprojHJ}, we will arrive at simpler and stronger versions of the above projective type result. In the following, among all possibilities, we only present two versions of such extensions.  

\subsubsection{Projected Hamiltonian operators and viscosity solutions -  Single-valued operators I}
 We introduce two single-valued Hamiltonian operators $H_0$ and $H_1$ as follow.
 \footnote{For notational simplicity, we slightly abuse notation by using the same notations $H_0, H_1$ as in the previous sub-section.}
  First, we define domains of the operators 
 \begin{align}\label{Domh01}
D(H_0) := \big\{  f= \sfp_{\inf} {\mathfrak f} : \exists {\mathfrak f} \in D({\mathcal H}_0)  \big\}, \quad
D(H_1) := \big\{   f =\sfp_{\sup} {\mathfrak f} : \exists {\mathfrak f} \in D({\mathcal H}_1) \big\}. 
\end{align}
Second, we define the operators
\begin{align}
H_0 f (x)  &:= \inf_{ \substack{ ({\mathfrak f},{\mathfrak g}) \in {\mathcal H}_0 \\ \text{ with } {\mathfrak f} \in (\sfp_{\inf})^{-1}(f)}}     
   \lim_{\substack{\epsilon \to 0\\ \epsilon>0}} 
        \sup_{y  \in E^-_\epsilon[\mathfrak f; x]}   {\mathfrak g}(y), 
 \quad \forall f \in D(H_0), \label{H0insuA} \\
H_1 f(x)  &:= \sup_{ \substack{ ({\mathfrak f},{\mathfrak g}) \in {\mathcal H}_1 \\\text{ with } {\mathfrak f} \in (\sfp_{\sup})^{-1}(f)}}
     \lim_{\substack{\epsilon \to 0 \\ \epsilon>0}}
        \inf_{y \in E^+_\epsilon[\mathfrak f; x]}  {\mathfrak g}(y),
\quad  \forall f  \in D(H_1). \label{H1suinA}
\end{align}
Noting 
\begin{align}\label{Sec3:Eincl}
E^{\pm}_{\epsilon^\prime}[\mathfrak f; x] \subset E^{\pm}_\epsilon[\mathfrak f; x],
 \quad \forall  0\leq \epsilon^\prime < \epsilon,
\end{align}
we have 
\begin{align*}
 \lim_{\substack{\epsilon \to 0 \\ \epsilon>0}} \sup_{E^-_\epsilon[\mathfrak f; x]}   {\mathfrak g}
 = \inf_{\epsilon > 0} \sup_{E^-_\epsilon[\mathfrak f; x]}   {\mathfrak g},
       \qquad  
\lim_{\substack{\epsilon \to 0 \\ \epsilon >0}} \inf_{E^+_\epsilon[\mathfrak f; x]}   {\mathfrak g}
  = \sup_{\epsilon > 0} \inf_{E^+_\epsilon[\mathfrak f; x]}   {\mathfrak g}.
\end{align*}

\begin{lemma}\label{spHJ1}
Let ${\mathfrak h}_i \in M(\sfY)$ be bounded for $i=0,1$.
Suppose that 
\begin{enumerate}
\item $\bar{\mathfrak f} \in M(\sfY; \bar{\R})$ is a sequential  {\bf strong} 
viscosity sub-solution to \eqref{subHY}.
\item $\bar{\mathfrak f}$ is invariant over $\sfp^{-1}(x)$ for every $x \in \sfX$ fixed:  
\begin{align}\label{fbarC}
\bar{\mathfrak f}(y) = \text{constant}=: \bar{f}(x), \quad \forall y \in \sfp^{-1}(x).
\end{align} 
\end{enumerate}
Then 
$H_0 :   D(H_0) \mapsto M(\sfX; \bar{\R})$,  and
the $\bar{f} \in M(\sfX; \bar{\R})$ is a   point-wise {\bf strong} viscosity sub-solution to \eqref{H0proj}
with $h_0 := \sfp_{\sup} {\mathfrak h}_0$. 

Similarly, suppose that 
\begin{enumerate}
\item $\underline{\mathfrak f} \in M(\sfY; \bar{\R})$ is a sequential  {\bf strong} 
viscosity super-solution to \eqref{supHY}.
\item $\underline{\mathfrak f}$ is invariant over $\sfp^{-1}(x)$ for every $x \in \sfX$ fixed:
\begin{align*}
 \underline{\mathfrak f}(y) = \text{constant} =: \underline{f}(x), \quad \forall y \in \sfp^{-1}(x).
\end{align*}
\end{enumerate}
Then  $H_1 :    D(H_1) \mapsto M(\sfX; \bar{\R})$ and 
the $\underline{f} \in M(\sfX; \bar{\R})$ is a {\bf strong}  point-wise viscosity super-solution to 
\eqref{H1proj} with $h_1 := \sfp_{\inf} {\mathfrak h}_1$. 
\end{lemma}
\begin{proof} 
Again, we only prove the sub-solution case. 

Let $f \in D(H_0)$ and $x_0 \in \sfX$ be such that $\sup_{\sfX}(\bar{f} - f) < \infty$ and that 
$(\bar{f} -  f)(x_0) = \sup_\sfX(\bar{f} - f)$.
From the definition of $H_0f(x_0)$ in \eqref{H0insuA}, we can select 
$({\mathfrak f}_n, {\mathfrak g}_n) \in {\mathcal H}_0$ with
$\sfp_{\inf} {\mathfrak f}_n =f$ and $\epsilon_n >0$, $\lim_{n \to \infty} \epsilon_n = 0$, such that
\begin{align}\label{appgfib}
  \sup_{E^-_{\epsilon_n}[{\mathfrak f}_n; x_0]}  {\mathfrak g}_n  \leq H_0 f(x_0) + \frac1n .
\end{align}

Next, we select 
$\{ y_{n,m} \}_{m=1,2,\ldots} \subset \sfp^{-1}(x_0)$ such that
\begin{align}\label{frakfsup}
 \lim_{m \to \infty} (\bar{\mathfrak f} - {\mathfrak f}_n)(y_{n,m}) = \sup_{\sfp^{-1}(x_0)}(\bar{\mathfrak f} - {\mathfrak f}_n).
\end{align}
Because that $\bar{\mathfrak f}$ is constant on $\sfp^{-1}(x)$ for each $x$ fixed, we have that for every $n$,
 \begin{align*}
  \sup_{\sfY} (\bar{\mathfrak f} - {\mathfrak f}_n) 
  &= \sup_{x \in \sfX} \sup_{y \in \sfp^{-1}(x)}  \big(\bar{\mathfrak f}(y)- {\mathfrak f}_n(y)\big)  =   \sup_{\sfX}   (\bar{f} - f) \\
&  =    (\bar{f} - f) (x_0) = \sup_{\sfp^{-1}(x_0)}(\bar{\mathfrak f} - {\mathfrak f}_n) 
= \lim_{m \to \infty} (\bar{\mathfrak f} - {\mathfrak f}_n) (y_{n,m}) .
\end{align*}
Since $\bar{\mathfrak f}$ is a sequential strong viscosity sub-solution to \eqref{subHY}, 
\begin{align}\label{YseqSub}
\limsup_{m \to \infty} \big( \bar{\mathfrak f}  - {\mathfrak h}_0 
  -  \alpha  {\mathfrak g}_n \big)(y_{n,m}) \leq 0.
\end{align}

From \eqref{frakfsup} and the assumption in \eqref{fbarC} (implying $\bar{\mathfrak f}(y_{n,m})=\bar{f}(x_0)$), 
we arrive at $y_{n,m} \in E^-_{\epsilon_n}[\mathfrak f_n; x_0]$ for 
$m$ sufficiently large while $n$ is fixed. Consequently,
\begin{align*}
   \big( \bar{f}(x_0) -  h_0(x_0) -  \alpha \sup_{E^-_{\epsilon_n}[\mathfrak f_n; x_0]} {\mathfrak g}_n \big)  
   \leq   \limsup_{m \to \infty} \big( \bar{\mathfrak f}  - {\mathfrak h}_0 
  -  \alpha  {\mathfrak g}_n \big)(y_{n,m}) \leq 0.
\end{align*}
In view of \eqref{appgfib}, 
\begin{align*}
 \bar{f}(x_0) -  h_0(x_0) \leq  \alpha H_0 f(x_0).
\end{align*} 
\end{proof}
Note again that, the above lemma and its proof does not guarantee the existence of the $x_0$. This is because that, in the definition of strong point-wise sub- super- solution, we only required the defining inequalities to hold when such extremal point $x_0$ exists. Therefore, when applying the above results, we need to explicitly construct the $x_0$ first.

\begin{lemma}\label{spHJ2}
Assume that ${\mathcal H}_0 \subset  \LSC(\sfX; \bar{\R}) \times \USC(\sfX; \bar{\R})$ 
and ${\mathcal H}_1 \subset \USC(\sfX; \bar{\R}) \times \LSC(\sfX; \bar{\R})$.
 
In addition, assuming each $\mathfrak f \in D(\mathcal H_0)$ has compact finite sub-levels on $\sfp^{-1}(x)$, then the $H_0$ in \eqref{H0insuA} admit simpler representation
\begin{align*}
H_0f(x)= \inf_{ \substack{ ({\mathfrak f},{\mathfrak g}) \in {\mathcal H}_0 \\{\mathfrak f} \in (\sfp_{\inf})^{-1}(f)}}  \sup_{E^-_0[\mathfrak f; x]}   {\mathfrak g},
\end{align*} 
Similarly, assuming additionally that, for each ${\mathfrak f} \in D(\mathcal H_1)$, $- \mathfrak f$ has compact finite sub-levels on $\sfp^{-1}(x)$, then the $H_1$ in \eqref{H1suinA} admit simpler representation
\begin{align*}
H_1 f(x)=  \sup_{ \substack{ ({\mathfrak f},{\mathfrak g}) \in {\mathcal H}_1 \\
  {\mathfrak f} \in (\sfp_{\sup})^{-1}(f)}} \inf_{ E^+_0[\mathfrak f;x]} \mathfrak g.
\end{align*}
\end{lemma}
\begin{proof}
 We note that, by inclusions $E^{\pm}_0[\mathfrak f;x]
 \subset E^{\pm}_\epsilon[\mathfrak f;x]$ (see \eqref{Sec3:Eincl}), it always holds that 
 \begin{align*}
 \lim_{\substack{\epsilon \to 0 \\ \epsilon>0}} \sup_{E^-_\epsilon[\mathfrak f; x]}   {\mathfrak g}
 = \inf_{\epsilon > 0} \sup_{E^-_\epsilon[\mathfrak f; x]}   {\mathfrak g} 
  \geq \sup_{E^-_0[\mathfrak f; x]}   {\mathfrak g},
       \qquad  
\lim_{\substack{\epsilon \to 0 \\ \epsilon >0}} \inf_{E^+_\epsilon[\mathfrak f; x]}   {\mathfrak g}
  = \sup_{\epsilon > 0} \inf_{E^+_\epsilon[\mathfrak f; x]}   {\mathfrak g}
   \leq \inf_{E^+_0[\mathfrak f; x]}   {\mathfrak g}.
\end{align*}
Hence we only need to verify the opposite sides of the inequalties. We only show the case for the first one.  For $\epsilon>0$, there exists $y_\epsilon \in E_\epsilon^-[\mathfrak f;x] \subset \sfp^{-1}(x)$ with 
\begin{align*}
 \sup_{E_\epsilon^-[\mathfrak f;x]} {\mathfrak g} < \epsilon + {\mathfrak g}(y_\epsilon).
\end{align*}
By the compact finite sub-levels assumption, we can find a subsequence so that $y_\epsilon \to y_0 \in \sfp^{-1}[x]$. By the lower semi-continuity of $\mathfrak f$, we have $y_0 \in E_0^-[\mathfrak f;x]$. By upper-semicontinuity of the $\mathfrak g$, 
\begin{align*}
 \inf_{\epsilon>0}\sup_{E_\epsilon^-[\mathfrak f;x]} {\mathfrak g} \leq \liminf_{\epsilon \to 0^+} {\mathfrak g}(y_\epsilon) \leq {\mathfrak g}(y_0) \leq \sup_{E^-_0[\mathfrak f; x]}   {\mathfrak g}.
\end{align*}
We conclude.
\end{proof}

\subsubsection{Projected Hamiltonian operators and viscosity solutions -  Single-valued operators II}
In fact, if we assume the conditions in Lemma~\ref{spHJ2}, we can have a much stronger result. Specifically,  compared to the definitions of $H_0, H_1$ in \eqref{H0insuA} and \eqref{H1suinA},  we can replace the 
 $\lim_{\epsilon \to 0^+} \sup_{  E^-_\epsilon[\mathfrak f; x]}$ in $H_0$ by $\inf_{E^-_0[\mathfrak f; x]}$, 
and $\lim_{\epsilon \to 0^+}\inf_{E^-_\epsilon[\mathfrak f; x]}$ in $H_1$ by $\sup_{E^-_0[\mathfrak f; x]}$.  Note that, regardless of any condition, by \eqref{Sec3:Eincl}, we always have
\begin{align*}
  \lim_{\substack{\epsilon \to 0 \\ \epsilon>0}} \inf_{E^-_\epsilon[\mathfrak f; x]}   {\mathfrak g}
 = \sup_{\epsilon > 0} \inf_{E^-_\epsilon[\mathfrak f; x]}{\mathfrak g} 
     \leq \inf_{E^-_0[\mathfrak f; x]}   {\mathfrak g},
       \qquad  
\lim_{\substack{\epsilon \to 0 \\ \epsilon >0}} \sup_{E^+_\epsilon[\mathfrak f; x]}   {\mathfrak g}
  = \inf_{\epsilon > 0} \sup_{E^+_\epsilon[\mathfrak f; x]}   {\mathfrak g}
     \geq \sup_{E^+_0[\mathfrak f; x]}   {\mathfrak g}.
\end{align*}
 
We define yet another set of $H_0, H_1$ operators by
\begin{align}
H_0 f (x)  &:= \inf_{ \substack{ ({\mathfrak f},{\mathfrak g}) \in {\mathcal H}_0 \\{\mathfrak f} \in (\sfp_{\inf})^{-1}(f)}}  \inf_{y  \in E^-_0[\mathfrak f; x]}   {\mathfrak g}(y), 
 \quad \forall f \in D(H_0), \label{H0infinf} \\
H_1 f(x)  &:= \sup_{ \substack{ ({\mathfrak f},{\mathfrak g}) \in {\mathcal H}_1 \\{\mathfrak f} \in (\sfp_{\sup})^{-1}(f)}} 
\sup_{y \in E^+_0[\mathfrak f; x]}  {\mathfrak g}(y),
\quad  \forall f  \in D(H_1). \label{H1supsup}
\end{align}
where the domains $D(H_0), D(H_1)$ are still defined as in \eqref{Domh01}.
 
 \begin{lemma}\label{spHJ3}
Assume that ${\mathcal H}_0 , {\mathcal H}_1 \subset M(\sfX; \bar{\R}) \times M(\sfX; \bar{\R})$, and that ${\mathfrak h}_i \in M(\sfY)$ is bounded, for $i=0,1$.

Suppose that 
\begin{enumerate}
\item $\bar{\mathfrak f} \in M(\sfY; \bar{\R})$ is a {\bf point-wise strong} 
viscosity sub-solution to \eqref{subHY}.
\item $\bar{\mathfrak f}$ is invariant over $\sfp^{-1}(x)$ for every $x \in \sfX$ fixed:  
\begin{align*}
\bar{\mathfrak f}(y) = \text{constant}=: \bar{f}(x), \quad \forall y \in \sfp^{-1}(x).
\end{align*} 
\item each $\mathfrak f \in D(\mathcal H_0)$ has compact finite-sub-levels over $\sfp^{-1}(x)$.
\end{enumerate}
Then the $H_0$ defined in \eqref{H0infinf} makes the $\bar{f} \in M(\sfX; \bar{\R})$ a point-wise {\bf strong} viscosity sub-solution to \eqref{H0proj} with $h_0 := \sfp_{\sup} {\mathfrak h}_0$. 

Similarly, suppose that 
\begin{enumerate}
\item $\underline{\mathfrak f} \in M(\sfY; \bar{\R})$ is a {\bf point-wise strong} 
viscosity super-solution to \eqref{supHY}.
\item $\underline{\mathfrak f}$ is invariant over $\sfp^{-1}(x)$ for every $x \in \sfX$ fixed:
\begin{align*}
 \underline{\mathfrak f}(y) = \text{constant} =: \underline{f}(x), \quad \forall y \in \sfp^{-1}(x).
\end{align*}
\item for each ${\mathfrak f} \in D(\mathcal H_1)$, $- \mathfrak f$ has compact finite-sub-levels 
over $\sfp^{-1}(x)$.
\end{enumerate}
Then  the $H_1$ defined in \eqref{H1supsup} makes the above $\underline{f} \in M(\sfX; \bar{\R})$ is  {\bf strong}  point-wise viscosity super-solution to \eqref{H1proj} with $h_1 := \sfp_{\inf} {\mathfrak h}_1$. 
\end{lemma}
\begin{proof} We only prove the sub-solution case, which is similar to that of Lemma~\ref{spHJ1}.

Let $f \in D(H_0)$ and $x_0 \in \sfX$ be such that $\sup_{\sfX}(\bar{f} - f) < \infty$ and that 
$(\bar{f} -  f)(x_0) = \sup_\sfX(\bar{f} - f)$.
By the definition of $H_0f(x_0)$ in \eqref{H0infinf}, we can select a sequence
$({\mathfrak f}_n, {\mathfrak g}_n) \in {\mathcal H}_0$ with
$\sfp_{\inf} {\mathfrak f}_n =f$, such that
\begin{align}\label{appgf2}
  \inf_{E^-_0[{\mathfrak f}_n; x_0]}  {\mathfrak g}_n  \leq H_0 f(x_0) + \frac1n .
\end{align}
Without lose of generality, we assume $H_0f(x_0)<+\infty$. This implies that $E^-_0[{\mathfrak f}_n; x_0] \neq \emptyset$, hence there exists $y_n \in E^-_0[{\mathfrak f}_n;x_0]$ with  
\begin{align*}
 {\mathfrak g}_n(y_n) \leq \inf_{E^-_0[{\mathfrak f}_n;x_0]} {\mathfrak g}_n + \frac1n .
\end{align*}

Because that $\bar{\mathfrak f}$ is 
constant over $\sfp^{-1}(x)$ for each $x$ fixed, we have that
\begin{align*}
  (\bar{\mathfrak f} - {\mathfrak f}_n)(y_n) = \sup_{\sfp^{-1}(x_0)}(\bar{\mathfrak f} - {\mathfrak f}_n),
\end{align*}
and that 
 \begin{align*}
  \sup_{\sfY} (\bar{\mathfrak f} - {\mathfrak f}_n) 
  &= \sup_{x \in \sfX} \sup_{y \in \sfp^{-1}(x)}  \big(\bar{\mathfrak f}(y)- {\mathfrak f}_n(y)\big)  =   \sup_{\sfX}   (\bar{f} - f) \\
&  =    (\bar{f} - f) (x_0) = \sup_{\sfp^{-1}(x_0)}(\bar{\mathfrak f} - {\mathfrak f}_n) 
= (\bar{\mathfrak f} - {\mathfrak f}_n) (y_n) .
\end{align*}
Since $\bar{\mathfrak f}$ is a point-wise strong viscosity sub-solution to \eqref{subHY}, the above implies that
\begin{align*}
 \big( \bar{\mathfrak f}  - {\mathfrak h}_0 \big) (y_n)
 \leq  \alpha  {\mathfrak g}_n (y_n).
\end{align*}
Therefore
\begin{align*}
 \alpha^{-1}\big( \bar{f} -  h_0\big) (x_0) \leq 
 \alpha^{-1} \big( \bar{\mathfrak f} - {\mathfrak h}_0\big) (y_n)
 \leq    {\mathfrak g}_n(y_n) \leq
   \inf_{E_0^-[\mathfrak f_n;x_0]} {\mathfrak g}_n + \frac1n
 \leq H_0 f(x_0) + \frac2n. 
 \end{align*}
We conclude by letting $n \to \infty$.
\end{proof}

\subsection{Projection of Hamiltonians defined with special test functions}\label{HJproj}
We continue by assuming $(\sfX, \sfd_\sfX)$ and $(\sfY, \sfd_\sfY)$ are metric spaces, and $\sfp: \sfY \mapsto \sfX$ a submetry map (see Appendix~\ref{App:MQ}). Throughout this subsection, we will make extensive use of simple smooth functions  
 ${\mathcal S}_\sfX^{\pm}$ and ${\mathcal S}_\sfY^{\pm}$ as defined in \eqref{SS+} and \eqref{SS-}. Here, subscripts $\sfX, \sfY$ are added to emphasize the metric space dependencies. 
 
 \subsubsection{Composition of distance functions as test functions}
 We consider operator  ${\mathcal H}$ whose domain $D({\mathcal H})$ is a subset of functions satisfying in particular ${\mathcal S}_\sfY^+ \cup  {\mathcal S}_\sfY^- \subset D({\mathcal H})$. We write 
\begin{align}
 {\mathfrak f}_{0; y_1, \ldots, y_K}(y)
 & :=\psi \big(\sfd_\sfY^2(y, y_1), \ldots, \sfd_\sfY^2(y, y_K)\big) 
\in {\mathcal S}^+_\sfY, \label{Adeffrakf} \\
  f_{0; x_1, \ldots, x_K}(x) 
  & := \psi \big(\sfd_\sfX^2(x, x_1), \ldots, \sfd_\sfX^2(x, x_K)\big) \in {\mathcal S}^+_\sfX, \label{Alexdeff0}
\end{align}
where $K \in \N, \psi \in {\bf \Psi}_K$, $y_1,\ldots,y_K \in \sfY$ and $x_1, \ldots, x_K \in \sfX$ (see \eqref{defPsiK} for definition of ${\bf \Psi}_K$).   Note that, for separable metric spaces, the collection of distance functions introduces a nice system of local coordinates. 

We assume that $\sfp$ is a {\em strong} submetry map (Definition~\ref{subMDef}). 
For every $y_0 \in \sfY$ and $x_0 = \sfp(y_0)$, by the $2$-point lifting property (Lemma~\ref{2pts} in Appendix),  there exists $y_k \in \sfp^{-1}(x_k)$ for $k=1,\ldots,K$ such that 
\begin{align*}
 \sfd_\sfY(y_k, y_0) = \sfd_\sfX (x_k, x_0),  \quad k=1,\ldots,K.
 \end{align*}
Therefore, the following definition of a section is non-empty at least in that case
\begin{align}\label{sfSdef}
\sfS(x_0; x_1,\ldots, x_K)& :=\Big\{(y_0,y_1,\ldots, y_K) \in \Pi_{k=0}^K \sfp^{-1}(x_k):   \\
&\qquad \qquad \text{ such that } 
  \sfd_\sfY(y_k, y_0) = \sfd_\sfX (x_k, x_0),  k=1,\ldots,K  \Big\}. \nonumber
\end{align}
See the following graph.
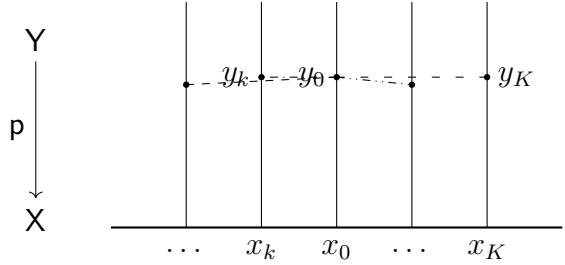
\begin{align*}
\begin{tikzpicture}
\node (Y) at (-4, 2.5) {$\sfY$};
\node (X) at (-4, 0.1) {$\sfX$};
\draw[black, thick] (-3,0) --(3,0);
\draw (-2,3) -- (-2,0);
\filldraw (-2,1.9) circle (1pt);
\draw [dashed] (0,2) --(-2,1.9);
\draw (-1,3) -- (-1,0);
\filldraw (-1,2) circle (1pt) node [anchor=east] {$y_k$};
\draw (0,3) -- (0,0);
\filldraw (0,2) circle (1pt) node [anchor=east] {$y_0$};
\draw [dash dot] (0,2) -- (-1,2);
\draw (1,3) -- (1,0);
\filldraw (1,1.9) circle (1pt)   ;
\draw [dash dot dot] (0,2) -- (1, 1.9);
\draw (2,3) -- (2,0);
\filldraw (2,2) circle (1pt)  node [anchor=west] {$y_K$};
\draw [loosely dashed] (0,2) --(2,2);
\node  at  (-2, -0.3) {\ldots};
\node     at (-1,-0.3) {$x_k$};
\node    at (0,-0.3) {$x_0$};
\node    at (1,-0.3) {\ldots};
\node    at (2,-0.3) {$x_K$};
\draw [->] (Y) -- (X)  node[midway, left] {$\sfp$};
 \end{tikzpicture}
\end{align*}
When the $y_0 \in \sfp^{-1}(x_0)$ is held fixed, we also write  
\begin{align}\label{Sdef}
S_{x_1,\ldots,x_K}(y_0) & := \Big\{(y_1,\ldots, y_K) \in \Pi_{k=1}^K \sfp^{-1}(x_k)
:  \text{ such that }   \\
& \qquad \qquad 
  \sfd_\sfY(y_k, y_0) = \sfd_\sfX (x_k, x_0), \text{ where } x_0 =p(y_0), k=1,\ldots,K  \Big\} .   \nonumber
\end{align}

 \begin{lemma}
We define, for $x_k \in \sfX$,
\begin{align}\label{Sec3:MPproj}
  F_{x_1,\ldots, x_K}(y):= \inf_{\substack{y_1 \in \sfp^{-1}(x_1)\\ \ldots\\ y_K \in \sfp^{-1}(x_K) } } {\mathfrak f}_{0;y_1,\ldots, y_K}(y), \quad \forall y \in \sfY.
\end{align}
With reference to the context and notations of Lemma~\ref{parextHJ}, we introduce index set
\begin{align}\label{LamDef}
\Lambda:= \Lambda(x_1, \ldots, x_K):= \Pi_{k=1}^K \sfp^{-1}(x_k) \subset \sfY \times \ldots \times \sfY,
\end{align}
then
\begin{align}\label{FDEF}
  F_{x_1,\ldots, x_K}(y) = \inf_{\substack{(y_1,\ldots,y_K) \\ \in \Lambda(x_1,\ldots,x_K)}} {\mathfrak f}_{0;y_1,\ldots, y_K}(y).
\end{align}
Assume that $\sfp$ is a strong submetry map. Then the following holds:
\begin{enumerate}
\item  for every $(y_1,\ldots, y_K)   \in S_{x_1,\ldots,x_K}(y)$ and 
$x = \sfp(y)$, we have
\begin{align}\label{projff}
 f_{0; x_1, \ldots, x_K}(x)={\mathfrak f}_{0; y_1, \ldots, y_K}(y) 
 = \inf_{ z \in \sfp^{-1}(x)} {\mathfrak f}_{0; y_1, \ldots, y_K}(z),
\end{align}
where the $f_{0;x_1,\ldots, x_K}$ defined by \eqref{Alexdeff0}.
\item  the $F_{x_1,\ldots, x_K}$ is constant along each fiber $\sfp^{-1}(x)$ for every $x \in \sfX$. Indeed,
\begin{align}\label{FCnfib}
 F_{x_1,\ldots,x_K}(y) = f_{0;x_1,\ldots,x_K}(\sfp(y)),
\end{align}
 which implies that
\begin{align}\label{finfinff}
  f_{0; x_1, \ldots, x_K}(x) 
   =  \inf_{y \in \sfp^{-1}(x)}F_{x_1,\ldots, x_K}(y) = 
    \inf_{\substack{(y,y_1,\ldots, y_K) \\ \in \sfS(x;x_1,\ldots,x_K)}}
   {\mathfrak f}_{0; y_1, \ldots, y_K}(y).
\end{align}
\item With reference to the notation in \eqref{E-def} appearing in Lemma~\ref{parextHJ},
\begin{align}\label{BeqE}
{\mathcal E}_{\Lambda(x_1,\ldots,x_K)}^-[{\mathfrak f}_{0;\cdot}(y)] = S_{x_1,\ldots, x_K}(y);
\end{align}
\item With reference to the notation in \eqref{defEe-}, 
\begin{align}\label{E0Fpx}
 E_0^-[F_{x_1,\ldots,x_K};x] = \sfp^{-1}(x).
\end{align}
\end{enumerate}
\end{lemma}
\begin{proof}
The $\mathfrak f_{0;y_1, \ldots, y_K}$ and $f_{0;x_1, \ldots, x_K}$ are both defined in terms of a function $\psi \in {\bf \Psi}_K$ in \eqref{Adeffrakf}-\eqref{Alexdeff0}.
By non-emptiness of the $S_{x_1,\ldots, x_K}(y)$ and by component-wise monotonicity in the $\psi$, \eqref{projff} follows.
Taking $\inf_{(y_1,\ldots, y_K) \in \Lambda(x_1,\ldots, x_K)}$ over both sides of \eqref{projff} gives \eqref{FCnfib}, hence \eqref{finfinff}.

To see \eqref{BeqE} holds, we notice that for $y_k \in \sfp^{-1}(x_k)$ to satisfy 
\begin{align*}
 {\mathfrak f}_{0;y_1,\ldots, y_K}(y) 
 = \inf_{\substack{  z_k\in \sfp^{-1}(x_k) \\ k=1,\ldots, K}}
  {\mathfrak f}_{0; z_1,\ldots, z_K}(y)
\end{align*}
is equivalent to 
\begin{align*}
 \sfd_\sfY(y, y_k) = \sfd_\sfX(x,x_k), \quad \forall k=1,\ldots, K.
\end{align*}

\eqref{E0Fpx} follows directly from the constant along fiber property \eqref{FCnfib}.

\end{proof}

Recall the $f_{0;x_1,\ldots,x_K}$s in \eqref{Alexdeff0}, 
we define
\begin{align}\label{projH0}
 H_0 f_{0;x_1,\ldots,x_K}(x)
  : = \inf_{\substack{(y_0, y_1,\ldots,y_k) \\ \in \sfS(x;x_1, \ldots,x_K)}}
   {\mathcal H} {\mathfrak f}_{0; y_1, \ldots, y_K}(y_0).
\end{align}
In the following result, without pursuing generality, we assume that $\sfp^{-1}(x)$ is compact for every $x \in \sfX$. There are a number of ways to relax this assumption. 
\begin{lemma}\label{PSubHJ}
Let $\overline{\mathfrak f} \in M(\sfY;\R)$ be a point-wise strong viscosity sub-solution to \eqref{subHY}. 
Suppose that,  for every $x \in \sfX$, $\sfp^{-1}(x)$ is compact in $\sfY$,
and that 
\begin{align*}
\overline{\mathfrak f}(y) = \text{constant},  \quad \forall y \in \sfp^{-1}(x).
\end{align*}
We define $\overline{f}(x) := \overline{\mathfrak f}(y)$, $ \forall y \in \sfp^{-1}(x)$.
Then $\overline{f} \in M(\sfX)$ is a strong point-wise viscosity sub-solution to \eqref{H0proj} with the  operator $H_0$ in \eqref{projH0}, and with the $h_0:= \sfp_{\sup} {\mathfrak h}_0$.
\end{lemma}
\begin{proof}
Under the assumption that $\sfp^{-1}(x)$ being compact for each $x \in \sfX$, the submetry $\sfp: \sfY \mapsto \sfX$ becomes a {\em strong} submetry.

Let $x_1,\ldots, x_K \in \sfX$, $y_1, \ldots, y_K \in \sfY$, ${\mathfrak f}_{0;y_1,\ldots, y_K} \in {\mathcal S}^+_\sfY$ and  $F_{x_1,\ldots, x_K}(y)$ be defined as above. By  viscosity extension Lemma~\ref{parextHJ} and in view of \eqref{FDEF} and \eqref{BeqE}, we define operator
\begin{align*}
{\tilde {\mathcal H}}F_{x_1,\ldots, x_K}(y):= 
\inf_{\substack{(y_1,\ldots, y_K) \in \\  
 {\mathcal E}^{-}_{\Lambda(x_1,\ldots,x_K)}[{\mathfrak f}_{0;\cdot}(y)]}} 
 {\mathcal H} {\mathfrak f}_{0; y_1, \ldots, y_K}(y)
  = \inf_{\substack{(y_1, \ldots, y_K) \\ \in S_{x_1,\ldots, x_K}(y)}}
 {\mathcal H} {\mathfrak f}_{0;y_1, \ldots, y_K}(y),
\end{align*}
then $\overline{\mathfrak f}$ is a strong point-wise sub-solution to 
\begin{align*}
\overline{\mathfrak f} - \alpha \tilde{\mathcal H} \overline{\mathfrak f}
  \leq {\mathfrak h}_0.
\end{align*}

Next, we apply the projected viscosity solution Lemma~\ref{spHJ3}.
We have by \eqref{FCnfib} and \eqref{finfinff},
\begin{align*}
f_{0;x_1,\ldots, x_K}  = \sfp_{\inf}  F_{x_1,\ldots, x_K} .
\end{align*} 
In addition,  in view of \eqref{E0Fpx}, 
\begin{align*}
 \inf_{y \in E_0^-[F_{x_1,\ldots,x_K};x]} \tilde{\mathcal H}F_{x_1,\ldots, x_K}(y) 
= \inf_{y \in \sfp^{-1}(x)}
    \inf_{\substack{(y_1, \ldots, y_K) \\ \in S_{x_1,\ldots, x_K}(y)}}
 {\mathcal H} {\mathfrak f}_{0;y_1, \ldots, y_K}(y)
     = H_0 f_0 (x).
\end{align*}
Hence the conclusion follows.
 \end{proof}

In the same vein, we consider 
\begin{align*}
{\mathfrak f}_1(y) &:= {\mathfrak f}_{1;y_1,\ldots, y_K}(y):=
 - \psi\big(\sfd_\sfY^2(y,y_1), \ldots, \sfd_\sfY^2(y,y_K)\big) \in {\mathcal S}^-_\sfY, \\
f_1(x)&:= f_{1;x_1,\ldots, x_K}(x):= - \psi\big(\sfd_\sfX^2(x,x_1), \ldots, \sfd_\sfX^2(x,x_K)\big) 
 \in {\mathcal S}^-_\sfX,
\end{align*}
and define operator
\begin{align}\label{projH1}
 H_1 f_1(x)
  : = \sup_{\substack{(y_0, y_1,\ldots,y_k) \\ \in  \sfS (x; x_1, \ldots,x_K)} }
   {\mathcal H} {\mathfrak f}_{1; y_1, \ldots, y_K}(y_0).
\end{align}
Then the following super-solution result holds.  
\begin{lemma}\label{PSupHJ}
Let $\underline{\mathfrak f} \in M(\sfY; \R)$ be a point-wise strong viscosity super-solution 
to \eqref{supHY}. Suppose that $\sfp^{-1}(x)$ is compact in $\sfY$ for each $x \in \sfX$,
and that the $\underline{\mathfrak f}$ is constant along each fiber $\sfp^{-1}(x)$, $x \in \sfX$. We define $\underline{f}(x) := \underline{\mathfrak f}(y)$, $\forall y \in \sfp^{-1}(x)$. Then this $\underline{f} \in M(\sfX;\R)$ is a point-wise strong viscosity super-solution to \eqref{H1proj} with the above defined 
single-valued operator $H_1$, and with $h_1:= \sfp_{\inf} {\mathfrak h}_1$.
\end{lemma}

 \subsection{Projected Hamiltonian with perturbed test functions}\label{PHJPert}
In multi-scale convergence applications, we usually need to introduce an extra perturbative term ${\mathfrak g}$ relative to those test functions appearing in Lemmas~\ref{PSubHJ} and \ref{PSupHJ}. 
For instance, in the sub-solution case, we consider test functions on $\sfY$ taking the form
\begin{align}\label{frakfg}
 {\mathfrak f}_{\mathfrak g; y_1,\ldots, y_K}  := {\mathfrak f}_{0;y_1,\ldots, y_K} + 
 {\mathfrak g}_{y_1,\ldots, y_K} \in D(\mathcal H)
\end{align}
with 
\begin{align*}
{\mathfrak f}_{0;y_1,\ldots, y_K} =
{\mathfrak f}_{0;y_1,\ldots, y_K} (y) = \psi\big( \sfd_\sfY^2(y,y_1), \ldots, \sfd_\sfY^2(y, y_K)\big) \in {\mathcal S}_\sfY^+,
\end{align*}
and a perturbative term
\begin{align*}
{\mathfrak g} := {\mathfrak g}_{y_1,\ldots, y_K}\in M(\sfY).  
\end{align*}

In certain class of asymptotic problems concerning sequence of Hamiltonian PDEs, there could be a separation of scale phenomenon. The perturbative term ${\mathfrak g}$ can be used to separate micro-scale structural information in the Hamiltonians from those of macro-scale.  The following  projective abstract viscosity solution theory are developed with such context in mind. 
 
\subsubsection{Simple perturbations} In this subsection, 
we consider a relatively simple scenario where the term $\mathfrak g$ can always be chosen to satisfy the following.  
\begin{condition}\label{gfreeDy}
For each ${\mathfrak f}_{0;y_1, \ldots, y_K}$, the perturbative term 
${\mathfrak g}:= {\mathfrak g}_{y_1,\ldots, y_K}$ in \eqref{frakfg} can only depend on the parameters  $y_1,\ldots, y_K$  through $x_1=\sfp(y_1), \ldots, x_K  = \sfp(y_K)\in \sfX$. That is,
\begin{align*}
  {\mathfrak g}_{x_1,\ldots, x_K}(y)
  := \inf_{\substack{y_k\in \sfp^{-1}(x_k),\\ k=1,\ldots,K}}{\mathfrak g}_{y_1,\ldots, y_K}(y) 
 = \sup_{\substack{y_k\in \sfp^{-1}(x_k),\\ k=1,\ldots,K}} {\mathfrak g}_{y_1,\ldots, y_K}(y), 
 \quad \forall y \in \sfY.
\end{align*}
\end{condition}
A trivial case satisfying the above is when the perturbative term has no dependence on the $y_1, \ldots, y_K$ at all.

Each of the ${\mathfrak g}:=  {\mathfrak g}_{x_1,\ldots, x_K}:={\mathfrak g}_{y_1,\ldots, y_K} : \sfY \mapsto \R$ satisfying Condition~\ref{gfreeDy} induces a 
$g:=g_{x_1,\ldots, x_K}: \sfX \mapsto \R$ defined by
\begin{align*}
g(x):= g_{x_1,\ldots,x_K}(x):= \inf_{y \in \sfp^{-1}(x)} {\mathfrak g}_{x_1,\ldots,x_K}(y).
\end{align*}

We recall the notations ${\mathfrak f}_{\mathfrak g; y_1, \ldots, y_K}$, which is defined in \eqref{frakfg} and $F_{x_1,\ldots, x_k}$ in \eqref{FDEF}. We introduce test functions on $\sfX$ of the form
\begin{align*}
f_{g; x_1,x_2,\ldots,x_K}(x)&:=  
 \inf_{y \in \sfp^{-1}(x)} \inf_{\substack{y_k\in \sfp^{-1}(x_k),\\ k=1,\ldots,K}} 
  {\mathfrak f}_{\mathfrak g; y_1,\ldots, y_K}(y) \\
  & =   \inf_{y \in \sfp^{-1}(x)} \big(F_{x_1,\ldots, x_K}(y) 
   + {\mathfrak g}_{x_1,\ldots,x_K}(y)\big) \\
  &=  f_{0;x_1,\ldots,x_K}(x) + g_{x_1,\ldots, x_K}(x).
\end{align*}
In the above, the second equality follows from Condition~\ref{gfreeDy} 
and the defining identity \eqref{Sec3:MPproj}, the third equality follows from  
\eqref{FCnfib}. Next, we define a Hamiltonian operator on this class of test functions by
\begin{align}\label{pH0Alt}
 H_0 f_{g; x_1,\ldots, x_K}(x)
  : =  \inf_{y \in E^-_0[{\mathfrak g}_{x_1,\ldots,x_K}; x]}  \inf_{\substack{(y_1,\ldots,y_k) \\ \in  S_{x_1,\ldots,x_K}(y)}}
 {\mathcal H} {\mathfrak f}_{\mathfrak g; y_1, \ldots, y_K} (y).
\end{align}

\begin{lemma}\label{PSubHJ2}[Projecting sub-solutions]
Let $\overline{\mathfrak f} \in M(\sfY;\R)$ be a point-wise strong viscosity sub-solution 
to \eqref{subHY}. Suppose that $\sfp^{-1}(x)$ is compact in $\sfY$ for every $x \in \sfX$,
and that the $\overline{\mathfrak f}$ satisfies the following for every $x$: 
\begin{align*}
 \overline{\mathfrak f}(y) = \text{constant}, \quad \forall y \in \sfp^{-1}(x).
\end{align*}
We define 
\begin{align*}
\overline{f}(x) := \overline{\mathfrak f}(y), \quad \forall y \in \sfp^{-1}(x).
\end{align*}
We consider a class of functions $\mathfrak g$ satisfying Condition~\ref{gfreeDy}, and use this class to define 
test functions of the form $f_{g;x_1,\ldots, x_K}$, and then define operator $H_0$ according to \eqref{pH0Alt}.  

Then $\overline{f} \in M(\sfX)$ and it is a point-wise strong viscosity sub-solution to \eqref{H0proj} with the $H_0$ defined by \eqref{pH0Alt} and with the $h_0:= \sfp_{\sup} {\mathfrak h}_0$.
\end{lemma}
\begin{proof}
We make two observations: First, we recall definitions of $\mathcal E^-_\Lambda$ in \eqref{E-def} and of $\Lambda:=\Lambda(x_1,\ldots, x_K)$ in \eqref{LamDef}. 
Since ${\mathfrak g}$ satisfies Condition~\ref{gfreeDy}, we have the first identity below; by \eqref{BeqE}, the second equality also holds:
\begin{align}\label{Sec3:calEid}
{\mathcal E}_{\Lambda(x_1,\ldots,x_K)}^-[ {\mathfrak f}_{\mathfrak g; \cdot}(y)] 
={\mathcal E}_{\Lambda(x_1,\ldots,x_K)}^-[{\mathfrak f}_{0;\cdot}(y)] 
=   S_{x_1,\ldots,x_K}(y).
\end{align}
Second, we introduce a new test function $F_{\mathfrak g;x_1,\ldots,x_K}: \sfY \mapsto \R$ by 
\begin{align*}
 F_{\mathfrak g;x_1,\ldots,x_K}(y) := \inf_{\substack{y_1 \in \sfp^{-1}(x_1) \\ \ldots \\  y_K \in \sfp^{-1}(x_K)}} {\mathfrak f}_{\mathfrak g; y_1,\ldots, y_K}(y)
  = F_{x_1,\ldots, x_K}(y)+ {\mathfrak g}_{x_1,\ldots, x_K}(y),
\end{align*}
where the $F_{x_1,\ldots, x_K}$ is defined in \eqref{FDEF}. In view of \eqref{FCnfib}, we have
\begin{align}\label{Sec3:Eid}
E_0^-[F_{\mathfrak g;x_1,\ldots, x_K};x ]
= E_0^-[{\mathfrak g}_{x_1,\ldots,x_K};x ] .
\end{align}

 Consequently, by \eqref{Sec3:calEid},
\begin{align*}
 \tilde{\mathcal H} F_{\mathfrak g; x_1,\ldots,x_K}(y):= 
  \inf_{\substack{(y_1,\ldots, y_k) \in \\ 
  {\mathcal E}_{\Lambda(x_1,\ldots,x_K)}^-[{\mathfrak f}_{{\mathfrak g}; \cdot}(y)] }} {\mathcal H} {\mathfrak f}_{\mathfrak g;y_1,\ldots,y_K}(y)  = \inf_{\substack{(y_1,\ldots, y_K) \\ \in S_{x_1,\ldots,x_K}(y)}} 
  {\mathcal H} {\mathfrak f}_{\mathfrak g;y_1,\ldots,y_K}(y),
\end{align*}
and by \eqref{Sec3:Eid},
\begin{align*}
 \inf_{y \in E_0^-[F_{\mathfrak g;x_1,\ldots, x_K};x ]}
 \tilde{\mathcal H}  F_{\mathfrak g; x_1,\ldots,x_K}(y) 
 = H_0 f_{\mathfrak g;x_1,\ldots,x_K}(x).
\end{align*}

Next, as in the proof of Lemma~\ref{PSubHJ},  we apply Lemmas~\ref{parextHJ} and ~\ref{spHJ3} to conclude.
 \end{proof}
 
Next, we consider the super-solution case. We define 
\begin{align}\label{frakfg1}
{\mathfrak f}_{\mathfrak g;y_1,\ldots,y_K}:=  {\mathfrak f}_{1;y_1,\ldots, y_K}+  
 {\mathfrak g}_{y_1,\ldots, y_K} \in D({\mathcal H}), 
 \end{align}
with 
\begin{align*}
{\mathfrak f}_{1;y_1,\ldots,y_K}(y) :=  
 - \psi\big(\sfd_\sfY^2(y,y_1), \ldots, \sfd_\sfY^2(y,y_K)\big) \in {\mathcal S}^-_\sfY,
 \end{align*}
 and ${\mathfrak g} \in M(\sfY)$. We assume the ${\mathfrak g}$ satisfies Condition~\ref{gfreeDy}.
Note that there is a slight abuse of notation as the above ${\mathfrak f}_{\mathfrak g;y_1,\ldots, y_K}$  is different than the  one in the sub-solution case \eqref{frakfg}. Denoting 
\begin{align*}
  f_{1;x_1,\ldots, x_K}(x):= - \psi\big(\sfd_\sfX^2(x,x_1), \ldots, \sfd_\sfX^2(x,x_K)\big) 
 \in {\mathcal S}^-_\sfX,
\end{align*}
and
\begin{align*}
 g_{x_1,\ldots, x_K}(x):= \sup_{y \in \sfp^{-1}(x)} {\mathfrak g}_{y_1,\ldots, y_K}(y), \quad \forall \sfp(y_k) =x_k, k=1,\ldots,K;
\end{align*}
we construct a new class of test functions on $\sfX$ by
 \begin{align*}
 f_{g;x_1,\ldots,x_K}(x):= f_{1;x_1,\ldots, x_K}(x) + g_{x_1,\ldots, x_K}(x),
\end{align*}
and introduce a Hamiltonian operator for functions on $\sfX$ by
\begin{align}\label{pH1Alt}
 H_1 f_{g;x_1,\ldots,x_K}(x)
  : = \sup_{y \in E^+_0[\mathfrak g_{x_1,\ldots, x_K};x]}
  \sup_{\substack{(y_1,\ldots,y_k) \\ \in  S_{x_1, \ldots,x_K}(y)} }
   \big({\mathcal H} {\mathfrak f}_{\mathfrak g; y_1, \ldots, y_K}\big)(y).
\end{align}
 
\begin{lemma}\label{PSupHJ2}[Projecting super-solutions]
The statements in Lemma~\ref{PSupHJ} still holds, when we replace the $H_1$ by a new one defined according to \eqref{pH1Alt}.
\end{lemma}

\subsubsection{A further simplifying situation}
Lemmas~\ref{PSubHJ2} and \ref{PSupHJ2} simplify significantly under the following.
\begin{condition}\label{Cg}
Condition~\ref{gfreeDy} holds, hence those ${\mathfrak g}_{y_1,\ldots, y_K}$s only appearing in \eqref{frakfg} and \eqref{frakfg1} only depend on $(y_1, \ldots y_K)$ through $(x_1, \ldots, x_K)$
\begin{align*}
  {\mathfrak g}_{x_1,\ldots, x_K}:= {\mathfrak g}_{y_1,\ldots,y_K} \in M(\sfY), \qquad \forall y_1 \in \sfp^{-1}(x_1), \ldots, y_K \in \sfp^{-1}(x_K).
\end{align*}

In addition, these ${\mathfrak g}_{y_1,\ldots, y_K}$s appearing 
in Lemmas~\ref{PSubHJ2} and \ref{PSupHJ2} 
are constant along fibers $\sfp^{-1}(x)$, for every $x \in \sfX$:
\begin{align}\label{gCfibr}
 {\mathfrak g}_{x_1,\ldots, x_K} (y) = \text{constant}_{x_1,\ldots,x_K}, \quad \forall y \in \sfp^{-1}(x).
\end{align}
\end{condition}
The above condition implies that
\begin{align*}
 E^-_0[{\mathfrak g}_{x_1,\ldots, x_K};x]  
 = E^+_0[{\mathfrak g}_{x_1,\ldots, x_K};x] = \sfp^{-1}(x),
\end{align*}
and that
\begin{align*}
 g_{x_1,\ldots, x_K}(x):= \inf_{y \in \sfp^{-1}(x)} {\mathfrak g}_{x_1,\ldots, x_K}(y) = \sup_{y \in \sfp^{-1}(x)} {\mathfrak g}_{x_1,\ldots, x_K}(y).
\end{align*}

\begin{lemma} \label{sPsubsup}
Suppose that Condition~\ref{Cg} holds.
 Replacing the $H_0$ defined in \eqref{pH0Alt}  by  
\begin{align*}
H_0 f_{g;x_1,\ldots, x_K}(x)   := \inf_{\substack{(y_0; y_1,\ldots, y_K)\\ \in \sfS(x;x_1,\ldots, x_K)}}
 {\mathcal H} {\mathfrak f}_{{\mathfrak g}; y_1, \ldots, y_K}(y),
 \end{align*}
 then the conclusions of Lemma~\ref{PSubHJ2} still holds.
Similarly, replacing the $H_1$ in \eqref{pH1Alt} by
\begin{align*}
 H_1 f_{g;x_1,\ldots, x_K}(x)   := \sup_{\substack{(y_0; y_1,\ldots, y_K)\\ \in \sfS(x;x_1,\ldots, x_K)}}
 {\mathcal H} {\mathfrak f}_{{\mathfrak g}; y_1, \ldots, y_K}(y),
\end{align*}
  the conclusions of Lemma~\ref{PSupHJ2} holds as well. 
\end{lemma}

\subsubsection{Parameter dependent perturbations, beyond simple situations}\label{Sec3:relax}
Our hydrodynamic limit application has a multi-scale averaging nature. In such setting, we will need to select the perturbative term $\mathfrak g_{y_1, \ldots, y_K}:= \mathfrak g_{y_1, \ldots, y_K}(y)$ in \eqref{frakfg} depending upon differential $d_y {\mathfrak f}_{0;y_1, \ldots, y_k}$, which makes Condition~\ref{gfreeDy} not satisfied. Next, we develop versions of Lemmas~\ref{PSubHJ2} and \ref{PSupHJ2} which are still applicable to such general situation, by using notions of $\delta$-Sections. These   approximate versions of the ${\sfS}$ and $S$ are defined in \eqref{sfSeDef} and \eqref{Sedef}.

We will also introduce an extra term $\mathfrak u$ in the test functions below, for a purpose different than mentioned above. Its usefulness will be clear once we combine results next with a Hamiltonian operator convergence theory in Section~\ref{CnvHJ} to verify convergence for sequence of solutions in Section~\ref{Sec6}.

We consider perturbed test functions on $\sfY$ taking the form
\begin{align}\label{Sec3:fug}
 {\mathfrak f}_{{\mathfrak u}, \epsilon {\mathfrak g}; y_1,\ldots, y_K}  := \big( {\mathfrak f}_{0;y_1,\ldots, y_K} + {\mathfrak u}  \big) + \epsilon {\mathfrak g}_{y_1,\ldots, y_K} \in D(\mathcal H), \quad \epsilon >0,
\end{align}
where
\begin{align}\label{Sec3:frakf0}
{\mathfrak f}_{0;y_1,\ldots, y_K} =
{\mathfrak f}_{0;y_1,\ldots, y_K} (y) = \psi\big( \sfd_\sfY^2(y,y_1), \ldots, \sfd_\sfY^2(y, y_K)\big) \in {\mathcal S}_\sfY^+,
\end{align}
with the ${\mathfrak u} \in M(\sfY)$ does not having dependence on any of the parameters $y_1, \ldots, y_K$; and being constant along $\sfp^{-1}(x)$ for each $x \in \sfX$.~\footnote{Hence Condition~\eqref{Cg} is satisfied with the $\mathfrak g$ replaced by $\mathfrak u$} We also require the perturbative term ${\mathfrak g} := {\mathfrak g}_{y_1,\ldots, y_K}  \in C_b(\sfY)$.  
 
The above assumptions on the $\mathfrak u$ implies in particular that, for every $x_k \in \sfX$ with $k=1,\ldots, K$,
\begin{align*}
u(x):= \inf_{y \in \sfp^{-1}(x)}
\inf_{\substack{y_1 \in \sfp^{-1}(x_1) \\ \ldots \\  y_K \in \sfp^{-1}(x_K)}}{\mathfrak u}(y)
=\inf_{y \in \sfp^{-1}(x)} {\mathfrak u}(y) 
=\sup_{y \in \sfp^{-1}(x)} {\mathfrak u}(y)
= \sup_{y \in \sfp^{-1}(x)} 
  \sup_{\substack{y_1 \in \sfp^{-1}(x_1) \\ \ldots \\  y_K \in \sfp^{-1}(x_K)}} 
   {\mathfrak u}(y)  
\end{align*}
is independent of the $(x_1, \ldots, x_K)$. It defines a function $u \in M(\sfX)$. Next, we define
\begin{align}\label{Sec3:projf3}
f_{u, \epsilon {\mathfrak g};x_1, \ldots, x_K}(x)& :=
\inf_{y \in \sfp^{-1}(x)} \inf_{\substack{y_1 \in \sfp^{-1}(x_1) \\ \ldots \\  y_K \in \sfp^{-1}(x_K)}}
 {\mathfrak f}_{{\mathfrak u}, \epsilon {\mathfrak g}; y_1,\ldots, y_K} (y) \\
 & = \inf_{y \in \sfp^{-1}(x)} \inf_{\substack{y_1 \in \sfp^{-1}(x_1) \\ \ldots \\  y_K \in \sfp^{-1}(x_K)}}
 \big( {\mathfrak f}_{0; y_1,\ldots, y_K} (y) + \epsilon {\mathfrak g}_{y_1,\ldots, y_K}(y)\big)+ u(x). 
 \nonumber
\end{align}
We have estimate
\begin{align*}
\sup_{x \in \sfX} \big| f_{u, \epsilon {\mathfrak g};x_1, \ldots, x_K}(x) 
 - \big( f_{0;x_1, \ldots, x_K}(x) + u(x)\big) \big|
 \leq  \epsilon \Vert {\mathfrak g}\Vert_\infty. 
\end{align*}
Here and below, we use notation
\begin{align*}
 \Vert {\mathfrak g}\Vert_\infty := 
  \sup_{y\in \sfY} \sup_{y_1, \ldots, y_K \in \sfY} |{\mathfrak g}_{y_1,\ldots, y_K}(y)|.
\end{align*}

To simplify, we impose the following.
\begin{condition}\label{Sec3:gCPar}
The map $(y, y_1, \ldots, y_K) \mapsto {\mathfrak g}_{y_1,\ldots, y_K}(y)$ is continuous.
\end{condition}

Next, we introduce notions of  $\delta$-approximate sections, and relate them with extremal parameter set defined in~\eqref{E-def}.
 
Let $x_0, x_1, \ldots, x_K \in \sfX$ and $K \in \N$ be given. We define, for each $\delta>0$, a notion of $\delta$-approximate section
\begin{align}\label{sfSeDef}
\sfS^\delta(x_0; x_1, \ldots, x_K) & :=\Big\{(y_0,y_1,\ldots, y_K) \in \Pi_{k=0}^K \sfp^{-1}(x_k):  
     \text{ such that } \\
&\qquad \qquad
 | \sfd_\sfY(y_k, y_0) - \sfd_\sfX (x_k, x_0)|<\delta,   k=1,\ldots,K  \Big\}.\nonumber
\end{align}
When $y_0 \in \sfp^{-1}(x_0) \subset \sfY$ is held fixed, we also write
\begin{align}\label{Sedef}
S^\delta_{x_1,\ldots,x_K}(y_0) & := \Big\{(y_1,\ldots, y_K) \in \Pi_{k=1}^K \sfp^{-1}(x_k)
:  \text{ such that }   \\
& \qquad \qquad 
  |\sfd_\sfY(y_k, y_0) - \sfd_\sfX (x_k, x_0)| < \delta, \text{ where } x_0 =p(y_0), k=1,\ldots,K  \Big\}.   \nonumber
\end{align}
See the following graph.
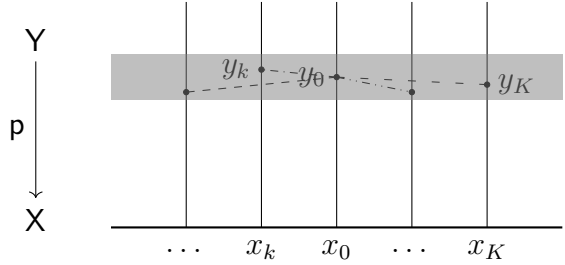
\begin{align*}
\begin{tikzpicture}
\node (Y) at (-4, 2.5) {$\sfY$};
\node (X) at (-4, 0.1) {$\sfX$};
\draw[black, thick] (-3,0) --(3,0);
\draw (-2,3) -- (-2,0);
\filldraw (-2,1.8) circle (1pt);
\draw [dashed] (0,2) --(-2,1.8);
\draw (-1,3) -- (-1,0);
\filldraw (-1,2.1) circle (1pt) node [anchor=east] {$y_k$};
\draw [dash dot] (0,2) -- (-1,2.1);
\draw (0,3) -- (0,0);
\filldraw (0,2.0) circle (1pt) node [anchor=east] {$y_0$};
\draw (1,3) -- (1,0);
\filldraw (1,1.8) circle (1pt)   ;
\draw [dash dot dot] (0,2) -- (1, 1.8);
\draw (2,3) -- (2,0);
\filldraw (2,1.9) circle (1pt)  node [anchor=west] {$y_K$};
\draw [loosely dashed] (0,2) --(2,1.9);
\node  at  (-2, -0.3) {\ldots};
\node     at (-1,-0.3) {$x_k$};
\node    at (0,-0.3) {$x_0$};
\node    at (1,-0.3) {\ldots};
\node    at (2,-0.3) {$x_K$};
\draw [->] (Y) -- (X)  node[midway, left] {$\sfp$};
\fill[gray, opacity=0.5] (-3,1.7 ) rectangle (3, 2.3);
 \end{tikzpicture}
\end{align*}

We have the following.
\begin{lemma}\label{Sec3:EcS}
Let $\mathfrak f_{0;y_1,\ldots, y_K}$ be as in \eqref{Sec3:frakf0} with the $\psi \in \Psi_K$. For each $\epsilon>0$, there exists $\delta := \delta(\epsilon; \Vert {\mathfrak g}\Vert_\infty, \psi)>0$, such that  
\begin{align*}
\mathcal E^-_{\Lambda(x_1,\ldots, x_K)}[ {\mathfrak f}_{{\mathfrak u}, \epsilon{\mathfrak g};\cdot}(y)] \subset S^\delta_{x_1,\ldots, x_K}(y),
\end{align*}
and both sets are non-empty.
Moreover, $\lim_{\epsilon \to 0^+} \delta (\epsilon; \Vert {\mathfrak g}\Vert_\infty, 
\psi) = 0$.
\end{lemma}
\begin{proof}
Let $y \in \sfY$. By assumption on the $u$,
\begin{align*}
\mathcal E^-_{\Lambda(x_1,\ldots, x_K)}[ {\mathfrak f}_{{\mathfrak u}, \epsilon{\mathfrak g};\cdot}(y)] =
 \mathcal E^-_{\Lambda(x_1,\ldots, x_K)}[ {\mathfrak f}_{0 ;\cdot}(y)+ \epsilon{\mathfrak g}_{\cdot}(y)].
\end{align*}
Hence we only need to prove the claim by setting ${\mathfrak u}=0$. 
Also, by Condition~\ref{Sec3:gCPar} and earlier assumption that $\sfp^{-1}(x_k)$ is compact in $\sfY$, the above set is non-empty.

Let 
\begin{align*}
(y_1, \ldots, y_K) \in   \mathcal E^-_{\Lambda(x_1,\ldots, x_K)}[ {\mathfrak f}_{0 ;\cdot}(y)+ \epsilon{\mathfrak g}_{\cdot}(y)] \subset \sfp^{-1}(x_1) \times \ldots \times \sfp^{-1}(x_K). 
\end{align*}
Then 
\begin{align*}
{\mathfrak f}_{0;y_1,\ldots, y_K}(y)+ \epsilon{\mathfrak g}_{y_1,\ldots,y_K}(y)
 &= \inf_{\substack{y_k^\prime \in \sfp^{-1}(x_k),\\ k=1,\ldots,K}}
 \big( {\mathfrak f}_{0;y_1^\prime,\ldots, y_K^\prime}(y)
  + \epsilon{\mathfrak g}_{y_1^\prime,\ldots,y_K^\prime}(y)\big) \\
  & \leq \inf_{\substack{y_k^\prime \in \sfp^{-1}(x_k),\\ k=1,\ldots,K}}
 \big( {\mathfrak f}_{0;y_1^\prime,\ldots, y_K^\prime}(y)\big)
  + \epsilon \Vert {\mathfrak g}\Vert_\infty 
  = f_{0;x_1,\ldots, x_K}(x) + \epsilon \Vert {\mathfrak g}\Vert_\infty,
\end{align*}
where the $x =\sfp(y)$. Consequently,
\begin{align*}
\big| \psi\big(\sfd^2_\sfY(y,y_1), \ldots, \sfd^2_\sfY(y,y_K) \big) 
- \psi\big(\sfd^2_\sfX(x,x_1), \ldots, \sfd^2_\sfX(x,x_K)\big) \big| < 2 \epsilon \Vert \mathfrak g\Vert_\infty.
\end{align*}
Since $\partial_k \psi >0$ for $k=1,\ldots, K$, we can find $\delta:= \delta(\epsilon; \Vert \mathfrak g\Vert_\infty, \psi)$ with desired property and
\begin{align*}
 \sup_{k=1,\ldots, K} |\sfd_\sfY(y,y_k) - \sfd_\sfX(x,x_k)| < \delta.
\end{align*}
Therefore $(y_1, \ldots , y_K) \in S^\delta_{x_1,\ldots, x_K}$.
\end{proof}

\begin{lemma} \label{sPrelax} 
In the context of Lemmas~\ref{PSubHJ2}, the conclusion still holds if we replace Condition~\ref{gfreeDy} by Condition~\ref{Sec3:gCPar},  and replace the $H_0$  
in \eqref{pH0Alt} by  
\begin{align*}
H_0 f_{u, \epsilon{\mathfrak g}; x_1,\ldots, x_K}(x)   := \sup_{\substack{(y; y_1,\ldots, y_K)\\ \in \sfS^\delta(x;x_1,\ldots, x_K)}}
 {\mathcal H}{\mathfrak f}_{{\mathfrak u}, \epsilon {\mathfrak g}; y_1,\ldots, y_K} (y),
 \end{align*}
for $f_{u, \epsilon{\mathfrak g}; x_1,\ldots, x_K}$ given by \eqref{Sec3:projf3}. The $\delta:= \delta(\epsilon; \Vert \mathfrak g\Vert_\infty, \psi)$ above can be any choice that is given by Lemma~\ref{Sec3:EcS}. 
\end{lemma}

\begin{proof}
Again, as in the proof of Lemma~\ref{PSubHJ}, we apply Lemma~\ref{parextHJ} and then Lemma~\ref{spHJ3} to conclude. Key details are given below.
 
Let a test function from $ \sfY \mapsto \R$ be defined as
\begin{align*}
 \tilde{F}_{x_1,\ldots,x_K}(y) :=\inf_{\substack{y_1 \in \sfp^{-1}(x_1) \\ \ldots \\  y_K \in \sfp^{-1}(x_K)}}
 {\mathfrak f}_{{\mathfrak u}, \epsilon {\mathfrak g}; y_1,\ldots, y_K} (y)  
=   \inf_{\substack{y_1 \in \sfp^{-1}(x_1) \\ \ldots \\  y_K \in \sfp^{-1}(x_K)}}
 \big( {\mathfrak f}_{0; y_1,\ldots, y_K} (y) + \epsilon {\mathfrak g}_{y_1,\ldots, y_K}(y)\big)+ u(x).
 \end{align*}
Let a new Hamiltonian operator $\tilde{\mathcal H}$ be defined on all such test functions by
\begin{align*}
\tilde{\mathcal H} \tilde{F}_{x_1, \ldots, x_K}(y)& :=\sup_{(y_1, \ldots, y_K) \in  S^\delta_{x_1,\ldots, x_K}(y)} {\mathcal H}  
  {\mathfrak f}_{{\mathfrak u}, \epsilon {\mathfrak g}; y_1,\ldots, y_K} (y) \\
  & \geq \inf_{(y_1, \ldots, y_K) \in \mathcal E^-_{\Lambda(x_1,\ldots, x_K)}[ {\mathfrak f}_{{\mathfrak u}, \epsilon{\mathfrak g};\cdot}(y)]}
 {\mathcal H}   {\mathfrak f}_{{\mathfrak u}, \epsilon {\mathfrak g}; y_1,\ldots, y_K}(y),
\end{align*}
where the $\delta>0$ is the one selected in Lemma~\ref{Sec3:EcS}, and the above inequality follows by that lemma.
Apply Lemma~\ref{parextHJ}, $\overline{\mathfrak f}$ is a strong point-wise sub-solution to 
$\overline{\mathfrak f} - \alpha \tilde{\mathcal H} \overline{\mathfrak f} \leq {\mathfrak h}_0$. 

We have shown that $f_{u, \epsilon{\mathfrak g}; x_1,\ldots, x_K}(x)=\inf_{y \in \sfp^{-1}(x)} \tilde{F}_{x_1,\ldots,x_K}(y)$ in \eqref{Sec3:projf3}. Next, we note
\begin{align*}
 \bigcup_{\substack{y \in E^-_0[\tilde{F}_{x_1,\ldots,x_K};x_0] \\
 \subset \sfp^{-1}(x_0)}}
\Big( \{ y\} \times \big( \bigcup_{y \in \sfp^{-1}(x_0)} S^\delta_{x_1,\ldots, x_K}(y) \big)
\Big)  \subset \sfS^\delta(x_0; x_1, \ldots, x_K).
\end{align*}
Therefore
\begin{align*}
 \inf_{y \in E^-_0[\tilde{F}_{x_1,\ldots,x_K};x]} \tilde{\mathcal H} \tilde{F}_{x_1, \ldots, x_K}(y)
  \leq      \sup_{\substack{(y; y_1,\ldots, y_K)\\ \in \sfS^\delta(x;x_1,\ldots, x_K)}}
 {\mathcal H} {\mathfrak f}_{{\mathfrak u}, \epsilon {\mathfrak g}; y_1,\ldots, y_K} (y).
\end{align*}
Apply Lemma~\ref{spHJ3} and we conclude.
 \end{proof}

In the same vein, a super-solution version of the result holds as well:
We consider test functions in $\sfY$ with the form
\begin{align*}
 {\mathfrak f}_{{\mathfrak u}, \epsilon {\mathfrak g}; y_1,\ldots, y_K}  
 := \big( {\mathfrak f}_{1; y_1,\ldots, y_K}  + {\mathfrak u}\big) 
 + \epsilon {\mathfrak g}_{y_1,\ldots, y_K} \in D({\mathcal H}).
 \end{align*}
Similar to \eqref{Sec3:projf3}, we have  
\begin{align}\label{Sec3:projf4}
f_{u, \epsilon {\mathfrak g};x_1, \ldots, x_K}(x)& :=
\sup_{y \in \sfp^{-1}(x)} \sup_{\substack{y_1 \in \sfp^{-1}(x_1) \\ \ldots \\  y_K \in \sfp^{-1}(x_K)}}
 {\mathfrak f}_{{\mathfrak u}, \epsilon {\mathfrak g}; y_1,\ldots, y_K} (y) \\
 & = \sup_{y \in \sfp^{-1}(x)} \sup_{\substack{y_1 \in \sfp^{-1}(x_1) \\ \ldots \\  y_K \in \sfp^{-1}(x_K)}}
 \big( {\mathfrak f}_{1; y_1,\ldots, y_K} (y) + \epsilon {\mathfrak g}_{y_1,\ldots, y_K}(y)\big)+ u(x); 
 \nonumber
\end{align}
with estimate
\begin{align*}
\sup_{x \in \sfX} \big| f_{u, \epsilon {\mathfrak g};x_1, \ldots, x_K}(x) 
 - \big( f_{1;x_1, \ldots, x_K}(x) + u(x)\big) \big|
 \leq  \epsilon \Vert {\mathfrak g}\Vert_\infty. 
\end{align*}
Moreover, we have the following estimate.
\begin{lemma}\label{Sec3:EcS2}
For each $\epsilon>0$, there exists $\delta := \delta(\epsilon; \Vert {\mathfrak g}\Vert_\infty, 
\psi)>0$ where the $\psi \in \Psi_K$ is the one defining $\mathfrak f_{1;y_1,\ldots, y_K}$, such that
\begin{align*}
\mathcal E^+_{\Lambda(x_1,\ldots, x_K)}[ {\mathfrak f}_{{\mathfrak u}, \epsilon{\mathfrak g};\cdot}(y)] \subset S^\delta_{x_1,\ldots, x_K}(y).
\end{align*}
Moreover, $\lim_{\epsilon \to 0^+} \delta (\epsilon; \Vert {\mathfrak g}\Vert_\infty, 
\psi) = 0$.
\end{lemma}

\begin{lemma} \label{sPrelax2} 
In the context of Lemmas~\ref{PSupHJ2}, the conclusions still hold, if we replace Condition~\ref{gfreeDy} by Condition~\ref{Sec3:gCPar},  and replace the $H_1$  
in \eqref{pH1Alt} by  
 \begin{align*}
 H_1 f_{u, \epsilon{\mathfrak g}; x_1,\ldots, x_K}(x)   := \inf_{\substack{(y; y_1,\ldots, y_K)\\ \in \sfS^\delta(x;x_1,\ldots, x_K)}}
 {\mathcal H} {\mathfrak f}_{{\mathfrak g}; y_1, \ldots, y_K}(y),
\end{align*}
 where the $\delta:= \delta(\epsilon; \Vert {\mathfrak g}\Vert_\infty, 
\psi)$ is the one from Lemma~\ref{Sec3:EcS2}.
\end{lemma}

\newpage

\section{A viscosity convergence theory in metric spaces} \label{CnvHJ}
To rigorously handle the hydrodynamic limit problem in our introduction, 
we need a convergence theory for viscosity solutions in space of probability measures. In this section, we build such a theory by generalizing the Barles-Perthame convergence scheme~\cites{BP87,BP88} for Hamilton-Jacobi equations to a general metric space setting. To a large extent, such work was developed in Feng and Kurtz~\cite{FK06}. The primary goal there was to apply the result to probabilistic large deviation theory, hence some estimates were formulated probabilistically.  Next, we adapt the same ideas and translate the arguments using only the language of classical analysis. There are more than one way to achieve this, here we choose an approach by formulating conditions on special test functions. Such formualtion is better suited for PDE applications with minimal structural assumptions on the equations or solutions.
 
Throughout this section, $\alpha>0$ is a fixed number, $(\sfX_n,\sfd_{\sfX_n})$ 
and $(\sfX, \sfd_\sfX)$ are complete metric spaces. We are given Hamiltonian operators 
\begin{align*}
H_{n,0} \subset M(\sfX_n;\R) \times M(\sfX_n; \R), \quad
H_{n,1} \subset  M(\sfX_n;\R ) \times M(\sfX_n, \R), 
\end{align*}
and functions $h_{n,0}, h_{n,1}: \sfX_n \mapsto \R$, with
$\overline{f}_n, \underline{f}_n: \sfX_n \mapsto \bar{\R}$ 
respectively viscosity sub- and super-solutions in the sequential sense to 
\begin{align}
\overline{f}_n - \alpha H_{n,0} \overline{f}_n & \leq h_{n,0},  \label{Hnsub}\\
\underline{f}_n - \alpha H_{n,1} \underline{f}_n & \geq h_{n,1}. \label{Hnsup}
\end{align}
We are also given two operators 
\begin{align*}
H_0 & \subset M(\sfX;\bar{\R} ) \times M(\sfX;  \bar{\R}), \\ 
 H_1 & \subset M(\sfX;\bar{\R} ) \times M(\sfX;  \bar{\R}). 
\end{align*}
Throughout this paper, we implicit assume that domains $D(H_0)$ and $D(H_1)$ of the operators consist of non-trivial functions (that is, $f \not \equiv \pm \infty$). 
One can think of them as playing roles of upper- and lower- bounds on limits of the $H_{n,i}$s, in a sense to be made precise next. We will define a kind of upper limit to the $\overline{f}_n$s by $\overline{f}$ in \eqref{fbar}, and a kind of lower limit to the $\underline{f}_n$s by $\underline{f}$ in \eqref{lbarf}.
We then show that (see Theorem~\ref{Sec4:BPThm}) they are respectively sub- and super-solutions to equations 
\begin{align}
 \overline{f} - \alpha H_0 \overline{f} \leq h_0, \label{Hsub} \\
 \underline{f} - \alpha H_1 \underline{f} \geq h_1. \label{Hsup}
\end{align}
In the special case when $f_n:= \overline{f}_n =\underline{f}_n$, assuming a comparison principle holds between the above two equations, we can conclude that $f:=\overline{f}=\underline{f}$ and $f_n$ convergences to $f$ in appropriately defined senses. 

The main technical difficulty here
is that we need to handle convergence of functions (and operators acting on such functions) in possibly non-locally compact metric spaces. There were three key ingredients introduced in Feng and Kurtz~\cite{FK06}. One, it used the sequential definition of viscosity solution (Definition~\ref{SeqVisDef}) which was motivated by maximum principle considerations. 
\footnote{Throughout this section, we don't explicitly assume our Hamilton-Jacobi operator satisfies a nonlinear maximum principle (e.g. Appendix A.3 in \cite{FK06}). However, from a functional analytic point of view, the whole generalized viscosity method is only natural when this is true.}
Two, it relied on uniform estimates on sequence of certain probability measures on compact sets. These are occupation measures arising from integral kernel representation for resolvents of Hamilton-Jacobi equations. In context of stochastic optimization problems, such representations always exist. See Lemma 5.9 and estimate (7.38) of \cite{FK06} which are consequences from a probabilistically formulated Condition 2.8  in that book. See also a special property in Lemma A.11 in the appendix of \cite{FK06}.
Third, the book~\cite{FK06} also introduced a multi-valued viscosity operator approach to handle technical difficulties arising from multi-scale convergence of Hamiltonians using  a variational approach.  In the following, we will present a purely analytic and somewhat different approach to the similar ideas mentioned above. In a much simpler setting involving only PDE in Euclidean spaces,  Feng, Fouque and Kumar described yet another similar approach in Section 4 of~\cite{FFK12}. 

We mention that the second ingredient in \cite{FK06} mentioned above, regardless of what language is used, requires a property that does not hold in our hydrodynamic limit example here. We could alter the growth condition on external potential term $U$ to enforce such property. But that will create complication involving semi-continuities for functions and operators which becomes difficult to handle. In this paper, we introduce another way to solve the issue by using multiple topologies. In addition to the abstract developments in this section, we also refer to concrete calculations and estimates in Section~\ref{Sec6} for details of applying this new technique.   
 
\subsection{Convergence of metric spaces -- generalized Gromov-Hausdorff convergence } 
Building on earlier works in semigroup convergence theory, Feng and Kurtz~\cite{FK06} introduced a notion of topological convergence of spaces to space, and subsequently, notion on convergence of functions on these spaces to functions in the limiting space. The idea can be traced back at least to Trotter~\cite{Tro58} with generalizations by Kurtz~\cites{Kurtz69,Kurtz70}. These formulations emphasize on {\em almost isomorphisms} of the approximating spaces.  Next, we strengthen these notions by requiring a kind of metric convergence, placing emphasize on {\em approximate isometry} (e.g.~\cite{Villani09}). It is meaningful to do this because that, for the applications we have in mind, the test functions are basically compositions of distance squared functions. See $\mathcal S^+$ and $\mathcal S^-$ defined in \eqref{SS+} and \eqref{SS-}. Such development can be more useful when applying to equations defined with a metric geometry nature. Additionally, using such more restrictive notion of convergence simplifies the method of \cite{FK06}, making the results more accessible. 

Consequently, we are lead to generalize the notion of Gromov-Hausdorff convergence to (possibly) non-locally compact metric spaces 
\begin{align*}
\{ (\sfX_n, \sfd_{\sfX_n})\}_{n \in\N} 
       \stackrel{\rm  gGH}{\longrightarrow}_{\mathcal Q} (\sfX, \sfd_\sfX)
\end{align*}
with respect to some pre-chosen index set $\mathcal Q$. See Definition~\ref{gGH} next.

Let closed subsets $A_n \subset \sfX_n$ and $A \subset \sfX$. We recall several equivalent definitions and properties of Gromov-Hausdorff convergence of  $A_n \stackrel{\rm GH}{\rightarrow} A$.  For definitions and generic properties, we refer to Chapter 27 of Villani~\cite{Villani09}, Chapter 7 of Burago, Burago and Ivanov~\cite{BBI01}, and pages 70-77 in Bridson and Haeflinger~\cite{BH99}.  In particular, we can define a metric $\sfd_{\rm GH}(K_1, K_2)$ as in (27.1) and (27.2) in~\cite{Villani09} to measure the distance between two metric spaces $K_1, K_2$. In fact, when the background metric spaces where these $K_1, K_2$s live in are compact, the Gromov-Hausdorff convergence is given by this metric topology.  For non-compact cases, modifications are needed. Villani~\cite{Villani09} gave several alternative definitions on pages 755-758. In the following, we adapt one of them into Definition~\ref{gGH}. 
For simplicity, we will use only the $\epsilon_n$-isometry version of definition of the Gromov-Hausdorff convergence as introduced through  properties (a') and (b') on top part of page 750 of~\cite{Villani09}. We will make explicit reference to these approximate isometries in our formulation next.

Our point of departure is the following basic setup: 
\begin{enumerate}
\item $(\sfX, \sfd_\sfX)$ and $(\sfX_n, \sfd_{\sfX_n})$, $n=1,2,\ldots$,  are metric spaces; 
\item $\mathcal Q$ is a given index set, with
$\{ K_n^q \subset \sfX_n : q \in \mathcal Q\}$ a family of  closed subsets in 
$(\sfX_n, \sfd_{\sfX_n})$ and $\{ K^q \subset \sfX : q\in \mathcal Q\}$  a family of closed subsets in $(\sfX, \sfd_\sfX)$; 
\item \label{Sec4:ApIso} there is a family of maps $\{ \eta^q_n : n=1,2,\ldots \}_{q \in \mathcal Q}$, such that each $\eta_n^q : K^q_n \mapsto K^q$ is an $\epsilon_n$-isometry with $\epsilon_n \to 0$; 
\item for every $q_1, q_2 \in \mathcal Q$, there exists $q_3 \in \mathcal Q$ such that 
$ K^{q_1} \cup K^{q_2} \subset K^{q_3}$ and that the approximate isometries are consistent in the sense that
\begin{align*}
 \eta_n^{q_3}\big|_{K_n^{q_1}} = \eta_n^{q_1}, \quad \eta_n^{q_3}\big|_{K_n^{q_2}} =\eta_n^{q_2}.
\end{align*}
\end{enumerate}
\begin{remark}
The last consistency condition is automatically satisfied, if there exists $\eta_n: \sfX_n \mapsto \sfX$ and the choice $\eta_n^q:= \eta_n \big|_{K_n^q}$ forms sequences of approximate isometries.
\end{remark}


\begin{definition}[Generalized Gromov-Hausdorff convergence]\label{gGH}
The sequence of spaces $(\sfX_n, \sfd_{\sfX_n})$ is said to converge to $(\sfX, \sfd_\sfX)$ in sense of {\em generalized Gromov-Hausdorff convergence with respect to $\mathcal Q$ and by means of approximate isometries $\{ \eta^q_n : n \in \N \}_{q \in \mathcal Q}$}, denoted by
\begin{align*}
(\sfX_n, \sfd_{\sfX_n})\stackrel{\rm  gGH}{\longrightarrow}_{\mathcal Q} (\sfX, \sfd_\sfX),
\end{align*}
if the followings are satisfied
\begin{enumerate}
 \item $\sfX_0:=\cup_{q \in \mathcal Q} K^q$ is dense in $(\sfX, \sfd_\sfX)$;
 \item for each $q \in \mathcal Q$, $K^q \subset \sfX$ is compact in $(\sfX, \sfd_\sfX)$;
\item  for each $q \in \mathcal Q$,  
\begin{align}\label{KnK}
\lim_{n \to \infty} \sfd_{\rm GH}(K_n^q, K^q) =0 
\end{align}
by means of the approximate isometries $\{ \eta^q_n \}_{n \in \N}$ as given above in the basic setup.
\end{enumerate}
\end{definition}
Unlike Villani's Definition 27.11 in \cite{Villani09}, the above definition does not require the $K_n^q$s to be compact in $\sfX_n$. However, they are necessarily ``asymptotically compact" by the requirement in \eqref{KnK}. For instance, even in the case when $\sfX_n = \sfX$, we can choose the $K_n^q$s to be closure of $\delta_n$-fattenings of some compact sets $K^q$s, with $\delta_n \to 0$.  

From now on, we require the following.
\begin{condition}\label{HConSpa}\
  For the given $\mathcal Q$ and $\{ \eta_n^q \}_{n \in \N, q \in \mathcal Q}$, we have
\begin{align}\label{XntoX}
 (\sfX_n, \sfd_{\sfX_n}) \stackrel{\rm gGH}{\longrightarrow}_{\mathcal Q} (\sfX, \sfd_\sfX).
\end{align} 
\end{condition}

It can be useful to identify a special point within each space $\sfX_n$ playing the role of ``origin" of the space.  We introduce the following notation.
\begin{condition}\label{Sec4:ptMS}[Pointed metric spaces]
There exists $x_0 \in K^{q_0} \subset \sfX$ for some $q_0 \in \mathcal Q$, and $x_{n,0} \in K_n^{q_0}$ with $\eta_n^{q_0}(x_{n,0}) \to x_0$.
\end{condition}

\subsection{A metric space version of the half-relaxed limit theory}

Let $f: K^q \mapsto \bar{\R}$ and $\eta_n^q: K_n^q \mapsto K^q$, we denote $(\eta_n^q f):= f\circ \eta_n^q : K_n^q \mapsto \bar{\R}$. 

\begin{definition}[Generalized $\Gamma$-convergence]
Let $f_n: \sfX_n \mapsto \bar{\R}$ and $f: \sfX \mapsto \bar{\R}$. 
We say that $f_n$ Gamma-converges to $f$  over sets indexed by $\mathcal Q$, denoted by 
$f_n \stackrel{\rm  \Gamma}{\longrightarrow}_{\mathcal Q} f$, if for every $q \in \mathcal Q$, $K^q$ and $K_n^q$ with approximate isometry $\eta_n^q : K_n^q \mapsto K^q$, the following properties hold:
\begin{enumerate}
\item (Liminf property:) for every $x_0 \in K^q$ and $x_n \in K_n^q$ satisfying $\eta_n^q(x_n) \to x_0$, we have
\begin{align*}
\liminf_{n \to \infty} f_n(x_n) \geq f(x_0);
\end{align*}
\item (Existence of recovering sequence:) for each $x_0 \in K^q$, there exists $q^\prime \in \mathcal Q$ such that $x \in K^{q^\prime}$, $x_n \in K_n^{q^\prime}$, and $\eta_n^{q^\prime}(x_n) \to x_0$, and that
\begin{align*}
\limsup_{n \to \infty} f_n(x_n) \leq f(x_0).
\end{align*}
\end{enumerate}
\end{definition}

\subsubsection{Conditions on convergence of functions and operators}

\begin{condition}\label{HJhCND}[Convergence of the $h_{n,0}$ and $h_{n,1}$]
The $h_{n,i}: \sfX_n \mapsto \R$ and $h_i \in C(\sfX)$, $i=0,1$, $n \in \N$ have following property:
  for every $q \in \mathcal Q$ and the associated $K_n^q, K^q$ with approximate isometry $\eta_n^q: K_n^q \mapsto K^q$, we have
\begin{align*}
  \limsup_{n \to \infty} \sup_{K_n^q} (h_{n,0}- \eta_n^q h_0) \leq 0, \text{ and } 
  \liminf_{n \to \infty} \inf_{K_n^q} (h_{n,1} -\eta_n^q h_1) \geq 0.
\end{align*}
 \end{condition}

Recall that, at this point,  $f_0 \in D(H_0)$ and $f_1 \in D(H_1)$ may be discontinuous functions. Indeed, the $f_0, f_1 \in M(\sfX; \bar{\R})$ may not even be finite on the whole $\sfX$.
\begin{condition}\label{HConSub}[Convergence of Hamiltonian operators, sub-solution case]
Every $f_0 \in D(H_0)$ has the following property: for every $x \in \sfX$, there exists $x_n \in \sfX_0$ such that $\lim_{n \to \infty} \sfd_\sfX(x_n, x_0)=0 $ and $\limsup_{n \to \infty} f_0(x_n) \leq f_0(x_0)$.

For each $(f_0, g_0) \in H_0$, there exists $(f_{n,0},  g_{n,0}) \in H_{n,0}$
satisfying the following:
\begin{enumerate}
\item \label{uno} [Operator convergence]  
\begin{align*}
f_{n,0} \stackrel{\rm  \Gamma}{\longrightarrow}_{\mathcal Q} f_0 ;
\end{align*}
and for every $q \in \mathcal Q$ with
 $K_n^q \stackrel{\rm GH}{\to} K^q$ by means of the approximate isometries $\eta_n^q$, we have 
\begin{align*}
 \limsup_{n \to \infty} \sup_{K_n^q}   \big( g_{n,0}  -  \eta_n g_0  \big) \leq 0;
 \end{align*}
 \item \label{due} [Solution growth properties]
 there exists a non-decreasing $\zeta\in C(\R; \R)$ with 
 super-linear growth at $+\infty$, 
 and sub-linear growth at $-\infty$:
 \begin{align*}
 \text{ i.e. } \quad \liminf_{r \to +\infty} r^{-1} \zeta(r) =+\infty, 
 \quad \text{ and } \quad
 \liminf_{r \to -\infty} |r|^{-1} \zeta(r) =0, 
\end{align*}
such that
 \begin{align}\label{Sec4:f0domf}
 \zeta \big( \overline{f}_n(x)\big) \leq f_{n,0}(x), 
 \quad  \zeta \big( h_{n,0}(x) \big) \leq f_{n,0}(x), 
 \quad \forall x \in \sfX_n, n\in \N;
\end{align}
\item\label{tre} [Almost compactness properties]
for each $L>0$, there exists $q:= q(L) \in \mathcal Q$, such that 
\begin{align*}
\{ x \in \sfX_n : f_{n,0}(x) \leq L \} \cap
\{ x \in \sfX_n:  g_{n,0}(x) \geq - L \} \subset K_n^q.
\end{align*}
\end{enumerate}
 \end{condition}
 \begin{remark}
 Condition~\ref{HConSub}.\ref{due} is trivially satisfied, if 
 \begin{align*}
\sup_n (\sup_{\sfX_n} \overline{f}_n +\sup_{\sfX_n} h_{n,0}) < +\infty, \quad
 \inf_n \inf_{\sfX_n} f_{n,0} > -\infty.
\end{align*}
\end{remark}
\begin{condition}\label{HConSup}[Convergence of Hamiltonian operators, super-solution case]
Every $f_1 \in D(H_1)$ has the following property: for each $x \in \sfX$, there exists $x_n \in \sfX_0$ such that $\lim_{n \to \infty} \sfd_\sfX(x_n, x) =0$ and that $\liminf_{n \to \infty} f_1(x_n) \geq f_0(x)$.

For each $(f_1, g_1) \in H_1$, there exists a sequence of $(f_{n,1},  g_{n,1}) \in H_{n,1}$
satisfying the following: 
\begin{enumerate}
\item \label{uno1}  
\begin{align*}
- f_{n,1} \stackrel{\rm  \Gamma}{\longrightarrow}_{\mathcal Q} (-f_1);
\end{align*}
for every $q \in \mathcal Q$ with $K_n^q \stackrel{\rm GH}{\to} K^q$ 
by means of  the approximate isometries $\eta_n^q$, we have
\begin{align*}
 \liminf_{n \to \infty} \inf_{K_n^q} \big(g_{n,1} -  g_1\big) \geq 0;
\end{align*}
\item \label{due1} there exists a non-decreasing function
 $\zeta \in C(\R; \R)$ with super-linear growth at $+\infty$ and 
  sub-linear growth at $-\infty$, just as in 
 Condition~\ref{HConSub}.\ref{due}, such that
 \begin{align*}
 - \zeta \big( - \underline{f}_n(x)\big) \geq f_{n,1}(x), \quad 
 - \zeta \big( - h_{n,1}(x) \big) \geq f_{n,1}(x), \quad \forall x \in \sfX_n, n\in \N;
\end{align*}
 \item\label{tre1} 
For each $L>0$, there exists $q:= q(L) \in \mathcal Q$, such that 
\begin{align*}
\{ x \in \sfX_n:  f_{n,1}(x)   \geq - L \} \cap
\{ x \in \sfX_n:    g_{n,1}(x) \leq  L \} \subset K_n^q.
\end{align*}
\end{enumerate}
 \end{condition}

\subsubsection{Construction of limiting sub- super-solutions}
First, we introduce two functions defined on the $\sfX_0 =\cup_{q \in \mathcal Q} K^q$.
Assuming Condition~\ref{HConSpa}, let $x_0 \in  \sfX_0$,  for each $q \in \mathcal Q$ such that $x_0 \in K^q$ and corresponding $K_n^q \subset \sfX_n$, we have $\lim_{n \to \infty} \sfd_{\rm GH}(K_n^q, K^q) =0$ by means of $\epsilon_n$-isometries $\eta_n^q : K_n^q \mapsto K^q$ with some $\epsilon_n \to 0$. In particular, because of the almost surjective property of $\epsilon_n$-isometry (property (b') on page 750  in \cite{Villani09}), there exists $x_n \in K_n^q$ such that $\sfd_\sfX(\eta_n^q(x_n), x_0) \leq \epsilon_n \to 0$.
We define, for $x_0 \in \sfX_0$,
\begin{align}
\hat{f}(x_0) & :=  \sup_{\substack{ q\in \mathcal Q \\ 
 \text{ s.t.} x_0 \in K^q}}   
\sup \big\{\limsup_{n \to \infty} \overline{f}_n(x_n) : \exists x_n  \in K_n^q  
 \text{ s.t.}\lim_{n \to \infty} \sfd_\sfX(\eta_n^q(x_n) ,x_0)=0 \big\},   \label{fhat} \\
 \hat{\hat{f}}(x_0)&:=  \inf_{\substack{ q\in \mathcal Q, \\ 
   \text{ s.t.} x_0 \in K^q}}  \inf \big\{ \liminf_{n \to \infty} \underline{f}_n(x_n) :   
    \exists x_n \in K_n^q 
  \text{ s.t.} \lim_{n \to \infty} \sfd_\sfX( \eta_n^q(x_n), x_0) =0 \big\}.  \label{fhathat} 
 \end{align}
 Second, we extend definitions of $\hat{f}$ and $\hat{\hat{f}}$ from domain 
 $\sfX_0$ to  $\sfX$: for each $x \in \sfX$, we define
 \begin{align}
\overline{f}(x)& :=    \lim_{\epsilon \to 0} \quad
 \sup \big\{ \hat{f}(x_0) :  x_0 \in \sfX_0,  
 \sfd_\sfX (x_0, x) <\epsilon \big\}, \label{fbar} \\
\underline{f}(x)&:=   \lim_{\epsilon \to 0}  \quad
 \inf \big\{ \hat{\hat{f}}(x_0) : x_0 \in \sfX_0,  
  \sfd_\sfX (x_0, x) <\epsilon \big\}.   \label{lbarf}
\end{align}
\begin{lemma}\label{barfP}
The above defined $\overline{f}, \underline{f} : \sfX \mapsto \R$ have the following properties:
\begin{enumerate}
\item  $\overline{f} \in \USC(\sfX; \bar{\R})$ and
 $\underline{f} \in \LSC(\sfX ; \bar{\R} )$;
\item Suppose that $\underline{f}_n \leq \overline{f}_n$, then 
$\underline{f} \leq \overline{f}$ in $\sfX$.
\item \label{barfP3} Let $q \in \mathcal Q$ and $K_n^q \stackrel{\rm GH}{\to} K^q$  by means of the $\epsilon_n$-isometry $\eta_n^q$.
Suppose that $x_n \in K_n^q$ and $x_0 \in K^q$ are such that 
$\lim_{k \to \infty} \sfd_\sfX \big(\eta_{n_k}^q(x_{n_k}),  x_0 \big)=0$ along some subsequence $\{ n_k : k=1,2,\ldots\}$. Then 
\begin{align*} 
 \limsup_{k \to \infty} \overline{f}_{n_k}(x_{n_k}) \leq \hat{f}(x_0) \leq \overline{f}(x_0), \quad  
 \limsup_{k \to \infty} \underline{f}_{n_k}(x_{n_k}) \geq \hat{\hat{f}}(x_0) \geq \underline{f}(x_0).
\end{align*}
\end{enumerate}
\end{lemma}
\begin{proof}
First, the semi-continuities of $\overline{f}, \underline{f}$ are consequences of their definitions. Specifically, let $x_n \to x \in \sfX$ in $(\sfX, \sfd_\sfX)$.
By \eqref{fbar}, there exists $x_{n,0} \in \sfX_0:=\cup_{q \in \mathcal Q}K^q$ such that    
$\lim_{n \to \infty} \sfd_\sfX(x_n, x_{n,0})=0$ and that $\limsup_{n \to \infty} \overline{f}(x_n)  \leq \limsup_{n \to \infty} \hat{f}(x_{n,0})$. On the other hand, by \eqref{fbar}, we also have 
\begin{align*}
 \limsup_{n \to \infty} \hat{f}(x_{n,0}) \leq \overline{f}(x).
\end{align*} 
Combine the above together gives $\overline{f} \in \USC(\sfX; \bar{\R})$. The case $\underline{f} \in \LSC(\sfX; \bar{\R})$ is similar.

Second, suppose that  $\underline{f}_n \leq \overline{f}_n$, then
 $\hat{\hat{f}} \leq \hat{f}$ in $\sfX_0$, consequently,
$\underline{f} \leq \overline{f}$ in $\sfX$.

Third, the last property of the Lemma follows from \eqref{fhat}-\eqref{fbar} and \eqref{fhathat}-\eqref{lbarf}.
 \end{proof}
 
Frequently, we can find {\it a priori} modulus of continuity estimates for the 
$\{\overline{f}_n, \underline{f}_n \}_{n \in \N}$.
 
 \begin{condition}\label{Sec4:Wmodfn}[Uniform modulus of continuity estimates]
For every $q \in \mathcal Q$, there exists $\omega_q \in C( \R_+ ; \R_+)$ with $\omega_q(0) =0$ such that 
\begin{align}\label{modd}
| \overline{f}_n (x) - \overline{f}_n(y) | 
+ |\underline{f}_n (x) - \underline{f}_n(y) | 
\leq \omega_q\circ \sfd_\sfX \big(\eta^q_n(x), \eta^q_n(y)\big), 
\quad \forall x, y \in K_n^q.
\end{align}
 \end{condition}
 
The above estimates usually hold in stronger forms.
 \begin{condition}\label{Sec4:modfn}[A strengthened form of modulus of continuity estimates]
 Condition~\ref{Sec4:Wmodfn} holds. Moreover, either one of the following holds:  
 \begin{enumerate}
 \item  for each compact $K \subset \sfX$ in $(\sfX, \sfd)$, there exists a $q \in \mathcal Q$ such that $K \subset K^q$;
 \item for each $x \in \sfX$, there exists a sufficiently small $\delta:=\delta(x)>0$, and
 a modulus $\omega_{x,\delta}:= \omega_{x,\delta}(r) \in C(\R_+; \R_+)$ with $\omega_{x,\delta}(0) =0$ such that 
 \begin{align*}
| \overline{f}_n (x) - \overline{f}_n(y) | 
+ |\underline{f}_n (x) - \underline{f}_n(y) | 
& \leq \omega_{x,\delta} \circ \sfd_\sfX \big(\eta^q_n(x), \eta^q_n(y)\big), \\
& \qquad   \quad \forall x, y \in K_n^q \cap \bar{B}(x;\delta), \quad \forall q \in {\mathcal Q}.
\end{align*}
where the $\bar{B}(x,\delta)$ denotes a closed $\sfd$-metric ball of size $\delta$ with center $x$.
 \end{enumerate}
\end{condition}
 
 \begin{lemma}\label{Sec4:fC}
Under Condition~\ref{Sec4:modfn}, $\overline{f}, \underline{f} \in C(\sfX)$.
\end{lemma}
\begin{proof}
We only verify $\overline{f} \in C(\sfX)$ next. The case of $\underline{f}$ can be handled similarly.  

 First, we assume Condition~\ref{Sec4:Wmodfn}. Let $\epsilon>0$ be given. For every $x^\prime, y^\prime \in \sfX_0$, by definition  of $\hat{f}$, there exists $q := q_\epsilon \in \mathcal Q$ and $x_n \in K_n^q$ with $\lim_{n \to \infty} \sfd(\eta_n(x_n) ,x^\prime) =0$ such that 
$\hat{f}(x^\prime) < \epsilon + \limsup_{n \to \infty}  \overline{f}_n(x_n)$. Next, we can re-choose the $q$ if necessary to make certain $y^\prime \in K^q$ as well.
Therefore, for every $y_n \in K_n^q$ with $\lim_{n \to \infty} \sfd(\eta_n^q(y_n) , y^\prime)=0$, we have
\begin{align*}
 \hat{f}(x^\prime) - \hat{f}(y^\prime) & < \epsilon 
 + \limsup_{n \to \infty} \big( \overline{f}_n(x_n) - \overline{f}_n(y_n) \big) \\
 & \leq \epsilon + \limsup_{n \to \infty} 
 \omega_q \circ \sfd( \eta^q_n(x_n), \eta^q_n(y_n)) \\
 & \leq \epsilon + \omega_q \circ \sfd(x^\prime, y^\prime),
\end{align*}
for some  $q \in \mathcal Q$ that only depends on $x^\prime, y^\prime$.  

Second, let $x, x_n \in \sfX$ with $\lim_{n \to \infty}\sfd(x_n, x) =0$. Condition~\ref{Sec4:modfn} enables us to conclude, using the above estimate,  
\begin{align*}
\liminf_{n \to \infty} \overline{f}(x_n) \geq \overline{f}(x).
\end{align*}
In view of Lemma~\ref{barfP}, we conclude.

\end{proof}

\begin{lemma}\label{f0f1max}
Suppose Condition~\ref{Sec4:Wmodfn} holds.
\begin{enumerate}
\item Let $f_0: \sfX \mapsto \bar{\R}$ and $f_{n,0} : \sfX_n \mapsto \R$ be such that 
\begin{align*}
f_{n,0} \stackrel{\rm  \Gamma}{\longrightarrow}_{\mathcal Q} f_0.
\end{align*}
Then for every $\delta>0$ and $K^q$ with $q \in \mathcal Q$, there exists another $q^\prime:=q^\prime(\delta, q) \in \mathcal Q$ such that 
\begin{align*}
 \sup_{K^q} (\hat{f}  - f_0) \leq \delta + 
 \limsup_{n \to \infty} \sup_{K_n^{q^\prime}}(\overline{f}_n - f_{n,0}). 
\end{align*}
Moreover, if Condition~\ref{Sec4:modfn} holds, then
\begin{align}\label{abHJf0m}
 \sup_\sfX (\overline{f}  - f_0) \leq \limsup_{n \to \infty} \sup_{\sfX_n}(\overline{f}_n - f_{n,0}). 
\end{align}
\item Similarly, let $f_1: \sfX \mapsto \bar{\R}$ and $f_{n,1} : \sfX_n \mapsto \bar{\R}$ be such that 
\begin{align*}
(- f_{n,1}) \stackrel{\rm  \Gamma}{\longrightarrow}_{\mathcal Q} (-f_1).
\end{align*}
Then for every $\delta>0$ and $K^q$ with $q \in \mathcal Q$, there exists another 
$q^\prime:=q^\prime(\delta, q)$ such that 
\begin{align*}
\sup_{K^q} (f_1-\underline{f}) \leq \delta + 
 \limsup_{n \to \infty} \sup_{K_n^{q^\prime}}(f_{n,1} -\underline{f}_n); 
\end{align*}
and, if Condition~\ref{Sec4:modfn} holds, then
\begin{align}\label{abHJf1m}
 \sup_\sfX (f_1 - \underline{f}) \leq   \limsup_{n \to \infty} \sup_{\sfX_n}(f_{n,1} - \underline{f}_n). 
\end{align}
\end{enumerate}
\end{lemma}
\begin{proof}
For the given $\delta$ and $K^q$, there exists a $x_0 \in K^q$ such that  
\begin{align*}
\sup_{K^q} (\hat{f} - f_0) \leq \delta/2 + (\hat{f} -f_0)(x_0).
\end{align*}
We assume without loss of generality that $f_0(x_0)<+\infty$.
By the definition of $\hat{f}$ in \eqref{fhat}, there exists $q^\prime \in \mathcal Q$ with the $x_0 \in K^{q^\prime}$ and $x_n \in K_n^{q^\prime} \stackrel{\rm GH}{\to} K^{q^\prime}$ by means of an approximate isometry $\eta_n^{q^\prime}$, such that $\lim_{n \to \infty} \sfd_\sfX(\eta_n^{q^\prime}(x_n),x_0)=0$, and that
\begin{align*}
 \hat{f}(x_0) < \frac{\delta}{2} + \limsup_{n \to \infty} \overline{f}_n(x_n).
\end{align*}
By $f_{n,0} \stackrel{\rm  \Gamma}{\longrightarrow}_{\mathcal Q} f_0$, there exists $\hat{x}_n \in K_n^q$ such that (re-choose the $q^\prime \in \mathcal Q$ if necessary) $\eta_n^{q^\prime}(\hat{x}_n) \to x_0$ and that $f_{n,0}(\hat{x}_n) \to f_0(x_0)$. 
Therefore, 
\begin{align*}
 \sup_{K^q} \big(\hat{f} -f_0\big) 
  & < \delta + \limsup_{n \to \infty} \big( \overline{f}_n(x_n)  - f_{n,0}(\hat{x}_n)\big) \\
  & < \delta + \limsup_{n \to \infty} \big( \overline{f}_n(x_n)  - \overline{f}_n(\hat{x}_n)\big) 
   +  \limsup_{n \to \infty} \sup_{K_n^{q^\prime}}(\overline{f}_n - f_{n,0})\\
 & \leq  \delta +
 \limsup_{n \to \infty} \sup_{K_n^{q^\prime}}(\overline{f}_n - f_{n,0}),
\end{align*}
where we used Condition~\ref{Sec4:Wmodfn} to get the last inequality.

Next, by density of $\sfX_0$ in $\sfX$, by $\overline{f} \in C(\sfX)$ (Lemma~\ref{Sec4:fC}), and by the property listed in the beginning of Condition~\ref{HConSub} for $f_0$,  it follows that  
\begin{align*}
 \sup_{\sfX} (\overline{f} - f_0) = \sup_{q \in \mathcal Q} \sup_{K^q} (\hat{f} - f_0),
\end{align*}
hence  \eqref{abHJf0m} follows.

The case of \eqref{abHJf1m} is verified similarly.
\end{proof}
\begin{remark}\label{Sec4:testfC}
With simple modifications in the proof, the above  results still hold if we remove the uniform modulus Conditions~\ref{Sec4:Wmodfn} and \ref{Sec4:modfn}, but strengthen the assumption $f_{n,0} \stackrel{\rm  \Gamma}{\longrightarrow}_{\mathcal Q} f_0$ by stronger conditions that $f_0 \in C(\sfX)$ and that $\lim_{n \to \infty} \sup_{K_n^q} | f_{n,0} - \eta_n f_0|=0$ for every $q \in \mathcal Q$. Similar remark applies also to the case involving $\underline{f}$ and $f_{n,1}$ and $f_1$. 
\end{remark}
  
\subsubsection{A half-relaxed limit theorem}
\begin{definition}[Comparison principle for a pair of Hamilton-Jacobi equations]\label{pairCMPdef}
We say that {\em comparison principle} holds between sequential (resp. point-wise, etc) sub-solutions of \eqref{Hsub} and sequential (resp. point-wise, etc) super-solutions of \eqref{Hsup}, if for every such sub-solution $f^* \in \USC(\sfX;\R)$ and every such super-solution $f_* \in \LSC(\sfX;\R)$, we have  
\begin{align*}
 \sup_\sfX (f^* - f_*) \leq \sup_\sfX (h_0 - h_1).
\end{align*}
\end{definition}

\begin{theorem}\label{Sec4:BPThm}
Let $\overline{f}$ and $\underline{f}$ be defined by \eqref{fbar} and \eqref{lbarf}. Assume both of them are finite functions (i.e. $\overline{f},\underline{f} : \sfX \mapsto \R$). Suppose that Conditions~\ref{HConSpa}, \ref{HJhCND}, \ref{HConSub} and \ref{HConSup} hold, and that $D(H_0) \subset \LSC(\sfX;\R \cup \{+\infty\})$ 
and $D(H_1) \subset \USC(\sfX;\R \cup \{-\infty\})$.  We also assume that either Condition~\ref{Sec4:modfn} holds, or the modified requirements in
Remark~\ref{Sec4:testfC} holds. Then 
\begin{enumerate}
\item $\overline{f} \in \USC(\sfX; \R)$ is a  sub-solution to \eqref{Hsub} in the sequential viscosity solution sense. 
Similarly,   $\underline{f} \in \LSC(\sfX; \R)$ is a  super-solution to \eqref{Hsup} in the sequential viscosity solution sense.
\item in the special case where $f_n:=\overline{f}_n=\underline{f}_n$ and $h =h_0=h_1$, if, in addition, we assume that the comparison principle holds between sequential sub-solutions of \eqref{Hsub} and sequential super-solutions of \eqref{Hsup}.   Then 
\begin{align*}
 f:=\overline{f} = \underline{f} \in C(\sfX)
\end{align*}
and  
\begin{align}\label{Sec4:cnvfn}
 \lim_{n \to \infty} \sup_{K_n^q} |f_n  - \eta_n^q f|=0,  
 \end{align}
 for every $q \in \mathcal Q$ and the associated $K_n^q, K^q$ with approximate isometries 
 $\eta_n^q: K_n^q \mapsto K^q$.
 \item In the above, if we strengthen requirements on limiting operator that $H_0 \subset \LSC(\sfX;\bar{\R}) \times \USC(\sfX;\bar{\R})$, then the $\overline{f}$ is a point-wise viscosity sub-solution.
 Similarly, assuming $H_1\subset \USC(\sfX;\bar{\R}) \times \LSC(\sfX;\bar{\R})$, then the $\underline{f}$ is a point-wise viscosity super-solution.
 \end{enumerate}
\end{theorem}
\begin{proof}
We only show that the $\overline{f}$ is a sub-solution, the proof for $\underline{f}$ being a super-solution can be done similarly.

Let $(f_0, g_0) \in H_0$ be such that $ \sup_\sfX (\overline{f} -f_0)<\infty$. 
Then there exists 
$(f_{n,0}, g_{n,0}) \in H_{n,0}$ satisfying Condition~\ref{HConSub}. Since $\overline{f}_n$ is a viscosity sub-solution in the sequential sense, we can find $\epsilon_n \to 0$ and $x_n \in \sfX_n$ such that  
\begin{align}\label{Sec4:BPseqV}
  \sup_{\sfX_n}(\overline{f}_n - f_{n,0}) \leq \epsilon_n+ (\overline{f}_n - f_{n,0})(x_n),   \text{ and }
    (\overline{f}_n - h_{n,0} - \alpha g_{n,0}) (x_n) \leq \epsilon_n.
\end{align}
By first part of the above estimates, and in view of \eqref{Sec4:f0domf} and \eqref{abHJf0m} (see also Remark~\ref{Sec4:testfC}),
\begin{align*}
-\infty< \sup_\sfX( \overline{f} -f_0 ) \leq \limsup_{n \to \infty} 
\sup_{\sfX_n} (\overline{f}_n -f_{n,0}) \leq \limsup_{n \to \infty} (\overline{f}_n - \zeta \circ \overline{f}_n)(x_n). 
 \end{align*}
Since $\lim_{r \to +\infty} r^{-1}\zeta(r) =+\infty$, selecting sub-sequence if necessary,
\begin{align*}
 \sup_n \overline{f}_n(x_n) <+\infty.
\end{align*}
Reapply the above estimate back to the first part of \eqref{Sec4:BPseqV},
\begin{align*}
\limsup_{n \to \infty} f_{n,0} (x_n) \leq - \sup_\sfX( \overline{f} -f_0 ) 
 + \sup_{n \in \N} \overline{f}_n(x_n) < +\infty.
\end{align*}
Next, from the second part of \eqref{Sec4:BPseqV},  again invoking \eqref{abHJf0m} and growth estimate \eqref{Sec4:f0domf},
\begin{align*}
\liminf_{n \to \infty} \alpha g_{n,0}(x_n) 
& \geq \liminf_{n \to \infty} (\overline{f}_n -f_{n,0})(x_n)
 +\liminf_{n \to \infty} (f_{n,0}-h_{n,0})(x_n) \\
&\geq \sup_{\sfX} (\overline{f} - f_0) 
 + \liminf_{n \to \infty} (\zeta \circ h_{n,0} - h_{n,0}) (x_n) \\
 & \geq  \sup_{\sfX} (\overline{f} - f_0) 
 + \inf_{r \in \R} \big(\zeta(r) -r\big)> -\infty.
\end{align*}
In summary, selecting subsequence if necessary, there is a large enough but finite $L >0$ such that 
\begin{align*}
 \sup_n f_{n,0}(x_n) \leq L, \text{ and } \inf_n g_{n,0}(x_n) \geq -L.
\end{align*}

In view of Condition~\ref{HConSub}.\ref{tre},  
there exists $q \in \mathcal Q$ such that $x_n \in K_n^q$ for all $n$. 
Since $K^q \subset \sfX$ is compact and  $K_n^q \stackrel{\rm GH}{\to}K^q$ by means of approximate isometry $\eta_n^q: K_n^q \mapsto K^q$, we can find a subsequence $n_k$ and a point $x_0 \in K^q$ such that 
\begin{align*}
\lim_{k \to \infty} \sfd_\sfX (\eta_{n_k}^q(x_{n_k}), x_0) =0, \text{ and } \limsup_{n \to \infty} \overline{f}_n(x_n) = \lim_{k \to \infty} \overline{f}_{n_k}(x_{n_k}).
\end{align*}
Therefore,
\begin{align*}
 \limsup_{n \to \infty} \sup_{\sfX_n} (\overline{f}_n - f_{n,0})  =  \limsup_{n \to \infty} (\overline{f}_n - f_{n,0})(x_n)
 \leq \hat{f}(x_0) - f_0(x_0)  \leq \overline{f}(x_0) - f_0(x_0).
\end{align*}
In the above, the first inequality follows from the $\liminf$ property of
 $f_{n,0} \stackrel{\rm  \Gamma}{\longrightarrow}_{\mathcal Q} f_0$, the second inequality from part~\ref{barfP3} of Lemma~\ref{barfP}. Combined with \eqref{abHJf0m}, we arrive at 
$(\overline{f} - f_0)(x_0) = \sup_\sfX (\overline{f}  - f_0)$. 
Since $f_{n_k,0}(x_{n_k}) \to f_0(x_0)$, the above also implies that 
\begin{align}\label{Sec4:barfCn}
\lim_{k \to \infty} \overline{f}_{n_k}(x_{n_k}) = \overline{f}(x_0).
\end{align}
Since $\overline{f} \in \USC(\sfX; \R)$ (Lemma~\ref{barfP}), 
\begin{align*}
\limsup_{k \to \infty} (\eta_{n_k}\overline{f}-\overline{f}_{n_k})(x_{n_k}) \leq 
 \overline{f}(x_0) - \lim_{k \to \infty} \overline{f}_{n_k}(x_{n_k}) = 0.
\end{align*}
Consequently, noting $x_n \in K_n^q$ and in view of the convergence assumptions in Conditions~\ref{HJhCND} and \ref{HConSub}, 
\begin{align*}
\limsup_{k \to \infty} \eta_{n_k}^q (\overline{f}- h_0 - \alpha g_0)(x_{n_k}) & \leq 
\limsup_{k \to \infty} (\overline{f}_{n_k}- h_{n_k,0} - \alpha g_{n_k,0})(x_{n_k}) \\
& \qquad  + \limsup_{k \to \infty} (\eta_{n_k}^q \overline{f} - \overline{f}_{n_k})(x_{n_k})\\
& \qquad \qquad 
  +\limsup_{n \to \infty} \big( \sup_{K_n^q} (h_{n, 0} -\eta_{n}^q h_0)
   + \alpha \sup_{K_n^q} (g_{n, 0} -\eta_{n}^q g_0)\big) \\
  & \leq 0.
\end{align*}
That is, the $\overline{f}$ is a sub-solution to \eqref{Hsub} in the sequential viscosity sense. 
 
Next, we assume that the comparison principle holds. In view of Lemma~\ref{barfP}, 
$f:= \overline{f}=\underline{f} \in C(\sfX)$. Especially, selecting subsequence if necessary and apply item \ref{barfP3} in that lemma give \eqref{Sec4:cnvfn}.

Finally, suppose $g_0 \in \USC(\sfX)$, then by \eqref{Sec4:barfCn} and the sequential viscosity sub-solution property,
\begin{align*}
 \overline{f}(x_0)  = \lim_{k \to \infty} \overline{f}_{n_k}(x_{n_k}) 
& \leq \limsup_{k \to \infty} \big( \alpha g_{n_k,0}(x_{n_k}) + h_{n_k,0}(x_{n_k}) \big) \\
& \leq \limsup_{k \to \infty} \big( \alpha \eta_{n_k}^q  g_0(x_{n_k}) + \eta_{n_k}^q h_0(x_{n_k}) \big)   \leq \alpha g_0(x_0) + h_0(x_0).
\end{align*}
That is, the $\overline{f}$ is a sub-solution to \eqref{Hsub} in the point-wise viscosity sense. 
\end{proof}

\subsection{Another version of  the half-relaxed limit theory - generalizations}\label{BPlimitA}
While Theorem~\ref{Sec4:BPThm} is readily applicable to the super-solution case of our hydrodynamic limit example, it does not apply to the sub-solution case directly. This is because that, provided we work with order $2$-Wasserstein space alone, we cannot construct test functions so that the almost compactness requirements in Condition~\ref{HConSub}.\ref{tre}, and the convergence requirements in 
Condition~\ref{HConSub}.\ref{uno} are satisfied simultaneously. Next, we introduce a variant of the previous arguments by using multiple topologies through embedding the original space $(\sfX, \sfd_\sfX)$ into a larger space $(\sfX^\prime, \sfd_{\sfX^\prime})$. The topology generated by $\sfd_{\sfX^\prime}$ is a weaker one, giving a larger family of neighborhood sets, hence helping some of the limit arguments. 
\footnote{This is inspired by the $B$-continuous solution idea introduced in Crandall and Lions~\cite{CL91}. The settings and structural properties we explore are very different though. As the example in Section~\ref{Sec6:Sub} shows, the Hamiltonian operator has no obvious ``coercive" term as critically used in \cite{CL91}. However, we will explore a perturbative argument \eqref{Sec6:fbarPert}, and a type of growth estimates on Hamiltonian operators acting on such perturbation, in the applications of this paper.}

\subsubsection{Basic setup on spaces}
We now work with the following:
\begin{enumerate}
\item $(\sfX_n, \sfd_{\sfX_n})$, $(\sfX, \sfd_\sfX)$ and $(\sfX^\prime, \sfd_{\sfX^\prime})$ are metric spaces;
\item $\sfX$ is a closed subset in $(\sfX^\prime, \sfd_{\sfX^\prime})$, where the topology generated by $\sfd_{\sfX^\prime}$ is weaker than that by $\sfd_\sfX$;  
\item $\mathcal Q$ is a prescribed index set; $\{ K_n^q \subset \sfX_n : q \in {\mathcal Q}\}$ is a family of closed subsets in $(\sfX_n, \sfd_{\sfX_n})$, and $\{ K^q \subset \sfX \subset \sfX^\prime : q \in {\mathcal Q}\}$ is a family of {\em closed} subsets in $(\sfX, \sfd_\sfX)$, and {\em compact} subsets in $(\sfX^\prime, \sfd_{\sfX^\prime})$;
\item  $\sfX = \overline{\sfX_0}^{\sfd\text{\rm-closure}}$ where
the $\sfX_0:= \cup_{q \in \mathcal Q} K^q$; 
\item \label{Sec4:gGH-C4} there is a family of maps $\{ \eta_n^q : n=1,2,\ldots\}_{q \in \mathcal Q}$ such that $\eta^q_n: K_n^q \mapsto K^q$  is an $\epsilon_n$-isometry, when the $K^q$ is considered as a compact set in metric space $(\sfX^\prime, \sfd_{\sfX^\prime})$, for some $\epsilon_n \to 0^+$.
\item for every $q_1, q_2 \in \mathcal Q$, there exists $q_3 \in \mathcal Q$ such that 
$K^{q_1} \cup K^{q_2} \subset K^{q_3}$ and that the following consistency holds:
\begin{align*}
 \eta_n^{q_3}|_{K_n^{q_1}} = \eta_n^{q_1}, \qquad 
 \eta_n^{q_3}|_{K_n^{q_2}} = \eta_n^{q_2}.
\end{align*}
 \end{enumerate}
The following is a counterpart of Condition~\ref{HConSpa} in current setup:
\begin{condition}\label{Sec4:Hausd}
The following holds:
\begin{align*}
 (\sfX_n, \sfd_{\sfX_n}) \stackrel{\rm gGH}{\longrightarrow}_{\mathcal Q} (\sfX^\prime, \sfd_{\sfX^\prime}),
\end{align*} 
where the $\mathcal Q$ here means that the approximate isometries in the convergence is given by the family of maps $\{ \eta_n^q \}_{n \in \N}, q \in \mathcal Q$ in item \ref{Sec4:gGH-C4} above.
\end{condition}
Note that the above implies in particular $\sfX^\prime= \overline{\sfX_0}^{\sfd^\prime\text{\rm-closure}}$.

\begin{example}
The result of this subsection is largely designed for the application in later Section~\ref{Sec6:Sub}.
For that example, we will take $\sfX:= {\mathcal P}_2(\R^d)$ with $\sfd_\sfX:= \sfd$ the $2$-Wasserstein metric, and $\sfX^\prime:= {\mathcal P}_{p_0}(\R^d)$ a $p_0$-Wasserstein space with $p_0 \in (1,2)$ fixed, and $p_0$-Wasserstein metric $\sfd_{\sfX^\prime}:= \sfd_{p_0}$. The $\sfX_n$ will be taken to be  space of $n$-points equally weighted empirical probability measures, with possibility of multiple identical points. We take $\sfd_{\sfX_n}:= \sfd_{p_0}\big|_{\sfX_n}$. Let $\eta_n (\rho) = \rho$ be the identity embedding map of $\sfX_n$ into $\sfX^\prime$. $\mathcal Q:=\R_+$ and $K^q:= \{ \rho \in \sfX^\prime : \sfd(\rho, \rho_0)\leq q\}$ and $K^q_n:= K^q \cap \sfX_n$; $\eta_n^q:= \eta_n \big|_{K_n^q}$.
 
One can verify that the above setup fits this situation and Condition~\ref{Sec4:Hausd} is satisfied.
\end{example}
\subsubsection{Upper limit gives a sub-solution}
Similar to the introduction of $\hat{f}$ and $\overline{f}$ in \eqref{fhat} and \eqref{fbar}, we define, for each $x_0 \in \sfX_0$,
\begin{align}\label{Sec4:ftil} 
\tilde{f}(x_0)  :=  \sup_{\substack{ q\in \mathcal Q \\ 
 \text{ s.t.} x_0 \in K^q}}   
\sup \big\{\limsup_{n \to \infty} \overline{f}_n(x_n) : \exists x_n  \in K_n^q  
 \text{ s.t.}\lim_{n \to \infty} \sfd_{\sfX^\prime}(\eta_n^q(x_n) ,x_0)=0 \big\};   
 \end{align}
 and for every $x \in \sfX$,
 \begin{align}\label{Sec4:fstarD}
 f^*(x) := \lim_{\epsilon \to 0^+} \sup \big\{ \tilde{f}(x_0) : x_0 \in \sfX_0, \sfd_\sfX(x_0, x) < \epsilon \big\}.
\end{align}
It follows then $\tilde{f} \leq f^*$. The $\hat{f}$ and $\overline{f}$ defined in \eqref{fhat} and \eqref{fbar} are all smaller than $f^*$. 
\begin{lemma}
The $f^* \in \USC(\sfX; \bar{\R})$. $\hat{f} \leq \tilde{f}$ on $\sfX_0$ and $\overline{f} \leq f^*$ on $\sfX$. 
\end{lemma}

\begin{lemma}\label{Sec4:supEst2}
Suppose that $f_0  \in C((\sfX, \sfd_\sfX))$. 
In addition, we assume that the $f_0$ is $\sfd_{\sfX^\prime}$-continuous in $K^q \subset \sfX$, for each $q \in \mathcal Q$ fixed: that is, for every $x_n, x_0 \in K^q$ with 
$\sfd_{\sfX^\prime}(x_n, x_0) \to 0$, we have $f_0(x_n) \to f_0(x_0)$. Let $f_{n,0} : \sfX_n \mapsto \R$ be such that 
\begin{align*}
\limsup_{n \to \infty} \sup_{K_n^q} |f_{n,0} - \eta_n^q f_0| =0, \quad \forall q \in \mathcal Q.
\end{align*} 
Then
\begin{align*}
 \sup_\sfX (f^* - f_0) \leq \limsup_{n \to \infty} \sup_{\sfX_n} (\overline{f}_n - f_{n,0}).
\end{align*}
\end{lemma}
\begin{proof}
We use arguments similar in the proof of Lemma~\ref{f0f1max}. 
For each $q \in \mathcal Q$ and $\delta >0$, by the defining relation \eqref{Sec4:ftil}, 
there exists $q^\prime \in \mathcal Q$, $x_0 \in K^q\subset K^{q^\prime}$  and $x_n \in K_n^{q^\prime}$ with $\lim_{n \to \infty}\sfd_{\sfX^\prime}\big(\eta^{q^\prime}_n(x_n), x_0\big) =0$ such that
\begin{align*}
 \sup_{K^q} (\tilde{f} -f_0) &< \frac{\delta}{2} +  (\tilde{f} -f_0)(x_0) \\
 &  < \delta + \limsup_{n \to \infty} \big(\overline{f}_n(x_n) - f_{n,0}(x_n)\big) 
 + \limsup_{n \to \infty} \big(f_{n,0}(x_n) - \eta_n^{q^\prime} f_0(x_n)\big) \\
 & \qquad \qquad +  \limsup_{n \to \infty} \big(\eta_n^{q^\prime} f_0(x_n) - f_0(x_0)\big).
\end{align*}
In the above, the second inequality follows from the definition of $\tilde{f}$ in \eqref{Sec4:ftil}.
By assumption on the convergence of $f_{n,0}$ to $f_0$, and by $\sfd_{\sfX^\prime}$-continuity of the $f_0$ in $K^{q^\prime}$, the last two limits on the right hand side above are both zero. In summary,
\begin{align*}
 \sup_{\sfX_0} (\tilde{f} -f_0)= \sup_{q \in \mathcal Q} \sup_{K^q} (\tilde{f} -f_0) \leq  \limsup_{n \to \infty} 
 \sup_{\sfX_n}(\overline{f}_n - f_{n,0}).
\end{align*}

We note that $f_0 \in C((\sfX, \sfd_\sfX))$.
In view of  \eqref{Sec4:fstarD}, by a density argument, 
\begin{align*}
\sup_\sfX ( {f}^* - f_0) = \sup_{\sfX_0} (\tilde{f} - f_0).
\end{align*}

Hence we conclude.
\end{proof}

\begin{condition}\label{Sec4:CnvCAlt}\
\begin{enumerate}
\item  $\overline{f}_n: \sfX_n \mapsto \bar{\R}$ is viscosity solution to \eqref{Hnsub} in the sequential sense;
\item \label{Sec4:dprmC} for each $q \in {\mathcal Q}$, $D(H_0)$ consists of functions $f_0$s which are $\sfd_{\sfX^\prime}$-continuous in $K^q \subset \sfX$:  Namely, for every $x_n, x_0 \in K^q$ such that 
$\lim_{n \to \infty} \sfd_{\sfX^\prime}(x_n, x_0) =0$, 
we have 
\begin{align*}
\lim_{n \to \infty} f_0(x_n) = f_0(x_0);
\end{align*}
\item \label{Sec4:f0C} 
each $f_0 \in D(H_0)$ is $\sfd_\sfX$-continuous in $\sfX$; 
\item   for each $(f_0, g_0) \in H_0$, there exists $(f_{n,0}, g_{n,0}) \in H_{n,0}$ 
satisfying the following
\begin{enumerate}
\item \label{Sec4:fnfAlt} for every $q \in \mathcal Q$,
\begin{align*}
 \lim_{n \to \infty} \sup_{K^q_n} |f_{n,0}- \eta^q_n f_0| =0;
 \end{align*}
\item \label{Sec4:EKpc} there exists $\epsilon_0>0$ and $q \in \mathcal Q$ such that 
\begin{align*}
 E_{\epsilon_0}^+[ \overline{f}_n - f_{n,0} ]
 := \big\{ x \in \sfX_n : \sup_{\sfX_n}(\overline{f}_n - f_{n,0}) 
    \leq \epsilon_0 +  (\overline{f}_n - f_{n,0})(x) \big\} \subset K_n^q.
\end{align*}
\item \label{Sec4:gngCnvAlt}
for each $q \in {\mathcal Q}$ and every $x_n  \in E_{\epsilon_0}^+[ \overline{f}_n - f_{n,0} ] \subset K_n^q$  and $x_0 \in K^q$  with 
$\lim_{n \to \infty} \sfd_{\sfX^\prime}(\eta_n^q(x_n), x_0)=0$, we have
\begin{align*}
\limsup_{n \to \infty}  g_{n,0}(x_n) \leq  g_0(x_0).
\end{align*}
\end{enumerate}
\end{enumerate}
\end{condition}
\begin{remark}
In non-locally compact metric space $(\sfX, \sfd_\sfX)$ situation, we usually need a coercive term in the $\overline{f}_n$, in order to verify Condition~\ref{Sec4:CnvCAlt}.\eqref{Sec4:EKpc}. This sometimes can be achieved through another perturbation to original sub-solution to equation $(I-\alpha H_n) \overline{f}_n \leq h_{n,0}$. See the arguments in Section~\ref{Sec6:Sub}.
\end{remark}

\begin{definition}[Property ${\mathscr P_N}$]
\label{Sec4:ProPN}
A sequence of functions $\{ f_n : \sfX_N  \mapsto \R\}_{n\in\N}$ and $f: \sfX \subset \sfX^\prime \mapsto \R$ is said to satisfy {\it Property ${\mathscr P_N}$} if the following holds: 

For each $q \in {\mathcal Q}$  and every $x_n \in K_n^q$ and $x_0 \in K^q$  with 
$\lim_{n \to \infty} \sfd_{\sfX^\prime}(\eta_n^q(x_n), x_0)=0$, 
we have $\limsup_{n \to \infty} f_n (x_n) \leq f(x_0)$.
\end{definition}

\begin{condition}\label{Sec4:CnvhnAlt}
The sequence $\{ h_{n,0} : n \in \N\}$ and $h_0$ satisfies Property ${\mathscr P_N}$.
\end{condition}

\begin{lemma}\label{Sec4:BPAlt}
Suppose that Conditions~\ref{Sec4:CnvCAlt} and \ref{Sec4:CnvhnAlt} hold.
Then the $f^* \in \USC(\sfX;\R)$ is a sub-solution to \eqref{Hsub} in the point-wise viscosity sense.  

Moreover, the sequence of functions $\{\overline{f}_n : n \in \N\}$ and $f^*$ satisfy the Property ${\mathscr P_N}$ (see Definition~\ref{Sec4:ProPN}) as well.
\end{lemma}
\begin{proof}
We only need to slightly modify arguments in the proof of Lemma~\ref{Sec4:BPThm}.
By the sequential viscosity solution assumption, we can find $\epsilon_n \to 0^+$ with $x_n \in \sfX_n$ such that \eqref{Sec4:BPseqV} still holds. 
By Condition~\ref{Sec4:CnvCAlt}.\ref{Sec4:EKpc}, there exists $q \in \mathcal Q$ with $x_n \in K_n^q$. Moreover, selecting subsequence if necessary, $\lim_{n \to \infty} \sfd_{\sfX^\prime}\big(\eta_n^q(x_n), x_0\big) =0$ for some $x_0 \in K^q \subset \sfX_0$. We claim that
\begin{align*}
 \lim_{n \to \infty} \sup_{\sfX_n}(\overline{f}_n - f_{n,0}) 
  = (f^* - f_0)(x_0) = \sup_\sfX(f^* - f_0).
\end{align*}
This is because that
\begin{align}
(f^*- f_0)(x_0) \leq \sup_\sfX(f^* -f_0) 
& \leq \limsup_{n \to \infty} \sup_{\sfX_n}(\overline{f}_n - f_{n,0}) \nonumber \\
&  \leq \limsup_{n \to \infty} (\overline{f}_n - f_{n,0})(x_n) \label{Sec4:BPff0Alt} \\
& \leq (\tilde{f}- f_0)(x_0) \leq (f^*-f_0)(x_0), \nonumber
\end{align}
where the second inequality follows from Lemma~\ref{Sec4:supEst2}, the third inequality from the first part of \eqref{Sec4:BPseqV}, the fourth inequality from \eqref{Sec4:ftil} and Condition~\ref{Sec4:CnvCAlt}.\ref{Sec4:fnfAlt}, and the last inequality from the definition of $f^*$ in \eqref{Sec4:fstarD}.

From Conditions~\ref{Sec4:CnvCAlt}.\ref{Sec4:dprmC} and 
~\ref{Sec4:CnvCAlt}.\ref{Sec4:fnfAlt}, $\lim_{n \to \infty} f_{n,0}(x_n) = f_0(x_0)$. Therefore, we conclude $\lim_{n \to \infty} \overline{f}_n(x_n) = f^*(x_0)$ from the above sequence of inequalities. Consequently, the second part of \eqref{Sec4:BPseqV} gives
\begin{align*}
 f^*(x_0) \leq \alpha g_0(x_0) + h_0(x_0).
\end{align*}

Finally, to verify Property ${\mathscr P_N}$, let $x_n \in K_n^q, x_0\in K^q$ be such that $\sfd_{\sfX^\prime}(\eta_n(x_n),x_0) \to 0$, then it follows from \eqref{fhat} that
\begin{align*}
\limsup_{n \to \infty} \overline{f}_n(x_n) \leq \hat{f}(x_0) \leq \overline{f}(x_0).
\end{align*}

\end{proof}

\vspace{1cm}

Following arguments similar to the above, there is a super-solution version to Lemma~\ref{Sec4:BPAlt}. We have no use of such result in this paper, hence will not write down details here.


\newpage
 
\section{Hamilton-Jacobi equation in space of empirical probability measures with finite number of point masses}\label{finPartH}
In Section~\ref{hydroIntro}, we described a hydrodynamic limit problem at a formal level.
In this section, we apply the abstract viscosity solution theories developed in Section~\ref{VisMetProj} to that problem and identify corresponding Hamiltonians defined in space of empirical probability measures with a fixed finite number of particles. We also prepare some estimates which will be useful in later sections. The issue of passing particle numbers to infinity will be discussed in detail in Section~\ref{Sec6}.

\subsection{Basic setup}\label{Bsetup}
We recall the setup and notations of Section~\ref{Mfield}.
We endow $\sfY_N:= (\R^d)^N$ with the usual Euclidean metric $\sfd_{Y_N}$, and call this  ``ordered-particle space". Let $\sfG_N$ denote discrete permutation group on $N$-indices, it acts on the $\sfY_N$ through relation
\begin{align*}
\tau {\bf x} := (x_{\tau(1)}, \ldots, x_{\tau(N)}), \quad \forall \tau \in \sfG_N,
{\bf x}:= (x_1, \ldots, x_N) \in \sfY_N.
\end{align*}
It follows then Condition~\ref{condG} is satisfied, and for each $N$ fixed, 
$\sfY_N$ is metrically foliated by $\sfG_N$ with a quotient structure $\sfX_N:= \sfY_N/~\sfG_N$. See Appendix~\ref{App:SubM} and \ref{App:MQ} for details. We denote the corresponding metric quotient space $(\sfX_N, \sfd_{\sfX_N})$.  This $\sfX_N$ models the ``space of un-ordered particles", we identify it with the space of empirical probability measures for $N$-particles \eqref{Sec5:pDef}.
Because of such identification, we denote a typical element in $\sfX_N$ by $\rho:=\rho(dy)$ and introduce a projection map $\sfp_N: \sfY_N \mapsto \sfX_N$ by
 \begin{align}\label{Sec5:pDef}
\rho_N = \sfp_N ({\bf x}) := \frac1N \sum_{i=1}^N \delta_{x_i}, \quad {\bf x}:=(x_1, \ldots, x_N)  \in (\R^d)^N.
\end{align}
Since the group action ${\bf x} \mapsto \tau {\bf x}$ is an isometric one, the above defined $\sfp_N$ is a submetry from $\sfY$ to $\sfX$ (Lemma~\ref{App:GsubMe}). Indeed, because that the $\sfG_N$ is a finite group, $\sfp_N$ is a strong submetry.

We denote $\sfX:= \mathcal P_2(\R^d)$ with the Wasserstein order-2 metric $\sfd$. The identity map ${\Id}_N: \sfX_N \mapsto \sfX$ induces a natural isometric embedding \eqref{isoddN}:
\begin{lemma}
For ${\bf x},{\bf y} \in \sfY_N$, we denote corresponding equivalence classes ${\bf x}^*, {\bf y}^* \in \sfX_N$, and identify them with
\begin{align*}
\rho(dx):= \frac1N \sum_{i=1}^N \delta_{x_i}(dx), 
\quad \gamma(dy):= \frac1N \sum_{i=1}^N \delta_{y_i}(dy).
\end{align*}
Then
\begin{align*}
\sfd_{\sfX_N}^2({\bf x}^*, {\bf y}^*)  := 
\inf_{{\bf x} \in {\bf x}^*, {\bf y} \in {\bf y}^*} \sfd_{Y_N}^2({\bf x}, {\bf y}) 
 =\inf_{\pi \in \sfG_N} \sfd_{Y_N}^2({\bf x}, \pi {\bf y}) 
=  \inf_{  \pi \in \sfG_N } \frac1N \sum_{i=1}^N |x_i -  y_{\pi(i)}|^2
=  \sfd^2(\rho, \gamma).
\end{align*}
\end{lemma}
In particular, one can also see the last identity as a direct consequence of Choquet's theorem on extremal points and Birkhoff theorem on stochastic matrices
(see the end of page 5 in Villani~\cite{Villani03}).


\subsection{Hamilton-Jacobi equations in $\sfY_N$}
We recall that, in the introduction section, a single particle level Hamiltonian function is defined as
\begin{align*}
\sfH:=\sfH(q,p) : \R^d \times \R^d \mapsto \R.
\end{align*}
A corresponding Lagrangian  $\sfL:=\sfL(q,\xi)$ is defined in \eqref{sfLdef} through the Legendre transform. We assume that Condition~\ref{PerCND} is satisfied. 
Following Section~\ref{hydroIntro}, we consider the hydro-dynamically rescaled $N$-particle Hamiltonian function $H_N: (\R^d \times \R^d)^N \mapsto \R$ given by \eqref{defHN}: 
\begin{align*}
H_N({\bf x},{\bf P}) := \frac1N \sum_{i=1}^N \Big( \sfH(\frac{x_i}{\epsilon}, P_i) - U(x_i)
 - \frac1N \sum_{i=1}^N V(\frac{x_i}{\epsilon_N}-\frac{x_j}{\epsilon_N}) \Big).
\end{align*} 
We define differential operator 
\begin{align*}
 {\mathcal H}_N {\mathfrak f}({\bf x}):= H_N\big({\bf x}, \nabla_N {\mathfrak f}({\bf x})\big), 
  \quad \forall {\mathfrak f} \in C^1\big((\R^d)^N\big);
\end{align*}
and consider equation
\begin{align}\label{HNHJB}
 {\mathfrak f}_N - \alpha {\mathcal H}_N {\mathfrak f}_N = {\mathfrak h}_N,
\end{align}
where the ${\mathfrak h}_N \in C\big((\R^d)^N\big)$ and $\alpha>0$ are given.
When $\sup_{\sfX_N} {\mathfrak h}_N <+\infty$, by classical PDE results, a candidate solution to \eqref{HNHJB} is given through the dynamical programming principle by
\begin{align}
{\mathfrak f}_N({\bf x}) & := \sup \Big\{ \int_0^\infty e^{-\frac{s}{\alpha}}
 \Big( \frac{{\mathfrak h}_N}{\alpha} \big({\bf z}(s)\big) - {\mathcal L}_N\big({\bf z}(s), {\bf \dot{z}}(s)\big)\Big)ds :  \label{PartfDef}\\
 & \qquad    {\bf z}(0) ={\bf x},   
  {\bf z}(\cdot) \in AC\big((0, \infty);(\R^d)^N\big) \cap C\big([0,\infty);(\R^d)^N\big) \Big\}, 
  \quad {\bf x} \in (\R^d)^N, \nonumber
\end{align}
where an $N$-particle level Lagrangian function 
${\mathcal L}_N: (\R^d)^N \times (\R^d)^N \mapsto \R\cup\{+\infty\}$ is introduced as
\begin{align}\label{Sec5:LN}
{\mathcal L}_N({\bf x, v})& := \sup_{{\bf P} \in (\R^d)^N} \big(\langle {\bf P, v} \rangle_N - H_N({\bf x, P}) \big) \\
& = \frac1N \sum_{i=1}^N \big( \sfL(\frac{x_i}{\epsilon}, v_i) + U(x_i) 
 + \frac1N\sum_{j=1}^N V(x_i-x_j) \big),  \nonumber
\end{align}
with the $\sfL$ given in \eqref{sfLdef}.
We observe a very rough estimate holds under Conditions~\ref{PerCND}, \ref{U0CND} and \ref{VCND}:
\begin{align*}
{\mathcal L}_N \geq \inf_{\R^d \times \R^d} \sfL + \inf_{\R^d} U + \inf_{\R^d} V \geq - \sup_{q \in \R^d} \sfH(q,0) + \inf_{\R^d} U + \inf_{\R^d} V>-\infty.
\end{align*} 
Hence an upper bound for the ${\mathfrak f}_N$ follows from \eqref{PartfDef}:
\begin{align}\label{Sec5:frafUbd}
 {\mathfrak f}_N \leq \sup_{(\R^d)^N} {\mathfrak h}_N + \alpha \big( \sup_{q \in \R^d} H(q, 0) - \inf_{\R^d} U - \inf_{\R^d} V \big) < +\infty.
\end{align}
Next, we give a rough estimate of the ${\mathfrak f}_N$ from below.
Taking a special path of ``resting" particles ${\bf z}(t) = {\bf x}$ for $t \geq 0$, then 
\begin{align}\label{frakfbl}
 {\mathfrak f}_N({\bf x}) \geq {\mathfrak h}_N({\bf x}) -\alpha {\mathcal L}_N({\bf x},0) 
    \geq {\mathfrak h}_N({\bf x}) + \alpha \inf_{q,p} \sfH(q,p) 
    - \alpha \langle U, \rho \rangle - \alpha \langle V*\rho, \rho \rangle
\end{align}
where $\rho:= \frac1N \sum_{i=1}^N \delta_{x_i}$.

We also note the following invariant property holds
\begin{align*}
{\mathcal L}_N\big(\tau {\bf z}, \dot{(\tau {\bf z})}\big) = {\mathcal L}_N({\bf z}, \dot{\bf z}), \quad \forall \tau \in \sfG_N
\end{align*} 
for every admissible curve ${\bf z}:={\bf z}(\cdot)$ with finite action in the definition of ${\mathfrak f}_N$. 

Given these observations, and in view of Lemmas~\ref{seq2pw}, \ref{Sec3:seq2spw} and their super-solution counterparts, we summarize some well know PDE results regarding  Hamilton-Jacobi equation in Euclidean spaces into the following. 
\begin{lemma}\label{parHJprop}
Suppose that Conditions~\ref{PerCND}, \ref{U0CND} and \ref{VCND} hold; and that ${\mathfrak h}_N \in C\big((\R^{d})^N\big)$ with $\sup {\mathfrak h}_N<\infty$. 
 Then $\mathfrak f_N \in C\big((\R^d)^N\big) \cap \Lip_{\loc}\big((\R^d)^N\big)$ with 
 upper bound \eqref{Sec5:frafUbd}. Moreover, 
\begin{enumerate}
\item  the ${\mathfrak f}_N$ is a point-wise strong viscosity solution to \eqref{HNHJB}. It also satisfies estimate \eqref{frakfbl} from below, hence
\begin{align}\label{Sec5:fNbdB}
 {\mathfrak f}_N (\bx) - {\mathfrak h}_N (\bx) \geq  -\beta\big(\sfd(\rho,\delta_0)\big),
 \quad  \rho = \frac1N \sum_{i=1}^N \delta_{x_i},
\end{align}
for some concave, increasing and sub-linear function $\beta: \R_+ \mapsto \R$.
\item 
The ${\mathfrak f}_N$ satisfies dynamic programming principle
\begin{align}\label{Sec5:DPP}
 {\mathfrak f}_N({\bf x}) & = \sup\big\{ \int_0^t e^{-\frac{s}{\alpha}} \big( \frac{{\mathfrak h}({\bf z}(s))}{\alpha} - {\mathcal L}_N({\bf z}, \dot{\bf z}) \big) ds + e^{-\frac{t}{\alpha}} {\mathfrak f}_N({\bf z}(t)) : \\
 & \qquad \qquad \qquad {\bf z}(\cdot) \in AC\big([0,t];(\R^d)^N\big), {\bf z}(0) ={\bf x}  \big\}. \nonumber
\end{align}
\item If additionally the ${\mathfrak h}_N$ is $\sfG_N$-invariant in that 
$ {\mathfrak h}_N(\tau {\bf x}) = {\mathfrak h}_N({\bf x})$ for every $\tau \in \sfG_N$,
then the ${\mathfrak f}_N$ is $\sfG_N$-invariant as well:
\begin{align}
 {\mathfrak f}_N(\tau {\bf x}) =  {\mathfrak f}_N({\bf x}), \quad \forall \tau \in \sfG_N.
\end{align}
That is, the ${\mathfrak f}_N({\bf x})$ is constant for ${\bf x} \in \sfp^{-1}(\rho)$, for each $\rho \in \sfX_N$ fixed.
\end{enumerate}
\end{lemma}

 \subsection{Submetry projection of Hamiltonians (from configuration spaces $\sfY_N$ to $\sfX_N$) - I, the sub-solution case}\label{Sec5:Sub}
Equation \eqref{HNHJB} is defined in $\sfY_N:=(\R^d)^N$. Next, using the abstract arguments in Section~\ref{PHJPert}, we derive sub-solution to a new equation \eqref{prHN0Eqn} defined in $\sfX_N$. The result is summarized in Lemma~\ref{SubQuotHJ}. 
 
There is a slight abuse of notations between this section and the abstract results in earlier sections. This is because that there are inconsistencies between  established notations in optimal transport theory which we use for the hydrodynamic example, and notations in metric space analysis which we used earlier for abstract development. To establish a clear notational correspondence, we begin with a graphical illustration in the current context. The translation of notations becomes apparent when one compare this graph with the one given in Section~\ref{projHJ}.
\begin{align*}
\begin{tikzpicture}
\draw[black, thick] (-4,0) --(4,0);
\draw (-3,4) -- (-3,0);
\filldraw (-3,1) circle (1pt) node [anchor= east] {${\bf y}_K$} ;
\draw (-2,4) -- (-2,0);
\filldraw (-2,1.7) circle (1pt);
\draw (-1,4) -- (-1,0);
\filldraw (-1,3) circle (1pt) node [anchor=east] {${\bf y}_1$};
\draw (0,4) -- (0,0);
\filldraw (0,2) circle (1pt) node [anchor=east] {${\bf x}$};
\draw (1,4) -- (1,0);
\filldraw (1,2.3) circle (1pt) node [anchor=east] {${\bf y}_k$};
\draw (2,4) -- (2,0);
\filldraw (2,0.5) circle (1pt) node [anchor=east] {${\bf y}_2$};
\draw (3,4) -- (3,0);
\filldraw (3,3.7) circle (1pt);
\node (a)  at (-3,-0.3) {$\gamma_K$};
\node (b) at  (-2, -0.3) {\ldots};
\node (c)  at (-1,-0.3) {$\gamma_1$};
\node (d)  at (0,-0.3) {$\rho$};
\node (e)  at (1,-0.3) {$\gamma_k$};
\node (f)  at (2, -0.3) {$\gamma_2$};
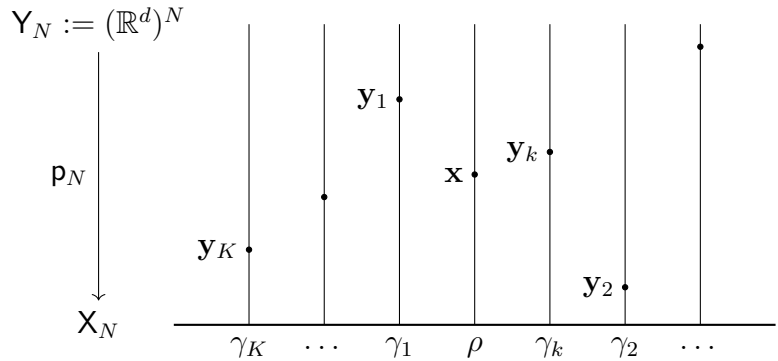
\node (g)  at (3, -0.3) {\ldots};
\node (Y) at (-5, 4.0) {$\sfY_N:= (\R^d)^N$};
\node (X) at (-5, 0.0) {$\sfX_N$};
\draw [->] (Y) -- (X)  node[midway, left] {$\sfp_N$};
\end{tikzpicture}
\end{align*}
In this section, we denote typical elements in $\sfY_N$ by 
\begin{align}\label{bfxy2xy}
{\bf x} := (x_1,\ldots, x_N), \quad {\bf y}_k &:=  (y^k_1, \ldots, y^k_N)\in \sfY_N=(\R^d)^N, \quad k=1,2,\ldots, K,
\end{align}
and typical elements in $\sfX_N$ by
\begin{align}\label{rhogamK}
\rho:=\sfp_N({\bf x}):= \frac1N \sum_{i=1}^N \delta_{x_i}, 
\quad \gamma_k:=\sfp_N({\bf y}_k)= \frac1N \sum_{i=1}^N \delta_{y^k_i}  \in \sfX_N.
\end{align}
We recall the isometric embedding $(\sfX_N, \sfd_{\sfX_N}) \subset (\sfX:={\mathcal P}_2(\R^d), \sfd)$, where $\sfd_{\sfX_N}$ is the natural quotient metric (as abstractly defined in Appendix~\ref{App:MQ}, and explicitly identified in \eqref{dNstar}). The $\sfd$ is order-$2$ Wasserstein metric. 

We will apply Lemma~\ref{sPrelax}. The test functions \eqref{Sec3:fug} used there has three terms, we write them explicitly in current context, one by one, next.

First, as in \eqref{mtranf}, \eqref{mtrang} and \eqref{Adeffrakf}, \eqref{Alexdeff0}, 
we denote $\mathcal S^{\pm}_{\sfX_N}$ 
two classes of simple smooth functions on $\sfX_N$. In particular, 
each $f_0 \in \mathcal S^+_{\sfX_N}$ can be written as
\begin{align}\label{f0rhodef}
f_0(\rho):= f_{0;\gamma_1,\ldots, \gamma_K}(\rho):= \psi \big( \sfd^2(\rho, \gamma_1), \ldots, \sfd^2(\rho, \gamma_K)\big), \quad  \partial_k \psi \geq 0,\forall K \in \N, 
\end{align}
where $\psi \in {\bf \Psi}_K$ (see \eqref{defPsiK} for definition). 
For each such $f_0$, we approximate it with another simple test function on $\sfY_N$ by
\begin{align*}
 {\mathfrak f}_{0; {\bf y}_1, \ldots, {\bf y}_K}({\bf x}) 
  & := \psi\big( \sfd_{\sfY_N}^2({\bf x}, {\bf y}_1), \ldots, \sfd_{\sfY_N}^2({\bf x}, {\bf y}_K) \big) \\
   &= \psi\big( \frac1N \sum_{i=1}^N |x_i-y_i^1|^2, \ldots, \frac1N \sum_{i=1}^N |x_i-y_i^K|^2 \big) \quad \in {\mathcal S}^+_{\sfY_N}.
\end{align*}
The ${\mathfrak f}_{0; \by_1, \ldots, \by_K}$ and $f_0:= f_{0;\gamma_1,\ldots, \gamma_K}$ are related in the same way as the relation between \eqref{Adeffrakf} and \eqref{Alexdeff0} in the abstract setting.

Second, we take
\begin{align}\label{Sec5:fraku}
{\mathfrak u}(\bx) :={\mathfrak u}_\zeta (\bx):= \frac1N \sum_{i=1}^N \zeta(|x_i|^2),
\end{align}
where the $\zeta \in C^2(\R; \R_+)$ is an arbitrary non-negative smooth function in this section.
 
Third, we specify perturbative test functions $\mathfrak g$ as appeared in \eqref{Sec3:fug}. We denote
\begin{align}\label{defF0}
\mathcal{F}_0 :=\Big\{ \phi:=\phi(x,P; q) := \sum_{l=1}^L \alpha_l(x,P) \varphi_l(q) \quad
  \big| \quad 0 \leq \alpha_l \in C^\infty_c(\R^d \times \R^d), \varphi_l \in C^\infty_\per(\R^d) \Big\} . 
\end{align}
For each $\phi:=\phi(x,P; q) \in {\mathcal F}_0$, we introduce 
\begin{align}
 {\mathfrak g}({\bf x})   := {\mathfrak g}_{\phi,\psi;{\bf y_1,\ldots, y_K}}({\bf x})
 & :=  \frac1N \sum_{i=1}^N \phi\big( x_i, P_{i;{\bf y_1, \ldots, y_K}}({\bf x}); \frac{x_i}{\epsilon}\big),  \nonumber \\
 \text{ with }  
 P_{i;{\bf y_1, \ldots, y_K}}({\bf x})&:= \nabla_{N,x_i} {\mathfrak f}_{0;{\bf y_1, \ldots, y_K}} = N \nabla_{x_i}{\mathfrak f}_{0;{\bf y_1, \ldots, y_K}}  \label{PiyDef} \\
 &  = \sum_{k=1}^K 2  \partial_k \psi (\star) (x_i-y_i^k) \in \R^d, \nonumber \\
 \text{ where the }  (\star)& := \big(\sfd_{\sfY_N}^2({\bf x}, {\bf y}_1), \ldots, \sfd_{\sfY_N}^2({\bf x}, {\bf y}_K)\big). \label{Sec5:star}
 \end{align}
The term $P_{i;{\bf y_1, \ldots, y_K}}({\bf x})$ depends on ${\bf y_1,\ldots, y_K}$ as well as on the $\psi$ (which defines the ${\mathfrak f}_{0; {\bf y}_1, \ldots, {\bf y}_K}$). Such perturbed test function does not have the kind of symmetry required by Condition~\ref{gfreeDy}.  Hence we turn to the relaxed version of results developed in Section~\ref{Sec3:relax}.

Finally, we assemble all the above components to get a class of test functions on $\sfY_N$:
\begin{align*}
 {\mathfrak f}_{{\mathfrak u}, \epsilon {\mathfrak g}; {\bf y_1,\ldots, y_K}}  
 := ({\mathfrak f}_{0;{\bf y_1, \ldots, y_K}}+ {\mathfrak u})
  +\epsilon {\mathfrak g}_{\phi, \psi; {\bf y_1,\ldots, y_K}}.
\end{align*}
The $\mathfrak u$ has no dependence on the parameters ${\bf y_1, \ldots, y_K}$. It
has invariance ${\mathfrak u}(\tau {\bf x}) ={\mathfrak u}({\bf x})$ for $\tau \in \sfG_N$. Hence, defining $u : \sfX_N \mapsto \R$ by
\begin{align}\label{Sec5:defu}
u(\rho):= \int_{\R^d} \zeta(|x|^2) \rho(dx),
\end{align}
the function in \eqref{Sec3:projf3} becomes
\begin{align}
& f_{u, \epsilon{\mathfrak g};\gamma_1,\ldots, \gamma_K}(\rho) \nonumber \\
&:= \inf_{   {\bf x} \in \sfp^{-1}_N(\rho) } 
   \inf_{ {\bf y}_1 \in \sfp_N^{-1}(\gamma_1)} \ldots
    \inf_{ {\bf y}_K \in \sfp^{-1}_N(\gamma_K) } 
\big({\mathfrak f}_{0;{\bf y_1, \ldots, y_K}}  + {\mathfrak u}
  +\epsilon {\mathfrak g}_{\phi, \psi; {\bf y_1,\ldots, y_K}}\big) ({\bf x}) \label{deffphi} \\
&= \inf_{   {\bf x} \in \sfp^{-1}_N(\rho) } 
   \inf_{ {\bf y}_1 \in \sfp_N^{-1}(\gamma_1)} \ldots
    \inf_{ {\bf y}_K \in \sfp^{-1}_N(\gamma_K) }       \big({\mathfrak f}_{0;{\bf y_1, \ldots, y_K}}  
  +\epsilon {\mathfrak g}_{\phi, \psi; {\bf y_1,\ldots, y_K}}\big) ({\bf x})   + u(\rho). \nonumber
\end{align}
Consequently, we have a quantitative estimate on what the $f_{u, \epsilon{\mathfrak g};\gamma_1,\ldots, \gamma_K}$ looks like in the $N \to \infty$ asymptotics:
\begin{align*}
\sup_N \sup_{\rho \in \sfX_N} \big| f_{u, \epsilon{\mathfrak g};\gamma_1,\ldots, \gamma_K}(\rho) - \big(f_{0;\gamma_1,\ldots, \gamma_K} + u\big)(\rho)  \big| \leq \epsilon \Vert \phi \Vert_\infty.
\end{align*}

To apply Lemma~\ref{sPrelax}, we also need a good estimate on ${\mathcal H}_N {\mathfrak f}_{{\mathfrak u}, \epsilon {\mathfrak g}; {\bf y_1,\ldots, y_K}}$. 
We denote
\begin{align}\label{Sec5:etaDef}
 \eta^\phi:=\eta^\phi(x,P):= \sup_{q \in \R^d} \sfH\big(q, P+ \nabla_q \phi(x,P;q)\big),
 \quad \forall \phi := \phi(x,P;q) \in {\mathcal F}_0.
\end{align}
Note that, in particular, 
\begin{align*}
\eta(x,P):=\eta^{\phi=0}(x,P) = \eta^{\phi=0}(P) = \sup_{q \in \R^d} \sfH(q,P). 
\end{align*}
We define
\begin{align*}
  \bM^{\bf x; y_1, \ldots, y_K} := \frac1N \sum_{i=1}^N \delta_{(x_i; y_i^1, \ldots, y_i^K)}   \in {\mathcal P}\big((\R^d)^N\big).
\end{align*}
and
\begin{align}\label{Sec5:mu}
 & \bmu^{\bf x;y_1,\ldots, y_K}(dx,dP)   \\
 & \qquad \qquad  := \int_{(y_1, \ldots, y_K) \in (\R^d)^K}  \delta_{\sum_{k=1}^K 2 \partial_k \psi(\star) (x-y_k)}(dP)  \bM^{\bf x;y_1,\ldots, y_K}(dx; dy_1, \ldots, d y_K). \nonumber
\end{align}
That is,
\begin{align*}
 \int_{(x,P) \in \R^{2d}} \varphi(x,P) \bmu^{\bf x;y_1,\ldots, y_K}(dx,dP)
  = \frac1N \sum_{i=1}^N \varphi \big(x_i, P_{i;{\bf y_1,\ldots, y_K}}({\bf x})\big), 
  \quad \forall \varphi \in C_c(\R^d \times \R^d),
\end{align*}
where the $P_{i; {\bf y_1, \ldots, y_K}}({\bf x})$ is defined according to \eqref{PiyDef}.

\begin{lemma}\label{Sec5:HNfest}
\begin{align*}
    {\mathcal H}_N {\mathfrak f}_{{\mathfrak u}, \epsilon {\mathfrak g}; {\bf y_1,\ldots, y_K}} ({\bf x})   
     & \leq \int_{\R^{2d}}  \Big\{ \eta^\phi\big(x,P + 2x \zeta^\prime(|x|^2) \big) -U(x) \\
     & \qquad  \qquad \qquad   
       - \big(V* \rho \big)(x) \Big\} 
    \bmu^{\bf x;y_1,\ldots, y_K}(dx,dP) + O(\epsilon),
 \end{align*}
where the $O(\epsilon)$ satisfies: for each $M>0$ finite, $\psi \in \Psi_K $ and $\phi \in {\mathcal F}_0$ fixed,
\begin{align*}
  \sup_N \sup_{\substack{{\bf x} \in \sfp^{-1}_N(\rho)\\ \int_{\R^d} |x|^2 \rho(dx) \leq M}} \frac{O(\epsilon)}{\epsilon} <+\infty.
\end{align*}
\end{lemma}
\begin{proof}
To simplify notations, we only prove the case $V=0$. The general case only differs slightly notationally.

We identify $\nabla_{N,x_i}   {\mathfrak f}_{\epsilon {\mathfrak g}; {\bf y_1,\ldots, y_K}}$ first. For such purpose, we compute $d \times d$-matrix
\begin{align*}
 D_{x_i} P_{j;{\bf y_1, \ldots, y_K}}({\bf x})  = 2 \sum_{k=1}^K \Big( \partial_k \psi(\star) \delta_{ij} I_{d \times d} + \sum_{l=1}^K \partial^2_{kl} \psi(\star) \frac1N 2 (x_i - y_i^l) \otimes (x_j- y_j^k)  \Big),
\end{align*}
where the shorthand notation $\star$ was defined in \eqref{Sec5:star}. 
We also write
\begin{align*}
 (\nabla_1 \phi)(x,P;q) & := \nabla_x \phi(x,P;q), \qquad (\nabla_2 \phi)(x,P;q) := \nabla_P \phi(x,P;q), \\
 & \qquad \qquad   (\nabla_3 \phi)(x,P;q) := \nabla_q \phi(x,P;q),
\end{align*}
and introduce another shorthand notation
 \begin{align}\label{Sec5:ast}
(\ast_j):= \big(x_j, P_{j;{\bf y_1, \ldots, y_K}}({\bf x}); \frac{x_j}{\epsilon}\big).
\end{align}
Then
\begin{align*}
\nabla_{N,x_i}   {\mathfrak f}_{\epsilon {\mathfrak g}; {\bf y_1,\ldots, y_K}}  
 & := P_{i;{\bf y_1, \ldots, y_K}}({\bf x})
  + \epsilon (\nabla_1 \phi)(\ast_i) 
  + (\nabla_3 \phi)(\ast_i) + \epsilon p_i({\bf x}) \in \R^d,
\end{align*}
where the  
\begin{align*}
p_i({\bf x})& :=   \sum_{j=1}^N \big(D_{x_i} P_{j;{\bf y_1,\ldots, y_K}}({\bf x})\big) (\nabla_2 \phi)(\ast_j) \\
& = 2   \sum_{k=1}^K (\partial_k \psi)(\star) \nabla_2 \phi(\ast_i) 
+ 4   \sum_{k,l=1}^K \partial^2_{kl} \psi(\star)
 \Big( \frac1N \sum_{j=1}^N \big( (x_j - y^k_j) \cdot \nabla_2 \phi (\ast_j)\big)\Big)  (x_i - y_i^l).
\end{align*}

Consequently,
 \begin{align*}
    {\mathcal H}_N {\mathfrak f}_{{\mathfrak u}, \epsilon {\mathfrak g}; {\bf y_1,\ldots, y_K}} ({\bf x})   
& = H_N \big({\bf x}, \nabla_N  {\mathfrak f}_{{\mathfrak u}, \epsilon {\mathfrak g}; {\bf y_1,\ldots, y_K}} ({\bf x})  \big) \\
&  =  \frac1N \sum_{i=1}^N \Big\{ \sfH\big(\frac{x_i}{\epsilon},  \nabla_{N,x_i} {\mathfrak f}_{{\mathfrak u}, {\mathfrak g};{\bf y_1,\ldots, y_K}}   \big) 
   - U(x_i)  \Big\}  \\
   & \leq \frac1N \sum_{i=1}^N \Big\{ 
    \sfH\big(\frac{x_i}{\epsilon}, P_{i;{\bf y_1,\ldots, y_K}}({\bf x}) +2 x_i\zeta^\prime(|x_i|^2)
    + \nabla_q \phi(x_i,P_{i;{\bf y_1,\ldots, y_K}}({\bf x}); \frac{x_i}{\epsilon} ) \big) \\
     & \qquad \qquad \qquad   - U(x_i)   \Big\}  + O(\epsilon) \\
     & \leq \int_{\R^{2d}}  \Big\{ \eta^\phi\big(x,P + 2x \zeta^\prime(|x|^2) \big) - U(x)\Big\} \bmu^{\bf x;y_1,\ldots, y_K}(dx,dP) + O(\epsilon).
 \end{align*}
 \end{proof}

We now define operator $H_{N,0}$ acting on test functions $f_{u, \epsilon{\mathfrak g};\gamma_1,\ldots, \gamma_K}$. Recall that, in Definition~\ref{MultOpt} and equation \eqref{bnu0Def}, we respectively introduced the notions of optimal multi-plans $\Gamma^{\opt}(\rho;\gamma_1, \ldots, \gamma_K)$ and the measure $\bnu^{\bM}_{f_{0;\gamma_1, \ldots, \gamma_K}}$ for a given test function $f_{0;\gamma_1, \ldots, \gamma_K}$. Next, we define $\delta$-approximate versions of both concepts. We write
\begin{align*}
 \Gamma^{\opt}_\delta (\rho;\gamma_1, \ldots, \gamma_K)& :=\Big\{ {\bf M} := {\bf M}(dx;dy_1, \ldots, dy_K) \in {\mathcal P}_2(\R^{(1+K) d}) \text{ such that } \\
 & \qquad \qquad \pi^{1+k}_\# {\bf M} = \gamma_k, k=1,\ldots, K; \text{ and } \\
 & \qquad \qquad \qquad
  \Big( \int_{\R^d \times \R^d} |x-y_k|^2 \big(\pi^{1,1+k}_\#{\bf M}\big)(dx,dy_k) \Big)^{1/2}
    \leq \sfd(\rho, \gamma_k) +\delta\Big\}.
\end{align*}
We still denote $\bnu^{\bM}_{f_{0;\gamma_1, \ldots, \gamma_K}}$ in the same way as \eqref{bnu0Def}, but with the above $\bM$.
 With these notations, we define
\begin{align}\label{HN0def}
  H_{N,0} f_{u, \epsilon{\mathfrak g};\gamma_1,\ldots, \gamma_K}(\rho) 
  & := \sup_{\bM \in \Gamma^{\opt}_\delta(\rho;\gamma_1,\ldots, \gamma_K)}  
 \int_{\R^{2d}} \Big\{  \eta^\phi\big(x,   P + 2 x \zeta^\prime(|x|^2) \big) \\
 & \qquad \qquad \qquad   - U(x) - \big(V* \rho\big)(x)
  \Big\}  \bnu^{\bM}_{f_{0;\gamma_1,\ldots, \gamma_K}}(dx,dP) + O(\epsilon),\nonumber 
 \end{align}
 where the $O(\epsilon)$ term is the same as in Lemma~\ref{Sec5:HNfest}. 

\begin{remark}  
In particular, considering a special case of the test functions $f_{0;\gamma_1,\ldots, \gamma_K}$ with
\begin{align*}
 K=1, \qquad \psi(r):= \frac{\alpha}{2} r, \quad \alpha >0,
\end{align*}
then have above expression reduces to
\begin{align}\label{Sec5:HN0f}
H_{N,0} f_{u, \epsilon{\mathfrak g}; \gamma_1}(\rho) & = \sup_{\bpi \in \Gamma^{\opt}_\delta(\rho,\gamma_1)} 
  \int_{\R^{2d}} \Big\{  \sup_{q \in \R^d}
  \eta^{\phi} \big(q, 2 \alpha(x-y) + 2 x\zeta^\prime(|x|^2) \big) \\
  & \qquad \qquad - U(x) - \big(V*\rho \big)(x) \Big\} \bpi(dx, dy) + O(\epsilon). \nonumber
\end{align}
 \end{remark}
 
 \begin{lemma}\label{Sec5:relxEst}
For each $\epsilon>0$ and $f_{u, \epsilon {\mathfrak g};\gamma_1,\ldots, \gamma_K}$, there exists a $\delta:=\delta(\epsilon; \Vert \mathfrak g\Vert_\infty, \psi)>0$ such that  
\footnote{Note that the $\psi$ is the one appearing in definition of $f_{0;\gamma_1,\ldots, \gamma_K}$. }
\begin{align}\label{Sec5:deleps}
 \lim_{\epsilon \to 0^+} \delta(\epsilon; \Vert \mathfrak g\Vert_\infty, \psi) =0,
\end{align}
and that  
\begin{align*}
    \sup_{({\bf x; y_1,\ldots, y_K}) \in \sfS^\delta(\rho;\gamma_1, \ldots, \gamma_K) } {\mathcal H}_N {\mathfrak f}_{{\mathfrak u}, \epsilon {\mathfrak g}; {\bf y_1,\ldots, y_K}} ({\bf x})   
  \leq  H_{N,0} f_{u, \epsilon {\mathfrak g};\gamma_1,\ldots, \gamma_K}(\rho).
\end{align*}
Recall that the notion of a $\delta$-section $\sfS^\delta$ is defined in \eqref{sfSeDef}.
\end{lemma}
\begin{proof}
The existence of $\delta>0$ satisfying \eqref{Sec5:deleps} follows from Lemma~\ref{Sec3:EcS}. For the selected $\delta>0$,  when $({\bf x; y_1,\ldots, y_K}) \in \sfS^\delta(\rho;\gamma_1, \ldots, \gamma_K)$, we have (by definitions)
\begin{align*}
  \sup_{k=1}^K |\sfd_{\sfY_N}({\bf x}, {\bf y}_k) - \sfd(\rho, \gamma_k)| <  \delta,  
\text{ and } \bM^{\bf x; y_1,\ldots, y_K} \in \Gamma^{\opt}_\delta(\rho;\gamma_1, \ldots, \gamma_K),
\end{align*}
 and 
\begin{align}\label{munueq}
  \bmu^{\bf x;y_1,\ldots, y_K}(dx,dP) = \bnu^{\bM^{\bf x;y_1,\ldots, y_K}}_{f_{0;\gamma_1,\ldots, \gamma_K}}(dx,dP)
\end{align}
where the left hand side notation refers to the one defined by \eqref{Sec5:mu} and the right hand side refers to the one defined by \eqref{bnu0Def}.

Therefore, the conclusion is just a re-statement of the result in Lemma~\ref{Sec5:HNfest}.
\end{proof}

\begin{lemma}\label{SubQuotHJ}
In the context of Lemma~\ref{parHJprop}, additionally assume that   
\begin{align}\label{Sec5:fhtx}
 {\mathfrak h}_N(\tau {\bf x}) ={\mathfrak h}_N({\bf x}), \quad \forall \tau \in \sfG_N. 
\end{align}
That is,
\begin{align*}
 h_{N,0}(\rho)  := \sup_{\substack{ (x_1,\ldots, x_N) \in (\R^d)^N \\
 \text{such that }  \rho= \frac1N \sum_{i=1}^N \delta_{x_i}  } } {\mathfrak h}_N(x_1, \ldots, x_N) 
 = \inf_{\substack{ (x_1,\ldots, x_N) \in (\R^d)^N \\
 \text{such that }  \rho= \frac1N \sum_{i=1}^N \delta_{x_i}  } }
 {\mathfrak h}_N(x_1, \ldots, x_N) .
\end{align*}
Then 
\begin{enumerate}
\item  the function ${\mathfrak f}_N$ defined in \eqref{PartfDef} is bounded from above and continuous.
It is also $\sfG_N$-invariant
\begin{align*}
\overline{f}_N(\rho) :=  \sup_{\substack{ (x_1,\ldots, x_N) \in (\R^d)^N \\
 \text{such that }  \rho= \frac1N \sum_{i=1}^N \delta_{x_i}  } } {\mathfrak f}_N(x_1, \ldots, x_N) 
 = \inf_{\substack{ (x_1,\ldots, x_N) \in (\R^d)^N \\
 \text{such that }  \rho= \frac1N \sum_{i=1}^N \delta_{x_i}  } }
 {\mathfrak f}_N(x_1, \ldots, x_N);
\end{align*}
and is a point-wise strong viscosity solution to \eqref{HNHJB}.
\item  $\overline{f}_N \in C(\sfX_N)$ is bounded from above and is a strong point-wise viscosity sub-solution to 
\begin{align}\label{prHN0Eqn}  
  ( I - \alpha H_{N,0} ) \overline{f}_N \leq h_{N,0}.
\end{align}
In particular, the $\delta$ appearing in  
$\Gamma^{\opt}_\delta(\rho;\gamma_1,\ldots, \gamma_K)$ 
in \eqref{HN0def} has the property \eqref{Sec5:deleps}. 
\end{enumerate}
\end{lemma}
\begin{proof}The first part of conclusion follows from Lemma~\ref{parHJprop}. The second part follows from the estimate in Lemma~\ref{Sec5:relxEst} applied to the abstract results in Lemma~\ref{sPrelax}.
\end{proof}

\subsection{Submetry-projection of Hamiltonians - II, the super-solution case}
The main result of this section is Lemma~\ref{SupQuotHJ}.  
As in the sub-solution case, the following diagram translates notations in this sub-section into those in Section~\ref{PHJPert} in a graphical way.  
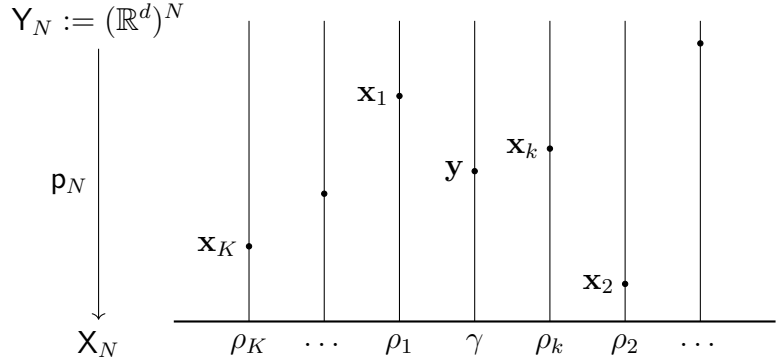
\begin{align*}
\begin{tikzpicture}
\draw[black, thick] (-4,0) --(4,0);
\draw (-3,4) -- (-3,0);
\filldraw (-3,1) circle (1pt) node [anchor= east] {${\bf x}_K$} ;
\draw (-2,4) -- (-2,0);
\filldraw (-2,1.7) circle (1pt);
\draw (-1,4) -- (-1,0);
\filldraw (-1,3) circle (1pt) node [anchor=east] {${\bf x}_1$};
\draw (0,4) -- (0,0);
\filldraw (0,2) circle (1pt) node [anchor=east] {${\bf y}$};
\draw (1,4) -- (1,0);
\filldraw (1,2.3) circle (1pt) node [anchor=east] {${\bf x}_k$};
\draw (2,4) -- (2,0);
\filldraw (2,0.5) circle (1pt) node [anchor=east] {${\bf x}_2$};
\draw (3,4) -- (3,0);
\filldraw (3,3.7) circle (1pt);
\node (a)  at (-3,-0.3) {$\rho_K$};
\node (b) at  (-2, -0.3) {\ldots};
\node (c)  at (-1,-0.3) {$\rho_1$};
\node (d)  at (0,-0.3) {$\gamma$};
\node (e)  at (1,-0.3) {$\rho_k$};
\node (f)  at (2, -0.3) {$\rho_2$};
\node (g)  at (3, -0.3) {\ldots};
\node (Y) at (-5, 4.0) {$\sfY_N:= (\R^d)^N$};
\node (X) at (-5, -0.3) {$\sfX_N$};
\draw [->] (Y) -- (X)  node[midway, left] {$\sfp_N$};
\end{tikzpicture}
\end{align*}

We denote
\begin{align*}
{\bf y} :=(y_1, \ldots, y_N), \quad {\bf x}_k:=  (x^k_1, \ldots, x^k_N)\in  (\R^d)^N, \quad k=1,2,\ldots, K,
\end{align*}
with $\rho_k :=\sfp_N({\bf x}_k)$ and $\gamma:= \sfp_N({\bf y})$.
We consider test function $f_1 \in {\mathcal S}^-_{\sfX_N}$ written as
\begin{align}\label{f1gamdef}
f_{1;\rho_1, \ldots, \rho_K}(\gamma):= - \psi \big( \sfd^2(\gamma, \rho_1), \ldots, \sfd^2(\gamma, \rho_K)\big),
\end{align}
and its counterpart defined on $\sfY_N$  
\begin{align*}
  {\mathfrak f}_{1;   \bx_1, \ldots, \bx_K}(\by) 
:= -\psi\big( \sfd_{\sfY_N}^2(\by,   \bx_1), \ldots, \sfd_{\sfY_N}^2(\by, \bx_K) \big).
\end{align*}
We also introduce the $\mathfrak u$ in \eqref{Sec5:fraku}, as in the sub-solution case.

For each $\phi :=\phi(y,P; q)\in  {\mathcal F}_0$ (recall definition in \eqref{defF0}),  we define 
\begin{align}\label{Sec5:g1}
{\mathfrak g}({\by}) := {\mathfrak g}_{\phi,\psi;{\bx_1}, \ldots, {\bx_K}}({\bf y})  
= \frac{1}{N} \sum_{i=1}^N  \phi(y_i, P_{i; \bx_1, \ldots, \bx_K}({\by}); 
 \frac{y_i}{\epsilon});
\end{align}
where the
\begin{align*}
P_{i; \bx_1, \ldots, \bx_K}({\bf y})  
 := \nabla_{N, y_i} {\mathfrak f}_{1; \bx_1, \ldots, \bx_K}
= N \nabla_{y_i} {\mathfrak f}_{1; \bx_1, \ldots, \bx_K}.
\end{align*}
We also take
\begin{align}\label{Sec5:uu1}
 {\mathfrak u}(\by) := {\mathfrak u}_{-\zeta}(\by) := \frac1N \sum_{i=1}^N (-\zeta)(y_i).
\end{align}
where the $\zeta \in C^2(\R; \R_+)$ is bounded from below.

%

We now consider perturbed test function on $\sfY_N$:
\begin{align}\label{Sec5:ff1}
 {\mathfrak f}_{{\mathfrak u}, \epsilon {\mathfrak g}; \bx_1, \ldots, \bx_K} 
 :=  ( {\mathfrak f}_{1; \bx_1,\ldots, \bx_K} + {\mathfrak u}) 
  + \epsilon {\mathfrak g}_{\phi,\psi; \bx_1,\ldots, x_K}.
 \end{align}
As in the sub-solution case, its counterpart on $\sfX_N$ is
\begin{align}\label{Sec5:f1}
 f_{u,\epsilon{\mathfrak g};\rho_1, \ldots, \rho_K}(\gamma) 
 := \sup_{\by \in \sfp^{-1}_N(\gamma)} \sup_{\bx_1 \in \sfp_N^{-1}(\rho_1)} 
  \ldots \sup_{\bx_K \in \sfp^{-1}_N(\rho_K)} \big({\mathfrak f}_{1;\bx_1, \ldots, \bx_K}(\by) + \epsilon {\mathfrak g}_{\phi,\psi;\bx_1, \ldots, \bx_K}(\by)\big) + u(\gamma). 
\end{align}
We have
\begin{align}\label{Sec5:fcnv}
 \sup_N \sup_{\gamma \in \sfX_N} |f_{u,\epsilon{\mathfrak g};\rho_1, \ldots, \rho_K}(\gamma) - (f_{1;\rho_1,\ldots, \rho_K}+u) (\gamma) | \leq \epsilon \Vert \phi \Vert_\infty.
\end{align}

For each given $\bM  \in \Gamma^{\opt}_\delta(\gamma;\rho_1, \ldots, \rho_K)$ and $\psi$ (defining the $f_1$ and appearing in the expression of   the $\beta_k$s in \eqref{betaDef}),
we denote
\begin{align}\label{bnu0f1def}
 \bnu^{\bM}_{f_{1;\rho_1, \ldots, \rho_K}}(dy, dP) :=   \int_{(x_1,\ldots, x_K) \in \R^{Kd}} 
  \delta_{\sum_{k=1}^K  \beta_k 2 (x_k-y)} (dP)\bM(dy, dx_1,\ldots, dx_K);
 \end{align}
and define
  \begin{align}\label{HN1def}
H_{N,1} f_{u, \epsilon{\mathfrak g}; \rho_1, \ldots, \rho_K}(\gamma) 
 & := \inf_{ \bM \in \Gamma^{\opt}_\delta(\gamma; \rho_1,\ldots, \rho_K)} 
   \int_{\R^{2d}} \Big\{\eta_\phi \big(y, P- 2 y \zeta^\prime(|y|^2)\big) \\
   & \qquad \qquad - U(y) - \big(V*\gamma\big) (y) \Big\}  \bnu^{\bf M}_{f_{1;\rho_1, \ldots, \rho_K}}(dy, dP) + O(\epsilon),   \nonumber
\end{align}
where the 
\begin{align}\label{Sec5:etainf}
\eta_\phi:=\eta_\phi(y,P):= \inf_{q \in \R^d}\sfH\big(q, P + \nabla_q \phi(y,P;q)\big).
\end{align}
The $\delta$ is chosen to satisfy Lemma~\ref{Sec3:EcS2} applied in such context. In particular, it satisfies \eqref{Sec5:deleps}.  Like in the sub-solution case, the $O(\epsilon)$ is meant to satisfy the estimate in Lemma~\ref{Sec5:HNfest}.
 
 \begin{lemma}\label{SupQuotHJ}
 In the context of Lemma~\ref{parHJprop}, assume
that Condition~\ref{PerCND} holds and that ${\mathfrak h}_N$ 
satisfies \eqref{Sec5:fhtx}. Consequently,
\begin{align}\label{Sec5:h1Def}
 h_{N,1}(\rho) & := \inf_{\{(x_1,\ldots, x_N) : \rho= \frac1N \sum_{i=1}^N \delta_{x_i}\}} 
 {\mathfrak h}_N(x_1, \ldots, x_N) \\
 & = \sup_{\{(x_1,\ldots, x_N) : \rho= \frac1N \sum_{i=1}^N \delta_{x_i}\}} {\mathfrak h}_N(x_1, \ldots, x_N). \nonumber
\end{align}
Then 
\begin{enumerate}
\item  the solution ${\mathfrak f}_N({\bf x})$ for \eqref{HNHJB} defined in \eqref{PartfDef} is $\sfG_N$ invariant
\begin{align}\label{Sec5:lfNDef}
\underline{f}_N(\rho) &:= \inf_{\{(x_1,\ldots, x_N) : \rho= \frac1N \sum_{i=1}^N \delta_{x_i}\}} {\mathfrak f}_N(x_1, \ldots, x_N) \\
 &=  \sup_{\{(x_1,\ldots, x_N) : \rho= \frac1N \sum_{i=1}^N \delta_{x_i}\}} {\mathfrak f}_N(x_1, \ldots, x_N). \nonumber
\end{align}
\item The $\underline{f}_N \in C(\sfX_N)$ with growth estimate 
\begin{align*}
 \underline{f}_N(\rho) \geq h_{N,1}(\rho) -\beta \circ \sfd(\rho, \delta_0) 
\end{align*}
for some concave, increasing and sub-linear function $\beta: \R_+ \mapsto \R$.  
There is a choice of the $O(\epsilon)$ in \eqref{HN1def} such that the $\underline{f}_N$ is a point-wise strong viscosity super-solution to 
 \begin{align}\label{prHN1Eqn}
  ( I - \alpha H_{N,1} ) \underline{f}_N \geq h_{N,1}.
\end{align}
\end{enumerate}
\end{lemma}
\begin{proof}
The proof follows from symmetric arguments as in the sub-solution case in Lemma~\ref{SubQuotHJ}. The growth estimate comes from \eqref{Sec5:fNbdB}.
\end{proof}
 
 \begin{remark}\label{Sec5:HN1Rmk}
 Note that in the special case of $K=1$, $\psi(r) = \frac{\alpha}{2} r$ with $\alpha>0$, 
 \begin{align*}
H_{N,1} f_{u, \epsilon{\mathfrak g};\rho_1}(\gamma) 
  = \inf_{ \bpi \in \Gamma^{\opt}_\delta(\rho_1,\gamma)} 
   \int_{\R^{2d}} \Big\{ \eta_\phi \big(y, \alpha(x-y)- 2y\zeta^\prime(|y|^2) \big)  
    - U(y) - V*\gamma(y) \Big\}  \bpi(dx, dy).  
\end{align*}
 \end{remark}

  \subsection{Uniform modulus of continuity estimate}
Let ${\mathfrak f}_N : (\R^d)^N \mapsto \R$ be defined according to \eqref{PartfDef}. We suppose that 
 the ${\mathfrak h}_N \in C\big((\R^d)^N\big)$ always has 
 the invariance property in \eqref{Sec5:fhtx}.
 Then, by Lemma~\ref{parHJprop}, the ${\mathfrak f}_N$ is $\sfG_N$-invariant as well. Hence it can be identified with a function in $\sfX_N:= (\R^d)^N/\sfG_N$. We note that, 
on one hand, $\sfX_N$ is a finite dimensional space; on the other, it can be identified with space of empirical probability measures with $N$ unit point masses. 
Denoting a typical element in $\sfX_N$ using empirical probability measure $\rho$, we write 
\begin{align*}
f_N(\rho):= {\mathfrak f}_N({\bf x}), \qquad \forall \rho:= \frac1N \sum_{i=1}^N \delta_{x_i}.
\end{align*}

Next, we provide an estimate regarding modulus of continuity for $f_N$.    
We begin with a technical lemma.
\begin{lemma}\label{LemMod}
Let $C_\delta>0$ for each $\delta>0$. We define $\omega(r):= \inf_{\delta>0} ( \delta + C_\delta r)$. Then such defined $\omega$ is a concave modulus in the sense that $\omega \in C(\R_+ ; \R_+)$ is non-decreasing, with $\omega(0) =0$ and $r \mapsto \omega(r)$ is concave.
\end{lemma}

We now state the main result of this subsection. In the following, $\sfd$ is the $2$-Wasserstein metric.
  \begin{lemma}\label{Sec5:fnmest}
 Suppose that $h_N$ is uniformly (in $N$) bounded from below in $\sfd$-balls of finite radius:
  \begin{align*}
  \inf_N \inf_{\substack{\sigma \in\sfX_N \\
    \sfd(\sigma, \delta_0)\leq R}} h_N (\sigma)> -\infty,
   \quad \forall R \in \R_+.
\end{align*}
Then for each $R>0$, there exists a modulus $\omega_R \in C(\R_+; \R_+)$ such that
\begin{align*}
f_N (\rho) - f_N(\gamma) \leq \omega_R (\sfd(\rho,\gamma)), \quad \quad
\forall \rho, \gamma \in \sfX_N
\text{ and } \sfd(\rho, \delta_0) + \sfd(\gamma, \delta_0) \leq R, 
\end{align*}
holds uniformly for all $N \in \N$.
 \end{lemma}
 \begin{proof}Again, to save space and notation, we only prove the case $V=0$.
 
For every $\rho, \gamma \in \sfX_N$, there exists 
\begin{align*}
\bpi:=\bpi_N(dx,dy) := \frac1N \sum_{i=1}^N \delta_{(x_i, y_i) \in \R^d \times \R^d} (dx, dy) 
\in {\mathcal P}_2(\R^{2d})
\end{align*}
  such that $\sfd^2(\rho, \gamma) = \int_{\R^d \times \R^d} |x-y|^2 \bpi(dx,dy)$. For any $\delta >0$,
  we define 
\begin{align*}
z_i(t):= x_i + t \frac{y_i-x_i}{\delta} \in C^1([0,\delta] ;\R^d),  \text{ and } 
\sigma(t): = \frac1N \sum_{i=1}^N \delta_{z_i(t)} \in AC([0,\delta];\sfX_N).
\end{align*}
Such $\sigma(0)=\rho$, $\sigma(\delta) =\gamma$ and the curve $t \mapsto \sigma(t)$ has constant speed with velocity
\begin{align*}
 \bnu(t):=\bnu_N(t;dx,dv) := \frac1N \sum_{i=1}^N \delta_{z_i(t), \frac{y_i-x_i}{\delta}}(dx, dv) = \dot{\sigma}(t) 
  \in \Tan_{\rho(t)} {\sfX_N}, \quad \forall t \in [0,\delta].
\end{align*}
Therefore, $\sfd(\sigma(t), \sigma(0)) = (t/\delta) \sfd(\rho, \gamma) \leq R$
 for $t \in [0,\delta]$.

According to the dynamical programming principle identity \eqref{Sec5:DPP} and in view of \eqref{Sec5:LN}, the following holds
\begin{align*}
 f_N(\rho)&  \geq \int_0^\delta \alpha^{-1} e^{-\frac{r}{\alpha}} h_N(\sigma(r)) dr 
   -\int_0^\delta  e^{-\frac{r}{\alpha}} \Big( \int_{\R^{2d}} \big( \sfL(\frac{x}{\epsilon}, v) + U(x)  \big) \bnu(r;dx,dv) \Big) dr \\
   & \qquad \qquad + e^{- \frac{\delta }{\alpha}} f_N(\gamma).
\end{align*}
Therefore, 
\begin{align*}
f_N(\gamma) - f_N(\rho) & \leq (1- e^{-\frac{\delta}{\alpha}})  f_N (\gamma)
+ \int_0^\delta    \Big( \int_{\R^{2d}} \big( \sfL(\frac{x}{\epsilon},v) + U(x) \big)\bnu(r;dx,dv) \Big) dr \\
& \qquad +( e^{-\frac{\delta}{\alpha}}-1)\inf_{\substack{\sigma \in\sfX_N \\
    \sfd(\sigma, \delta_0)\leq 2R}} h_N (\sigma) \\
& \leq    (1- e^{-\frac{\delta}{\alpha}}) \sup_{\sfX_N} f_N 
+ C_0 \delta  + \frac{C_1}{\delta} \sfd^2(\rho, \gamma) + C_2\delta \big(1 +
\sfd(\rho, \delta_0)+  \sfd(\rho,\gamma) \big) \\
& \qquad +( e^{-\frac{\delta}{\alpha}}-1) 
\inf_{\substack{\sigma \in\sfX_N \\
    \sfd(\sigma, \delta_0)\leq 2R}} h_N (\sigma).
\end{align*}
For the above estimates,  we used the facts that (recall \eqref{Sec1:sfLlb})  
 \begin{align*}
 \sfL(q,v) \leq C_0+ C_1 |v|^2 \quad  \text{ for some } C_1 >0, C_0 \in \R ;
\end{align*}
and that (by Condition~\ref{U0CND})
\begin{align*}
 \langle U, \sigma(r) \rangle \leq C_2 \int_{\R^d} \big(1+|z|) \sigma(r,dz) 
    \leq C_2 \big(1+  \int_{\R^d} |x| \rho(dx) + \frac{t}{\delta} \sfd(\rho,\gamma)\big),
     \quad r \in [0, \delta]
\end{align*}
for some $C_2 >0$.

 The conclusion now follows from Lemma~\ref{LemMod}.
 \end{proof}

\subsection{Submetry projection of Hamiltonians - III, revisiting the sub-solution case.}\label{Sec5:RevSub}
For reasons which will be clear when we develop limit theorems in Section~\ref{Sec6:Sub}, we need to generalize the result of Lemma~\ref{SubQuotHJ}. We still take $\sfY_N:= (\R^d)^N$, but with a new metric $\sfd_{\sfY_N}$ corresponding to $p$-norm ($p \in (1,\infty)$) when the $(\R^d)^N$ is viewed as a Banach space. Specifically, instead of using
\begin{align*}
 \sfd_{\sfY_N}(\bx, \by) := \big(\frac1N \sum_{i=1}^N |x_i - y_i|^2\big)^{1/2}, 
\end{align*}
we now use a new one  
\begin{align*}
 \sfd_{\sfY_N}(\bx, \by):= \big(\frac1N \sum_{i=1}^N |x_i-y_i|^p\big)^{1/p}.
\end{align*}
The corresponding quotient (with respect the to permutation group $\sfG_N$) space $\sfX_N$ is still the space of empirical probability measures for $N$ equally weighted particles. However, the quotient metric now can be identified with the $p$-Wasserstein metric $\sfd_p$ as follows: 
\begin{align*}
 \sfd_{\sfX_N}^p (\rho, \gamma) 
  := \inf_{\pi \in \sfG_N} \frac1N \sum_{i=1}^N |x_i - y_{\pi(i)}|^p
   = \inf_{\bm \in \Gamma(\rho, \gamma)} \int_{\R^d \times \R^d} |x-y|^p \bm(dx, dy) =: (\sfd_p)^p(\rho, \gamma)
\end{align*}
for 
\begin{align*}
 \rho(dx) := \frac1N \sum_{i=1}^N \delta_{x_i}(dx),
     \quad \gamma(dy):=\frac1N \sum_{i=1}^N \delta_{y_i}(dy).
\end{align*}

Since all the arguments are completely in parallel with those in Section~\ref{Sec5:Sub}, we only highlight differing details. 

\subsubsection{A $p$-Wasserstein version of the sub-solution Lemma~\ref{SubQuotHJ}} We revisit the arguments in Section~\ref{Sec5:Sub}. 

First of all, we replace the $2$-Wasserstein metric $\sfd$ that were used everywhere, with the $p$-Wasserstein metric $\sfd_p$. In particular, the $f_0 \in {\mathcal S}_{\sfX_N}^+$ in \eqref{f0rhodef} now becomes 
\begin{align}\label{Sec5:f0Alt}
 f_0(\rho):= f_{0;\gamma_1, \ldots, \gamma_K}(\rho):=\psi\big( \sfd_p^2(\rho,\gamma_1), \ldots, \sfd_p^2(\rho, \gamma_K) \big);
\end{align}
and 
\begin{align*}
{\mathfrak f}_{0;{\bf y_1,\ldots, y_K}}({\bf x}) :=
 \psi\big( \frac1N  (\sum_{i=1}^N |x_i-y_i^1|^p )^{\frac{2}{p}}, 
  \ldots, \frac1N  (\sum_{i=1}^N |x_i-y_i^K|^p )^{\frac{2}{p}} \big).
\end{align*}
The new version of $P_{i;\by_1, \ldots, \by_K}$ in \eqref{PiyDef} becomes 
~\footnote{We take convention $|z|^{p-1} \frac{z}{|z|} =0$ when $z=0$, $p>1$.} 
\begin{align}\label{Sec5:newP}
P_{i;\by_1, \ldots, \by_K}(\bx) 
 & := N \nabla_{x_i} {\mathfrak f}_{0;{\bf y_1,\ldots, y_K}} \\
& =   \sum_{k=1}^K  2 
 \Big( |x_i - y_i^k|^{p-1} \frac{x_i - y_i^k}{|x_i - y_i^k|}\Big)
 \sfd_p^{2-p}(\rho,\gamma_k)
 \partial_k \psi\big(\sfd_p^2(\rho,\gamma_1), \ldots, \sfd_p^2(\rho, \gamma_K)\big) .
 \nonumber
\end{align}
We also introduce $p$-Wasserstein version of the collection of optimal transport measure
\begin{align*}
\Gamma^{\opt}_p(\rho; \gamma)&:=  \Big\{ \bmu:=\bmu(dx;dy)\in {\mathcal P}_p(\R^{2d}) \text{ such that }  \pi^1_\# \bmu =\rho, \pi^2_\# \bmu= \gamma \\
& \qquad \qquad \text{ and } \int_{\R^d \times \R^d} |x-y|^p  \bmu(dx,dy) 
   =\sfd_p^p(\rho, \gamma), \Big\};
\end{align*}
and its multi-marginal analogue: 
\begin{align}\label{pMultOpt}
 \Gamma_p^{\opt}(\rho;\gamma_1, \ldots, \gamma_K) & :=
  \Big\{ \bM:=\bM(dx;dy_1,\ldots, dy_K)\in {\mathcal P}_p(\R^{(1+K)d}) \text{ such that } \\
   & \qquad \qquad 
   \pi_\#^{1,1+k} \bM \in \Gamma^{\opt}_p(\rho, \gamma_k)  \Big\}. \nonumber 
\end{align}
For the above $f_0$ and $\bM \in \Gamma_p^{\opt}(\rho;\gamma_1, \ldots, \gamma_K)$, we now extend the measure in \eqref{bnu0Def} to the $p$-Wasserstein setting by
\begin{align}\label{Sec5:p-nu}
&  \bnu_{f_0}^{\bM}(dx, dP) \\
& \qquad := \int_{(y_1, \ldots, y_k) \in \R^{Kd}}
  \delta_{\sum_{k=1}^K  2 \sfd_p^{2-p}(\rho,\gamma_k) 
  |x-y_k|^{p-1}\frac{x-y_k}{|x-y_k|}\alpha_k}(dP) \bM(dx;dy_1, \ldots, dy_K), \nonumber
\end{align}
with  
\begin{align*}
\alpha_k:=\alpha_k(\rho;\gamma_1, \ldots, \gamma_K): = \partial_k \psi\big( \sfd^2_p(\rho, \gamma_1), \ldots, \sfd^2_p(\rho, \gamma_K) \big).
\end{align*}

Secondly, we introduce a $p$-Wasserstein version of the operator $H_{N,0}$ in \eqref{HN0def}. By replacing the $2$-Wasserstein distance $\sfd$ by the $p$-Wasserstein version $\sfd_p$, we define a counterpart for the $\Gamma_\delta^{\opt}$ which we denote $\Gamma_{\delta,p}^{\opt}$. For the $f_0$ in \eqref{Sec5:f0Alt} and any given $\phi :=\phi(x,P;q) \in {\mathcal F}_0$ in \eqref{defF0}, we define a perturbative function $\mathfrak g$ just as in \eqref{PiyDef} but with 
the new $P_{i;\bf y_1, \ldots, y_K}$s as given by \eqref{Sec5:newP}, and with the $\star$ in \eqref{Sec5:star} replaced by squares of $p$-Wasserstein distance functions. For each $\zeta \in C^2$, we define $u$ in \eqref{Sec5:defu} and then $f_{u,\epsilon {\mathfrak g}; \gamma_1, \ldots, \gamma_K}$ according to \eqref{deffphi}, and
\begin{align}\label{p-HN0def}
  H_{N,0} f_{u, \epsilon {\mathfrak g}; \gamma_1, \ldots, \gamma_K}(\rho)& :=  
   \sup_{\bM \in \Gamma^{\opt}_{\delta, p}(\rho; \gamma_1, \ldots, \gamma_K)}
  \int_{\R^{2d}} \Big\{ \eta^\phi\big(x, P+ 2x \zeta^\prime(|x|^2) \big) \\
  & \qquad \qquad \qquad
   - U(x) - \big(V*\rho\big)(x) \Big\}\bnu^{\bM}_{f_{0;\gamma_1,\ldots, \gamma_K}}(dx, dP) + O (\epsilon).  
   \nonumber
\end{align}

\begin{lemma}\label{Sec5:relxEstA}
With the above notational changes, the statements in Lemma~\ref{SubQuotHJ} still hold when the space of empirical probability measures with $N$ equal mass particles, still denoted  $\sfX_N$, is identified as a closed sub-space of the $p$-Wasserstein space with metric $\sfd_p$, $p > 1$.
\end{lemma}

\subsubsection{A perturbative version of the $p$-Wasserstein formulation}
In this subsection, we establish a perturbative variant of Lemma~\ref{Sec5:relxEstA}.
See Remark~\ref{Sec5:PerFbar} for necessity of considering such perturbation. 

Let $\overline{f}_N(\rho)$ be defined according to Lemma~\ref{SubQuotHJ}, and
$\theta>0$. We take $\zeta(r):= \theta (r/2)$  in \eqref{Sec5:defu}, hence 
$u(\rho):=u_\theta(\rho):= \int_{\R^d} \zeta(|x|^2) \rho(dx)$. We denote
\begin{align}\label{Sec5:barfNthe}
 \overline{f}_{N, \theta}(\rho):=  \overline{f}_N(\rho) -   u_\theta(\rho).
\end{align}
With reference to the $f_{u,\epsilon {\mathfrak g}; \gamma_1, \ldots, \gamma_K}$
given by \eqref{deffphi}, we write 
\begin{align}\label{Sec5:perTestf0}
 f_{0,\epsilon {\mathfrak g}; \gamma_1, \ldots, \gamma_K}
 :=f_{u=0,\epsilon {\mathfrak g}; \gamma_1, \ldots, \gamma_K}.
\end{align}
Next, recall the constants $c, C$ in Condition~\ref{Sec1:Tone2} for the $\sfH$, we define 
\begin{align}\label{Sec5:HN0the}
 H_{N,0}^\theta  f_{0, \epsilon {\mathfrak g}; \gamma_1, \ldots, \gamma_K}(\rho) 
 &:=  H_{N,0} f_{0, \epsilon {\mathfrak g}; \gamma_1, \ldots, \gamma_K}(\rho)\\
 &   \qquad + (1-\frac{1}{\lambda})\Big\{ c- \inf{\sfH}
   + 4C  \sup_{q,P,x}|\nabla_q \phi|^2  \nonumber \\
   & \qquad \qquad + 4C
  \sup_{\bM \in \Gamma^{\opt}_{\delta, p}(\rho; \gamma_1, \ldots, \gamma_K)}
    \int_{\R^{2d}} |P|^2 \bnu^{\bM}_{f_{0;\gamma_1,\ldots, \gamma_K}}(dx, dP)\Big\} 
      \nonumber \\
 & \qquad \qquad \qquad
   + \frac{2 C \lambda}{\lambda -1} \theta^2 \int |x|^2 \rho(dx). \nonumber
\end{align}
The $H_{N,0}^\theta f_{0, \epsilon {\mathfrak g}; \gamma_1, \ldots, \gamma_K}$ should be viewed as multi-valued, with a free varying parameter $\lambda>1$. 

Note that for fixed finite $N$, $(\sfX_N, \sfd_{p>1})$ and $(\sfX_N, \sfd_{p=2})$ are topologically and metrically equivalent. Moreover, the above constructed $\overline{f}_N$ and $\overline{f}_{N,\theta}$ are independent of the $p > 1$. 
$D(H_{N,0}^\theta)$ consists of test functions of the form \eqref{Sec5:perTestf0} with {\em every } $p \in (1,2)$.

\begin{lemma}\label{Sec5:SubQuoA}
Under the assumptions of Lemma~\ref{SubQuotHJ}, for each $\theta >0$,
the above $\overline{f}_{N,\theta} \in C(\sfX_N)$ is bounded from above and is a strong viscosity sub-solution in the point-wise sense to 
\begin{align}\label{Sec5:HNtheta}
 ( I - \alpha H_{N,0}^\theta) \overline{f}_{N,\theta} \leq h_{N,0}.
\end{align}
\end{lemma}
 \begin{proof}
The proof is based upon one observation: We may consider $\overline{f}_N$ as a viscosity sub-solution with $u(\cdot)$ as part of the test functions for the Hamiltonian operator $H_{N,0}$ in \eqref{p-HN0def}; we may also consider $\overline{f}_{N,\theta}$ as a viscosity sub-solution with another Hamiltonian operator $H_{N,0}^\theta$. With proper error estimates, the first one implies the second. 
 
We establish some estimates regarding the $\eta^\phi$ in \eqref{Sec5:etaDef} in Lemma~\ref{Sec5:etaphi} next. In particular, estimate \eqref{Sec5:etaP} implies that  
\begin{align*}
  H_{N,0} f_{u, \epsilon {\mathfrak g}; \gamma_1, \ldots, \gamma_K}(\rho) 
 & \leq  H_{N,0}^\theta  f_{0, \epsilon {\mathfrak g}; \gamma_1, \ldots, \gamma_K}(\rho).
\end{align*}
Here, the $f_{u,\epsilon {\mathfrak g}; \gamma_1, \ldots, \gamma_K}
 =f_{0,\epsilon {\mathfrak g}; \gamma_1, \ldots, \gamma_K}  + u$.
Again, the right hand side above means a multi-valued function with free varying parameter $\lambda>1$.

The conclusion now follows from Lemma~\ref{Sec5:relxEstA}.
\end{proof}

We establish some estimates which will be useful in Section~\ref{Sec6:Sub}. One of them was also used in the proof of previous lemma.
\begin{lemma}\label{Sec5:etaphi}
Suppose Conditions~\ref{PerCND} and~\ref{Sec1:Tone2} hold. For each $\phi \in \mathcal F_0$ fixed,  
\begin{enumerate}
\item  the map in \eqref{Sec5:etaDef} is continuous $\eta^\phi:=\eta^\phi(x,P) \in C(\R^{2d})$.
\item  there exists a finite constant $C_\phi>0$ such that
\begin{align*}
 \nabla_q \phi(x,P;q) =0, \quad \forall |x| > C_\phi, P,q \in \R^d.
\end{align*}
Hence
\begin{align*}
 \eta^\phi(x,P) = \eta^{\phi=0}(P), \quad \forall P \in \R^d, \text{ whenver } |x| > C_\phi.
\end{align*}
By \eqref{amQforH}, the $\eta^{\phi}$ has no more than quadratic growth in $P$ at infinity uniformly in $x$:
\begin{align}\label{Sec5:etapest}
  \eta^{\phi}(x,P) \leq c_{\phi}(1+ |P|^2),\quad \exists c_{\phi}>0.
\end{align}
\item for each $0<p<2$ fixed, the map (see definition of $\bnu_{f_0}^{\bM}$ in \eqref{Sec5:p-nu})
\begin{align*}
{\mathcal P}_p(\R^{(1+K)d}) \ni \bM \mapsto \int_{\R^{2d}}
  \eta^{\phi}(x, P) \bnu_{f_0}^{\bM}(dx,dP), 
\end{align*}
is continuous in the topology given by $p$-Wasserstein metric on ${\mathcal P}_p(\R^{(1+K)d})$.
\item  let $c,C>0$ be the constants in \eqref{amQforH}, and $L_\phi$ be defined as in 
\begin{align}\label{Sec5:LphiDef}
 L_\phi:= \sup_{q,P} \sup_{x \neq y} 
  \frac{| \nabla_q \phi(y,P;q) - \nabla_q\phi(x,P;q)|}{|x-y|}.
\end{align}
Then for every $\lambda>1$,
\begin{align}\label{Sec5:etaDiff}
\lambda \eta^\phi(y,\frac{P}{\lambda}) - \eta^\phi(x, P) 
\leq c(\lambda -1) + \frac{CL_\phi^2}{\lambda-1}  |x-y|^2.
\end{align}
\item  we have
\begin{align}\label{Sec5:etaP}
\eta^\phi(x, P+ \xi) & \leq \eta^\phi(x,P) +  (1-\frac{1}{\lambda})
  \Big(c - \inf{\sfH} + 4C\sup_{q,P,x}|\nabla_q \phi|^2+ 4C|P|^2 \Big) \\
  & \qquad \qquad \qquad + \frac{2 C\lambda }{\lambda-1}|\xi|^2. \nonumber
\end{align}
 \end{enumerate} 
\end{lemma}
\begin{proof}
 We begin by recalling
 $\eta^\phi(x,P) = \sup_{q \in \T^d} \sfH(q,P+ \nabla_q \phi(x,P;q))$, 
 and that $(q, x,P) \mapsto \sfH(q, P + \nabla_q \phi(x,P;q) )$ is continuous. 
 Continuity of the $\eta^\phi$ follows from Lemma~\ref{infsupSC} in Appendix:
First, by first part of that lemma, we have $\eta^\phi \in \LSC(\R^{2d})$. 
 Second, by compactness of $\T^d$, the second part of the lemma implies $\eta^\phi \in \USC(\R^{2d})$. 

From convexity of $p \mapsto \sfH(q,p)$ and by Condition~\ref{Sec1:Tone2}, we have
\begin{align} \label{Sec5:Hlam}
 \lambda \sfH(q,\frac{p}{\lambda}) - \sfH(q, p^\prime) 
   \leq (\lambda-1) \sfH\big(q, \frac{p-p^\prime}{\lambda-1}\big) 
      \leq c(\lambda -1) + \frac{C}{\lambda-1} |p-p^\prime|^2.   
\end{align}
Therefore, by definition of $\eta^\phi$ in \eqref{Sec5:etaDef}, 
\begin{align*}
\lambda \eta^\phi(y,\frac{P}{\lambda}) - \eta^\phi(x, P) & \leq  
\sup_{q \in \R^d} \Big( \lambda 
 \sfH\big(q,\frac{P+ \nabla_q \phi(y,P;q)}{\lambda}\big) 
             - \sfH\big(q, P+ \nabla_q\phi(x,P;q)\big) \Big) \\
             & \leq c(\lambda -1) + \frac{C}{\lambda-1} 
  |\nabla_q \phi(y,P;q) - \nabla_q\phi(x,P;q)|^2,
\end{align*}
giving \eqref{Sec5:etaDiff}.

Given $p^{\prime \prime} \in \R^d$, take $p=\lambda p^{\prime \prime}$ in \eqref{Sec5:Hlam}, using the at most quadratic growth of $p \mapsto \sfH$ requirement in Condition~\eqref{amQforH}, we also get 
\begin{align*}
\sfH(q,p^{\prime \prime}) & \leq \frac{1}{\lambda} \sfH(q, p^\prime) + c(1-\frac{1}{\lambda}) + \frac{C}{\lambda(\lambda-1)} \Big(\lambda |p^{\prime \prime}-p^\prime|
+ (\lambda -1) |p^\prime| \Big)^2 \\
& \leq  \Big( \sfH(q, p^\prime) - (1-\frac{1}{\lambda}) \inf \sfH\Big) +  c (1-\frac{1}{\lambda}) + \frac{2C\lambda}{\lambda -1}
|p^{\prime \prime}-p^\prime|^2 + 
 \frac{2C(\lambda-1)}{\lambda} |p^\prime|^2  \\
  & = \sfH(q, p^\prime) + (1-\frac{1}{\lambda})
  \Big(c - \inf{\sfH} + 2C |p^\prime|^2\Big) 
   +   \frac{2 C\lambda }{\lambda-1}|p^{\prime \prime}-p^\prime|^2.
\end{align*}
Therefore, 
\begin{align*}
\eta^\phi(x,P+\xi) -\eta^\phi(x,P) 
&\leq \sup_{q \in \R^d} \Big( \sfH\big(q,P+\xi+ \nabla_q\phi(x,P;q)\big) 
    -  \sfH\big(q,P+ \nabla_q\phi(x,P;q)\big)\Big) \\
    & \leq (1-\frac{1}{\lambda})
  \Big(c - \inf{\sfH} + 4C |P|^2 + 4C \sup_{q,P,x}|\nabla_q \phi|^2\Big) 
   +   \frac{2 C\lambda }{\lambda-1}|\xi|^2.
\end{align*}
\end{proof}
  
 \begin{remark}\label{Sec5:PerFbar}
 In Section~\ref{CnvHJ}, we extended the Barles-Perthame half relaxed limit theory to metric space settings. There are two versions of such extension: a ``simpler-minded" version was given as Theorem~\ref{Sec4:BPThm}, a more subtler version was also developed in Section~\ref{BPlimitA}.    

We would like to apply these abstract results to establish limiting behaviors for the sub- super- solutions from $H_{N,0}$ and the $H_{N,1}$ respectively. In the case of $H_{N,0}$,  we have a sub-solution result from Lemma~\ref{Sec5:relxEstA} (namely, the $\overline{f}_N$ in Lemma~\ref{SubQuotHJ} is a viscosity sub-solution to \eqref{prHN0Eqn}). However, such result is not compatible with the ``simpler-minded" metric space version of the half-relaxed limit theory in Theorem~\ref{Sec4:BPThm}, for the purpose of deriving limit. This is because that the conditions required are not satisfied. To apply the subtler version, we need to consider the perturbed problem in Lemma~\ref{Sec5:SubQuoA}. For more detailed explanation, see opening paragraph
in Section~\ref{Sec6:Sub}.

We note that, because of the $u_\theta$ (with $\theta>0$) term, 
the $\overline{f}_{N,\theta}$ has a compact sub-levels in $p$-Wasserstein space with $0<p<2$. Moreover, test functions \eqref{Sec5:perTestf0} are continuous in 
${\mathcal P}_p(\R^d)$. These properties will play important roles in Section~\ref{Sec6:Sub}, when we derive limit properties for sub-solutions of 
 \eqref{Sec5:HNtheta} as $N \to \infty$. 
 
%
\end{remark}

\newpage 
 
\section{Convergence of Hamiltonians in the hydrodynamic limit for infinite particles}\label{Sec6}
In this section, we apply the abstract viscosity solution theories developed in Section \ref{CnvHJ} and the explicit estimates of finite particle Hamiltonians in Section~\ref{finPartH} to our hydrodynamic limit problem. We introduce a pair of Hamiltonian operators $H_0, H_1$ and show that they are respectively upper- and lower- limits of the $H_{N,0}$ and $H_{N,1}$ in \eqref{HN0def}, \eqref{HN1def} in proper senses. We only discuss convergence of the Hamiltonian operators in this section. We leave for later sections about comparison principle and issues on convergence of solutions to associated Hamilton-Jacobi equations.
 
 \subsection{Convergence of Hamiltonians, the super-solution case} 
Unlike other parts of the paper, we discuss the super-solution case first. This is because that, in this case, Condition~\ref{HConSup}.\ref{tre1} in the abstract viscosity convergence Theorem~\ref{Sec4:BPThm} can be readily verified through constructing relatively simple test functions. The case of sub-solution does not follow by symmetric arguments. More complicated arguments are needed to verify the counterpart Condition~\ref{HConSub}.\ref{tre}. Hence we delay its developments until the super-solution case is cleared.

Recall that we assumed both $U,V$ have sub-linear growth at infinity (Conditions~\ref{U0CND} and \ref{VCND}).

\subsubsection{Convergence of spaces}
Let metric spaces $(\sfX_N, \sfd_{\sfX_N})$ and $(\sfX, \sfd)$ be as in Section~\ref{Bsetup}. That is, $\sfX:= \mathcal P_2(\R^d)$ and $\sfd$ are the order-2 Wasserstein space and metric respectively, and the $\sfX_N$ is identified as space of empirical measures for $N$ number of points with equal mass.  
We denote $\eta_N:=\Id$ the identity map that embeds $\sfX_N$ into $\sfX$. Namely, the $\eta_N$ maps empirical probability measure to itself identified as a probability measure with finite second moment.  We introduce index set  
\begin{align}\label{indQ}
{\mathcal Q}:= \{ q= (\zeta , M) : \zeta \in C^1(\R; \R_+)  \text{ such that }   
 \lim_{C \to +\infty} \inf_{|r|\geq C}\frac{\zeta(r)}{1+|r|} =+\infty, 
 M \in \R_+ \}.
\end{align}
This $\mathcal Q$ induces a family of compact subsets in $\sfX$   (e.g. Proposition 7.1.5 in \cite{AGS08})
\begin{align*}
 K^q:=K^{\zeta, M}:= \big\{ \rho \in \sfX :   \int_{\R^d} \zeta(|x|^2) \rho(dx) \leq M \big\}, \quad  q:=(\zeta, M) \in {\mathcal Q}, 
\end{align*}
and similarly, a family of compact subset in $\sfX_N$ 
\begin{align*}
 K^q_N:=K^{\zeta, M}_N:= \big\{ \rho \in \sfX_N :   \int_{\R^d} \zeta(|x|^2) \rho(dx) \leq M \big\}, \quad  q:=(\varphi, M) \in {\mathcal Q}.
\end{align*}
We choose $\eta^q_N := \eta_N\big|_{K_N^q} : K_N^q \mapsto K^q$.

We recall that, by a uniform integrability characterization of compact set of Wasserstein order-2 space (e.g. Proposition 7.1.5 in \cite{AGS08}),
 every $K^{\zeta,M}$ is a  compact sets in $\sfX$. Moreover,  by proper choice of the $\zeta$,  every compact set $K \subset\subset \sfX$ can be contained in one of the $K^{\zeta,M}$s.  
By a density argument of empirical probability measures in space of probability measures, we have
\begin{align*}
 \lim_{N \to \infty} \sfd_{\rm GH} (K_N^{\zeta, M}, K^{\zeta, M}) =0.
\end{align*} 
Therefore, the following holds. 
\begin{lemma}  $\sfX_N$ converges to $\sfX$ in the sense of generalized Gromov-Hausdorff convergence with respect to index set $\mathcal Q$:
\begin{align*}
 (\sfX_N, \sfd_{\sfX_N}) \stackrel{\rm gGH}{\longrightarrow}_{\mathcal Q}   (\sfX, \sfd).
\end{align*}
\end{lemma}
 
In Section~\ref{finPartH}, we established results projecting sub- and super-solutions of Hamilton-Jacobi equations in ordered-particle space $\sfY_N:=(\R^d)^N$ to un-ordered particle space $\sfX_N:= \sfY_N/\sfG_N$. We note that definition of the projected Hamiltonian operators $H_{N,1}f_{\phi,\epsilon}$ in \eqref{HN1def} involves a parameter $\phi:= \phi(y, P; q) \in {\mathcal F}_0$ (defined in \eqref{defF0}). This $\phi$ records information regarding highly oscillating  microstructures of the dynamic through the variable $q$. It also records dependency between highly oscillating structure and slowly oscillating macrostructures of the dynamic through the $(y,P)$ variable. In the $N \to \infty$ limit, this parameter $\phi$ disappears in the limiting test function $f_1$. However, it remains in the multi-valued $H_1 f_1$ as an extra index. Then, through the defining inequality property of viscosity super-solution, we can optimize over such $\phi$ to tighten up the estimates, giving a variational structure of the limiting effective Hamiltonian. While implementing this procedure, there are subtle technical twists, we develop these details next.

\subsubsection{Limiting Hamiltonian operator}
We consider a special class of the $\zeta$s as appeared in \eqref{indQ}. Following \eqref{Sec5:defu}, the $\zeta \in C^1(\R; \R_+)$ needs to have super-linear growth at infinity. Indeed, we require something even more:  the $\zeta$ is non-decreasing and 
\begin{align}\label{Sec6:zpriGro}
 \liminf_{r \to +\infty} \zeta^\prime (r)=+\infty.
\end{align}
We write
\begin{align*}
u(\gamma):=u_{-\zeta}(\gamma) := - \int_{\R^d} \zeta(|y|^2) \gamma(dy), \quad \forall \gamma \in \sfX.
\end{align*}
Such class of $\zeta$s is large enough to have the following property: 
for each $\gamma_0 \in \sfX$, we can find a $\zeta$ with the above property and satisfying $\int_{\R^d} \zeta(|y|^2) d\gamma_0(dy) <+\infty$.

Let $\rho_k \in \sfX$ for $k =1,2,\ldots, K$. Following \eqref{f1gamdef}, we write a class of simple test functions
\begin{align*}
f_1(\gamma):= f_{1;\rho_1,\ldots, \rho_K}(\gamma) = - \psi\big(\sfdist^2_{\rho_1}(\gamma), \ldots, \sfdist^2_{\rho_K}(\gamma) \big) \in {\mathcal S}^-,
\end{align*}
and  perturbed test functions
\begin{align}\label{Sec6:furho}
 f_{u; \rho_1, \ldots, \rho_K}:=  f_{1;\rho_1,\ldots, \rho_K} + u_{-\zeta}.
\end{align}
We denote $\widehat{\mathcal S}^-$ the collection of perturbed test functions defined in the last line.
With $\eta_\phi$ introduced in \eqref{Sec5:etainf}:  
\begin{align*}
 \eta_\phi:= \eta_\phi(y,P) := \inf_{q \in \R^d} \sfH\big( q, P+ \nabla_q \phi(y,P;q)\big), 
 \quad \forall \phi:=\phi(y,P;q) \in {\mathcal F}_0,
\end{align*}
we write
\begin{align}\label{Sec6:GinfDef}
 G_{f_{u;\rho_1, \ldots, \rho_K}}^\phi (\gamma) 
 &:= \inf_{\bM \in \Gamma^{\opt}(\gamma;\rho_1, \ldots, \rho_K)} 
 \int_{\R^{2d}}   \eta_\phi \big( y, P -2y\zeta^\prime(|y|^2) \big) \bnu^{\bM}_{f_{1;\rho_1,\ldots, \rho_K}}(dy, dP) \\
 & \qquad \qquad - \langle U, \gamma\rangle - \langle V*\gamma, \gamma\rangle, \nonumber
\end{align}
where the  $\bnu^{\bM}_{f_{1;\rho_1,\ldots, \rho_K}}$ is defined just as in \eqref{Sec2:bnu1Def}, and the notation $\Gamma^{\opt}(\gamma;\rho_1, \ldots, \rho_K)$ defined in \eqref{MultOpt}. 

We now define an operator $H_1$, identified through its graph, by
\begin{align}\label{Sec6:H1Def}
H_1:=\Big\{ (f_{u;\rho_1, \ldots, \rho_K}, G_{f_{u;\rho_1, \ldots, \rho_K}}^\phi) : f_{u,\rho_1,\ldots, \rho_K} \in \widehat{\mathcal S}^-, \phi \in {\mathcal F}_0 \Big\}.
\end{align}

\begin{remark}\label{Sec6:H1Rmk}
Expression for the above operator becomes more explicit when it acts on a special class of the test functions. Namely, for those $f_{u;\rho_1}=(f_{1,\rho_1} + u_{-\zeta}) \in D(H_1)$ with
 \begin{align*}
 f_{1;\rho_1}(\gamma) =-\frac{\alpha}{2} \sfd^2(\rho_1, \gamma), \quad \alpha>0, \rho_1 \in \sfX,
\end{align*}
we have
\begin{align}\label{Sec6:GAlt}
 G_{f_{u; \rho_1}}^\phi(\gamma) = \inf_{\bpi \in \Gamma^{\opt}(\rho_1;\gamma)} 
 \int_{\R^{2d}}  \eta_\phi \big( y, \alpha(x-y) -2y\zeta^\prime(|y|^2) \big)   \bpi(dx,dy)
  - \langle (U+ V*\gamma), \gamma\rangle.
\end{align}
It is easier to see this by following a notational convention which we practice throughout this paper. We always write $d\gamma:=\gamma(dy)$ and $d\rho=\rho(dx)$ to associate the $x$ with the $\rho$, and $y$ with the $\gamma$.
\end{remark}

Next, in two steps, we show that this $H_1$ is a lower limit to the $H_{N,1}$s in \eqref{HN1def} in the sense as required by Condition~\ref{HConSup}. 

\begin{lemma}\label{Sec6:eta}
  Assuming Condition~\ref{PerCND}, then for each $\phi \in \mathcal F_0$, 
\begin{enumerate}
\item  the map $\eta_\phi:=\eta_\phi(y,P) \in C(\R^{2d})$.
\item  there exists a finite constant $C_\phi>0$ such that
\begin{align*}
 \nabla_q \phi(y,P;q) =0, \quad \forall |y| > C_\phi, P,q \in \R^d.
\end{align*}
Hence
\begin{align*}
 \eta_\phi(y,P) = \eta_{\phi=0}(P), \quad \forall P \in \R^d, \text{ whenver } |y| > C_\phi.
\end{align*}
By \eqref{amQforH}, the $\eta_{\phi}$ has exactly quadratic growth in $P$ at infinity uniformly in $y$:
\begin{align}\label{Sec6:etaQ}
 - c_\phi + C_\phi^{-1} |P|^2 \leq \eta_{\phi}(y,P) \leq c_{\phi}+ C_\phi |P|^2,\quad \exists c_{\phi}, C_\phi>0.
\end{align}
 \end{enumerate} 
\end{lemma}
\begin{proof}
We prove the continuity of $\eta_\phi$ only. The other properties follow directly.
We note that $(q,y,P) \mapsto \sfH\big(q, P+ \nabla_q \phi(y,P;q)\big)$ is continuous. By second part of Lemma~\ref{infsupSC}, $\eta_\phi \in \USC(\R^{2d})$.  Next, by compactness of $\T^d$ 
(noting periodicity of $q \mapsto \phi(y,P;q)$ and $q \mapsto \sfH(q,p)$),  the first part of Lemma~\ref{infsupSC} gives $\eta_\phi \in \LSC(\R^{2d})$.
\end{proof}

\begin{lemma}\label{Sec6:Kflev}
Suppose that Conditions~\ref{PerCND}, ~\ref{U0CND}, ~\ref{VCND} and~\ref{Sec1:Tone2} hold.
 Then  function $\gamma \mapsto G_{f_{u; \rho_1, \ldots, \rho_K}}^\phi(\gamma)$ is lower semicontinuous in $(\sfX, \sfd)$. Indeed, every finite sub-level set of the function is compact in $(\sfX, \sfd)$. 
\end{lemma}
\begin{proof}
From the estimate \eqref{Sec6:etaQ},  $(y,P) \mapsto \eta_\phi \big( y; P- 2y \zeta^\prime(|y|^2) \big)$ is bounded below by constant $-c_\phi$;  moreover, 
from \eqref{Sec6:zpriGro}, the dominating term of growth estimate from above is 
$|P|^2+ |y|^2 (\zeta^\prime)^2(|y|^2)$.  Consequently, the dominating term in $G_{f_{u;\rho_1, \ldots, \rho_K}}^\phi(\gamma)$, as it blows up, is $\gamma \mapsto \langle |y|^2 \big(\zeta^\prime (|y|^2)\big)^2, \gamma \rangle$. Therefore, for every finite $L>0$,
\begin{align*}
 \{ \gamma \in \sfX : G_{f_{u; \rho_1, \ldots, \rho_K}}^\phi(\gamma) \leq L \}
\end{align*}
 is a relatively compact subset in $(\sfX, \sfd)$. 
 
 Next, we verify that $G_{f_{u; \rho_1, \ldots, \rho_K}}^\phi  \in \LSC(\sfX;\R \cup\{+\infty\})$ to conclude. To simplify, we only treat the case of special test functions in Remark~\ref{Sec6:H1Rmk}. Proof in the general situation only requires notational changes. We recall the expression of $G_{f_{u; \rho_1}}^\phi$ in \eqref{Sec6:GAlt}. 
Let $\gamma_n \to \gamma_0$ in $\sfd$ and choose a $\bpi_n \in \Gamma^{\opt}(\rho_1;\gamma_n)$ to be such that
\begin{align*}
 G_{f_{u; \rho_1}}^\phi(\gamma_n) \geq \int_{\R^{2d}} \eta_\phi \big(x,   \alpha(x-y)  -2y\zeta^\prime(|y|^2)   \big) \bpi_n(dx,dy) -\langle U, \gamma_n\rangle - \langle V*\gamma_n, \gamma_n\rangle - \frac1n .
\end{align*}
Then by a variant of tightness argument,
\footnote{The proof of Lemma~\ref{Sec6:H1SupUlt} will use such argument again. See there for more explanation.}
 at least along subsequence, $\bpi_n \to \bpi_0$ for some 
$\bpi_0 \in \Gamma^{\opt}(\rho_1;\gamma_0)$ in order-$2$ Wasserstein metric in $\mathcal P_2(\R^{2d})$. Hence, by a version of the Fatou's lemma, 
 \begin{align*}
 \liminf_{n \to \infty}  G_{f_{u; \rho_1}}^\phi(\gamma_n) & \geq  
           \int_{\R^{2d}}  \eta_\phi \big(y,   \alpha(x-y)  -2y\zeta^\prime(|y|^2)  \big) \bpi_0(dx,dy) - \langle U, \gamma_0\rangle - \langle V*\gamma_0, \gamma_0\rangle  \\
           &  \geq G_{f_{u; \rho_1}}^\phi(\gamma_0) .
\end{align*}
We conclude.
 \end{proof}

\begin{lemma}\label{Sec6:SupCnv}
Suppose that Conditions~\ref{PerCND}, ~\ref{U0CND} and~\ref{Sec1:Tone2} hold.
Then the $H_{N,1}$s in \eqref{HN1def} and $H_1$ satisfy Condition~\ref{HConSup}.   
\end{lemma}
\begin{proof}
Again, we only  write out details for the case when $f_{u;\rho_1} \in D(H_1)$ is the special test function in Remark~\ref{Sec6:H1Rmk}. Proof for general case only requires notational changes. 
  
First, we approximate the $f_{u;\rho_1}$. For the given $\rho_1$, we construct empirical probability measures $\rho_1^N:= \frac1N \sum_{i=1}^N \delta_{x_i^N}$ satisfying
 \begin{align*}
 \lim_{N \to \infty} \sfd(\rho_1^N, \rho_1) =0,
\end{align*}
with the $x_i^N \in \R^d$.  We denote $\bx_1^N:= (x_1^N, \ldots, x_N^N)$.
As in \eqref{Sec5:g1} -\eqref{Sec5:f1}, we introduce 
\begin{align*}
{\mathfrak f}_{1;\bx_1^N}(\by):= -\frac{\alpha}{2} \frac1N \sum_{i=1}^N |y_i - x_i^N|^2,  
\end{align*}
and ${\mathfrak u}:={\mathfrak u}_{-\zeta}$. For every $\phi :=\phi(y,P;q) \in {\mathcal F}_0$, we define $\mathfrak g:= {\mathfrak g}_{\phi; \bx_1^N}$  and 
\begin{align*}
{\mathfrak f}_{{\mathfrak u}, \epsilon {\mathfrak g}; \bx_1^N}
 := ({\mathfrak f}_{1;\bx_1^N} + {\mathfrak u}) + \epsilon {\mathfrak g},
\end{align*}
and 
\begin{align*}
f_N(\gamma) := \sup\Big\{ {\mathfrak f}_{{\mathfrak u}, \epsilon {\mathfrak g}; \bx_1^N}(\by) :   \by =(y_1,\ldots, y_N)\in (\R^d)^N  \text{such that } 
 \gamma= \frac1N\sum_{i=1}^N \delta_{y_i} \Big\},
\end{align*}
Then according to \eqref{Sec5:fcnv}, as $\epsilon:= \epsilon_N \to 0$, we have~\footnote{ Note that $\eta_N : \sfX_N \mapsto \sfX$ is the identity embedding map, we may simply ignore it in the  expression.}
\begin{align*}
 \lim_{N \to \infty} \sup_{K_N^q} | f_N - \eta_N f_{u;\rho_1}  |   = 0, \quad \forall q \in {\mathcal Q}.
\end{align*}
Since $f_{u;\rho_1} \in \USC(\sfX; \bar{\R})$, we have $(-f_N) 
 \stackrel{\rm  \Gamma}{\longrightarrow}_{\mathcal Q} (-f_{u;\rho_1})$.

Second, we show that, for each $q \in \mathcal Q$ fixed and $\gamma_N \in K_N^q$, $\gamma_0 \in K^q$ such that $\sfd(\gamma_N, \gamma_0) \to 0$, we have
 \begin{align*}
 \liminf_{N \to \infty}   \big( H_{N,1} f_N \big)(\gamma_N) \geq  G_{f_{u,\rho_1}}^\phi(\gamma_0).
\end{align*}
Following Remark~\ref{Sec5:HN1Rmk} for the expression of $H_{N,1} f_N$, we can find $\bpi_N \in  \Gamma^{\opt}_\delta(\rho_1^N;\gamma_N)$  with $\delta$ satisfying property \eqref{Sec5:deleps},  and some $\bpi_0 \in \Gamma^{\opt}(\rho_1;\gamma_0)$ such that $\bpi_N \Rightarrow \bpi_0$ in the narrow convergence sense, and such that (by Fatou's lemma)
\begin{align*}
 \liminf_{N \to \infty} H_{N,1} f_N(\gamma_N)   & \geq 
  \liminf_{N \to \infty}   \int_{\R^{2d}} \big( \eta_\phi(\star \star) - U(y) - V*\gamma_N(y) \big)   \bpi_N(dx, dy) \\
 & \geq  \int_{\R^{2d}} \big( \eta_\phi(\star \star) - U(y) - V*\gamma_0(y) \big)  \bpi_0 (dx, dy) 
        \geq  G_{f_{u;\rho_1}}^\phi(\gamma_0),
\end{align*} 
where the notation
\begin{align*}
 (\star \star):= \Big(y, \alpha(x-y) - 2y \zeta^\prime(|y|^2)\big) .
\end{align*}

Finally, we note that Condition~\ref{HConSup}.\ref{tre1} is already verified in Lemma~\ref{Sec6:Kflev}.
\end{proof}

\subsubsection{A lower limit estimate using super-solution for a limiting equation}
Invoking Theorem~\ref{Sec4:BPThm}, we obtain the following one-sided limit result regarding the solution ${\mathfrak f}_N$ to \eqref{HNHJB}.

Let the $\mathfrak h_N$ satisfy invariance \eqref{Sec5:fhtx} and $h_{N,1} $ be defined according to \eqref{Sec5:h1Def}. 
\begin{condition}\label{Sec6:h1Cvg}
$h_1 \in C(\sfX)$ satisfies the following properties: 
\begin{enumerate}
\item for every $\rho_N := \frac1N \sum_{k=1}^N \delta_{x_k^N} \in \sfX_N$ and $\rho_0 \in \sfX$ with $\sfd(\rho_N, \rho_0)=0$, we have 
\begin{align*}
\liminf_{N \to \infty} h_{N,1}(\rho_N) \geq h_1(\rho_0);
\end{align*}
\item moreover, there exists an increasing sub-linear function $\beta: \R_+ \mapsto \R$ such that 
\begin{align*}
\inf_N  h_{N,1}(\rho) \geq   - \beta\circ \sfd(\rho, \delta_0), 
\quad \forall \rho \in \sfX_N,
\quad  \text{ and } \sup_N \sup_{\sfX_N} h_{N,1} < +\infty.
\end{align*}
\end{enumerate}
\end{condition}

 \begin{lemma}\label{Sec6:SupLem}
Suppose that Conditions~\ref{PerCND}, ~\ref{U0CND} and~\ref{Sec1:Tone2} hold. Let ${\mathfrak f}_N: (\R^d)^N \mapsto \R$ and $\underline{f}_N: \sfX_N \mapsto \R$ be defined according to \eqref{PartfDef} and \eqref{Sec5:lfNDef} respectively. We also introduce $\underline{f}: \sfX \mapsto \R$ as in \eqref{fhathat}-\eqref{lbarf}.

We assume that the $h_{N,1}$s and $h_1$ are related by Condition~\ref{Sec6:h1Cvg}. Then the $\underline{f}$ is a super-solution to 
\begin{align}\label{Sec6:H1sup}
 (I - \alpha H_1) \underline{f} \geq h_1,
\end{align}
in the point-wise viscosity solution sense. 
\end{lemma}
 \begin{proof}
The results from Lemmas~\ref{SupQuotHJ}, \ref{Sec5:fnmest} and \ref{Sec6:SupCnv} verify conditions required for applying Theorem~\ref{Sec4:BPThm}, hence the conclusion follows.
 \end{proof}

Up to this point, our definition of the $H_1$ in \eqref{Sec6:H1Def} always has a non-zero $\zeta$ term -- see the $f_{u;\rho_1,\ldots, \rho_K}$ in \eqref{Sec6:furho} which is the origin of the hat on notation $\hat{\mathcal S}^-$. The $\zeta$ played a significant role in producing compactness type arguments in previous proofs 
(e.g. Lemma~\ref{Sec6:Kflev}). Next, we get rid of this term. 
We explore localness of the operator $H_1$, and conclude by a variant of Lemma~\ref{Sec3:seq2spw}.
To reduce an already long list of notations, with a slight abuse of notation,  we still use $H_1$ to denote such reduction. That is, the new $H_1$ is notationally defined as the old one by setting $\zeta=0$, hence the the $f_{u;\rho_1, \ldots, \rho_K} = f_{1;\rho_1, \ldots, \rho_K} \in D(H_1)$ for the new one.
  
 \begin{lemma}\label{Sec6:H1SupUlt}
Let $\underline{f} \in \LSC(\sfX)$ be an at most linear growth super solution to 
\eqref{Sec6:H1sup} in the point-wise viscosity sense, with the $H_1$ defined as in \eqref{Sec6:H1Def} with non-zero $\zeta$ term. Then it is also a super-solution in the strong point-wise viscosity sense, with the $\zeta$ term removed from both the test functions and the operator $H_1$. 
\end{lemma}
\begin{proof}
Let $f_1:=f_{1,\rho_1, \ldots, \rho_K} \in D(H_1)$ for the new $H_1$, and let $\gamma_0 \in \sfX$ be such that
\begin{align*}
 (f_1 - \underline{f})(\gamma_0) =\sup_\sfX (f_1 - \underline{f}).
\end{align*}
Then by adding an extra term as in Lemma~\ref{Sec3:locF}
\begin{align*}
 f_{1,\theta}:=f_{1;\theta; \rho_1,\ldots, \rho_K, \gamma_0}
 := f_{1; \rho_1, \ldots, \rho_K} - \theta \sfd^2(\cdot, \gamma_0),
\end{align*}
the $\gamma_0$ becomes a strict global maxima of $f_{1,\theta} - \underline{f}$.
Let $Y_0$ be a $\R^d$-valued random variable with probability distribution $\gamma_0$, then $E[|Y_0|^2]<\infty$. By a uniform integrability result due to La Vall\'ee Poussin (e.g. see the constructive proof of (1) implies (2) for Theorem T22 in Chapter II, on page 19 of Meyer~\cite{Meyer66}), there exists a non-negative $\zeta \in C^1$ in the class of test functions that we considered earlier, such that
\begin{align*}
\int_{\R^d} \zeta(|y|^2) \gamma_0(dy)=E[\zeta(|Y_0|^2)]<+\infty.
\end{align*}
 We introduce 
 \begin{align*}
f_{1,\theta, n}:=f_{1,\theta, n}(\gamma):= f_{1,\theta}(\gamma)- \int_{\R^d} \frac1n \zeta(|y|^2)  \gamma(dy) \in \hat{\mathcal S}^-.
\end{align*}

By point-wise viscosity super-solution property for \eqref{Sec6:H1sup}, there exists $\gamma_n:=\gamma_{\theta,n} \in \sfX$ with 
\begin{align*}
 (f_{1,\theta,n} - \underline{f})(\gamma_n) = \sup_\sfX (f_{1,\theta, n} - \underline{f}),
\end{align*}
and  
\begin{align}\label{Sec6:AltSupEq}
\alpha^{-1} (\underline{f}- h_1)(\gamma_n) & \geq  
\inf_{\bM \in \Gamma^{\opt}(\gamma_n;\rho_1, \ldots, \rho_K, \gamma_0)} 
 \int_{\R^{2d}}   \eta_\phi \big( y, P -\frac{2}{n}y\zeta^\prime(|y|^2) \big) \bnu^{\bM}_{f_{1,\theta;\rho_1,\ldots, \rho_K;\gamma_0}}(dy, dP) \\
 & \qquad \qquad - \langle (U + V*\gamma_n), \gamma_n\rangle. \nonumber
 \end{align}

Summarizing the above, we have 
\begin{align*}
(f_{1,\theta} - \underline{f})(\gamma_n) & \geq  (f_{1,\theta,n} - \underline{f})(\gamma_n) = \sup_\sfX (f_{1,\theta,n} - \underline{f})\\
&\qquad \qquad \geq (f_{1,\theta,n} - \underline{f})(\gamma_0) = \sup_\sfX (f_{1,\theta} - \underline{f}) - \frac1n \int_{\R^d} \zeta(|y|^2) \gamma_0(dy).
\end{align*}
That is, the $\gamma_n:= \gamma_{\theta, n}$ has a super-solution version of the property as in \eqref{Sec3:seqmax}. Hence by a super-solution version of Lemma~\ref{Sec3:locF}, we have (for every fixed $\theta$)
\begin{align*}
 \lim_{n \to \infty} \sfd(\gamma_n, \gamma_0) =0, \quad 
 \lim_{n \to \infty} \underline{f}(\gamma_n) = \underline{f}(\gamma_0), \quad
 \lim_{n \to \infty} f_{1,\theta}(\gamma_n) = f_{1,\theta}(\gamma_0).
\end{align*}
By Fatou's lemma, passing $n \to \infty$ in \eqref{Sec6:AltSupEq} gives
\begin{align*}
\alpha^{-1} (\underline{f}- h_1)(\gamma_0)  \geq  (H_1 f_{1,\theta}) (\gamma_0).
\end{align*}

Invoking similar arguments as in Lemma~\ref{Sec3:seq2spw}, 
we send $\theta \to 0^+$ and conclude. 
 \end{proof}
 \subsection{Convergence of Hamiltonians  - the sub-solution case}\label{Sec6:Sub}
As mentioned in the opening paragraph of Section~\ref{BPlimitA}, we cannot directly apply the half-relaxed-limit Theorem~\ref{Sec4:BPThm} to the $\overline{f}_N$s and $H_{N,0}$s in Lemma~\ref{SubQuotHJ}, in the sub-solution case. We will apply the generalized version of results given by Lemma~\ref{Sec4:BPAlt} to those perturbed Hamilton-Jacobi equations \eqref{Sec5:HNtheta} in Lemma~\ref{Sec5:SubQuoA}. In this step, the Hamiltonians act on test functions defined using $p$-Wasserstein distance for $1< p <2$. We use sequence of sub-solutions which are those perturbed ones given by \eqref{Sec5:barfNthe} with a small parameter $\theta>0$. Upon getting limiting equation using the abstract results of Lemma~\ref{Sec4:BPAlt}, we will then let the $p \to 2$, followed by letting the $\theta \to 0^+$, to arrive at another limiting Hamiltonian defined on functions over $2$-Wasserstein space. This last step relies upon a type of viscosity extension method which was first introduced in Feng and Kurtz~\cite{FK06}.
 
We present details of the above procedure step by step next. 
\subsubsection{Convergence of spaces}
Let $p_0 \in (1,2)$ be arbitrary but fixed. We take $\sfX:= {\mathcal P}_2(\R^d)$ with 
$\sfd_\sfX:=\sfd$ the $2$-Wasserstein metric, and 
$\sfX^\prime:= {\mathcal P}_{p_0}(\R^d)$ with 
$\sfd_{\sfX^\prime}:=\sfd_{p_0}$ the $p_0$-Wasserstein metric. 
 The space $(\sfX_N, \sfd_{\sfX_N})$ is taken according to the definition in Section~\ref{Sec5:RevSub}. In particular, we use $\sfd_{\sfX_N}$ here the $p_0$-Wasserstein metric restricted to the space of empirical probability measures with $N$ equally-weighted point masses. 
 
We take $\mathcal Q := \R_+$ and for each $q \in \mathcal Q$,
 \begin{align*}
 K_N^q:= \{ \rho \in \sfX_N : \int_{\R^d} |x|^2 \rho(dx) \leq q\}, \quad
 K^q:= \{ \rho \in \sfX : \int_{\R^d} |x|^2 \rho(dx) \leq q\}.
\end{align*}
Let $\eta_N : \sfX_N \mapsto \sfX \subset \sfX^\prime$ be the identity map and 
$\eta_N^q:= \eta_N \big|_{K^q_N} : K_N^q \mapsto K^q$.
The $K^q$ is compact in $(\sfX^\prime, \sfd_{\sfX^\prime})$ for each $q \in \mathcal Q$; $\cup_{q \in \mathcal Q} K^q = \sfX$ and 
\begin{align*}
\sfX =\overline{\cup_{q \in \mathcal Q} K^q}^{\sfd_{\sfX}-{\rm closure}} = \cup_{q \in \mathcal Q} K^q, 
\qquad \overline{\sfX}^{\sfd_{\sfX^\prime}-{\rm closure}} = \sfX^\prime.
\end{align*}
Moreover, for each $\rho_N \in K_N^q$ and $\rho_0 \in K^q$ with 
$\lim_{N \to \infty} \sfd_{\sfX^\prime}(\rho_N, \rho_0) =0$, by uniformly bounded second moments property of the $\rho_N$s, we have indeed a stronger convergence
\begin{align*}
 \lim_{N \to \infty} \sfd_p(\rho_N, \rho_0) =0, \quad \forall p \in [1,2).
\end{align*}

\begin{lemma}
\begin{align*}
 (\sfX_N, \sfd_{\sfX_N}) \stackrel{\rm gGH}{\longrightarrow}_{\mathcal Q}   (\sfX^\prime, \sfd_{\sfX^\prime}).
\end{align*}
\end{lemma}
 
\subsubsection{A limiting Hamiltonian operator and a limiting sub-solution}\label{Sec6:limHsub}
  Take $p \in [p_0, 2)$, and let test functions $f_0: \sfX \mapsto \R$ be defined according to 
\begin{align}\label{Sec6:f0p}
f_0(\rho):=f_{0;\gamma_1,\ldots, \gamma_K} (\rho)
 :=f_{0,p;\gamma_1,\ldots, \gamma_K} (\rho)
   := \psi\big(\sfd^2_p(\rho,\gamma_1),\ldots, \sfd^2_p(\rho,\gamma_K)\big),
\end{align}
just as in \eqref{Sec5:f0Alt}. We also define measure 
$\bnu^{\bM}_{f_0}:= \bnu^{\bM}_{f_0; \gamma_1, \ldots, \gamma_K}$ by \eqref{Sec5:p-nu}, representing a type of differential of the $f_0$.  Let $\phi:=\phi(x,P;q) \in {\mathcal F}_0$. We recall the $\eta^\phi$ defined in \eqref{Sec5:etaDef}:
\begin{align*}
 \eta^\phi(x,P):= \sup_{q \in \R^d} \sfH\big( q, P+ \nabla_q \phi(x,P;q)\big),
\end{align*}
and the collection of optimal multi-plans $\Gamma^{\opt}_p(\rho; \gamma_1, \ldots, \gamma_K)$ in \eqref{pMultOpt}. With all these notations, we denote
\begin{align}\label{Sec6:GGDef}
 G^{f_{0;\gamma_1, \ldots, \gamma_K}}_\phi(\rho)  
  & := \sup_{\bM \in \Gamma_p^{\opt}(\rho;\gamma_1,\ldots, \gamma_K)}  
 \int_{\R^{2d}}    \eta^\phi (x, P)    \bnu^{\bM}_{f_{0;\gamma_1,\ldots, \gamma_K}}(dx,dP)  \\
 & \qquad \qquad - \langle U, \rho \rangle - \langle V*\rho, \rho \rangle.    \nonumber 
 \end{align}

Next, we introduce a multi-valued Hamiltonian operator $H_0 \subset C\big((\sfX, \sfd); \R\big) \times M\big((\sfX, \sfd);\bar{\R}\big)$ identified through its graph by
 \begin{align}\label{Sec6:H0Def}
H_0  := \Big\{ \big(f_0, G^{f_0}_{\phi}\big): f_0:=f_{0,p; \gamma_1, \ldots, \gamma_K}  
  \text{ as above}, p \in [p_0,2), \phi \in {\mathcal F}_0 \Big\}.
\end{align}
Expression of the $H_0$ in the following special situation simplifies: let $K=1$ and $\psi(r):= \frac{\alpha}{2} r$ with $\alpha>0$, then
\begin{align*}
 f_0:=f_{0;\gamma_1}(\rho)= \frac{\alpha}{2} \sfd_p^2(\rho,\gamma_1).
\end{align*}
and 
\begin{align*}
 G^{f_0}_{\phi}(\rho) &=  \sup_{\bpi \in \Gamma^{\opt}_p(\rho;\gamma_1)}
     \int_{\R^{2d}} \eta^\phi \Big(x,  \alpha \sfd_p^{2-p}(\rho,\gamma_1) |x-y|^{p-1}\frac{x-y}{|x-y|} \Big) \bpi(dx,dy) \\
     & \qquad \qquad \qquad \qquad - \langle (U+ V* \rho), \rho \rangle.
\end{align*}
That is, for each fixed $f_0$ of the above type, $H_0 f_0$ is a set of functions
 $\{ G^{f_0}_{\phi}: \phi  \in {\mathcal F}_0 \}$. 
 
We recall regularity estimates for the $\eta^\phi$ in Lemma~\ref{Sec5:etaphi}, which give the following. 
\begin{lemma}
 $ \rho \mapsto G^{f_0}_\phi(\rho)$ is upper semi-continuous in the topology given by $p$-Wasserstein metric, for each $p \in (1,2)$.  
\end{lemma}
%
%
  
Let $\overline{f}_N$ be those constructed in Lemma~\ref{SubQuotHJ}. They solve \eqref{prHN0Eqn} in sub-solution sense as have been made precise in that lemma. Let $\theta >0$ be fixed. 
Recall that, through \eqref{Sec5:barfNthe}, we defined
\begin{align}\label{Sec6:fbarPert}
\overline{f}_{N,\theta}(\rho) := \overline{f}_N(\rho) -\frac{\theta}{2} \int_{\R^d} |x|^2 \rho(dx).
\end{align}
By Lemma~\ref{Sec5:SubQuoA}, the above function is a strong viscosity sub-solution in the point-wise sense to an equation \eqref{Sec5:HNtheta} given by the perturbed 
Hamiltonian operator $H_{N,0}^\theta$, defined in that lemma. Following the notations there about the constants $c, C>0$, we denote for every $\lambda>1$,
\begin{align}\label{Sec6:GSubthelam}
  G^{f_0;\theta, \lambda}_\phi(\rho)  
&   := G_\phi^{f_0}(\rho)  +  (1-\frac{1}{\lambda})\Big\{ c- \inf\sfH 
+ 4C \sup_{\bM \in \Gamma^{\opt}_p(\rho; \gamma_1, \ldots, \gamma_K)}
  \int_{\R^{2d}} |P|^2 \bnu^{\bM}_{f_{0;\gamma_1,\ldots, \gamma_K}}(dx, dP)  \\
 & \qquad  \qquad  + 4C \sup_{q,x,P} |\nabla_q \phi|^2 \Big\}
  + \frac{4C \lambda}{\lambda -1} \theta\Big\{
 \sup_N \sup_{\sfX_N} \overline{f}_N -\psi(0,\ldots, 0) +1 \nonumber \\
& \qquad \qquad \qquad \qquad 
 - \inf_N \overline{f}_N(\delta_0) +  f_0(\delta_0)   \Big\} \nonumber \\
 & \quad =: G_\phi^{f_0}(\rho)  + \text{Error}_{1,\lambda} (\rho) 
  +  \text{Error}_{2, \theta,\lambda}. \nonumber
 \end{align}
 We introduce a perturbed operator (multi-valued) 
\begin{align}\label{Sec6:H0the}
 H_0^\theta   := \Big\{ \big(f_0, G^{f_0;\theta, \lambda}_{\phi}\big): f_0:=f_{0,p; \gamma_1, \ldots, \gamma_K}, p \in (1,2), \phi \in {\mathcal F}_0, \lambda >1 \Big\}.
\end{align}

\begin{lemma}\label{Sec6:Err12}
For every finite $L>0$, we have 
\begin{align*}
 \lim_{\lambda \to 1^+} \sup_{\rho \in \sfX, \int_{\R^d} |x|^2 \rho(dx) \leq L} \text{\rm Error}_{1,\lambda} (\rho) =0,  \qquad 
  \lim_{\lambda \to 1^+} \lim_{\theta \to 0^+}  \text{\rm Error}_{2, \theta,\lambda}=0.
\end{align*}
\end{lemma}

We recall that the ${\mathfrak h}_N : (\R^d)^N \mapsto \R$ and the $h_{N,0}: \sfX_N \mapsto \R$ are defined as in Lemma~\ref{SubQuotHJ}. 
We assume the following sub-solution counterpart of the Condition~\ref{Sec6:h1Cvg}.
 \begin{condition}\label{Sec6:h0hN}
 $h_0 \in C(\sfX)$ and $h_{N,0} : \sfX_N \mapsto \sfX$ satisfy the following:  
\begin{enumerate}
\item the $h_{N,0}$s are uniformly bounded:
\begin{align*}
\sup_N \sup_{\sfX_N} h_{N,0} <+\infty.
\end{align*}
\item the sequence $\{ h_{N,0} : N \in \N\}$ and the $h_0$ satisfies Property 
${\mathscr P_N}$ (Definition~\ref{Sec4:ProPN}). Specifically in current context, this becomes: for every $\rho_N := \frac1N \sum_{k=1}^N \delta_{x_k^N} \in \sfX_N$ and $\rho_0 \in \sfX$ with $\sfd_{p=1}(\rho_N, \rho_0)=0$, and with $\sup_N \int_{\R^d} |x|^2 \rho_N (dx) <+\infty$, we have $\limsup_{N \to \infty} h_{N,0}(\rho_N) \leq h_0(\rho_0)$.
\end{enumerate}
\end{condition}
 
  \begin{lemma}
  Suppose that Conditions~\ref{PerCND}, \ref{U0CND} and \ref{Sec6:h0hN} hold.
  Then Condition~\ref{Sec4:CnvCAlt} is satisfied for the operators $H_{N,0}^\theta$ and $H_0^\theta$.
  \end{lemma}
  \begin{proof}
Since the $K^q$ is $2$-Wasserstein closed balls of radius $q \in \R_+$, it is compact in the $p$-Wasserstein space for any $p \in (1,2)$.
Recall that $\eta_N^q : K_N^q \mapsto K^q$ is simply the identity embedding map from $\sfX_N$ to $\sfX$. Therefore, if $\rho_N \in K^q_N$, $\rho_0 \in K^q$ and $\lim_{N \to \infty} \sfd_{\sfX^\prime}( \rho_N, \rho_0) =0$,  then $\rho_N \to \rho_0$ in any $p$-Wasserstein metric with $1<p<2$. In particular, this implies that each $f_0 \in D(H_0^\theta)$ is $\sfd_{\sfX^\prime}$-continuous in the closed subset $K^q \subset \sfX$, verifying Condition~\ref{Sec4:CnvCAlt}.\ref{Sec4:dprmC}.

Let $(f_0, G^{f_0;\theta,\lambda}_\phi) \in H^\theta_0$, where the $f_0 := f_{0,p;\gamma_1,\ldots, \gamma_K}$ for some $p \in [p_0,2)$ as in \eqref{Sec6:f0p}. We can find $\gamma_k^N := \frac1N \sum_{i=1}^N \delta_{y_i^k} \in \sfX_N$ such that $\lim_{N \to \infty} \sfd_\sfX(\gamma_k^N, \gamma_k) =0$ for $k=1,\ldots,K$. We define $\mathfrak g$ according to \eqref{PiyDef} and 
$f_{0,N}:=f_{0, \epsilon_N \mathfrak g; \gamma_1^N,\ldots, \gamma_K^N}$ 
according to \eqref{Sec5:perTestf0}.

Take $\epsilon_0 \in (0,1)$. If, additionally, we assume that $\rho_N \in E_{\epsilon_0}^+[ \overline{f}_{N,\theta} - f_{0,N} ] \subset K_N^q$ as in Condition~\ref{Sec4:CnvCAlt}.\ref{Sec4:EKpc}, then  
\begin{align*}
\frac{\theta}{2} \int_{\R^d} |x|^2 \rho_N(dx) 
& \leq (\overline{f}_N  - f_{0,N})(\rho_N) +1  
   -  (\overline{f}_{N,\theta} - f_{0,N})(\delta_0)\\
& \leq \big( \sup_N \sup_{\sfX_N} \overline{f}_N 
       -\psi(0,\ldots, 0)  - \epsilon_N \inf \phi \big) +1 \\
       & \qquad \qquad
 - \inf_N \overline{f}_N(\delta_0) 
  + \big( f_0(\delta_0) + \epsilon_N \sup\phi \big)  <+\infty,
\end{align*}
where the $\psi$ is the one appearing in definition of the $f_0$.
The above estimate implies $\{ \rho_N\}_N$ is relatively compact in any topology given by $r$-Wasserstein metric with $r \in (1,2)$.

Recall the definition of $H_{N,0}^\theta f_{0,N}$ in \eqref{Sec5:HN0the}, 
we have
\begin{align*}
  \limsup_{N \to \infty} H_{N,0}^\theta f_{0,N}(\rho_N) \leq G_\phi^{f_{0;\theta,\lambda}}(\rho_0). 
\end{align*}
  \end{proof}
  
We recall the definition of $\overline{f}_{N,\theta}$ in \eqref{Sec5:barfNthe} and the result of Lemma~\ref{Sec5:SubQuoA}. Following \eqref{Sec4:ftil} in Section~\ref{BPlimitA}, we define, for each $\theta >0$ and $\rho_0 \in \sfX =\overline{\cup_{q \in \mathcal Q} K^q}^{\sfd_{\sfX}-{\rm closure}} = \cup_{q \in \mathcal Q} K^q$, 
\begin{align*}
\tilde{f}_\theta(\rho_0)& := \sup_{L >0} \sup \Big\{ 
   \limsup_{N \to \infty}   \overline{f}_{N,\theta} (\rho_N)  : \exists \rho_N \in \sfX_N, 
  \int_{\R^d} |x|^2 \rho_N(dx) \leq L,
 \lim_{N \to \infty} \sfd_{\sfX^\prime}(\rho_N,\rho_0) =0 \Big\}  \\
& = \sup_{L >0} \sup \Big\{ \limsup_{N \to \infty} \Big( \overline{f}_N(\rho_N) 
  - \frac{\theta}{2} \int_{\R^d} |x|^2 \rho_N(dx) \Big) : \exists \rho_N \in \sfX_N, \\
 & \qquad \qquad \qquad \qquad \int_{\R^d} |x|^2 \rho_N(dx) \leq L,
 \lim_{N \to \infty} \sfd_{\sfX^\prime}(\rho_N,\rho_0) =0 \Big\},
\end{align*}
and for every $\rho \in \sfX$,
\begin{align}\label{Sec6:fstarThe}
f_\theta^*(\rho) := \lim_{\epsilon \to 0^+}\sup \{ \tilde{f}_\theta(\rho_0) :
 \rho_0 \in {\mathcal P}_2(\R^d),  \sfd_\sfX(\rho_0, \rho) < \epsilon \}. 
\end{align}
Then $f^*_\theta \in \USC\big((\sfX, \sfd) ;\R \big)$ is bounded from above~\footnote{We note that the $\overline{f}_N$ is bounded from above by Lemma~\ref{parHJprop}}. A little thinking reveals an even stronger result: the map $(\rho, \theta) \mapsto f^*_\theta(\rho)$ belongs to $\USC\big(\sfX \times [0,\theta_0] ;\R \big)$ for every $\theta_0>0$. Also, $f^*_\theta \leq f^*_{\theta^\prime}$ for every $0<\theta^\prime < \theta$. Let
\begin{align}\label{Sec6:fstarDef}
f^*(\rho):= \limsup_{\theta \to 0^+} f^*_\theta(\rho) =   \sup_{\theta \in [0,\theta_0]} f^*_\theta(\rho), \quad \forall \rho \in \sfX, \theta_0>0.
\end{align}
Then by Lemma~\ref{infsupSC}, $f^* \in \USC\big((\sfX, \sfd);\R\big)$ with the following property. We will use this property in proof of Theorem~\ref{Sec8:MainThm1} later.

\begin{lemma}\label{Sec6:Lfstar}
For every $\rho_N \in \sfX_N$ and $\rho_0 \in \sfX$ such that 
$\lim_{N \to \infty} \sfd_\sfX(\eta_N(\rho_N) , \rho_0)=0$, we have 
\begin{align*}
\limsup_{N \to \infty} \overline{f}_N(\rho_N) \leq f^*(\rho_0).
\end{align*}
In fact, the following stronger result holds: The sequence $\{\overline{f}_N \}_{N \in \N}$ and $f^*$ satisfy Property ${\mathscr P_N}$  in Definition~\ref{Sec4:ProPN}.
\end{lemma}
\begin{proof}
Suppose the $\rho_N, \rho_0$ are such that 
\begin{align*}
\sup_N \int_{\R^d} |x|^2 \rho_N(dx) <+\infty, \quad 
\lim_{N \to \infty} \sfd_{\sfX^\prime}(\eta_N(\rho_N),\rho_0) = 0.
\end{align*}
Then
\begin{align*}
\limsup_{N \to \infty} f_N(\rho_N) &= \limsup_{\theta \to 0^+}  \limsup_{N \to \infty} 
\Big( f_N(\rho_N) - \frac{\theta}{2} \int_{\R^d} |x|^2 \rho_N(dx) \Big) \\
& \leq   \limsup_{N \to \infty} f_N(\rho_N)  \leq \tilde{f}(\rho_0) \leq  f^*(\rho_0). 
\end{align*}
\end{proof}

An application of Lemma~\ref{Sec4:BPAlt}, to the sub-solution result of Lemma~\ref{Sec5:SubQuoA},  gives the next result. 
\begin{lemma}\label{Sec6:subCnvEqn}
Suppose that Conditions~\ref{PerCND}, \ref{U0CND} and \ref{Sec6:h0hN} hold. 
Then the above defined $f_\theta^*$ is a sub-solution to 
\begin{align}\label{Sec6:H0eqn}
(I - \alpha H_0^\theta) f_\theta^* \leq h_0,
\end{align}
 in the point-wise viscosity solution sense, in the $2$-Wasserstein space $(\sfX, \sfd)$.
 The $h_0$ is the one in Condition~\ref{Sec6:h0hN}. 
\end{lemma} 

The $f_\theta^*$ is more than just $\sfd_\sfX$-upper semi-continuous in the $\sfX$. In fact, the following $\sfd_{\sfX^\prime}$-upper semi-continuity in $\sfX$ property holds.
\begin{lemma}\label{Sec6:fwUSC}
For every $\rho_N, \rho_0 \in \sfX$ with $\lim_{N \to \infty} \sfd_{\sfX^\prime}(\rho_N, \rho_0) =0$, we have
\begin{align*}
  \limsup_{N \to \infty} f_\theta^*(\rho_N) \leq f_\theta^*(\rho_0).
\end{align*}
\end{lemma}
\begin{proof}
For the $\rho_N, \rho_0$, by definition of $f^*_\theta$, 
we can always find $\rho_N^\prime \in \sfX$ 
with $\lim_{N \to \infty}\sfd_\sfX(\rho_N^\prime, \rho_N) =0$ and 
$f^*_\theta(\rho_N) \leq \tilde{f}_\theta(\rho_N^\prime)+ \frac1N$. 
Consequently, $\lim_{N \to \infty} \sfd_{\sfX^\prime}(\rho_N^\prime, \rho_0) =0$, 
$\sup_N \int_{\R^d} |x|^2 \rho^\prime_N(dx) <+\infty$ and 
\begin{align*}
\limsup_{N \to \infty} f^*_\theta(\rho_N) \leq 
   \limsup_{N \to \infty} \tilde{f}_\theta(\rho_N^\prime) 
\leq \tilde{f}_\theta(\rho_0) \leq f^*_\theta(\rho_0).
\end{align*}
\end{proof}

\subsubsection{Viscosity extension for limiting Hamiltonian operators, the sub-solution case - I }
We established, in Lemma~\ref{Sec6:subCnvEqn}, that $f^*_\theta$ is a sub-solution to \eqref{Sec6:H0eqn}, for test functions in $D(H_0^\theta)$. Next, we enlarge the domain of test functions to include those of following type:
\begin{align}\label{Sec6:f0e}
 f_0^\epsilon(\rho) & := f_0(\rho) + \epsilon \sfd^2(\rho, \gamma_0), \quad f_0 \in D(H_0^\theta), \epsilon>0, \gamma_0 \in \sfX,\\
 & =  \psi(\sfd_p^2(\rho, \gamma_1), \ldots, \sfd_p^2(\rho,\gamma_K)) + \epsilon \sfd^2(\rho, \gamma_0), \quad \forall p \in (1,2).   \nonumber
\end{align}
We want to show that $f^*_\theta$ is still a sub-solution, by correspondingly extending the Hamiltonian operator. 
We note that, in this and subsequent subsections, $\sfd:= \sfd_\sfX$ denotes the $2$-Wasserstein metric. The $\gamma_1, \ldots, \gamma_K, \gamma_0 \in \sfX$ all have finite second moment.

The function $f^\epsilon_0$ can be approximated in monotone point-wise convergence sense, 
in the $p_n \to 2^-$ limit (assuming $1 < p_n < p_{n+1}   <2$), by 
\begin{align}\label{Sec6:fne}
f_{n}^\epsilon:= f_{n}^\epsilon(\rho) &:=  f_0(\rho) + \epsilon \sfd_{p_n}^2(\rho, \gamma_0) \\
& =  \psi(\sfd^2_p(\rho, \gamma_1), \ldots, \sfd^2_p(\rho,\gamma_K)) + \epsilon \sfd_{p_n}^2(\rho, \gamma_0) \in D(H_0^\theta). \nonumber
\end{align}
Let  $c, C$ be the constants in \eqref{amQforH}. Next, we define a perturbative version of the $G^{f_0;\theta,\lambda}_\phi$ in \eqref{Sec6:GSubthelam}. For each $\theta>0, \lambda >1$ and $\phi \in \mathcal F_0$, $\epsilon>0$, we define
\begin{align}\label{Sec6:Gtil}
\tilde{G}^{f_0^\epsilon;\theta,\lambda}_\phi(\rho) &:= G^{f_0;\theta,\lambda}_\phi(\rho) 
    + (1-\frac{1}{\lambda})\Big\{ c - \inf{\sfH} \\
    & \qquad  \qquad \qquad 
    + 4C \sup_{\bM \in \Gamma_p^{\opt}(\rho; \gamma_1, \ldots, \gamma_K)} \int_{\R^{2d}} |P|^2 \bnu^{\bM}_{f_{0;\gamma_1, \ldots, \gamma_K}}(dx,dP) 
     \nonumber  \\
    & \qquad \qquad \qquad \qquad 
  + 4C \sup_{x,P,q} |\nabla_q \phi|^2 \Big\} 
   + \frac{2C\lambda}{\lambda-1} 4 \epsilon^2  
      \sfd^2(\rho, \gamma_0). \nonumber
\end{align}
We also define operator
\begin{align*}
\tilde{H}_0^\theta:= H_0^\theta \cup \Big\{ (f_0^\epsilon, \tilde{G}^{f_0^\epsilon;\theta,\lambda}_\phi) :
 \forall f_0^\epsilon \text{ as above with } p \in [p_0,2),  \epsilon>0, \lambda >1, \phi \in \mathcal F_0 \Big\}.
\end{align*}

Our main result for this subsection is Lemma~\ref{Sec6:Htild}. Before stating it,  let us quote the following property regarding a special type of $\Gamma$-convergence. ~\footnote{The type referred to here is monotone pointwise convergence.}  
\begin{lemma}\label{Sec6:Attou}
Let $\sfX^\prime$ be a general metric space
and $F_n : \sfX^\prime \mapsto \R$ be such that $\{ F_n \}_{n \in \N}$ is a non-increasing sequence of upper semi-continuous functions with limit function $F$ (in point-wise convergence sense). Let $x_n, x_0 \in \sfX^\prime$ and $\epsilon_n>0$ be such that $\sup_{\sfX^\prime} F_n \leq F_n(x_n) + \epsilon_n$, $\epsilon_n \to 0^+$ and $x_n \to x_0$. Then
\begin{align*}
\lim_{n \to \infty} \sup_{\sfX^\prime} F_n = \lim_{n \to \infty} F_n(x_n) = F(x_0) = \sup_{\sfX^\prime} F. 
\end{align*}
\end{lemma}
\begin{proof}
See Lemma A.4 of Feng and Kurtz~\cite{FK06}, or more generally,  Proposition 2.42 of Attouch~\cite{Att84}.
\end{proof}

\begin{definition}\label{Sec6:weakUSC}
Let $r \in [1,2)$. A function $h : \sfX \mapsto \R$ is said to be $\sfd_r$-upper semi-continuous in $\sfX$ if the following holds: for every $\rho_n, \rho_0 \in \sfX$ with $\lim_{n \to \infty} \sfd_r(\rho_n,\rho_0) =0$, we have
\begin{align*}
 \limsup_{n \to \infty} h(\rho_n) \leq h(\rho_0).
\end{align*}
\end{definition}

\begin{lemma}\label{Sec6:Htild}
Let $h_0 \in C(\sfX)$ be obtained as in Lemma~\ref{Sec6:subCnvEqn}, and is $\sfd_{p=1}$ upper semicontinuous in $\sfX$ (see Definition~\ref{Sec6:weakUSC}).
Then the $f_\theta^*$ is a sub-solution in the point-wise viscosity sense to 
\begin{align*}
 (I - \alpha \tilde{H}_0^\theta) f_\theta^* \leq h_0.
\end{align*}
\end{lemma}
\begin{proof}  
The proof follows from a variational convergence method introduced in Lemmas 7.7 and 13.21 in Feng and Kurtz~\cite{FK06}. Some modifications are needed in order to be adapted here. 

Fix a $f^\epsilon_0$ in \eqref{Sec6:f0e}, we construct $f^\epsilon_n$ as in \eqref{Sec6:fne}.
By Lemma~\ref{Sec6:subCnvEqn}, there exists $\rho_n \in \sfX$ such that 
\begin{align}\label{Sec6:fthemax}
(f_\theta^* - f_n^\epsilon)(\rho_n)=\sup_\sfX(f_\theta^* - f_n^\epsilon) 
\end{align}
and that
\begin{align}\label{Sec6:ExtEQ1}
 \alpha^{-1} (f_\theta^* - h_0)(\rho_n) \leq G_\phi^{f_n^\epsilon;\theta, \lambda} (\rho_n), \quad \forall \theta>0, \lambda>1, \phi \in \mathcal F_0.
\end{align}
Since $\sup_\sfX f_\theta^* <\infty$ and the $\epsilon>0$ is fixed, we have 
\begin{align}\label{Sec6:rhonpB}
\sup_n \int |x|^{p_n}\rho_n(dx)<\infty.
\end{align}
Hence $\{ \rho_n \}_{n \in \N}$ is relatively compact in topology given by $r$-Wasserstein metric, for every $r \in (1,2)$ fixed. We take one satisfying $0<2(p-1)<r<2$ (this $p$ is the fixed parameter in the $f_0^\epsilon$), and label convergence subsequence still using $\{ \rho_n \}_{n \in \N}$. 
By Fatou's lemma applied to estimate \eqref{Sec6:rhonpB}, then there exists 
$\rho_0 \in \sfX$ such that
\begin{align*}
 \lim_{n \to \infty}\sfd_r(\rho_n, \rho_0) =0. 
\end{align*}
Note that
\begin{align*}
 f_n^\epsilon \leq f_{n+1}^\epsilon \leq \ldots \leq f_0^\epsilon, 
  \quad \lim_{n \to \infty} f_n^\epsilon(\rho) = f_0^\epsilon(\rho).
\end{align*}
Applying Lemmas~\ref{Sec6:fwUSC} and \ref{Sec6:Attou} to finite upper-level sets of $F_n := f_\theta^* - f_n^\epsilon$, we obtain
\begin{align}\label{Sec6:ftheeps}
\lim_{n \to \infty} \sup_{\sfX} (f_\theta^* - f_n^\epsilon)
=\lim_{n \to \infty} (f_\theta^* - f_n^\epsilon)(\rho_n) 
=   (f_\theta^* - f_0^\epsilon)(\rho_0) =  \sup_{\sfX} (f_\theta^* - f_0^\epsilon).
\end{align}
Therefore,
\begin{align}\label{Sec6:dpnEst}
\limsup_{n \to \infty} \epsilon \sfd_{p_n}^2(\rho_n, \gamma_0) & 
= \limsup_{n \to \infty} \big( (f^*_\theta - f_0)(\rho_n)
 - \sup_\sfX(f_\theta^* - f_n^\epsilon) \big) \\
& \leq  (f^*_\theta - f_0)(\rho_0)- \sup_\sfX(f_\theta^* - f_0^\epsilon) =\epsilon \sfd^2(\rho_0, \gamma_0), \nonumber
\end{align}
where the first identity above follows from \eqref{Sec6:fthemax} and 
the inequality from \eqref{Sec6:ftheeps} and Lemma~\ref{Sec6:fwUSC}.

A little thinking also reveals that the above implies
\begin{align}\label{Sec6:fTheStaCv}
\lim_{n \to \infty} f_\theta^*(\rho_n) = f_\theta^*(\rho_0), \quad 
 \lim_{n \to \infty} f_n^\epsilon(\rho_n) = f_0^\epsilon(\rho_0).
\end{align}

From the estimate of $\eta^\phi$ in \eqref{Sec5:etaP}, recall quantities $G^{f_0}_\phi$ and $G^{f_n^\epsilon}_\phi$ as defined by \eqref{Sec6:GGDef}, the $f_n^\epsilon$ and $f_0$ as appeared in \eqref{Sec6:fne}, as well as the expression in \eqref{Sec5:p-nu}, we have 
\begin{align*}
 G^{f_n^\epsilon}_\phi(\rho) & 
 \leq G^{f_0}_\phi(\rho) 
    + (1-\frac{1}{\lambda})\Big\{ c - \inf{\sfH} + 4C \Big(  \sup_{\bM \in \Gamma^{\opt}_p(\rho; \gamma_1, \ldots, \gamma_K)} \int_{\R^{2d}} |P|^2 \bnu^{\bM}_{f_{0;\gamma_1, \ldots, \gamma_K}}(dx,dP)  \\
    & \qquad \qquad 
    + \sup_{q,x,P} |\nabla_q \phi|^2 \Big) \Big\}
       + \frac{2C\lambda}{\lambda-1}  4 \epsilon^2 
        \big(\sfd_{p_n}^{2(2-p_n)}(\rho,\gamma_0) \big)
         \big( \sfd_{2(p_n-1)}^{2(p_n-1)}(\rho, \gamma_0)\big),
           \qquad \forall \rho \in \sfX.
\end{align*}
Note that the $f_0$ has dependency on the parameter $p$, and so are the $G^{f_0}_\phi(\rho)$ and $\bnu^{\bM}_{f_{0;\gamma_1, \ldots, \gamma_K}}$. Furthermore, in view of \eqref{Sec6:dpnEst}, 
\begin{align*}
\limsup_{n \to \infty}   \sfd_{p_n}^2(\rho_n, \gamma_0) \leq   \sfd^2(\rho_0,\gamma_0). 
\end{align*}
Hence
\begin{align*}
\limsup_{n \to \infty} \big(\sfd_{p_n}^{2(2-p_n)} 
 \times \sfd_{2(p_n-1)}^{2(p_n-1)}\big)(\rho_n, \gamma_0) 
 & \leq \limsup_{n \to \infty} \big(\sfd_2^{2(2-p_n)} 
 \times \sfd_2^{2(p_n-1)}\big)(\rho_n, \gamma_0) \\
&=  \limsup_{n \to \infty}   \sfd_2^2(\rho_n, \gamma_0) \leq  
 \sfd^2(\rho_0,\gamma_0).
\end{align*}
We note that convergence in the $r$-Wasserstein metric implies convergence in $2(p-1)$-Wasserstein metric. In fact,  since the convergence holds for arbitrary $r \in (1,2)$, the $p$-Wasserstein convergence holds and 
\begin{align*}
& \limsup_{n \to \infty} G^{f_n^\epsilon}_\phi(\rho_n) \\
 &\qquad \leq G^{f_0}_\phi(\rho_0) 
    + (1-\frac{1}{\lambda})\Big\{ c - \inf{\sfH} + 4C \sup_{\bM \in \Gamma_p^{\opt}(\rho_0; \gamma_1, \ldots, \gamma_K)} \int_{\R^{2d}} |P|^2 \bnu^{\bM}_{f_{0;\gamma_1, \ldots, \gamma_K}}(dx,dP) \\
    & \qquad \qquad \qquad \qquad \qquad \qquad
    + 4C \sup_{x,P,q} |\nabla_q \phi |^2\Big\}   + \frac{2C\lambda}{\lambda-1} 4 \epsilon^2 \sfd^2(\rho_0, \gamma_0).
\end{align*}

Through \eqref{Sec6:GSubthelam}, we see that the 
$G^{f_n^\epsilon;\theta,\lambda}_\phi$ and 
$G^{f_0;\theta,\lambda}_\phi$ are respectively perturbative 
versions of the $G^{f_n^\epsilon}_\phi$ and $G^{f_0}_\phi$.
  Consequently, the above estimate also gives
\begin{align*}
 \limsup_{n \to \infty} G^{f_n^\epsilon;\theta,\lambda}_\phi(\rho_n) 
   \leq \tilde{G}^{f_0^\epsilon;\theta,\lambda}_\phi(\rho_0),
\end{align*}
with the last quantity defined by \eqref{Sec6:Gtil}.

We now conclude the lemma by applying the above estimate and \eqref{Sec6:fTheStaCv} to \eqref{Sec6:ExtEQ1}, and by noting 
$\limsup_{n \to \infty} h_0(\rho_n) \leq h_0(\rho_0)$ 
(we assumed that $h_0$ satisfies the property in Definition~\ref{Sec6:weakUSC}).  
\end{proof}

\subsubsection{Viscosity extension for limiting Hamiltonian operators, the sub-solution case - II }
Next, we extend the $\tilde{H}_0^\theta$ to another slightly simplified new operator:
\footnote{We eliminated the $\epsilon$-dependence in the test functions \eqref{Sec6:f0e} in this step.}
\begin{align*}
\tilde{\tilde{H}}_0^\theta   := \tilde{H}_0^\theta \cup \Big\{ \big(f_0, \tilde{\tilde{G}}^{f_0;\theta, \lambda}_{\phi}\big): f_0:=f_{0,p=2}:=f_{0; \gamma_1, \ldots, \gamma_K} \in {\mathcal S}^+_\sfX, \phi \in {\mathcal F}_0, \lambda >1 \Big\},
\end{align*}
where the
\begin{align}\label{Sec6:GTTilde}
 \tilde{\tilde{G}}^{f_0;\theta, \lambda}_{\phi}(\rho) &  :=
 G^{f_0;\theta, \lambda}_{\phi}(\rho) +  (1-\frac{1}{\lambda})
\Big\{ c - \inf{\sfH} +4C \sup_{x,P,q} |\nabla_q \phi|^2 \\
& \qquad \qquad \qquad \qquad 
+  4C \sup_{\bM \in \Gamma_p^{\opt}(\rho_0; \gamma_1, \ldots, \gamma_K)} \int_{\R^{2d}} |P|^2 \bnu^{\bM}_{f_{0;\gamma_1, \ldots, \gamma_K}}(dx,dP) \Big\}.
 \nonumber
\end{align}

\begin{lemma}\label{Sec6:DtildH}
In context of Lemma~\ref{Sec6:Htild},
the above constructed $f_\theta^* \in \USC\big((\sfX, \sfd);\R\big)$ is a sub-solution to 
\begin{align*}
(I - \alpha \tilde{\tilde{H}}_0^\theta) f_\theta^* \leq h_0
\end{align*}
in the point-wise viscosity sense.  
\end{lemma}
\begin{proof}
Again, the proof follows lines of the method introduced in Lemmas 7.7 and 13.21 in~\cite{FK06}. 
There is an added twist. Using a variant of the argument in Lemma~\ref{Sec3:locF}, we will improve convergence in a weaker sense (in $r$-Wasserstein with $1<r<2$) of extremal points (in definition of the viscosity sub-solutions) to a stronger convergence (in $2$-Wasserstein metric). 

Let
\begin{align*}
f_0:=f_0(\rho) = \psi\big(\sfd^2(\rho, \gamma_1), \ldots, \sfd^2(\rho, \gamma_K)  \big) 
\in  {\mathcal S}^+_\sfX= D(\tilde{\tilde{H}}_0^\theta).
\end{align*}
Let $1< \ldots< p_n < p_{n+1} < \ldots <2$ be such that $\lim_{n \to \infty} p_n =2$.
We approximate the above $f_0$ by
\begin{align*}
f_{n,0}:=f_{n,0}(\rho) = \psi\big(\sfd_{p_n}^2(\rho, \gamma_1), \ldots, \sfd_{p_n}^2(\rho, \gamma_K) \big) \in D(H_0^\theta) \subset D(\tilde{H}_0^\theta).
\end{align*}
It follows that, 
\begin{align*}
f_{n,0} \leq f_{n+1,0} \leq \ldots \leq f_0, \quad \lim_{n \to \infty} f_{n,0} (\rho) = f_0(\rho), \quad \forall \rho \in \sfX.
\end{align*}
By Lemma~\ref{Sec6:subCnvEqn}, there exists $\rho_n \in \sfX$ such that 
$(f_\theta^* - f_{n,0})(\rho_n)=\sup_\sfX(f_\theta^* - f_{n,0})$ and that
\begin{align*}
 \alpha^{-1} (f_\theta^* - h_0)(\rho_n) \leq G_\phi^{f_{n,0};\theta, \lambda} (\rho_n), \quad \forall \theta>0, \lambda>1, \phi \in \mathcal F_0.
\end{align*}

{\bf Part A:} To simplify, we first consider cases where the $\psi:=\psi(r_1, \ldots, r_K)$   satisfies 
\begin{align}\label{Sec6:psiKgro}  
\lim_{r_k \to +\infty} \psi(r_1, \ldots, r_k, \ldots, r_K) =+\infty, \quad \exists k \in \{1,2,\ldots,K\}.
\end{align}

Since $\sup_\sfX f_\theta^*<\infty$, we conclude $\sup_n \int_{\R^d} |x|^{p_n} d\rho_n(dx) <\infty$. By Lemmas~\ref{Sec6:fwUSC} and \ref{Sec6:Attou},  there exists a $\rho_0:=\rho_{0,\theta} \in \sfX \subset \sfX^\prime$
\footnote{The $\rho_0 \in \sfX$ because that, by Fatou's lemma, $\int_{\R^d} |x|^2 d \rho_0 \leq \liminf_n \int_{\R^d} |x|^{p_n} d\rho_n(dx) <\infty$.}
 with $\sfd_{\sfX^\prime}(\rho_n, \rho_0) \to 0$ at least along subsequences, and 
\begin{align*}
 (f_\theta^* - f_0)(\rho_0)=\sup_\sfX(f_\theta^* - f_0).
\end{align*}
If we can also derive that
\begin{align*}
\liminf_{n \to \infty} G_\phi^{f_{n,0};\theta,\lambda} (\rho_n) \leq G_\phi^{f_0;\theta,\lambda}(\rho_0),
\end{align*}
and that 
\begin{align*}
f_\theta^*(\rho_n) \to f_\theta^*(\rho_0), \quad \text{ and }
 f_0(\rho_n) \to  f_0(\rho_0),
\end{align*}
then we can conclude. However, the convergence of $\rho_n$ to $\rho_0$ in $\sfd_{\sfX^\prime}$ is too weak for us to achieve these directly. Next, we introduce yet another perturbation to the above test functions, for such purpose.

As in Lemma~\ref{Sec3:locF}, we introduce 
\begin{align*}
f_0^\epsilon(\rho):= f_0(\rho) + \epsilon \sfd^2(\rho, \rho_0), \quad
 f_{n,0}^\epsilon(\rho):=f_{n,0}(\rho) + \epsilon \sfd^2(\rho,\rho_0) \in D(\tilde{H}_0^\theta).
\end{align*}
This makes $\rho_0$ the unique global strict maximizer of $f_\theta^* - f_0^\epsilon$. By Lemma~\ref{Sec6:Htild}, there exists $\rho_{n,\epsilon}  \in \sfX$ with
\begin{align*}
(f_\theta^* - f_{n,0}^\epsilon)(\rho_{n,\epsilon}) = \sup_{\sfX} (f_\theta^* - f_{n,0}^\epsilon),
\end{align*}
and
\begin{align}\label{Sec6:tildGpreLim}
\alpha^{-1} (f_\theta^* - h_0) (\rho_{n,\epsilon})
  \leq \tilde{G}_\phi^{f_{n,0}^\epsilon; \theta, \lambda}(\rho_{n,\epsilon}),
\quad \forall \phi \in {\mathcal F}_0, \lambda >1, \theta >0,
\end{align}
with the $\tilde{G}_\phi^{f_{n,0}^\epsilon; \theta, \lambda}$ defined in \eqref{Sec6:Gtil}.
Again, from monotone point-wise convergence of $f_{n,0}^\epsilon$ to $f_0^\epsilon$ in the $n \to \infty$ limit, using properties of Gamma convergence (Lemma~\ref{Sec6:Attou}), we have (at least along a subsequence)  
\begin{align*}
\lim_{n \to \infty}\sfd_{\sfX^\prime}(\rho_{n,\epsilon}, \rho_{0,\epsilon}) =0, \quad \exists \rho_{0,\epsilon} \in \sfX;
\end{align*}
and 
\begin{align*}
\lim_{n \to \infty} \sup_{\sfX} (f_\theta^* - f_{n,0}^\epsilon)
=\lim_{n \to \infty} (f_\theta^* - f_{n,0}^\epsilon)(\rho_{n,\epsilon}) 
=   (f_\theta^* - f_0^\epsilon)(\rho_{0,\epsilon}) =  \sup_{\sfX} (f_\theta^* - f_0^\epsilon).
\end{align*}
Since $\rho_0$ is the only one global maximizer of $(f_\theta^* - f_0^\epsilon)$, we conclude that $\rho_{0,\epsilon}=\rho_0$. Indeed, from the above equalities, it follows that
\begin{align}\label{Sec6:2tilHineq}
   \sup_{\sfX} (f_\theta^* - f_0)    =(f_\theta^* - f_0)(\rho_0) 
   & = (f_\theta^* - f_0^\epsilon)(\rho_0) = \sup_{\sfX} (f_\theta^* - f_0^\epsilon) \\
   & = \lim_{n \to \infty} (f_\theta^* - f_{n,0}^\epsilon)(\rho_{n,\epsilon}) \nonumber \\
&  \leq \limsup_{n \to \infty} (f_\theta^* - f_{n,0})(\rho_{n,\epsilon}) 
  - \liminf_{n \to \infty} \epsilon \sfd^2(\rho_{n,\epsilon}, \rho_0) \nonumber \\
&  \leq  (f_\theta^* - f_0)(\rho_0) 
  - \liminf_{p \to 2^-} \epsilon \sfd^2(\rho_{n,\epsilon}, \rho_0), \nonumber
\end{align}
(recall that $\sfd:= \sfd_\sfX$ here). Hence, the $\epsilon$-perturbation created a strong enough coercive effect, we have improved convergence result to
\begin{align*}
 \lim_{n \to \infty} \sfd_\sfX(\rho_{n,\epsilon}, \rho_0)=0, \quad \forall \epsilon >0.
\end{align*}
We have now, for every fixed $\phi \in {\mathcal F}_0, \lambda >1, \theta >0$, that
\begin{align*}
\lim_{n \to \infty} \tilde{G}_\phi^{f_{n,0}^\epsilon; \theta, \lambda}(\rho_{n,\epsilon}) 
\leq  \tilde{G}_\phi^{f_0^\epsilon; \theta, \lambda}(\rho_0). 
\end{align*}
Moreover, note that the $\rho_0$ is chosen independent of the $\epsilon>0$, we also have
\begin{align*}
 \limsup_{\epsilon \to 0^+} \tilde{G}_\phi^{f_0^\epsilon; \theta, \lambda}(\rho_0)
\leq \tilde{\tilde{G}}_\phi^{f_0; \theta,\lambda}(\rho_0).
\end{align*}

From \eqref{Sec6:2tilHineq}, we also obtain 
\begin{align*}
\lim_{n \to \infty} f^*_\theta(\rho_{n,\epsilon}) = f_\theta^*(\rho_0), \quad 
 \lim_{n \to \infty} f_{n,0}^\epsilon(\rho_{n,\epsilon}) = f_0(\rho_0),
\end{align*}
at least along subsequences. Taking limit on \eqref{Sec6:tildGpreLim}, consequently
\begin{align*}
\alpha^{-1}(f^*_\theta - h_0) (\rho_0)  \leq \tilde{\tilde{G}}_\phi^{f_0; \theta, \lambda}(\rho_0).
\end{align*}
We conclude.

{\bf Part B:} Next, we consider the case of general $\psi$ in the test function $f_0$.
We note that the definition of $\tilde{\tilde{G}}^{f_0;\theta, \lambda}_\phi(\rho)$ only involves localness of the $f_0$ at $\rho_0$, in neighborhood induced by the $\sfd$-metric. Hence, standard localization arguments can reduce current situation to that of Part A. 
\end{proof}

\subsubsection{Viscosity extension for limiting Hamiltonian operators, the sub-solution case - III }
Next, in the context of Lemma~\ref{Sec6:DtildH}, we would like to take $\theta \to 0$, so that the $f^*_\theta$ can be replaced by the $f^*$, and the Hamiltonian operator gets further simplified. 
 
 We recall the definitions of $f_{0,p}:=f_0$ in \eqref{Sec6:f0p} and $G^{f_0;p}_\phi:= G^{f_0}_\phi$ in \eqref{Sec6:GGDef} -- we added the parameter $p$ here to emphasize its explicit dependency. We introduced an operator $H_0$ in \eqref{Sec6:H0Def} which is defined on those test functions $f_{0,p}$ with $1<p<2$. Next, we consider the $p=2$ case. To reduce the amount of (already many) notations, 
with a slight abuse of notation, we will still use the $H_0$ by writing
\begin{align}\label{Sec6:H0Rev}
  H_0 := \big\{ (f_0, G_\phi^{f_0; p=2}) : f_0:=f_{0,p=2} \in {\mathcal S}^+_\sfX , 
 \phi \in {\mathcal F}_0 \big\}.
\end{align}

\begin{lemma}\label{Sec6:H0Ext}
Following the context of Lemma~\ref{Sec6:Htild},
the $f^* \in \USC(\sfX;\R)$ defined in \eqref{Sec6:fstarDef} is a viscosity sub-solution in the point-wise sense to 
\begin{align*}
 (I - \alpha H_0) f^* \leq h_0,
\end{align*}
with the $H_0$ given by \eqref{Sec6:H0Rev}.
\end{lemma}
\begin{proof}
Let
\begin{align*}
f_0:=f_0(\rho) = \psi\big(\sfd^2(\rho, \gamma_1), \ldots, \sfd^2(\rho, \gamma_K) \big) \in  D(H_0).
\end{align*}
As in Part B of the proof of Lemma~\ref{Sec6:DtildH}, by localization argument if needed, we proceed next by assuming \eqref{Sec6:psiKgro} holds. 
 
By Lemma~\ref{Sec6:DtildH}, there exists $\rho_\theta \in \sfX$ with
$(f_\theta^* - f_0)(\rho_\theta) = \sup_{\sfX} (f_\theta^* - f_0)$, and
\begin{align}\label{Sec6:Sub3Eq}
\alpha^{-1} (f_\theta^* - h) (\rho_\theta) \leq \tilde{\tilde{G}}_\phi^{f_0; \theta, \lambda}(\rho_\theta),
\quad \forall \phi \in {\mathcal F}_0, \lambda >1,
\end{align}
with the $\tilde{\tilde{G}}_\phi^{f_{K,0;p}; \theta, \lambda}$ defined by~\eqref{Sec6:GTTilde}. 

Since $\sup_\theta \sup_\sfX f^*_\theta <+\infty$, following arguments in the proof of Lemma~\ref{Sec6:DtildH}, we may assume without loss of generality that 
$\{ \rho_\theta \}_\theta$ is relatively compact in $\sfX^\prime$. Selecting sub-sequence if necessary, there exists $\rho_0  \in \sfX$  (note that \eqref{Sec6:psiKgro} holds) with $\lim_{\theta \to 0^+} \sfd_{p_0}(\rho_\theta , \rho_0) =0$. 
We note that the $f^*_\theta$ is $\sfd_{\sfX^\prime}$-upper semi-continuous in $\sfX$  (Lemma~\ref{Sec6:fwUSC}). Also, the definitions of $f^*_\theta$ and $f^*$ implies that 
\begin{align*}
f^*_\theta(\rho) \leq f^*_{\theta^\prime}(\rho) \text{ whenever } 0< \theta^\prime < \theta, 
\quad \lim_{\theta \to 0^+} f_\theta^*(\rho) = f^*(\rho), \qquad \forall \rho \in \sfX.
\end{align*}
Invoking Lemma~\ref{Sec6:Attou}, therefore
 \begin{align*}
\lim_{\theta \to 0^+}(f_\theta^* - f_0)(\rho_\theta) = (f^* - f_0)(\rho_0) = \sup_\sfX (f^* - f_0).
\end{align*}
The above also implies $f_\theta^*(\rho_K) \to f^*(\rho_0)$ and $f_0(\rho_\theta)  \to f_0(\rho_0)$.

We would be able to conclude the proof if 
\begin{align*}
 \limsup_{\theta \to 0^+}   \tilde{\tilde{G}}_\phi^{f_0; \theta, \lambda}(\rho_\theta) \leq  
 \tilde{\tilde{G}}_\phi^{f_0; \theta, \lambda}(\rho_0).
\end{align*}
However, we only have $\rho_\theta \to \rho_0$ in $\sfd_{\sfX^\prime}$, which is too weak for the above to hold.  Again, we go through the perturbative arguments as in Lemma~\ref{Sec6:DtildH} to improve the convergence.

Let  
\begin{align*}
f_0^\epsilon(\rho):= f_0(\rho) + \epsilon \sfd^2(\rho, \rho_0).
 \end{align*}
Invoking Lemma~\ref{Sec3:locF}, the above derived $\rho_0$ is also a global strict maximizer of $f^*_\theta - f_0^\epsilon$ and all of the following hold as consequences: There exists  
$\rho_\theta^\epsilon \in \sfX$ with $(f_\theta^* - f_0^\epsilon)(\rho_\theta^\epsilon) 
= \sup_{\sfX} (f_\theta^* - f_0^\epsilon)$, and
\begin{align*}
\alpha^{-1} (f_\theta^* - h) (\rho_\theta^\epsilon) \leq \tilde{\tilde{G}}_\phi^{f_0^\epsilon; \theta, \lambda}(\rho_\theta^\epsilon), \quad \forall \phi \in {\mathcal F}_0, \lambda >1;
\end{align*}
The sequence $\{ \rho^\epsilon_\theta \}_\theta$ is relatively compact in the $\sfd_{\sfX^\prime}$ metric topology,  with limiting point has to be the $\rho_0$. Consequently, using similar arguments in the proof of Lemma~\ref{Sec6:DtildH} 
(after estimates in \eqref{Sec6:2tilHineq}), we have
\begin{align*}
 \lim_{\theta \to 0^+} \sfd(\rho_\theta^\epsilon, \rho_0) = 0, \quad \forall \epsilon >0,
\end{align*}
and $(f^* - f_0)(\rho_0) = \sup_{\sfX} (f^* - f_0)$.

We recall definitions of $G^{f_0^\epsilon;\theta,\lambda}_\phi$ in \eqref{Sec6:GSubthelam}, $\tilde{\tilde{G}}^{f_0^\epsilon;\theta,\lambda}_\phi$ in \eqref{Sec6:GTTilde}, and the estimates in Lemma~\ref{Sec6:Err12}. Next, taking
 $\limsup_{\lambda \to 1^+} \limsup_{\epsilon \to 0^+} \limsup_{\theta \to 0^+}$ on both sides of the $\epsilon$-perturbed version of \eqref{Sec6:Sub3Eq} (when the $f_0$ is replaced by $f_0^\epsilon$),  we arrive at
\begin{align*}
 \alpha^{-1} (f^* - h) (\rho_0) \leq G_\phi^{f_0;p=2}(\rho_0),
 \quad \forall \phi \in {\mathcal F}_0, \lambda >1.
\end{align*}
We note that the $\rho_0$ is independent of the $\lambda>1$ during the above process. 

We conclude. 
\end{proof}

\subsubsection{Viscosity extension for limiting Hamiltonian operators, the sub-solution case - A summary}
The operator $H_0$ in \eqref{Sec6:H0Rev} can be viewed as a (multi-valued) first order differential operator acting on $\mathcal S^+_\sfX$. In that sense, the $H_0$ is a local operator. Because of this, together with upper semi-continuity regularity of $\rho \mapsto G^{f_0}_\phi(\rho)$ in the $2$-Wasserstein metric, we can verify  \eqref{Sec3:gapp} in Lemma~\ref{Sec3:seq2spw}.  This leads to further strengthening on the notation of sub-solution as obtained in Lemma~\ref{Sec6:H0Ext}, to become strong point-wise viscosity sub-solution.
 
The key to such strengthening is verification of the following property. 
\begin{lemma}\label{Sec6:locF}
 Let $f_0 \in {\mathcal S}_{\sfX}^+$ and $\rho_0 \in \sfX$ satisfies $(\overline{f} - f_0)(\rho_0) = \sup_{\sfX} (\overline{f} - f_0)$. We introduce a perturbation of the $f_0$ by
\begin{align*}
f^\epsilon_0(\rho) := f_0(\rho) + \epsilon \sfd^2(\rho, \rho_0), \quad \forall \epsilon >0.
\end{align*}
 Then $f^\epsilon_0 \in {\mathcal S}_\sfX^+$, and  for each $\epsilon >0$, $\phi \in {\mathcal F}_0$ and $\rho_{n,\epsilon} \to \rho_0$ in $\sfd$ as $n \to \infty$,  we have
\begin{align}\label{Sec6:Gapp}
  \limsup_{n \to 0^+} G^{f^\epsilon_0}_\phi(\rho_{n,\epsilon}) \leq G^{f_0}_\phi(\rho_0).
\end{align}
In the above, the $G^{f^\epsilon_0}_\phi$ is defined by \eqref{Sec6:GGDef}.
\end{lemma}
\begin{proof}
To fix notations, we denote $f_0:= f_{0;\gamma_1,\ldots, \gamma_K} \in {\mathcal S}_\sfX^+$. Then
\begin{align*}
 G^{f^\epsilon_0}_\phi(\rho)& = \sup \Big\{
 \int \Big( \eta^\phi \big(x, \sum_k \alpha_k (x -y_k) + 2\epsilon(x-x_0)\big) \\
 & \qquad \qquad  \qquad - U(x) - (V*\rho)(x) \Big)    \bM( dx; dy_1, \ldots, dy_K; dx_0) :  \bM \in \Gamma^{\opt}(\rho,\gamma_k); \\
 & \qquad \qquad  \qquad \qquad 
 \pi^{1,1+k}_\# \bM \in \Gamma^{\opt}(\rho,\gamma_k), k=1,\ldots,K; \pi^{1,K+2}_\# \bM \in \Gamma^{\opt}(\rho,\rho_0)  \Big\}  
\end{align*}
where the $\alpha_k:=\alpha_k(\rho, \gamma_1, \ldots, \gamma_K)$s are defined as in \eqref{alphaDef}.
In particular, 
\begin{align*}
 G^{f^\epsilon_0}_\phi(\rho_0)&=  \sup \Big\{
 \int \Big( \eta^\phi \big(x, \sum_k \alpha_k (x -y_k)  \big) - U(x) - (V*\rho)(x) \Big)    \bM( dx; dy_1, \ldots, dy_K) :  \\
 & \qquad \qquad  \qquad \bM \in \Gamma^{\opt}(\rho,\gamma_k); \pi^{1,1+k}_\# \bM \in \Gamma^{\opt}(\rho,\gamma_k), k=1,\ldots, K    \Big\} \nonumber \\
 & = G^{f_0^\epsilon}_\phi(\rho_0).
\end{align*}

We note the following property of optimal mass transport: for every $\rho_n, \gamma_n \in \sfX$ and 
 $\bpi_n(dx,d y) \in \Gamma^{\opt}(\rho_n,\gamma_n)$ with
\begin{align*}
 \lim_{n \to \infty} \sfd(\rho_n, \rho_0)+ \sfd(\gamma_n, \gamma_0) =0,
\end{align*}
the sequence $\{ \bpi_n : n=1,2,\ldots\}$  is relatively compact in $\mathcal P_2(\R^{2d})$ in the $2$-Wasserstein topology and any limiting point satisfies $\bpi_0 \in \Gamma^{\opt}(\rho_0, \gamma_0)$.  
In a similarly way, selecting subsequence and relabel if necessary, there exists $\bM_n \in \Gamma(\rho_n; \gamma_1,\ldots, \gamma_K; \rho_0)$ such that $\bM_n \to \bM_0 \in  \Gamma(\rho_0; \gamma_1,\ldots, \gamma_K; \rho_0)$ in the sense of convergence of joint distributions in the order-2 Wasserstein topology in $\mathcal P_2(\R^{(K+2)d})$. Moreover, we can choose the $\bM_n$ such that
\begin{align*}
 G^{f^\epsilon_0}_\phi(\rho_n)&\leq \frac1n + 
 \int \Big( \eta^\phi \big(x, \sum_k \alpha_k (x -y_k) + 2 \epsilon(x-x_0) \big) \\
 & \qquad \qquad  \qquad 
 - U(x) - (V*\rho)(x) \Big)    \bM_n( dx; dy_1, \ldots, dy_K; dx_0).
\end{align*}
By Lemma~\ref{Sec5:etaphi}, $\eta^\phi:=\eta^\phi(x,P) \in C(\R^{2d})$ and has at most quadratic growth at infinity. Consequently, by Fatou's lemma, we have
\begin{align*}
  \limsup_{n \to \infty} G^{f^\epsilon_0}_\phi(\rho_n) 
   \leq G^{f^\epsilon_0}_{\phi}(\rho_0) = G^{f_0}_\phi(\rho_0). 
\end{align*}
\end{proof}

We conclude, in view of all the above extension results for sub-solutions, with the following consolidated result, which we will use later.  
 \begin{lemma}\label{Sec6:seq2spw}
The $f^*$ in \eqref{Sec6:fstarDef} is bounded above and $f^* \in \USC(\sfX; \R)$. Let the $h_0$ satisfy Condition~\ref{Sec6:h0hN} and the requirement in Lemma~\ref{Sec6:Htild}, and has at most $\sfd$-sub-quadratic growth in $\sfX$. Then the $f^*$ is a  sub-solution to  
\begin{align}\label{Sec6:H0sub}
 (I - \alpha H_0)f^* \leq h_0
\end{align}
in the point-wise strong viscosity solution sense, with an extra property that:  for each $f_0 \in D(H_0)$, at least one maximizer $\rho_0 \in \sfX$ is guaranteed to exist  $(f^* - f_0)(\rho_0) = \sup_\sfX(f^* -f_0)$.

Moreover, the sequence $\{\overline{f}_N \}_{N \in \N}$ and $f^*$ satisfy 
Property ${\mathscr P_N}$  in Definition~\ref{Sec4:ProPN}.
\end{lemma}  
\begin{proof}
The result follows by combining Lemmas~\ref{Sec6:subCnvEqn}, \ref{Sec6:Htild}, \ref{Sec6:DtildH}, \ref{Sec6:H0Ext},   and by applying the property verified in Lemma~\ref{Sec6:locF} to Lemma~\ref{Sec3:seq2spw}.  

The Property ${\mathscr P_N}$ for $\{\overline{f}_N \}_{N \in \N}$ and $f^*$ follows from Lemma~\ref{Sec6:Lfstar}.
\end{proof}

\subsection{From multi-valued Hamiltonian operators $H_0$ and $H_1$ to single-valued ones ${\mathbb H}_0$ and ${\mathbb H}_1$}
\label{Sec6:bbHams}
We recall the $\bar{\sfH}(x,P;\rho)$ defined in \eqref{barHxP}:
\begin{align*}
\bar{\sfH}(x,P;\rho):= \bar{\sfH}(P) -U(x) - V*\rho(x).
\end{align*}
 Next, we introduce Hamiltonian operators (noting $\bnu^{\bM}_{f_0}$ is defined as in \eqref{bnu0Def})
 \begin{align}\label{Sec6:bbH0}
 {\mathbb H}_0 f_0(\rho) & := \sup_{\bM \in \Gamma^{\opt}(\rho;\gamma_1,\ldots, \gamma_K, \ldots)}  
 \int_{\R^{2d}} \bar{\sfH}(x, P;\rho)  \bnu^{\bM}_{f_{0;\gamma_1,\ldots, \gamma_K, \ldots}}(dx,dP), \\
 & \qquad \qquad \forall f_0:=f_{0;\gamma_1,\ldots, \gamma_K, \ldots} \in {\mathcal S}^{+,\infty} \text{ as defined in } \eqref{extSS+}, \nonumber 
 \end{align}
and (noting $\bnu^{\bM}_{f_1}$ is defined as in \eqref{Sec2:bnu1Def})
\begin{align}\label{Sec6:bbH1}
 {\mathbb H}_1 f_1(\gamma) & := \inf_{\bM \in \Gamma^{\opt}(\gamma;\rho_1,\ldots, \rho_K, \ldots)}  
 \int_{\R^{2d}} \bar{\sfH}(y, P; \gamma)  \bnu^{\bM}_{f_{1;\rho_1,\ldots, \rho_K, \ldots}}(dy,dP), \\
 & \qquad \qquad \forall f_1:=f_{1;\rho_1,\ldots, \rho_K, \ldots} \in {\mathcal S}^{-,\infty} \text{ as defined in } \eqref{extSS-}. \nonumber
 \end{align}

Both the ${\mathbb H}_0$ and ${\mathbb H}_1$ are single valued operators. By Lemma~\ref{S+SCC}, $\mathcal S^{+, \infty} \subset \SCC(\sfX;\R)$. Using notations and results in Section~\ref{SSmass}, these operators can also be equivalently expressed as follow.   Introducing 
the notation $\oplus$ as in Definition~\ref{LinPTan}, and $\beta_k$ as in \eqref{betaDef},  
for each $\bmu_k \in \exp_\gamma^{-1} (\rho_k)$, we have
$\oplus_{k=1}^\infty (2 \beta_k) \cdot \bmu_k \subset \bpartial_\gamma^{s,-} f_1$.
Note that (by Definition~\ref{Emap}) we have representation
\begin{align}\label{Sec6:pi2nu}
 \bmu_k(dy, dP) := \int_{x \in \R^d} \delta_{x-y}(dP) \bpi_{1,k}(dx,dy), 
 \quad \exists \bpi_{1,k} \in \Gamma^{\opt}(\rho_k,\gamma).
\end{align}
Also, conversely, each optimal plan $\bpi_{1,k}$ defines an inverse exponential map $\bmu_k$ through such identity. 
Therefore
\begin{align}\label{Sec6:bbH1A}
{\mathbb H_1} f_1(\gamma) &= \inf_{\substack{\bmu \in \oplus_{k=1}^\infty 2 \beta_k \cdot \bmu_k \\
   \bmu_k \in \exp^{-1}_{\gamma}(\rho_k)}} \Big\{ \int_{\R^{2d}}  \big( \bar{\sfH}(y,P;\gamma)\big) \bmu(dy,dP) \Big\}.  
\end{align} 
For the case of ${\mathbb H}_0f_0$, we introduce $\alpha_k$s according to \eqref{alphaDef}, then the $\bnu^{\bM}_{f_0} \in \bpartial_\gamma^{s,+} f_0$ (Lemma~\ref{Sec2:nupar}) and
\begin{align}\label{Sec6:bbH0A}
 {\mathbb H}_0 f_0(\rho) & = \sup_{\substack{\bmu \in \oplus_{k=1}^\infty 2 \alpha_k \cdot \bmu_k\\
    \bmu_k \in \exp_{\rho}^{-1}(\gamma_k)}}  
 \int_{\R^{2d}} \bar{\sfH}(x, -P;\rho)  \bmu(dx,dP) \\
 & = \sup_{\substack{\bmu \in \oplus_{k=1}^\infty 2 \alpha_k \cdot \bmu_k\\
    \bmu_k \in \exp_{\rho}^{-1}(\gamma_k)}}  
 \int_{\R^{2d}} \bar{\sfH}(x, P;\rho)  \big((-1) \cdot \bmu\big)(dx,dP). \nonumber
 \end{align}
See Definition~\ref{LinPTan} for definition of $(-1) \cdot \bmu$ and Lemma~\ref{negnu} for some of its properties.

\begin{lemma}\label{Sec6:bbH0Eqn}
Let $h_0 \in C(\sfX)$. 
Suppose that bounded above function $\overline{f} \in \USC(\sfX;\R)$ is a viscosity sub-solution in the point-wise (respectively, strong) sense to 
\begin{align}\label{Sec6:bbH0ext}
 (I - \alpha H_0) \overline{f} \leq h_0,
\end{align}
with the multi-valued operator $H_0$ defined in \eqref{Sec6:H0Rev}.  Then such $\overline{f}$ is also a viscosity sub-solution in the point-wise (respectively, strong) sense to 
\begin{align}\label{Sec6:bbH0Sub}
 (I - \alpha {\mathbb H}_0) \overline{f} \leq h_0.
\end{align}
\end{lemma}
\begin{proof}
We only prove the point-wise viscosity sub-solution sense. The case of strong point-wise viscosity sense can be done in a parallel way.

Let $f_0:=f_{0;\gamma_1,\ldots, \gamma_K} \in {\mathcal S}^{+,\infty}$. By definition of viscosity sub-solution in the point-wise sense for \eqref{Sec6:bbH0ext}, there exists a $\rho_0 \in \sfX$ (by definition of $H_0$, such $\rho_0$ is chosen independently of the $\phi$ below) with 
\begin{align*}
 (\overline{f} - f_0)(\rho_0) = \sup_\sfX( \overline{f} -f_0), \quad 
 (\overline{f} - h_0)(\rho_0) 
 \leq \inf_{\phi \in {\mathcal F}_0}G_\phi^{f_0}(\rho_0);
 \end{align*}
 where the $G_\phi^{f_0}$ is defined in \eqref{Sec6:GGDef} but with the $p=2$.
 Then by Lemma~\ref{App:supintH0} in Appendix,
\begin{align*}
 \inf_{\phi \in {\mathcal F}_0} G_\phi^{f_0}(\rho)
 &= \inf_{\phi \in {\mathcal F}_0} \sup_{\substack{\bmu \in \oplus_{k=1}^\infty 2 \alpha_k \cdot \bmu_k \\
  \bmu_k \in \exp^{-1}_{\rho}(\gamma_k) }}  \int_{\R^d \times \R^d} 
  \big(\eta^\phi(x,P) - U(x) - V*\rho(x)\big) \big((-1)\cdot \bmu\big)(dx,dP) \\
    & = 
 \sup_{\substack{\bmu \in \oplus_{k=1}^\infty 2 \alpha_k \cdot \bmu_k \\
  \bmu_k \in \exp^{-1}_{\rho}(\gamma_k) }}  \int_{\R^d \times \R^d} 
  \bar{\sfH}(x,P;\rho) \big((-1)\cdot \bmu\big)(dx,dP) = {\mathbb H}_0 f_0(\rho).
\end{align*}
 
Let 
\begin{align*}
 {\mathcal N}:= \big\{\bmu \in \oplus_{k=1}^\infty 2 \alpha_k \cdot \bmu_k,
  \bmu_k \in \exp^{-1}_{\rho}(\gamma_k) \big\}.
  \end{align*}
When applying Lemma~\ref{App:supintH0}, we need to show compactness of the $\mathcal N$ as a subset in ${\mathcal P}_2(\R^{2d})$. It is sufficient to show compactness of  $\exp^{-1}_{\rho}(\gamma_k)$ for each $k$ fixed, which is verified as follow. Since the $\rho, \gamma_k \in \sfX$ are fixed, for each $k$, we can find an increasing and convex $\beta \in C(\R_+ ; \R_+)$ with super-linear growth at infinity 
(e.g. Theorem T22 on page 19 of Meyer~\cite{Meyer66}) such that 
\begin{align*}
 \int_{\R^d} \beta(|x|^2) \rho(dx) + \int_{\R^d} \beta(|y|^2) \gamma_k(dy) <\infty,
\end{align*}
hence
\begin{align*}
 \sup_{\bmu \in \exp^{-1}_{\rho}(\gamma_k)} \int_{\R^d \times \R^d} 
 \big(\beta(|x|^2) + \beta(|P|^2) \big) \bmu(dx,dP) <\infty,
\end{align*}
verifying compactness of the $\exp^{-1}_{\rho}(\gamma_k)$. 
\end{proof}

 Next, we consider the super-solution case.

\begin{lemma}\label{Sec6:bbH1Eqn}
Suppose that $\underline{f} \in \LSC(\sfX;\R)$ is a viscosity super-solution in the point-wise sense to 
 \begin{align}\label{Sec6:bbH1ext}
 (I - \alpha H_1) \underline{f} \geq h_1,
\end{align}
with the multi-valued operator $H_1$ defined in \eqref{Sec6:H1Def}, and $h_1 \in C(\sfX)$.  Then such $\underline{f}$ is also a strong viscosity super-solution in the point-wise sense to 
\begin{align}\label{Sec6:bbH1res}
 (I - \alpha {\mathbb H}_1) \underline{f} \geq h_1.
\end{align}
\end{lemma}
\begin{proof}
The proof is  similar to that of the sub-solution case in Lemma~\ref{Sec6:bbH0Eqn}, with the minimax part having some subtle differences. 

By Lemma~\ref{Sec6:H1SupUlt}, such $\underline{f}$ is a point-wise strong super-solution to \eqref{Sec6:bbH1ext} with the $H_1$ defined without those $\zeta$s. Following \eqref{f1gamdef}, we write 
\begin{align*}
f_1:= f_{1;\rho_1,\ldots, \rho_K} \in {\mathcal S}^{-,\infty}.
\end{align*}
Then for every $\gamma \in \sfX$ such that 
\begin{align*}
(f_{1;\rho_1,\ldots, \rho_K}- \underline{f})(\gamma) 
= \sup_\sfX ( f_{1;\rho_1,\ldots, \rho_K}  - \underline{f}), 
\end{align*}
we have
\begin{align*}
( \underline{f} - h_1)(\gamma) \geq \sup_{\phi \in {\mathcal F}_0} G^\phi_{f_{1;\rho_1,\ldots, \rho_K}}(\gamma),
\end{align*}
where the $G^\phi_{f_1}$ is defined in \eqref{Sec6:GinfDef} by setting the $\zeta=0$. 

Next, by Lemma~\ref{App:infintH1} in Appendix,
\begin{align*}
   \sup_{\phi \in {\mathcal F}_0} G^\phi_{f_1}(\gamma)  
 &= \sup_{\phi \in {\mathcal F}_0} \inf_{\substack{\bmu \in \oplus_{k=1}^\infty 2 \beta_k \cdot \bmu_k \\
   \bmu_k \in \exp^{-1}_{\gamma}(\rho_k)}}  \int_{\R^{2d}} \eta_\phi\big(y, P\big)  \bmu(dy,dP)
    - \langle (U + V*\gamma), \gamma\rangle  \\
 & =
 \inf_{\substack{\bmu \in \oplus_{k=1}^\infty 2 \beta_k \cdot \bmu_k \\
   \bmu_k \in \exp^{-1}_{\gamma}(\rho_k)}}
    \int_{\R^{2d}} \sfH(y,P;\gamma) \bmu(dy,dP) = {\mathbb H}_1 f_1(\gamma).
\end{align*}
Hence we conclude the lemma.
\end{proof}

We close this section by stating the following result.
\begin{lemma}\label{Sec6:fsbbH0}
In the context of Lemma~\ref{Sec6:seq2spw}, the $f^*$ is also a point-wise strong viscosity sub-solution to \eqref{Sec6:bbH0Sub}. The $\underline{f}$ in Lemma~\ref{Sec6:SupLem} is also a point-wise strong super-solution to \eqref{Sec6:bbH1res}.
\end{lemma}
\begin{proof}
The sub-solution case follows by combining results of Lemmas~\ref{Sec6:seq2spw} 
and \ref{Sec6:bbH0Eqn}. The super-solution case follows by Lemmas~\ref{Sec6:SupLem} and \ref{Sec6:bbH1Eqn}.
\end{proof}
   
\newpage

\section{Comparison principles for Hamilton-Jacobi equations in space of probability measures}\label{Sec7}
As in the previous section, we denote $\sfX:= {\mathcal P}_2(\R^d)$ and $\sfd$ the Wasserstein order-2 metric on $\sfX$. The map $\rho \mapsto \sfd^2 (\rho,\gamma)$ is a semi-concave function in the sense of Definition~\ref{Sec2:SCCDef} (see Theorem 7.3.2 of \cite{AGS08}). The pair $(\sfX, \sfd)$ forms an Alexandrov metric space with non-negative curvature. Following \eqref{effsfLdef}, we define $\bar{\sfL}(v)$ and introduce $\bar{\sfH}(P)$ as Legendre transform of the convex function $v \mapsto \bar{\sfL}(v)$ (see \eqref{Sec1:barHdef}).  By Proposition~\ref{avgHrep}, under Condition~\ref{PerCND}, we have
\begin{align*}
\bar{\sfH}(P) = \inf_{\varphi \in C^\infty_c(\R^d)} \sup_{q \in \R^d} \sfH(q, P+ \nabla_q \varphi) = \sup_{\varphi \in C^\infty_c(\R^d)} \inf_{q \in \R^d} \sfH(q, P+ \nabla_q \varphi).
\end{align*}
We also recall that, with a slight abuse of notation, we introduced yet another $\bar{\sfH}$-notation through \eqref{barHxP}, 
 \begin{align*}
 \bar{\sfH}(x,P;\rho) := \bar{\sfH}(P) - U(x) - (V* \rho)(x).
\end{align*}

Next, similar to the ${\mathbb H}_0$ and ${\mathbb H}_1$ in Section~\ref{Sec6:bbHams},
we introduce yet another pair of single valued Hamiltonian operators
(noting $\bnu^{\bM}_{f_0}$ is defined as in \eqref{bnu0Def}):
 \begin{align}\label{Sec7:H0}
 {\bf H}_0 f_0(\rho) & := \inf_{\bM \in \Gamma^{\opt}(\rho;\gamma_1,\ldots, \gamma_K, \ldots)}  
 \int_{\R^{2d}} \bar{\sfH}(x, P;\rho)  \bnu^{\bM}_{f_{0;\gamma_1,\ldots, \gamma_K, \ldots}}(dx,dP), \\
 & = \inf_{\substack{\bmu \in \oplus_{k=1}^\infty 2 \alpha_k \cdot \bmu_k\\
    \bmu_k \in \exp_{\rho}^{-1}(\gamma_k)}}  
 \int_{\R^{2d}} \bar{\sfH}(x, -P;\rho)  \bmu(dx,dP)  \label{Sec7:H0Alt}
\\
 & = \inf_{\substack{\bmu \in \oplus_{k=1}^\infty 2 \alpha_k \cdot \bmu_k\\
    \bmu_k \in \exp_{\rho}^{-1}(\gamma_k)}}  
 \int_{\R^{2d}} \bar{\sfH}(x, P;\rho)  \big((-1) \cdot \bmu\big)(dx,dP). \nonumber \\
  & \qquad \qquad \forall f_0:=f_{0;\gamma_1,\ldots, \gamma_K, \ldots} \in {\mathcal S}^{+,\infty} \text{ as defined in } \eqref{extSS+}, \nonumber
 \end{align}
and (noting $\bnu^{\bM}_{f_1}$ is defined as in \eqref{Sec2:bnu1Def})
\begin{align}\label{Sec7:H1}
 {\bf H}_1 f_1(\gamma) & := \sup_{\bM \in \Gamma^{\opt}(\gamma;\rho_1,\ldots, \rho_K, \ldots)}  
 \int_{\R^{2d}} \bar{\sfH}(y, P; \gamma)  \bnu^{\bM}_{f_{1;\rho_1,\ldots, \rho_K, \ldots}}(dy,dP), \\
  &= \sup_{\substack{\bmu \in \oplus_{k=1}^\infty 2 \beta_k \cdot \bmu_k \\
   \bmu_k \in \exp^{-1}_{\gamma}(\rho_k)}} \Big\{ \int_{\R^{2d}}  \big( \bar{\sfH}(y,P;\gamma)\big) \bmu(dy,dP) \Big\} \label{Sec7:H1Alt} \\
    & \qquad \qquad \forall f_1:=f_{1;\rho_1,\ldots, \rho_K, \ldots} \in {\mathcal S}^{-,\infty} \text{ as defined in } \eqref{extSS-}. \nonumber
 \end{align}
 
We have two main results in this section. First, we prove the following comparison principle. Then, in subsection~\ref{Sec7:GradHam}, we  relate ${\bf H}_0$ and ${\bf H}_1$ with a variety of other pairs of Hamiltonian operators (including the ${\mathbb H}_0$ and ${\mathbb H}_1$ introduced in Section~\ref{Sec6:bbHams}). 
 \begin{theorem}\label{Sec7:CMP} 
 Let $\alpha >0$, $h_0, h_1 \in C(\sfX)$ satisfy
  \begin{align*}
 \sup_{\substack{\rho \in \sfX, \\ \sfd(\rho,\delta_0) < R}} h_0 (\rho)
 + \sup_{\substack{\rho \in \sfX, \\ \sfd(\rho,\delta_0) < R}} h_1 (\rho) <+\infty, \quad \forall R \in \R_+.
\end{align*}
Moreover, we assume that at least one of the $h_0$ and $h_1$ has modulus of continuity on each $\sfd$-balls of finite radius. That is, in the case of $h_0$, it holds that, for each $R \in \R_+$, there exists a modulus of continuity $\omega_{h_0;R}$ such that
\begin{align}\label{Sec7:h0mod}
 h_0(\rho) - h_0(\gamma) \leq \omega_{h_0;R}(\sfd(\rho, \gamma)), \quad \forall \rho, \gamma \text{ satisfying } \sfd(\rho, \delta_0) +\sfd(\gamma, \delta_0)\leq R.
\end{align}

Suppose that both
$\overline{f}$ and $\underline{f}$ has sub-linear  growth with respect to the metric $\sfd$, and that Conditions~\ref{PerCND}, \ref{U0CND}, \ref{VCND} hold. Let $\overline{f}$ be a sub-solution, in the point-wise strong viscosity sense, to equation 
\begin{align}\label{Sec7:subHJB}
 \overline{f} - \alpha {\bf H_0} \overline{f} \leq h_0,
\end{align}
and let $\underline{f}$ be a super-solution, in the point-wise strong viscosity sense, to equation  
\begin{align}\label{Sec7:supHJB}
\underline{f} - \alpha {\bf H_1} \underline{f} \geq  h_1.
\end{align}
Then, allowing possibility on right hand side of the following to be $+\infty$, we have
\begin{align}\label{Sec7:CMPineq}
\sup_{\sfX} ( \overline{f} - \underline{f} ) \leq \sup_{\sfX} (h_0- h_1).
\end{align}
\end{theorem}

\subsection{The comparison principle for $\bH_0$ and $\bH_1$ }\label{Sec7:CMP1}
We divide proof of Theorem~\ref{Sec7:CMP} into several parts in this section. 

\subsubsection{A two variable barrier function and its estimates }
Let $\bar{\rho}$ be a fixed probability measure with bounded support (for instance, take 
$\bar{\rho}:=\delta_0$). Let $\epsilon,\kappa, \delta \in (0,1)$, and $\zeta(r):= \sqrt{1+r}$. 
We note $\sup_{r\in \R_+} r(\zeta^\prime(r))^2<\infty$.
 We define a function on $\sfX \times \sfX$ by
 \begin{align}\label{Sec7:Psi}
 \Psi_\delta (\rho, \gamma)  := \Psi_{\epsilon, \kappa, \delta}(\rho,\gamma) 
  & := \frac{1}{1-\epsilon-\kappa} \overline{f}(\rho) - \frac{1}{1+\epsilon} \underline{f}(\gamma) 
  - \frac{1}{2\delta} \sfd^2(\rho,\gamma) 
   -  \frac{\kappa}{1-\epsilon-\kappa} \big(\zeta \circ \sfd^2(\rho, \bar{\rho}) \big).
 \end{align}
By the sub-linear growth condition on $\overline{f}$ and $\underline{f}$, $\sup_{\sfX \times \sfX} \Psi_{\delta}  <\infty$. Let $\rho_\delta, \gamma_\delta \in \sfX$ be such that 
\begin{align}\label{Sec7:appmax}
 \Psi_\delta (\rho_\delta, \gamma_\delta) > \sup_{\sfX \times \sfX} \Psi_\delta -\delta.
\end{align}
Invoke the Borwein-Preiss perturbed optimization principle (Lemma~\ref{BorPre} in Appendix), 
and noting the semi-continuity assumptions on the $\overline{f} \in \USC(\sfX; \R)$ 
and $\underline{f} \in \LSC(\sfX;\R)$, we have the following.
\begin{lemma}\label{BoPrDelta}
There exists $(\hat{\rho}_{\delta,k},  \hat{\gamma}_{\delta,k})$ and
$(\hat{\rho}_\delta, \hat{\gamma}_\delta)$ in the product space $\sfX \times \sfX$ such that 
\begin{align*}
\lim_{k \to \infty} \big( \sfd(\hat{\rho}_{\delta,k}, \hat{\rho}_\delta) + \sfd(\hat{\gamma}_{\delta,k}, \hat{\gamma}_\delta)\big) =0;
\end{align*}
and that the following hold: If we denote
 \begin{align}
 \Delta_\delta(\rho, \gamma)&:= \sum_{k=0}^\infty \frac{1}{2^{k+1}} (\sfd^2(\rho,\hat{\rho}_{\delta,k}) 
 + \sfd^2(\gamma, \hat{\gamma}_{\delta,k})), \label{Sec7:D} \\
 \Psi_{\delta, \Delta}(\rho, \gamma) & := \Psi_\delta(\rho, \gamma) - \sqrt{\delta} \Delta_{\delta}(\rho, \gamma),  \label{Sec7:PdD}
\end{align}
then 
\begin{align}
\Psi_{\delta, \Delta}(\hat{\rho}_\delta, \hat{\gamma}_\delta)
   &= \sup_{\sfX \times \sfX}\Psi_{\delta, \Delta},\label{PsiGmax} \\
\big( \sfd^2(\rho_\delta, \hat{\rho}_\delta) + \sfd^2(\gamma_\delta, \hat{\rho}_\delta) \big)&  \bigvee\sup_k \big(\sfd^2(\hat{\rho}_\delta, \hat{\rho}_{\delta, k})
    + \sfd^2(\hat{\gamma}_\delta, \hat{\gamma}_{\delta,k})\big)    \leq \sqrt{\delta}, \label{ddelt} \\
\Delta_\delta (\hat{\rho}_\delta, \hat{\gamma}_\delta) & \leq \sqrt{\delta},\label{Derr} \\
| D_{\hat{\rho}_\delta} \Delta_\delta(\cdot, \hat{\gamma}_\delta) | 
& = | \sum_{k=0}^\infty\frac{1}{2^{k+1}} D_{\hat{\rho}_\delta} \sfd^2(\cdot, \hat{\rho}_{\delta, k}) |
 \leq 2 \sum_{k=0}^\infty  \frac{1}{2^{k+1}} \sfd(\hat{\rho}_\delta, \hat{\rho}_{\delta, k}) 
  \leq 2 \delta^{1/4}. \label{gradDel}
\end{align}
\end{lemma}
The following is an approximate version of Proposition 3.7 in Crandall, Ishii and Lions~\cite{CIL92}. 
See also Lemma 9.2 of Feng and Kurtz~\cite{FK06}.
\begin{lemma}
For each $\epsilon, \kappa>0$  fixed, we have
\begin{align}\label{dzero}
 \lim_{\delta \to 0^+} \frac{1}{\delta} \sfd^2(\hat{\rho}_{\delta}, \hat{\gamma}_{\delta}) =0.
\end{align}
\end{lemma}
Combine the above with the definition of $\Psi_{\delta, \Delta}$ and sub-linear growth assumption on the $\overline{f}$ and $\underline{f}$, for $\kappa \in (0,1)$ fixed, we have
\begin{align}\label{Sec7:dbdest}
 \limsup_{\epsilon \to 0^+} \limsup_{\delta \to 0^+} \big( \sfd(\hat{\rho}_\delta, \bar{\rho}) 
  + \sfd(\hat{\gamma}_\delta, \bar{\rho}) \big)<+\infty.
\end{align}

Using convexity and at most quadratic growth of $P \mapsto \bar{\sfH}(P)$, we observe the following useful estimates.
\begin{lemma}\label{Sec7:Est1}
There exists a finite constant $C:=C_\sfH >0$, such that
\begin{enumerate}
\item for $f_{0,1} (\rho):=    \zeta (\sfd^2(\rho, \gamma))$ and every $\bM \in \Gamma^{\opt}(\rho;\gamma)$,
\begin{align*}
\int_{\R^{2d}} \bar{\sfH}(P) \bnu^{\bM}_{f_{0,1}}(dx,dP) \leq C(1+ |D_\rho f_{0,1}|^2) \leq 
  C \Big(1+ 4    \sfd^2(\rho,\gamma)\big( \zeta^\prime \circ \sfd^2(\rho,\gamma)\big)^2 )\Big)
      \leq 2 C  ;
\end{align*}
\item for  $f_{0,2}(\rho):= \frac{1-\epsilon-\kappa}{\epsilon} \sqrt{\delta} \Delta_\delta(\rho, \gamma)$ and every $\bM \in \Gamma^{\opt}(\rho;\hat{\rho}_{\delta,1}, \ldots, \hat{\rho}_{\delta,k}, \ldots)$, we have
\begin{align*}
\int_{\R^{2d}} \bar{\sfH}(P)  \bnu^{\bM}_{f_{0,2}}(dx, dP) \leq C(1+ |D_\rho f_{0,2} |^2) 
  \leq C\Big(1+ \sum_{k=0}^\infty \frac{1}{2^{k+1}}  (\frac{1-\epsilon-\kappa}{\epsilon})^2\delta 4 \sfd^2(\rho,\hat{\rho}_{\delta,k}) \Big).  
\end{align*}
\end{enumerate}
\end{lemma}
 
\subsubsection{Estimate on a coupling between ${\bf H_0}$ 
and ${\bf H_1}$}
We construct the $\Psi_\delta$ according to \eqref{Sec7:Psi}, then the $\Psi_{\delta, \Delta}$ as in \eqref{Sec7:PdD}.  We take
\begin{align*}
f_0(\rho) &:= (1-\epsilon-\kappa) \frac{\sfd^2(\rho, \hat{\gamma}_\delta)}{2 \delta}  
 +     \kappa  \zeta\circ \sfd^2(\rho, \bar{\rho})  
   + \epsilon \big(\frac{1-\epsilon-\kappa}{\epsilon} \sqrt{\delta} \Delta_\delta (\rho, \hat{\gamma}_\delta)\big).
 \end{align*}
Then \eqref{PsiGmax} implies that
\begin{align}\label{Sec7:ff0max}
 (\overline{f} - f_0)(\hat{\rho}_\delta) = \sup_\sfX (\overline{f} - f_0).
\end{align}
 By convexity of $P \mapsto \bar{\sfH}(P)$ and the estimates in Lemma~\ref{Sec7:Est1} (combined with \eqref{ddelt}), we have
\begin{align*}
 {\bf H_0}  f_0 (\hat{\rho}_\delta) &  \leq (1-\epsilon-\kappa)  \int_{\R^d \times \R^d} 
  \bar{\sfH}\Big(\frac{x-y}{\delta} \Big){\bpi}(dx; dy) + \kappa (2C)   \\
  & \qquad \qquad  +   \epsilon C \Big(1+ (\frac{1-\epsilon-\kappa}{\epsilon})^2 4 \delta^{\frac{3}{2}} \Big) 
 - \langle U+ V* \hat{\rho}_\delta, \hat{\rho}_\delta \rangle \\
 & \qquad \qquad \qquad 
 \forall \bpi \in \Gamma^{\opt}(\hat{\rho}_\delta, \hat{\gamma}_\delta).
\end{align*}

Similarly,  we consider
 \begin{align*}
 f_1(\gamma) &:= (1+ \epsilon) 
 \big(-\frac{\sfd^2(\hat{\rho}_\delta, \gamma)}{2\delta}  \big) 
    + \epsilon \big( - \frac{1+\epsilon}{\epsilon} \sqrt{\delta} 
     \Delta_\delta (\hat{\rho}_\delta, \gamma)\big) \\
  & = (1+ \epsilon) \big(-\frac{\sfd^2(\hat{\rho}_\delta, \gamma)}{2\delta}  \big) 
  + \epsilon \big( - \frac{1+\epsilon}{\epsilon} \sqrt{\delta} \sum_{k=0}^\infty \frac{1}{2^{k+1}} \sfd^2(\hat{\gamma}_{\delta, k}, \gamma)+ \text{Constant})\big).
\end{align*}
Then 
\begin{align}\label{Sec7:ff1min}
(f_1- \underline{f})(\hat{\gamma}_\delta) = \sup_\sfX ( f_1 - \underline{f}).
\end{align}
Denoting 
\begin{align*}
\rho_1: = \hat{\rho}_\delta, \quad \rho_{2+k}:= \hat{\gamma}_{\delta, k}, \quad k =0,1,2,\ldots,
\end{align*}
we consider each given choice of 
\begin{align*}
\bM :=\bM(dy; dx_1, \ldots, dx_k, \ldots)  \in \Gamma^{\opt}(\gamma; \rho_1, \ldots, \rho_k, \ldots)
 :=\Gamma^{\opt}(\gamma; \hat{\rho}_\delta,  \hat{\gamma}_{\delta, 0}, 
 \hat{\gamma}_{\delta, 1}, \ldots, \hat{\gamma}_{\delta,k}, \ldots).
\end{align*}
Then optimal plans $\bpi_j:=\bpi_j(dx_j,dy) := \pi^{1+j,1}_\# \bM \in \Gamma^{\opt}(\rho_j, \gamma)$ for $j=1,2,\ldots$, and
\begin{align*}
\bmu_k(dy,dP):= \int_{x_k \in\R^d} \delta_{(x_k-y)}(dP)\bpi_k(dx_k,dy) \in \exp_\gamma^{-1}(\rho_k).
\end{align*}
Moreover, letting (Lemma~\ref{supsubD})
\begin{align*}
\tilde{\bmu}(dy,dP) &:= \int_{(x_1, \ldots, x_k, \ldots) \in \R^d \times \ldots \times \R^d \times \ldots} 
 \delta_{\sum_{k=0}^\infty \frac{1}{2^{k+1}} 2 (x_{k+2} - y) }(dP) \bM(dy; dx_1,\ldots, dx_k, \ldots) \\
 & \qquad \qquad \in  \oplus_{k=0}^\infty \frac{1}{2^{k+1}} 2 \cdot \bmu_{k+2} 
  \subset \bpartial_{\hat{\gamma}_\delta}^{s,-} (-\Delta_{\delta})(\hat{\rho}_\delta,\cdot),
\end{align*}
by Remark~\ref{dggamma} and \eqref{dflest} in Lemma~\ref{dxfCCV} (see Definition~\ref{Sec2:Frech}), we have
\begin{align*}
\big((1+\epsilon)\frac{1}{\delta} \cdot \bmu_1\big) 
  \oplus \big(\epsilon (\frac{1+\epsilon}{\epsilon}\sqrt{\delta}) \cdot \tilde{\bmu}\big) \subset \bpartial_\gamma^{s,-} f_1.
\end{align*}
Also note that, by convexity of $\bar{\sfH}$,
\begin{align*}
 \bar{\sfH}(P) & \leq \frac{1}{1+\epsilon}\bar{\sfH}\Big( (1+\epsilon )P + \epsilon Q \Big)  
  + \frac{\epsilon}{1+\epsilon} \bar{\sfH}\big(-Q\big), \quad \forall P,Q \in \R^d.
\end{align*}
Taking
\begin{align*}
P:= \frac{x_1-y}{\delta},  \quad
Q:= \frac{(1+\epsilon)\sqrt{\delta}}{\epsilon} \sum_{k=0}^\infty \frac{1}{2^{k+1}} 2(x_{k+2} -y),
\end{align*} 
therefore,
\begin{align*}
 {\bf H_1}  f_1 (\hat{\gamma}_\delta) & \geq 
 \int \bar{\sfH}\big((1+\epsilon) P +\epsilon Q)\big) d \bM(dy; dx_1, \ldots, dx_k, \ldots)
     - \langle U + V* \hat{\gamma}_\delta, \hat{\gamma}_\delta \rangle \\
 &\geq   (1+ \epsilon) \int_{\R^{2d}}  
 \bar{\sfH}( P )  \big(\frac{1}{\delta} \cdot \bmu_1 \big)(dy,dP)   
  - \epsilon \int_{\R^{2d}} \bar{\sfH}(Q)
  \big(  (\frac{1+\epsilon}{\epsilon}\sqrt{\delta}) \cdot \tilde{\bmu}\big)(dy, dQ)  \\ 
 & \qquad \qquad  - \langle U+ V* \hat{\gamma}_\delta, \hat{\gamma}_\delta \rangle  \\
 & \geq (1+ \epsilon) \int_{\R^{2d}}  
 \bar{\sfH}\big( \frac{x-y}{\delta}\big)\bpi(dx, dy) 
   - \epsilon C \big(1+ (\frac{1+\epsilon}{\epsilon})^2  \delta^{3/2} \big)
     - \langle U+V* \hat{\gamma}_\delta, \hat{\gamma}_\delta \rangle \\
 & \qquad \qquad    \qquad \qquad \qquad \forall 
  \bpi \in \Gamma^{\opt}(\hat{\rho}_\delta, \hat{\gamma}_\delta),
\end{align*}
where, in the last step above, we used estimate \eqref{ddelt}.

Consequently, the following estimate follows
\begin{lemma}\label{H0DH1}
  \begin{align*}
& \frac{1}{1-\epsilon-\kappa} \bH_0 f_0(\hat{\rho}_\delta) 
 - \frac{1}{1+\epsilon} \bH_1 f_1(\hat{\gamma}_\delta) \\
& \leq \frac{2C\kappa}{1-\epsilon -\kappa}
 + \epsilon \Big( \frac{1}{1-\epsilon-\kappa} 
 + \frac{1}{1+\epsilon} \Big)  
C \Big(1 + 4 (\frac{1+\epsilon}{\epsilon})^2  \delta^{3/2}   \Big) \\
& \qquad \qquad 
 +\frac{1}{1+\epsilon}\langle U + V*\hat{\gamma}_\delta, \hat{\gamma}_\delta \rangle
  -  \frac{1}{1-\epsilon-\kappa} \langle U + V*\hat{\rho}_\delta, \hat{\rho}_\delta \rangle.
\end{align*}
\end{lemma}

\subsubsection{The comparison principle}
 \begin{lemma}
The comparison principle stated in Theorem~\ref{Sec7:CMP} holds.
\end{lemma}
\begin{proof}
Following the above constructions, because of \eqref{Sec7:ff0max},  by the strong viscosity sub-solution property (in the point-wise sense),
\begin{align*}
 \alpha^{-1} \big( \overline{f} - h_0 \big)(\hat{\rho}_\delta) \leq {\bf H_0} f_0(\hat{\rho}_\delta).
\end{align*}
Similarly, because of \eqref{Sec7:ff1min},
\begin{align*}
\alpha^{-1} \big( \underline{f} - h_1\big)(\hat{\gamma}_\delta) 
\geq {\bf H_1} f_1(\hat{\gamma}_\delta).
\end{align*}
Consequently
\begin{align*}
 {\rm I}:=  \frac{1}{1-\epsilon -\kappa}\overline{f}(\hat{\rho}_\delta) 
  - \frac{1}{1+\epsilon} \underline{f}(\hat{\gamma}_\delta)
  &  \leq \alpha \Big( \frac{1}{1-\epsilon -\kappa}(\bH_0 f_0)(\hat{\rho}_\delta)
     - \frac{1}{1+\epsilon} (\bH_1 f_1) (\hat{\gamma}_\delta)\Big) \\
     & \qquad \qquad +\Big( \frac{1}{1-\epsilon -\kappa}h_0 (\hat{\rho}_\delta) 
  - \frac{1}{1+\epsilon} h_1(\hat{\gamma}_\delta)\Big) \\
   & =: {\rm II} + {\rm III}.
\end{align*}
We note that right hand side above can be estimated by Lemma~\ref{H0DH1}.

On one hand, by definition of the $\Psi_{\delta, \Delta}$ and \eqref{PsiGmax}, 
for every $\rho\in \sfX$ fixed,
\begin{align*}
&  \frac{1}{1-\epsilon -\kappa}\overline{f}(\rho) - \frac{1}{1+\epsilon} \underline{f}(\rho)
  -\frac{\kappa}{1-\epsilon - \kappa} \zeta\circ \sfd^2(\rho,\bar{\rho})
  - \sqrt{\delta} \Delta_\delta(\rho,\rho) \\
  & \qquad \leq \Psi_{\delta,\Delta}(\rho, \rho)  
     \leq \Psi_{\delta,\Delta}(\hat{\rho}_\delta, \hat{\gamma}_\delta) 
    \leq \Psi_\delta (\hat{\rho}_\delta, \hat{\gamma}_\delta) \leq \rm{I}.
\end{align*}
On the other hand, in view of \eqref{Sec7:dbdest}, \eqref{Sec7:h0mod} and \eqref{dzero}, there exists $R:=R(\kappa) \in \R_+$ such that 
\begin{align*}
{\rm III} & \leq  \frac{1}{1-\kappa} \limsup_{\epsilon \to 0^+}\limsup_{\delta \to 0^+}
  \big( h_0(\hat{\rho}_\delta) - h_1(\hat{\gamma}_\delta)\big) \\
   & \leq \frac{1}{1-\kappa} \Big(  \limsup_{\epsilon \to 0^+}\limsup_{\delta \to 0^+}
  \omega_{h_0;R} \big(\sfd(\hat{\rho}_\delta, \hat{\gamma}_\delta)\big) 
  + \sup_{\gamma \in \sfX : \sfd(\gamma, \delta_0) \leq R} (h_0 - h_1)(\gamma)\Big)  \\
& \leq \frac{1}{1-\kappa}  \sup_\sfX(h_0-h_1).
\end{align*}
In deriving the second inequality above, we used \eqref{Sec7:h0mod}. If the modulus of continuity assumption was on the $h_1$, similar proof still holds, giving the third inequality.
Next, we estimate the term ${\rm II}$.
Again, in view of \eqref{Sec7:dbdest} and \eqref{dzero}, there exists a limiting point $\tilde{\rho}_\kappa \in \sfX$ such that 
\begin{align*}
 \hat{\rho}_\delta:= \hat{\rho}_{\delta,\epsilon,\kappa} \to \tilde{\rho}, 
 \quad \hat{\gamma}_\delta:=\hat{\gamma}_{\delta,\epsilon, \kappa} \to \tilde{\rho}_\kappa, \quad \text{ in $r$-Wasserstein metric, with any $r \in (0,2)$}. 
\end{align*}
Invoking the estimate in Lemma~\ref{H0DH1}, 
 \begin{align*}
 \limsup_{\epsilon \to 0^+}\limsup_{\delta \to 0^+} {\rm II} \leq 
 (1- \frac{1}{1-\kappa})  \langle U+V* \tilde{\rho}_\kappa, \tilde{\rho}_\kappa\rangle
 \leq (1- \frac{1}{1-\kappa}) (\inf_{\R^d} U + \inf_{\R^d} V).
\end{align*}

Combine the above estimates on  ${\rm I}, {\rm II}$ and ${\rm III}$ together, 
\begin{align*}
\frac{1}{1  -\kappa}\overline{f}(\rho) -   \underline{f}(\rho) 
 - \frac{\kappa}{1-\kappa} \zeta\circ \sfd^2(\rho,\bar{\rho}) 
  \leq (1- \frac{1}{1-\kappa}) (\inf_{\R^d} U + \inf_{\R^d} V)
   + \frac{1}{1-\kappa}  \sup_\sfX(h_0 - h_1).
 \end{align*}
We take $\lim_{\kappa \to 0^+}$ and conclude by the arbitrariness of $\rho \in \sfX$.
\end{proof}
 
 \subsection{Other forms of Hamiltonians}
 It will be useful to introduce several other related Hamiltonian operators and study their relations.
 \subsubsection{Hamiltonians as dual from effective Lagrangians}
We recall the definitions of $\bar{\sfL}$ and $\bar{\sfH}$ as in \eqref{effsfLdef} and \eqref{Sec1:barHdef}. By duality of the Legendre transforms, we also have 
\eqref{Sec1:barsfL}. With a slight abuse of notation, we wrote in \eqref{Sec1:LUV} and 
\eqref{Sec1:Lbnu}:
\begin{align*}
 \bar{\sfL}(x,v;\rho):= \bar{\sfL}_{U,V}(x,v ; \rho) : =  \bar{\sfL}(v) + U(x) + (V* \rho)(x), \quad \forall v \in \R^d, \rho \in \sfX,
\end{align*}
and
\begin{align*}
 L(\bnu) := \int_{\R^{2d}} \bar{\sfL}_{U,V}(x,v; \pi_\#^1 \bnu)   \bnu(dx, dv), \quad \forall \bnu \in {\mathcal P}_2(\R^d \times \R^d).
 \end{align*}
From Lemma~\ref{PbarH}, it follows that the $\bar{\sfL}$ has at most quadratic growth at infinity. Consequently, $L(\bnu)$ is always finite for the above choice of $\bnu$.
Let $f \in {\mathcal S}^{+, \infty} \cup {\mathcal S}^{-,\infty}$, by Lemma~\ref{Dmdist} and Remark~\ref{Sec2:Dmdext}, Lemma~\ref{WassDSS} and Remark~\ref{dggamma}, 
$d_\rho f$ exists and can be explicitly expressed. 
With these in mind, we introduce yet anther single-valued Hamiltonian operator in \eqref{Sec1:bHDef}:
\begin{align}\label{Sec7:bHDef}
{\bH} f(\rho):=  \sup \Big\{ (d_{\rho} f)(\bnu) - L(\bnu) : \bnu \in \Tan_\rho \Big\},
\quad \forall f \in{\mathcal S}^{+, \infty} \cup {\mathcal S}^{-,\infty}.
\end{align}
Following Section~\ref{appmass},   we recall the definition of $G(\rho)$ and that $\Tan_\rho:= \overline{G(\rho)}^{\sfd_\rho}$. We have the following.
\begin{lemma}\label{Sec7:HfG} 
For every $f \in {\mathcal S}^{+,\infty} \cup {\mathcal S}^{-,\infty}$,
\begin{align*}
{\bf H} f(\rho) =   \sup \Big\{ \big(d_\rho f\big) (\bmu) - L(\bmu) : \bmu \in G(\rho) \Big\},
\quad \forall \rho \in \sfX .
\end{align*}
\end{lemma}
\begin{proof}
 We note that, by Lemma~\ref{PbarH}, $\bar{\sfL} \in C(\R^d)$ and has at most quadratic growth at infinity. Hence, for each $\rho :=\pi^1_\# \bmu$ fixed, 
\begin{align*}
 \Tan_\rho \ni \bmu \mapsto \int_{\R^d \times \R^d} \bar{\sfL}(v) \bmu(dx,dv)
\end{align*}
is continuous under topology generated by the tangent cone metric $\sfd_\rho(\cdot, \cdot)$. 

Therefore, the map $\bmu \mapsto L(\bmu)$ is continuous in the cone space $\Tan_\rho$. Moreover, by Lemma~\ref{dDpm}, $\bmu \mapsto (d_\rho f_1)(\bmu)$ is also continuous in $\Tan_\rho$. Therefore, the conclusion follows by density of the $G(\rho)$ in $\Tan_\rho$.
\end{proof}

In the following proofs, to simplify notation, we only write the operators as if the test functions are in ${\mathcal S}^+ \cup {\mathcal S}^-$. The general cases follow by notationally allowing the $K=+\infty$. 
 
\begin{lemma}\label{Sec7:HleqH0}
We have
\begin{align*}
\bH f_0 \leq \bH_0 f_0, \qquad \forall f_0 \in \mathcal S^{+,\infty}.
\end{align*}
\end{lemma}
\begin{proof}
Let $f_0$ be as in \eqref{mtranf}, we recall the expression of $\bH_0 f_0$ in \eqref{Sec7:H0Alt} with  notation $\alpha_k$ defined in \eqref{alphaDef},
\begin{align}\label{Sec7:bH0Alt}
{\bH}_0 f_0(\rho)& :=  \inf_{\substack{\bmu \in \oplus_{k=1}^K 2 \alpha_k  \cdot \bmu_k \\ \text{ with } \bmu_k \in \exp_\rho^{-1}(\gamma_k)}}
  \int_{\R^d \times \R^d} \bar{\sfH}(x,-P;\rho) \bmu(dx,dP)   \\
& = \inf_{\bM_0 \in \Gamma^{\opt}(\rho; \gamma_1, \ldots, \gamma_K)} 
 \int_{\R^d \times \R^d} \bar{\sfH}\Big(x,  \sum_{k=1}^K 2 \alpha_k  (x-y_k);\rho \Big) \bM_0(dx; dy_1, \ldots, dy_K). \nonumber
 \end{align}

Writing 
\begin{align*}
\text{Adm}_\bnu:=\Big\{ \bM \in {\mathcal P}_2(\R^{(K+2)d}): \pi^{1,2, \ldots, K+1}_\# \bM \in \Gamma^{\opt}(\rho; \gamma_1, \ldots, \gamma_K),  \pi^{1,K+2}_\# \bM = \bnu \Big\}.
\end{align*}
Then, by Lemma~\ref{WassDSS},
\begin{align*}
& \bH f_0(\rho) \\
& = \sup_{\bnu \in \Tan_\rho} \inf_{\bM \in \text{Adm}_\bnu}\int_{\R^{(K+2)d}} \Big(\sum_{k=1}^K 2 \alpha_k (x-y_k) \cdot v - \bar{\sfL}_{U,V}(x,v;\rho)\Big) \bM(dx, dy_1, \ldots, dy_K; dv) \\
 & \leq   \inf_{ \substack{ \bM_0 =  \pi^{1,2,\ldots, K+1}_\#\bM, \\  \bM \in \text{Adm}_\bnu} }
 \Big\{   \int_{\R^{(K+2)d}} \bar{\sfH}\Big(x, \sum_{k=1}^{K} 2 \alpha_k (x-y_k); \rho\Big) 
 \bM_0(dx, dy_1, \ldots, dy_K) \Big\}   \\
 & = \bH_0 f_0 (\rho).
\end{align*} 
\end{proof}

\begin{lemma}\label{Sec7:HgeqH1}
We have
\begin{align*}
\bH f_1 \geq \bH_1 f_1, \qquad \forall f_1 \in \mathcal S^{-,\infty}.
\end{align*}
\end{lemma}
\begin{proof}
For the given $f_1:=f_1(\gamma)$ (see the $g$ in \eqref{mtrang} and the $\beta_k$s in \eqref{betaDef}), we denote
\begin{align*}
 P:=P(y, x_1, \ldots, x_K; \gamma, \rho_1, \ldots, \rho_K) 
  := \sum_{k=1}^K \beta_k 2 (x_k-y).
\end{align*}
For every $\epsilon>0$, there exists a measurable $v_\epsilon:=v_\epsilon(y,x_1, \ldots, x_K; \gamma, \rho_1, \ldots, \rho_K) \in \R^d$, with proper integrability as needed below, such that 
\begin{align*}
 \bar{\sfH}\big(y, P; \gamma \big) \leq \epsilon + P v_\epsilon -\bar{\sfL}(y, v_\epsilon; \gamma).
\end{align*}

By definition of ${\bf H}_1$ in \eqref{Sec7:H1}, therefore
\begin{align*}
& \bH_1 f_1(\gamma) \\
& = \sup_{\bM_0 \in \Gamma^{\opt}(\gamma; \rho_1, \ldots, \rho_K)}
\int_{\R^{(2+K)d}} \bar{\sfH}\big(y, P(y,x_1, \ldots, x_k;\gamma, \rho_1, \ldots, \rho_K); \gamma\big) \bM_0(dy;dx_1, \ldots, dx_K) \\
& \leq \epsilon + \sup_{\bM_0 \in \Gamma^{\opt}(\gamma; \rho_1, \ldots, \rho_K)}
\int_{\R^{(2+K)d}} \big( P v_\epsilon - \bar{\sfL}(y, v_\epsilon; \gamma)  \big) \bM_0(dy;dx_1, \ldots, dx_K).
\end{align*}

Let $\xi_\epsilon:=\xi_\epsilon(y,x_1,\ldots,x_K; \gamma, \rho_1, \ldots, \rho_K)
:=\Pi_\gamma(v_\epsilon)$ (which as a function of $y$ belongs to $L^2_{\nabla, \gamma}(\R^d)$) be the projected vector field as obtained from Lemma~\ref{defPiproj}. Then by Lemmas~\ref{Bproj} and \ref{p4coupling},  for each 
$\bM_0 \in \Gamma^{\opt}(\gamma; \rho_1, \ldots, \rho_K)$, we have
\begin{align*}
 \int P v_\epsilon d \bM_0 = \int P \xi_\epsilon d \bM_0
\end{align*}
Moreover, by Lemma~\ref{JensenLike},
\begin{align*}
\int \bar{\sfL}(y, \xi_\epsilon; \gamma) d \bM_0  \leq 
\int  \bar{\sfL}(y, v_\epsilon; \gamma) d \bM_0 .
\end{align*}
Consequently, defining 
\begin{align*}
\bM(dy;dx_1, \ldots dx_K; d\xi):= \bM_0(dy; dx_1, \ldots, dx_K) 
\delta_{\xi_\epsilon(y,x_1,\ldots,x_K; \gamma, \rho_1, \ldots, \rho_K)}(d\xi).
\end{align*}
Then $\bnu_\epsilon:=\pi_\#^{1,K+2} \bM \in \Tan_\gamma$ according to Lemma~\ref{Bproj}, and by Remark~\ref{dggamma},
\begin{align*}
 (d_\gamma f_1)(\bnu_\epsilon) \geq \int P v_\epsilon d \bM_0 = \int P \xi_\epsilon d \bM_0, \quad
 \int  \bar{\sfL}(y, \xi_\epsilon; \gamma) d \bM_0 = \bar{\sfL}(\bnu_\epsilon).
\end{align*}

Combine the above, we have 
\begin{align*}
& \bH_1 f_1(\gamma) \\
& \leq \epsilon + \sup_{\bM_0 \in \Gamma^{\opt}(\gamma; \rho_1, \ldots, \rho_K)}
\int_{\R^{(2+K)d}} \big( P \xi_\epsilon - \bar{\sfL}(y, \xi_\epsilon; \gamma)  \big) \bM_0(dy;dx_1, \ldots, dx_K) \\
& \leq \epsilon +  (d_\gamma f_1)(\bnu_\epsilon) - \bar{\sfL}(\bnu_\epsilon)
 \leq \epsilon + \bH f_1(\gamma).
\end{align*}
\end{proof}

\begin{lemma}\label{Sec7:CMP2}
In the context of Theorem~\ref{Sec7:CMP},  let $\overline{f}$ be a sub-solution, in the strong point-wise sense, for \eqref{Sec7:subHJB} with the $\bH_0$ replaced by $\bH$ with domain consisting of  {\em only} test functions in ${\mathcal S}^{+,\infty}$. Let $\underline{f}$ be a super-solution, in the strong point-wise sense, for \eqref{Sec7:supHJB} with the $\bH_1$ replaced by $\bH$ with domain consisting of  {\em only} test functions in ${\mathcal S}^{-,\infty}$.
Then the comparison principle \eqref{Sec7:CMPineq} still holds. 
\end{lemma}
\begin{proof}
Conclusion follows from combining results in Lemmas~\ref{Sec7:HleqH0}, \ref{Sec7:HgeqH1} and Theorem~\ref{Sec7:CMP}.
\end{proof}

\subsubsection{Hamiltonian operators expressed using gradients}\label{Sec7:GradHam}
Let $f_0 \in \mathcal S^{+,\infty}$, by Lemma~\ref{UniGrad},  $\bmu := \grad_\rho f_0 \in \Tan_\rho$ in the sense of Definition~\ref{GradDef} exists and is unique, it can be explicitly identified through Lemma~\ref{WassGf}. 
We define   
\begin{align}\label{Sec7:Hcirc}
 {\bf H}^\circ_0 f_0 (\rho):= \int_{\R^{2d}}   \bar{\sfH}(x, P;\rho) \bmu(dx, dP),   \text{ where } \bmu = \grad_\rho f_0;
\end{align}
and 
\begin{align}\label{Sec7:H2circ}
 {\bf H}^{\circ \circ}_0 f_0 (\rho)
 := \sup \Big\{ \langle \grad_\rho f_0, \bnu \rangle_\rho - L(\bnu) : 
    \bnu \in \Tan_\rho \Big\}.
\end{align}
Noting $\grad_\rho f_0$ admits representation \eqref{Sec2:gradf}, it follows from Lemma~\ref{WassNisio} that ${\bf H}_0^\circ = {\bf H}_0^{\circ \circ}$.
More generally, we have the following.
\begin{lemma}\label{Sec7:Hcom}
\begin{align*}
 {\bf H} f_0  \leq {\bf H}_0 f_0  \leq {\bf H}_0^{\circ \circ} f_0
  = {\bf H}_0^\circ f_0 \leq {\mathbb H}_0 f_0 , \quad \forall f_0 \in {\mathcal S}^{+,\infty}.
\end{align*}
\end{lemma}
\begin{proof}
The inequality ${\bf H} f_0  \leq {\bf H}_0 f_0$ was already established in Lemma~\ref{Sec7:HleqH0}. Next, we prove 
\begin{align*}
{\bf H}_0 f_0 \leq {\bf H}_0^\circ f_0 \leq {\mathbb H}_0 f_0,
\end{align*}
which is a consequence of \eqref{defnuu0}.
 
For the given $f_0:=f_{0;\gamma_1, \ldots, \gamma_K}$, 
let $\bpi_{0,k} \in \Gamma^{\opt}(\rho;\gamma_k)$ be those minimizers uniquely defined by \eqref{uniMini}. We denote
\begin{align}\label{Sec7:bM0}
\bM_0(dx; dy_1, \ldots, dy_K) := 
\bpi_{0,1}(dy_1|x) \times  \ldots \times \bpi_{0,K}(dy_K|x) \rho(dx) 
\in \Gamma^{\opt}(\rho;\gamma_1, \ldots, \gamma_K),
\end{align}
and following \eqref{bnu0Def}, we also denote
\begin{align*}
\bnu^{\bM_0}_{f_0}(dx,dP)& := \int_{(y_1, \ldots, y_K) \in (\R^d)^K} \delta_{\sum_{k=1}^K 2 \alpha_k(x-y_k)}(dP) \bM_0(dx; dy_1, \ldots, dy_K) \\
& = \delta_{u_{f_0}(x)}(dP)\rho(dx),
\end{align*}
with
\begin{align*}
 u_{f_0}:=u_{f_0}(x):= -\sum_{k=1}^K 2 \alpha_k v_k(x), 
  \quad \text{ where } v_k(x):= \int_{\R^d} (y-x) \bpi_{0,k}(dy|x).
\end{align*}
Then, by \eqref{defnuu0} in Lemma~\ref{WassGf}, we identify 
\begin{align*}
\grad_\rho f_0 =\bnu_{f_0}^{\bM_0}.
\end{align*}

Consequently,
\begin{align*} 
 {\bf H}_0^\circ f_0(\rho)  & = \int_{\R^d} \bar{\sfH}\Big(x, P ;\rho \Big)\bnu_{f_0}^{\bM_0}(dx, dP) \\
& \leq \sup_{\bM \in \Gamma^{\opt}(\rho;\gamma_1, \ldots, \gamma_K)}  \int_{\R^d \times \R^d} \bar{\sfH}\big(x,P;\rho \big) \bnu_{f_0; \gamma_1, \ldots, \gamma_K}^{\bM}(dx,dP) = \mathbb{H}_0 f_0(\rho).
\end{align*}
and  
\begin{align*}
{\bf H}_0^\circ f_0(\rho)   \geq \inf_{\bM \in \Gamma^{\opt}(\rho;\gamma_1, \ldots, \gamma_K)}  \int_{\R^d \times \R^d} \bar{\sfH}\big(x,P;\rho \big) \bnu_{f_0; \gamma_1, \ldots, \gamma_K}^{\bM}(dx,dP) = {\bf H}_0 f_0(\rho).
\end{align*}
\end{proof}

 Next, we consider the super-solution case. 
Recall that the definition of $\grad f_0$ is based upon semi-concavity of the function $f_0$. However, $f_1 \in {\mathcal S}^{-,\infty}$ is semi-convex. Hence, if we were to use properties that come out of Definition~\ref{GradDef}, we may use a different definition 
\begin{align*}
 \grad_\gamma f_1 := (-1) \cdot \big( \grad_\gamma (-f_1)\big).
\end{align*}
See Definition~\ref{LinPTan} for $(-1) \cdot \bM$ for $M \in {\mathcal P}_2(\R^{2d})$.
Then
\begin{align*}
\big( \grad_\gamma f_1\big)(dy,dv) = \delta_{\sum_{k=1}^\infty 2 \beta_k u_k(y)} (dv) \gamma(dy), \quad u_k(y) := \int_{\R^d} (x-y) \bpi_{0,k}(dx|y), 
\end{align*}
where the $\bpi_{0,k} \in \Gamma^{\opt}(\rho_k,\gamma)$ is the unique minimizer to 
\begin{align*}
\tilde{s}_k^2& := \int_{\R^d} |\int_{\R^d} x \bpi_{0,k}(dx|y) - y|^2 \gamma(dy) \\
& = \inf \Big\{ \int_{\R^d} |\int_{\R^d} x \bpi(dx|y) - y|^2 \gamma(dy) :
  \bpi \in  \Gamma^{\opt}(\rho_k,\gamma)\Big\}.
\end{align*}
We see that there is an asymmetry between the expressions of $\grad f_0$ and so defined $\grad f_1$.
 
Similar to \eqref{Sec7:Hcirc} and \eqref{Sec7:H2circ}, for $f_1 \in {\mathcal S}^{-,\infty}$, we introduce
\begin{align*}
{\bf H}_1^\circ f_1(\gamma) := \int_{\R^{2d}} \bar{\sfH}(x,P;\gamma) \bmu(dx,dP), \quad \bmu = \grad_\gamma f_1,
\end{align*}
and
\begin{align*}
{\bf H}_1^{\circ \circ} f_1(\gamma):=
 \sup \Big\{ (\langle \grad_\gamma f_1, \bmu\rangle_\gamma - L(\bmu) : \bmu \in \Tan_\gamma \Big\}.
\end{align*}
Then ${\bf H}_1^\circ = {\bf H}_1^{\circ \circ}$ by Lemma~\ref{WassNisio}. We also have
\begin{align}\label{Sec7:IneqH1s}
{\bf H} f_1 \geq {\bf H}_1 f_1 \geq {\bf H}_1^{\circ \circ} f_1 = {\bf H}_1^\circ f_1
  \geq {\mathbb H}_1 f_1.
\end{align}

\subsubsection{Further generalized notions of Hamiltonian using sub- super-gradients}
In a similar vein to the above arguments in this section, one can introduce even more Hamiltonians defined through Fr\'echet super- and sub-gradients (Definition~\ref{Sec2:Frech}). Lemmas~\ref{sizesD}, \ref{supsubD} and \ref{GradS} etc offer insights as to how to relate these definitions together. We do not pursue these anymore in this paper.

\newpage 
\section{Convergence of solutions for Hamilton-Jacobi equations arising from the hydrodynamic limit}
\label{Sec8}
  
In Section~\ref{finPartH}, we considered value functions of minimal action finite particle Lagrangian dynamics, with particle permutation symmetry. They are solutions to Hamilton-Jacobi PDEs in finite dimensions.  In Section~\ref{Sec6}, we derived hydrodynamic limit of the corresponding Hamiltonian operators (on functions defined over space of probability measures). We proved that upper- and lower-limits of the value functions as sub- and super-solutions of PDE in space of probability measures given by a pair of Hamiltonian operators $H_0,H_1$ (and by another pair ${\mathbb H}_0, {\mathbb H}_1$). The operators are defined by \eqref{Sec6:H0Rev} and \eqref{Sec6:H1Def} (respectively, by \eqref{Sec6:bbH0} and \eqref{Sec6:bbH1}). For summary of these results, see Lemmas~\ref{Sec6:seq2spw}, \ref{Sec6:fsbbH0}, and \ref{Sec6:SupLem}. In Section~\ref{Sec7}, we proved comparison principle for sub- and super- solutions of respective equations \eqref{Sec7:subHJB} and \eqref{Sec7:supHJB} with a different, yet another, pair of Hamiltonian operators ${\bf H}_0$ and ${\bf H}_1$. See also comparison result on the operator $\bH$ in Lemma~\ref{Sec7:CMP2}. In Section~\ref{Sec7:GradHam}, we introduced even more pairs of Hamiltonian operators which are natural in such context, and compared with the above ones. 

We note that, however, these results do not allow us to conclude any comparison between viscosity solutions for equations given by the pair of operators $H_0$ and $H_1$. 
In this section, we develop a technique on viscosity extension (first introduced in ~\cite{FK06}) for such purpose.
 
\subsection{A technical problem, intuitively explained}
To explain the intricacies among the above mentioned operators, we take a look at the $H_0 f_0(\rho)$ with a simple choice of test function $f_0(\rho):= \frac{\alpha}{2} \sfd^2(\rho,\gamma)$, where the $\gamma \in {\mathcal P}_2(\R^d)$ and $\alpha>0$ are fixed.  $H_0f_0$ is a set consisting of $G^{f_0}_\phi$ as elements, with arbitrary $\phi  \in \mathcal F_0$, as defined in \eqref{Sec6:GGDef} and \eqref{Sec6:H0Rev}. The $G^{f_0}_\phi$ involves a term $\sup_{\bpi \in \Gamma^\opt(\rho, \gamma)}$ (as given by \eqref{Sec6:GGDef} with the $p=2$). In a similar way, the $H_1 f_1(\gamma)$, for 
$f_1 (\gamma):= -\frac{\alpha}{2} \sfd^2(\rho, \gamma)$, 
 involves a term $\inf_{\bpi \in \Gamma^\opt(\rho, \gamma)}$. If we were directly to prove comparison principle using the $H_0$ and $H_1$, these terms are the source of difficulties, unless the set $\Gamma^\opt(\rho, \gamma)$ consists of a single element only. 
We recall that, when $\rho:=\rho(dx)$ (or the $\gamma$) does not give mass to ``small" sets, then the  $\Gamma^\opt(\rho, \gamma)$ contains only a single element, which is given by the Brenier optimal transport map. See Theorem 2.12 in Villani~\cite{Villani03} for details and precise statements. However, when the $\rho$ concentrates positive mass on a small set, such $\rho$ becomes a ``singular" point in the $\sfX \in \CBB(0)$. Then multiple connecting geodesics may appear, no matter how small the distance $\sfd(\rho,\gamma)$ is.  That is, the set $\Gamma^\opt(\rho, \gamma)$ may contain more than one element.  If we want to infer defining inequalities of sub-solutions for ${\bf H}_0f_0$ from those for $H_0f_0$,  we need to improve some inequalities from 
\begin{align*}
\ldots \leq \sup_{\bpi \in \Gamma^\opt(\rho, \gamma)} \ldots
\end{align*}
 into something like
 \begin{align*}
\ldots \leq \inf_{\bpi \in \Gamma^\opt(\rho, \gamma)} \ldots.
\end{align*}
Through a regularization method, we will show that (Lemma~\ref{Sec8:bbH0bfH0}) the above is indeed possible in a perturbative sense, for those $\rho$s appearing as maximizers of certain functions in definition of viscosity sub-solution.  We will be using special properties of the Wasserstein space, as stated in Lemma~\ref{Sec2:fdtouch}, to prove this.
 
Before we begin, it is also useful to trace origin of the $\sup_{\bpi \in \Gamma^\opt(\rho, \gamma)}$ term in the $H_0 f_0(\rho)$, and explain why we couldn't derive the $\inf_{\bpi \in \Gamma^\opt(\rho, \gamma)}$ term directly through our works on the hydrodynamic limit equations. This is because that, during the hydrodynamic limit, we lost ability to be  precise at recording which geodesic direction is relevant giving the viscosity sub-solution property, when making sense of relevant cotangent elements corresponding to derivatives of test functions in the Hamiltonian operator. Note that the ${\it \sup}$ in $H_0f_0$  simply means that ``{\em there exists some} geodesic direction" such that the defining inequality for viscosity sub-solution holds, while the ${\it \inf }$ in ${\bf H}_0f_0$ means the inequality needs to hold for  ``{\em for every} geodesic direction". 

We mentioned the loss of ability to track relevant geodesic directions. This happened during each of the following two steps in earlier derivations: 
\begin{enumerate}
\item the submetry projection of Hamiltonian operator with general non-symmetric perturbative term (compare the $H_0$ in Lemma~\ref{sPrelax} with that in Lemma~\ref{PSubHJ2} where the Hamiltonian is given by \eqref{pH0Alt}).
\item the passage of limit in derivation of a limit Hamiltonian for the sub-solution case. 
\end{enumerate}

\subsection{Viscosity extensions -- the sub-solution case, extensions from ${\mathbb H}_0$ to ${\bf H}_0$} \label{Sec8:HImSub}
We revisit equation \eqref{Sec6:bbH0Sub}. We recall that the $U,V$ satisfying Conditions~\ref{U0CND} and \ref{VCND} are globally Lipschitz. Therefore, there exists finite constant $C_{U,V}>0$ such that
\begin{align}\label{Sec8:UVmod}
 & \int_{\R^{2d}} \Big(\big(U(x) - U(y)\big) + \big(V*\rho(x) - 
 V*\gamma(y)  \big)\Big) \bpi(dx,dy)  \\
  & \qquad  \qquad   \qquad
  \leq C_{U,V} \sfd(\rho, \gamma), \quad \forall \bpi \in \Gamma^{\opt}(\rho,\gamma).
  \nonumber
 \end{align}
Motivated by  the estimate in Lemma~\ref{parHJprop},
we also consider the following.
\begin{condition}\label{Sec8:fGroEst}
 \begin{align}
 \overline{f}(\rho) \geq - \beta\circ \sfd(\rho, \delta_0), \quad \forall \rho \in \sfX,
\end{align}
for some concave, increasing $\beta:=\beta(r): \R_+ \mapsto \R_+$ growing at sub-linear rate to infinity as $r \to +\infty$. 
\end{condition}
 
 \begin{lemma}\label{Sec8:bbH0bfH0}
Let $\overline{f} \in \USC(\sfX;\R)$ with  $\sup_\sfX \overline{f} <+\infty$ and Condition~\ref{Sec8:fGroEst} be satisfied.
Suppose that $\overline{f}$ is a sub-solution to \eqref{Sec6:bbH0Sub} in the point-wise viscosity solution sense, with the operator ${\mathbb H}_0$ defined in \eqref{Sec6:bbH0}. 
We define 
\begin{align*}
 \overline{f}_\epsilon(\rho):= 
 \sup_{\gamma \in \sfX} \big( \overline{f}(\gamma)- \frac{\sfd^2(\rho,\gamma)}{2 \epsilon} \big), \quad \forall \epsilon >0,  \rho \in \sfX.
\end{align*}

We assume that $h_0 \in C(\sfX)$.
Let $C:=C_{U,V}>0$ be the constant in \eqref{Sec8:UVmod}.
 We introduce another function
\begin{align*}
 h_{0;\epsilon} &:=h_{0;\epsilon}(\rho)  \\
 & :=  \frac{\alpha C_{U,V}^2}{2}\epsilon + 
   \sup \Big\{ h_0(\gamma) : \forall \gamma \in \sfX \text{ s.t. } 
    \sfd(\gamma,\rho) 
   \leq   \sqrt{2 \epsilon \big(\sup_\sfX \overline{f}  
   +\beta (\sfd(\rho,\delta_0))\big)} \Big\}.
 \end{align*}

Then    
\begin{enumerate}
\item $\overline{f}_\epsilon \in  \Lip_{\loc}(\sfX;\R)$, and $\overline{f}\leq  \overline{f}_{\epsilon^\prime} \leq \overline{f}_\epsilon$ for every $0< \epsilon^\prime <\epsilon$ with $\lim_{\epsilon \to 0^+} \overline{f}_\epsilon(\rho) =\overline{f}(\rho)$ for each $\rho \in \sfX$ fixed.
\item for each $\epsilon>0$ fixed, the $\overline{f}_\epsilon$ is a strong point wise viscosity sub-solution to 
\begin{align}\label{Sec8:bH0Sub}
  (I  - \alpha {\bf H}_0 ) \overline{f}_\epsilon \leq h_{0;\epsilon}.
\end{align} 
\end{enumerate}
\end{lemma}
\begin{proof}
The local Lipschitz and other regularity properties of the $\overline{f}_\epsilon$ follow from its definition. See for instance, Chapter 3 in \cite{AGS08}.  Next, we prove the strong point-wise viscosity sub-solution property. To simplify notations, we only verify this for $f_0 \in {\mathcal S^+}$. The general $f_0 \in \mathcal S^{+,\infty}$ case only requires notational modification by allowing $K=+\infty$.
To summarize, for each given $f_0:= f_{0; \gamma_1, \ldots, \gamma_K} \in {\mathcal S}^+$ and $\rho_\epsilon \in \sfX$ satisfying
\begin{align*}
 (\overline{f}_\epsilon - f_0) (\rho_\epsilon) 
    = \sup_{\sfX} (\overline{f}_\epsilon - f_0),
\end{align*}
we only need to prove that
\begin{align*}
 \big( \overline{f}_\epsilon - h_{0,\epsilon} \big) (\rho_\epsilon) 
  \leq \alpha {\bf H}_0 f_0(\rho_\epsilon).
\end{align*}
The proof is divided into several steps.

{\bf Step one:} We assumed that the $\overline{f}$ is a point-wise viscosity sub-solution to \eqref{Sec6:bbH0Sub}. Noting an equivalent expression of ${\mathbb H}_0$ in \eqref{Sec6:bbH0A}, for the above $\rho_\epsilon$,  there exists a maximizer $\gamma_\epsilon:=\gamma_{\epsilon}(dy)\in \sfX$ in the definition of $\overline{f}_\epsilon$ such that
\begin{align}\label{Sec8:Egamma}
\overline{f}_\epsilon(\rho_\epsilon) = \overline{f}(\gamma_\epsilon)
      - \frac{\sfd^2(\rho_\epsilon, \gamma_\epsilon)}{2 \epsilon},
\end{align}
and that
\footnote{Comparing here with the definition of ${\mathbb H}_0$ in \eqref{Sec6:bbH0}, we note that  the roles of $\rho:=\rho(dx), \gamma:=\gamma(dy)$ (hence the $x,y$) are reversed.}
\begin{align}\label{sec8:subH0}
 \alpha^{-1}  \big( \overline{f} - h_0 \big) (\gamma_\epsilon)  
 & \leq {\mathbb H}_0 \Big(  \frac{\sfd^2(\rho_\epsilon, \cdot)}{2\epsilon}\Big)(\gamma_\epsilon)   \\
 &= \sup_{\bpi \in \Gamma^{\opt}(\rho_\epsilon,\gamma_\epsilon)} \int_{\R^{2d}} 
   \bar{\sfH} \big(y, \frac{y-x}{\epsilon}; \gamma_\epsilon\big) \bpi(dx, dy),  \nonumber \\
    & =  \int_{\R^{2d}}   \bar{\sfH} \big(y,\frac{y-x}{\epsilon}; \gamma_\epsilon \big) \bpi_\epsilon(dx, dy), \nonumber
 \end{align}
for some probability measure $\bpi_\epsilon \in \Gamma^{\opt}(\rho_\epsilon, \gamma_\epsilon)$. 
 
Existence of the above $\rho_\epsilon, \gamma_\epsilon$ is equivalent to
\begin{align}\label{sec8:f0df}
\overline{f}(\gamma_\epsilon)  - \frac{\sfd^2(\rho_\epsilon, \gamma_\epsilon)}{2\epsilon}  -  f_0(\rho_\epsilon)   = \big(\overline{f}_\epsilon -  f_0\big)(\rho_\epsilon) 
 =\sup_{\rho,\gamma \in \sfX}
 \Big( \overline{f} (\gamma) - \frac{\sfd^2(\rho, \gamma)}{2\epsilon} -  f_0(\rho)  \Big).
\end{align}
We note, from definition of the $\overline{f}_\epsilon$, that
\begin{align}\label{sec8:fbare}
 \overline{f}(\rho_\epsilon) \leq \overline{f}_\epsilon(\rho_\epsilon) = \overline{f}(\gamma_\epsilon) - \frac{\sfd^2(\rho_\epsilon, \gamma_\epsilon)}{2 \epsilon} \leq \overline{f}(\gamma_\epsilon).
\end{align}

In addition, from \eqref{sec8:f0df}, we also get 
\begin{align*}
 f_0(\rho_\epsilon) \leq \overline{f}(\gamma_\epsilon) - \sup_\sfX(\overline{f} - f_0)
 \leq \sup_\sfX \overline{f} - \sup_\sfX(\overline{f} - f_0) < +\infty.
\end{align*}
The $f_0$ is defined in terms of $\psi \in \Psi$ in \eqref{Sec2:defPsi}.
Since  $\lim_{r_k \to \infty} \psi(r_1, \ldots, r_k,\ldots)=+\infty$ for at least one of the variables $r_k$, we conclude that  
\begin{align*}
 \limsup_{\epsilon \to 0^+} \sfd(\rho_\epsilon,\delta_0) <+\infty.
\end{align*}

{\bf Step two:} Taking an {\it arbitrary} $\bM \in \Gamma^{\opt}(\rho_\epsilon; \gamma_1, \ldots, \gamma_K) \subset {\mathcal P}_2(\R^{(1+K)d})$ (See Definition~\ref{MultOpt}), for the $f_0:=f_{0;\gamma_1, \ldots, \gamma_K} \in {\mathcal S}^+$, we define \begin{align*}
\bnu^{\bM}_{f_{0;\gamma_1,\ldots, \gamma_K}}
:=\bnu^{\bM}_{f_{0;\gamma_1,\ldots, \gamma_K}}(dx,dP)
:=\int_{(y_1,\ldots, y_K) \in \R^{Kd}} \delta_{\sum_k 2 \alpha_k(x-y_k)}(dP) \bM(dx; dy_1,\ldots, dy_K), 
\end{align*}
  as in \eqref{bnu0Def}, where the $\alpha_k:= \alpha_k(\rho;\gamma_1,\ldots, \gamma_K)$ are defined according to \eqref{alphaDef}. The $\bpi_\epsilon \in \Gamma^{\opt}(\rho_\epsilon, \gamma_\epsilon)$ in \eqref{sec8:subH0} admits a measurable slicing decomposition
\begin{align*}
 \bpi_\epsilon(dx, dy) = \bpi_\epsilon(dy|x) \rho_\epsilon(dx).
\end{align*}
This allows us to construct a lifting of the $\bM$ by introducing
\begin{align*}
 \bN(dx, dy; dy_1, \ldots, d y_K) := \bpi_\epsilon(dy | x) \bM(dx; dy_1, \ldots, d y_K) \in {\mathcal P}_2(\R^{d(2+K)}),
\end{align*}
and a further lifting of the $\bN$ by introducing
\begin{align*}
\widehat{\bN}(dx,dy;d y_1,\ldots, d y_K, dP):= \delta_{\sum_k 2 \alpha_k (x-y_k)}(dP) \bN(dx,dy;dy_1,\ldots, dy_K) \in {\mathcal P}_2(\R^{d(3+K)}).
\end{align*}
In particular, $\bN$ and $\bM$ can be obtained as projections from the $\widehat{\bN}$:
\begin{align*}
\pi^{1,\ldots, K+2}_\# \widehat{\bN} =\bN, \quad
\pi^{1,3,\ldots, K+2}_\# \widehat{\bN} = \pi^{1,3,\ldots, K+2}_\# \bN = \bM, \quad
 \pi^{1, K+3}_\# \widehat{\bN} = \bnu^{\bM}.
\end{align*} 
Moreover,
\begin{align*}
 \pi^{1,2}_\# \bN = \bpi_\epsilon \in \Gamma^{\opt}(\rho_\epsilon, \gamma_\epsilon), \quad \pi^{2, k+2}_\# \bN =\pi^{1,k+1}_\# \bM \in \Gamma^{\opt}(\rho_\epsilon, \gamma_k), k=1,\ldots, K.
\end{align*}
Using standard probability arguments, one can even construct random variables making the above probability measures as respective joint distributions. This is illustrated below using a graph: 

  {\scriptsize {\em An informal graphical representation of the marginal probability measures $\gamma,\rho_\epsilon, \gamma_1, \ldots, \gamma_K$ as submetry projections of random variables $Y,X, Y_1, \ldots, Y_K$ defined in one canonical probability space 
 $([0,1], {\mathcal B}_{[0,1]}, {\text Leb})$:}}
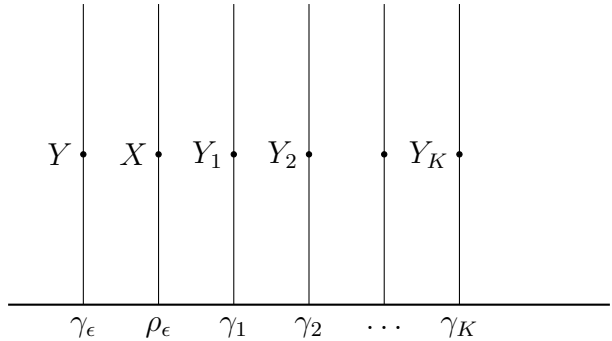
\begin{align*}
\begin{tikzpicture}
\draw[black, thick] (-4,0) --(4,0);
\draw (-3,4) -- (-3,0);
\filldraw (-3,2) circle (1pt) node [anchor= east] {$Y$};
\draw (-2,4) -- (-2,0);
\filldraw (-2,2) circle (1pt) node [anchor=east] {$X$};  
\draw (-1,4) -- (-1,0);
\filldraw (-1,2) circle (1pt) node [anchor=east] {$Y_1$};
\draw (0,4) -- (0,0);
\filldraw (0,2) circle (1pt) node [anchor=east] {$Y_2$};
\draw (1,4) -- (1,0);
\filldraw (1,2) circle (1pt) ;
\draw (2,4) -- (2,0);
\filldraw (2,2) circle (1pt) node [anchor=east] {$Y_K$};
\node (a)  at (-3,-0.3) {$\gamma_\epsilon$};
\node (b) at  (-2, -0.3) {$\rho_\epsilon$};
\node (c)  at (-1,-0.3) {$\gamma_1$};
\node (d)  at (0,-0.3) {$\gamma_2$};
\node (e)  at (1,-0.3) {$\ldots$};
\node (f)  at (2, -0.3) {$\gamma_K$};
\end{tikzpicture}
\end{align*}
 
{\bf Step three:} From \eqref{sec8:f0df}, we see that $\rho_\epsilon$ is a maximizer of function
\begin{align*}
 \rho \mapsto   -\frac{\sfd^2(\rho,\gamma_\epsilon)}{2\epsilon}-  f_0(\rho).
\end{align*}
In view of the results in Lemma~\ref{Sec2:fdtouch}, we have  
\begin{align*}
  \frac{x-y}{\epsilon} =   \sum_k 2\alpha_k(\rho_\epsilon; \gamma_1, \ldots, \gamma_K)(y_k-x), \quad  \bN \text{\rm - almost everywhere}.
\end{align*}
 Consequently,
\begin{align}\label{matchDer0}
 \big(y, \frac{y-x}{\epsilon} \big)= (y, P), \quad  \widehat{\bN} \text{\rm - almost everywhere}.
\end{align}
Therefore,  
 \begin{align*}
   &  \int_{\R^{2d}}  \bar{\sfH}\big( y, \frac{y-x}{\epsilon};\gamma_\epsilon\big) 
                            \bpi_\epsilon(dx,dy)  \\
   & = \int_{\R^{(3+K)d}}  \bar{\sfH} \big(y,  \frac{y-x}{\epsilon}; \gamma_\epsilon \big) 
   \widehat{\bN}(dx, dy; d y_1,\ldots, d y_K; dP) \\
    &  =  \int_{\R^{(3+K)d}}   \bar{\sfH}\big(y, P;\gamma_\epsilon \big) 
    d \widehat{\bN}   \\
  & =  \int_{\R^{2d}}  \bar{\sfH}\big(x,P;\rho_\epsilon \big) 
   \bnu^{\bM}_{f_{0;\gamma_1,\ldots, \gamma_K}}(dx,dP)
   + \int_{\R^{2d}} \big( U(x) -U(y)\big) \bpi_\epsilon(dx,dy) \\
& \qquad \qquad   + \int_{\R^{2d}}   \big(V*\rho_\epsilon(x) - 
 V*\gamma_\epsilon(y)  \big) \bpi_\epsilon(dx,dy).
 \end{align*}
 In the above, the first equality follows because 
 $\bpi_\epsilon= \pi^{1,2}_\# \bN = \pi^{1,2}_\# \widehat{\bN}$, the second equality follows from \eqref{matchDer0}. 
  
Combining the above with \eqref{sec8:subH0} and the equality part of \eqref{sec8:fbare}, and in view of estimate \eqref{Sec8:UVmod}, we have
  \begin{align*}
   \alpha^{-1} \Big( \overline{f}_\epsilon (\rho_\epsilon) 
       + \frac{\sfd^2(\rho_\epsilon, \gamma_\epsilon)}{2 \epsilon} 
           -   h_0 (\gamma_\epsilon)   \Big)  
   \leq  \int_{\R^{2d}}  \bar{\sfH}(x,P;\rho_\epsilon)  
      \bnu^{\bM}_{f_0; \gamma_1, \ldots, \gamma_K}(dx,dP) 
       + C_{U,V} \sfd(\rho_\epsilon,\gamma_\epsilon).
\end{align*}
 By arbitrariness of the $\bM$, and in view of \eqref{Sec7:H0},
we have
\begin{align*}
  \overline{f}_\epsilon  (\rho_\epsilon) 
       &   \leq \alpha {\bf H}_0 f_0(\rho_\epsilon) 
          + \Big(  C_{U,V} \alpha \sfd(\rho_\epsilon, \gamma_\epsilon)
           - \frac{\sfd^2(\rho_\epsilon, \gamma_\epsilon)}{2  \epsilon} \Big)
      +     h_0(\gamma_\epsilon)   \\
      &   \leq  \alpha {\bf H}_0 f_0(\rho_\epsilon) 
          +   \frac{\alpha C_{U,V}^2}{2}\epsilon + h_0(\gamma_\epsilon)  .
\end{align*}

From \eqref{sec8:f0df},
\begin{align*}
 \frac{\sfd^2(\rho_\epsilon, \gamma_\epsilon)}{2\epsilon} 
 \leq \overline{f}(\gamma_\epsilon)- \overline{f}(\rho_\epsilon)
 \leq \sup_\sfX \overline{f} + \beta\big(\sfd(\rho_\epsilon, \delta_0)\big).  \end{align*}
Consequently
\begin{align*}
 h_0(\gamma_\epsilon) & \leq \sup_{\gamma \in \sfX} 
  \Big\{ h_0(\gamma) :  \sfd^2(\gamma,\rho_\epsilon) 
   \leq   2 \epsilon \big(\sup_\sfX \overline{f}  
    +\beta(\sfd(\rho_\epsilon, \delta_0))\big) \Big\} \\
   & = h_{0;\epsilon}(\rho_\epsilon) - \frac{\alpha C_{U,V}^2}{2}\epsilon. 
\end{align*}
Therefore,  
\begin{align*}
  \overline{f}_\epsilon  (\rho_\epsilon) 
    \leq  \alpha {\bf H}_0 f_0(\rho_\epsilon) 
          +   h_{0;\epsilon}(\rho_\epsilon).
\end{align*}
We conclude.
\end{proof}

\begin{remark}
At beginning of the above lemma, we required that $\overline{f}$ is a sub-solution to \eqref{Sec6:bbH0Sub} in the point-wise viscosity sense. In particular, this implicitly means that maximum of $\overline{f} - f_0$ always exists for each $f_0 \in D({\mathbb H}_0)$. This guaranteed the existence of $\gamma_\epsilon$ in \eqref{Sec8:Egamma} in the above proof. We recall that the combined results of Lemmas~\ref{Sec6:seq2spw} and~\ref{Sec6:bbH0Eqn} ensured such assumption is not vacuous, and is useful in our context.  However, later application of a super-solution version of above result (proof of Lemma~\ref{Sec9:fReg}) will not have such property {\em a priori}. Consequently, we would like a version of the above lemma by not assuming existence of such extremal point. Indeed, because of the whole development in Section~\ref{Sec6} with results summarized in Lemma~\ref{Sec6:fsbbH0}, we only need to work with {\em strong} viscosity solutions.  
We have the following results.
\end{remark}

\begin{lemma}\label{Sec8:H0bfH0Alt}
In the context of Lemma~\ref{Sec8:bbH0bfH0}, if the $\overline{f}$ is a sub-solution to \eqref{Sec6:bbH0Sub} in the point-wise {\rm strong-} viscosity solution sense, then conclusions of the lemma still hold the same.
\end{lemma}
\begin{proof}
In the proof of comparison principle in Lemma~\ref{Sec7:CMP}, we used a perturbation method by invoking the Borwein-Preiss Lemma to produce maximum point. We use that argument here in similar ways to create maximizer satisfying a perturbed version of \eqref{Sec8:Egamma}. 

We note that operator ${\mathbb H}_0$ has the following property.
Let $\delta>0$ be a small parameter, and let a convergent sequence of 
$\{ \gamma_{\delta,k} \}_{k \in \N} \subset \sfX$ with limiting point 
$\gamma_\delta \in \sfX$ such that 
\begin{align*}
\sup_{\delta>0} \sfd(\gamma_\delta, \bar{\rho})<+\infty, \quad \exists \bar{\rho} \in \sfX.
\end{align*}
 We define
\begin{align*}
\Delta(\gamma) :=\Delta_\delta(\gamma)
:= \sum_{k=0}^\infty\frac{1}{2^{k+1}} \sfd^2(\gamma, \gamma_{k,\delta})
\end{align*}
Since 
\begin{align*}
\lim_{\delta \to 0^+} \sqrt{\delta} |D^+_{\gamma_\delta} \Delta| = 0,
\end{align*}
we have 
 \begin{align*}
\lim_{\delta \to 0^+} |{\mathbb H}_0 \big(f_0 + \sqrt{\delta} \Delta_\delta\big)(\gamma_\delta) 
 - {\mathbb H}_0 f_0(\gamma_\delta) | =0.
\end{align*}

The conclusion follows by adding such additional layer of approximation. 
\end{proof}

We see that the result in Lemma~\ref{Sec8:bbH0bfH0} is not perfect. There is a parameter $R$ in the definition of $h_{0;\epsilon, R}$. When such $R$ is fixed, and 
the $\sfd(\rho, \delta_0)$ becomes larger than the $R$, it is not apparent how to get useful information from the equation. 
Next, we introduce a technique to recover such information by exploring two features: one,  such $R$ can be chosen arbitrarily; two, sub-solution is stable with respect to another type of perturbation that reflects the growth estimates of the sub-solution.  

Let 
\begin{align*}
 \overline{f}_{\lambda,\theta}(\rho):= \lambda \overline{f}(\rho) 
 -  \theta  \sqrt{1+\sfd^2(\rho,\delta_0)},\quad \forall \theta>0, \lambda > 1. 
\end{align*}
Let $c, C \in \R_+$ be the constants in Lemma~\ref{PbarH}.  We define
\begin{align}\label{Sec8:epsrho}
{\rm Err}_{\lambda,\theta}(\rho):=(\lambda-1) c+ C\frac{\theta^2}{\lambda-1} 
 -(\lambda -1) (\langle U + V*\rho, \rho \rangle).
\end{align}

\begin{lemma}\label{Sec8:bfH0Pert}
Let $h_0 :\sfX \mapsto \R$ with $\sup_\sfX h_0 <+\infty$, and $\overline{f}$ be a strong point-wise viscosity sub-solution to 
\begin{align*}
 (I -\alpha {\bf H}_0) \overline{f} \leq h_0.
\end{align*}
We define
\begin{align*}
h_{0,\lambda, \theta}(\rho) := \lambda h_0(\rho) -\theta\sqrt{1+\sfd^2(\rho,\delta_0)} + \alpha  {\rm Err}_{\lambda,\theta}(\rho).
\end{align*}

Then the $\overline{f}_{\lambda,\theta}$ is a strong point-wise viscosity sub-solution to 
\begin{align*}
 (I - \alpha {\bf H}_0) \overline{f}_{\lambda,\theta} \leq h_{0,\lambda, \theta}.
\end{align*}
\end{lemma}
\begin{proof}
Let $f_0 :=f_{0;\gamma_1,\ldots, \gamma_K, \ldots} \in {\mathcal S}^{+,\infty}$ and 
$\rho_0 \in \sfX$ be such that 
\begin{align*}
 (\overline{f}_{\lambda,\theta} - f_0)(\rho_0) 
 = \sup_\sfX  (\overline{f}_{\lambda,\theta} - f_0).
\end{align*}
We write
\begin{align*}
f_{0;\lambda,\theta}:= \frac{1}{\lambda} \big( f_0 
 + \theta  \sqrt{1+\sfd^2(\cdot,\delta_0)}\big)
 \in {\mathcal S}^{+,\infty}.
\end{align*}
Then
\begin{align*}
 (\overline{f} - f_{0;\lambda,\theta})(\rho_0) 
 = \sup_\sfX (\overline{f} - f_{0;\lambda,\theta}).
 \end{align*}
By the strong point-wise sub-solution assumption,
\begin{align*}
 \overline{f}(\rho_0) \leq \alpha {\bf H}_0 f_{0;\lambda,\theta}(\rho_0) + h_0(\rho_0);
\end{align*}
or, equivalently,
\begin{align*}
\overline{f}_{\lambda,\theta}(\rho_0) 
 \leq \alpha \lambda {\bf H}_0 f_{0;\lambda,\theta}(\rho_0) 
  + \big(\lambda h_0(\rho_0) - \theta \sqrt{1+\sfd^2(\rho_0,\delta_0)}\big).
\end{align*}

Next, we have estimates
\begin{align*}
\lambda \big({\bf H}_0 f_{0;\lambda,\theta}\big) (\rho_0)
 &= \lambda {\bf H_0} \big( \frac{1}{\lambda}   f_0 
 + (1-  \frac{1}{\lambda}) \frac{\theta}{\lambda-1}
   \sqrt{1+\sfd^2(\cdot,\delta_0)}\big)(\rho_0) \\
 & \leq   {\bf H_0}   f_0  (\rho_0)
 + ( \lambda-1) {\bf H}_0\big( \frac{\theta}{\lambda-1}
   \sqrt{1+\sfd^2(\cdot,\delta_0)} \big)(\rho_0) \\
 & \leq {\bf H_0}   f_0(\rho_0) + (\lambda-1) \big( c+ 
 C\big|D_{\rho_0}^+\frac{\theta}{\lambda-1}
   \sqrt{1+\sfd^2(\cdot,\delta_0)} \big|^2 \\
   & \qquad \qquad \qquad -(\langle U  + V*\rho_0, \rho_0\rangle)\big) \\
 & \leq {\bf H_0}   f_0(\rho_0) + (\lambda-1) c+ 
 C\frac{\theta^2}{\lambda-1} 
 -(\lambda -1) (\langle U  + V* \rho_0, \rho_0 \rangle).
\end{align*}
In the second inequality above, we used an estimate on $\bar{\sfH}$ in Lemma~\ref{PbarH}. The constants $c,C \in \R_+$ are the ones there. 

Hence we conclude.
\end{proof}

\subsection{Viscosity extensions -- the super-solution case, extension from ${\mathbb H}_1$ to ${\bf H}_1$} \label{Sec8:HImSup}
Next, we revisit equation \eqref{Sec6:bbH1res}. 
Similar to arguments used in the proof of Lemma~\ref{Sec8:bbH0bfH0}, we  establish the following.

\begin{lemma}\label{Sec8:bbH1bfH1}
Let $\underline{f} \in \LSC(\sfX; \R)$ be such that $\sup_\sfX \underline{f}<+\infty$ and 
$\underline{f}(\gamma) \geq - \beta\circ \sfd(\gamma, \delta_0)$
for some concave, increasing $\beta:=\beta(r): \R_+ \mapsto \R$ growing at sub-linear rate to infinity as $r \to +\infty$.
Suppose that the $\underline{f}$ is a super-solution to \eqref{Sec6:bbH1res} in the point-wise viscosity solution sense. We define
\begin{align*}
 \underline{f}_\epsilon(\gamma) 
  := \inf_{\rho \in \sfX} \big( \underline{f}(\rho) + \frac{\sfd^2(\rho,\gamma)}{2 \epsilon} \big).
\end{align*}

We assume $h_1 \in C(\sfX)$. Let $C:=C_{U,V}>0$ be the constant in \eqref{Sec8:UVmod}. We introduce, for every $\epsilon>0$,
\begin{align}\label{Sec8:h1espD}
 h_{1;\epsilon}   :=h_{1;\epsilon}(\gamma) 
  :=  -\frac{\alpha C^2_{U,V}}{2}\epsilon 
   + \inf \Big\{ h_1(\rho) : \rho \in {\mathcal N}_\epsilon(\gamma) \Big\},
\end{align}
with  
\begin{align}\label{Sec8:NeDef}
{\mathcal N}_\epsilon(\gamma) := \Big\{ \rho \in \sfX :  \sfd^2(\rho,\gamma) 
   \leq 2 \epsilon\big(\sup_\sfX \underline{f} + \beta(\sfd(\rho,\delta_0)) \big) \Big\}
\end{align}
Note that, because of the sub-linear growth of $\lim_{r \to +\infty} \beta(r)=+\infty$, the set $N_\epsilon(\gamma)$ is a $\sfd$-bounded set in $\sfX$.

Then
\begin{enumerate}
\item $\underline{f}_\epsilon \in \Lip_{\loc}(\sfX; \R)$, and 
$\underline{f} \geq \underline{f}_{\epsilon^\prime} \geq \underline{f}_\epsilon$ 
for every $0<\epsilon^\prime < \epsilon$, with 
$\lim_{\epsilon \to 0^+} \underline{f}_\epsilon(\gamma) = \underline{f}(\gamma)$ for each $\gamma \in \sfX$ fixed. 
\item for each $\epsilon>0$, the $\underline{f}_\epsilon$ is a strong point-wise super-solution to 
\begin{align}\label{Sec8:bH1Sup}
 (I - \alpha {\bf H}_1) \underline{f}_\epsilon \geq h_{1; \epsilon}.
\end{align}
\end{enumerate}
\end{lemma}
\begin{proof}
We only highlight steps which are different than the sub-solution proof.
Let $\gamma_\epsilon \in \sfX$ be such that 
\begin{align*}
 (f_1-\underline{f}_\epsilon)(\gamma_\epsilon)
  = \sup_\sfX(f_1-\underline{f}_\epsilon).
\end{align*}
That is,
\begin{align*}
 \sup_{\rho \in \sfX} \big(f_1(\gamma_\epsilon) 
  -\frac{\sfd^2(\rho,\gamma_\epsilon)}{2 \epsilon} - \underline{f}(\rho)\big)
   \geq f_1(\gamma) 
  -\frac{\sfd^2(\rho,\gamma_\epsilon)}{2 \epsilon} - \underline{f}(\gamma), 
  \quad \forall \rho,\gamma \in \sfX.
\end{align*}
By point-wise viscosity solution property of the $\underline{f}$, there exists 
$\rho_\epsilon \in \sfX$ attaining the maximum on left of the above inequality. Moreover, 
\begin{align*}
 \underline{f}(\rho_\epsilon) \geq \alpha {\bf H}_1 \Big(-\frac{\sfd^2(\cdot, \gamma_\epsilon)}{2 \epsilon}\Big)(\rho_\epsilon) + h_1(\rho_\epsilon). 
\end{align*}

From the above, we obtain estimate
\begin{align*}
\sfd^2(\rho_\epsilon, \gamma_\epsilon) \leq 2 \epsilon \Big(\sup_\sfX \underline{f} + \beta \big(\sfd^2(\rho_\epsilon,\delta_0) \big) \Big).
\end{align*}
Hence
\begin{align*}
h_1(\rho_\epsilon) \geq \inf \{ h_1(\rho): \rho \in N_\epsilon(\gamma_\epsilon) \}
 = h_{1;\epsilon}(\gamma_\epsilon) +\frac{\alpha C_{U,V}^2}{2}\epsilon.
\end{align*}
\end{proof}

Similar to Lemma~\ref{Sec8:H0bfH0Alt}, by introducing an additional layer of approximation using the Borwein-Preiss perturbed optimization lemma, we have the following super-solution version. 
\begin{lemma}\label{Sec8:H1bfH1Alt}
In the context of Lemma~\ref{Sec8:bbH1bfH1}, if the $\underline{f}$ is a super-solution to \eqref{Sec6:bbH1res} in the point-wise {\rm strong-} viscosity sense (instead of point-wise viscosity sense), then the conclusion of that Lemma still holds the same.
\end{lemma}

Similar to Lemma~\ref{Sec8:bfH0Pert}, we have the following. 
Given $\underline{f} : \sfX \mapsto \R$, we introduce
\begin{align*}
\underline{f}_{\lambda, \theta}(\gamma) := \lambda^{-1} \underline{f} 
+ \theta \sqrt{1+ \sfd^2(\gamma, \delta_0)}, \quad \forall \theta >0, \lambda >1.
\end{align*}

\begin{lemma}\label{Sec8:bfH1Pert}
Let $h_1 :\sfX \mapsto \R$, and $\underline{f}$ be a strong point-wise viscosity sub-solution to 
\begin{align*}
 (I -\alpha {\bf H}_1) \underline{f} \geq h_1.
\end{align*}
We define
\begin{align*}
h_{1,\lambda, \theta}(\gamma) 
:= \lambda^{-1} h_1(\gamma) + \theta\sqrt{1+\sfd^2(\gamma,\delta_0)}
- \alpha {\rm Err}_{\lambda,\theta}(\gamma).
\end{align*}
where the ${\rm Err}_{\lambda,\theta}$ term is defined in \eqref{Sec8:epsrho}.

Then the $\underline{f}_{\lambda,\theta}$ is a strong point-wise viscosity super-solution to 
\begin{align*}
 (I - \alpha {\bf H}_1) \underline{f}_{\lambda,\theta} \geq h_{1,\lambda, \theta}.
\end{align*}
\end{lemma}
\begin{proof}
We only highlight changes in some of the key estimates.
First, by a convexity argument and in view of Lemma~\ref{PbarH}, 
for every $f_1 \in {\mathcal S}^{-,\infty}$, $\lambda >1$ and $\theta >0$,
\begin{align*}
& {\bf H}_1 \Big(\lambda \big(f_1 
          + \theta \sqrt{1+ \sfd^2(\cdot, \delta_0)}\big)\Big)(\gamma)\\
 & = {\bf H}_1 \Big( \lambda  f_1 - (\lambda -1) \frac{\theta}{\lambda-1}
 \big(- \sqrt{1+ \sfd^2(\cdot, \delta_0)}\big)\Big)(\gamma) \\
 & \geq \lambda ({\bf H}_1 f_1)(\gamma) - (\lambda -1) 
 {\bf H}_1 \big( - \frac{\theta}{\lambda-1} 
    \sqrt{1+ \sfd^2(\cdot, \delta_0)}\big)(\gamma) \\
 & \geq \lambda ({\bf H}_1 f_1)(\gamma) - (\lambda -1) \Big\{ c+ C
  \big| D_\gamma \big(- \frac{\theta}{\lambda-1} 
   \sqrt{1+ \sfd^2(\cdot, \delta_0)}\big)\big|^2
   - \langle U + V*\gamma, \gamma \rangle \Big\}. 
\end{align*}

Second, for $\gamma_0 \in \sfX$ such that 
\begin{align*}
(f_1 - \underline{f}_{\lambda, \theta})(\gamma_0) 
= \sup_\sfX (f_1 - \underline{f}_{\lambda, \theta}),
\end{align*}
we have
\begin{align*}
\underline{f}_\lambda(\gamma_0) & \geq 
 \alpha  \lambda^{-1} {\bf H}_1 \Big( \lambda \big(f_1 
  + \theta \sqrt{1+\sfd^2(\cdot,\delta_0)}\big)\Big) (\gamma_0) 
   + \lambda^{-1} h_1(\gamma_0) 
     +\theta \sqrt{1+ \sfd^2(\gamma,\delta_0)} \\
     & \geq \alpha {\bf H}_1 f_1 (\gamma_0) + h_{1,\lambda,\theta}(\gamma_0).
\end{align*}

The conclusion follows.
\end{proof}

\subsection{Convergence of viscosity solutions, from particle to continuum}
Let an appropriate sequence of ${\mathfrak h}_N \in C\big((\R^d)^N\big)$ be given.
We define ${\mathfrak f}_N$ through \eqref{PartfDef}. By Lemma~\ref{parHJprop},  such ${\mathfrak f}_N \in C\big((\R^d)^N\big)$ is the unique viscosity solution to Hamilton-Jacobi equation \eqref{HNHJB}.  Next, we study convergence of the ${\mathfrak f}_N$s and characterize the limit $f$ as viscosity solution of Hamilton-Jacobi equation in space of probability measures in proper senses -- See Theorem~\ref{Sec8:MainThm1}. Later, in Theorem~\ref{Sec9:MainThm2}, we will further improve the characterization of limit solution $f$.

\begin{definition} [Class $\mathcal C$] \label{Sec8:ClassC}
For sequence of functions with ${\mathfrak h}_N \in C ((\R^d)^N )$ and $h \in C(\sfX)$,
we define a special collection:
\begin{align*}
{\mathcal C} := \big\{ (\{{\mathfrak h}_N \}_{N \in \N}, h) \text{ satisfy the following properties } \big\} \subset C\big((\R^d)^N\big) \times \ldots \times C(\sfX),
\end{align*}
with 
\begin{enumerate}
\item \label{Sec8:CC1} 
${\mathfrak h}_N(\tau {\bf x}) = {\mathfrak h}_N({\bf x})$ for every $\tau \in \sfG_N$;
\item \label{Sec8:CC2}
for every $\rho_N:= \frac1N \sum_{i=1}^N \delta_{x_i}$ and every $\rho_0 \in \sfX$ such that $\sfd(\rho_N ,\rho_0) \to 0$, we have
\begin{align*}
\lim_{N \to \infty} {\mathfrak h}_N(x_1, \ldots, x_N) = h(\rho_0);
\end{align*}
 \item \label{Sec8:CC3}
 uniform growth estimates for ${\mathfrak h}_N$ and $h$:  
 \begin{align*}
\sup_N \sup_{(\R^d)^N} {\mathfrak h}_N + \sup_\sfX h <+\infty, 
\end{align*}
and there exists a concave, increasing and sub-linear function $\beta: \R_+\mapsto \R$ such that
\begin{align*}
 {\mathfrak h}_N(x_1,\ldots, x_N) & \geq -\beta \big( \sfd(\rho_N,\delta_0) \big), 
 \quad \forall \rho_N := \frac1N \sum_{i=1}^N \delta_{x_i};
\end{align*}
\item \label{Sec8:CC4}
for every $\rho_N:= \frac1N \sum_{i=1}^N \delta_{x_i}$ and every $\rho_0 \in \sfX$ such that $\sfd_{r=1}(\rho_N ,\rho_0) \to 0$ and $\sup_N \int_{\R^d} |x|^2 \rho_N(dx) <\infty$, we have
\begin{align*}
\lim_{N \to \infty} {\mathfrak h}_N(x_1, \ldots, x_N)  \leq h(\rho_0);
\end{align*}
(i.e. Property ${\mathscr P_N}$ as given by Definition~\ref{Sec4:ProPN} is satisfied);
\item \label{Sec8:CC5}
the $h$ is $\sfd_{p=1}$-upper semicontinuous in $\sfX$ (see Definition~\ref{Sec6:weakUSC});
\item \label{Sec8:CC6}
the $h$ has modulus of continuity with respect to $\sfd:=\sfd_{p=2}$-metric, on every $\sfd$-balls with finite radius.
\end{enumerate}
\end{definition}

\begin{example}
Let $h \in {\mathcal S}^{-}$, and ${\mathfrak h}_N$ be the empirical measure versions of the $h$. Then the $(\{ {\mathfrak h}_N \}_{N \in \N}, h ) \in {\mathcal C}$. Such class of functions can be used to identify closed sets $A \subset \sfX$ by approximating the function:
\begin{align*}
\chi_A (x) := 
 \begin{cases}
  0, & \text{ when } x \in A; \\
  -\infty, & \text{ when } x \not\in A.
 \end{cases}
\end{align*}
\end{example}

\begin{theorem}\label{Sec8:MainThm1}
Suppose that $(\{ {\mathfrak h}_N \}_{N \in \N}, h ) \in {\mathcal C}$.
We define ${\mathfrak f}_N$ according to \eqref{PartfDef}. 
By Lemma~\ref{parHJprop}, such ${\mathfrak f}_N \in C\big((\R^d)^N\big)$ is the unique viscosity solution to Hamilton-Jacobi equation \eqref{HNHJB} with at most linear growth. We define
\begin{align*}
\overline{f}_N(\rho) :=  \sup_{\substack{ (x_1,\ldots, x_N) \in (\R^d)^N \\
 \text{such that }  \rho= \frac1N \sum_{i=1}^N \delta_{x_i}  } } {\mathfrak f}_N(x_1, \ldots, x_N) 
 = \inf_{\substack{ (x_1,\ldots, x_N) \in (\R^d)^N \\
 \text{such that }  \rho= \frac1N \sum_{i=1}^N \delta_{x_i}  } }
 {\mathfrak f}_N(x_1, \ldots, x_N).
\end{align*}
See Lemma~\ref{SubQuotHJ} for validity of the above definition.
As in \eqref{Sec5:barfNthe}, we introduce small perturbation
\begin{align*}
\overline{f}_{N,\theta}(\rho):=\overline{f}_N(\rho) - \theta \int_{\R^d} |x|^2 \rho(dx), 
\quad \theta >0,
\end{align*}
and define \ $f_\theta^*$ according to \eqref{Sec6:fstarThe} and $f^*$ as in \eqref{Sec6:fstarDef}.
We also define $\underline{f}$ as in the context of Lemma~\ref{Sec6:SupLem}.

Then
\begin{enumerate}
\item we have relation
\begin{align*}
f:= f^* = \underline{f} \in C(\sfX).
\end{align*}	
\item the $f$ is both a viscosity sub-solution in the point-wise strong sense to 
\eqref{Sec6:bbH0Sub}, as well as a super-solution in the point-wise strong sense to \eqref{Sec6:bbH1res}, with the $h_0=h_1=h$. 
\item The sequence $(\{ {\mathfrak f}_N \}_{N \in \N}, f)$ satisfies all the properties of being in class ${\mathcal C}$, except \eqref{Sec8:CC5} and \eqref{Sec8:CC6} 
regarding $f$ being $\sfd_{p=1}$-upper semicontinuous in $\sfX$, and with 
$\sfd_{p=2}$-modulus of continuities in $\sfd_{p=2}$-balls of finite radius. In particular, we have that 
\begin{align}\label{Sec6:fNCvgf}
\lim_{N \to \infty} {\mathfrak f}_N(x_1, \ldots, x_N) = f(\rho_0), 
\quad \forall \rho_N:= \frac1N \sum_{i=1}^N \delta_{x_i}, \rho_0 \in \sfX \text{ with } \sfd(\rho_N ,\rho_0) \to 0.
\end{align}
\end{enumerate}
\end{theorem}
\begin{remark}
Indeed, properties \eqref{Sec8:CC5} and \eqref{Sec8:CC6} in Definition~\ref{Sec8:ClassC} also hold. That is, we have 
 $(\{ {\mathfrak f}_N \}_{N \in \N}, f) \in {\mathcal C}$. However, we won't prove this claim until Theorem~\ref{Sec9:MainThm2}, after introducing additional variational characterization for the $f$.
\end{remark}
\begin{proof}
By Lemma~\ref{Sec6:seq2spw}, the $f^* \in \USC(\sfX)$ is a viscosity sub-solution in point-wise strong sense to \eqref{Sec6:H0sub} given by Hamiltonian operator $H_0$. By Lemma~\ref{Sec6:bbH0Eqn}, it is also a point-wise strong sub-solution to \eqref{Sec6:bbH0Sub} with Hamiltonian operator ${\mathbb H}_0$.
By Lemma~\ref{Sec8:H0bfH0Alt}, its Yosida regularization
 \begin{align*}
 \big(f^*\big)_\epsilon(\rho) := \sup_{\gamma \in \sfX} \Big( f^*(\gamma) 
   -\frac{\sfd^2(\rho,\gamma)}{2\epsilon}\Big),
\end{align*}
is a strong point-wise sub-solution to \eqref{Sec8:bH0Sub} given by Hamiltonian ${\bf H}_0$. By Lemma~\ref{Sec8:bfH0Pert}, 
\begin{align*}
\big(f^*\big)_{\epsilon; \lambda, \theta}
:= \big(f^*\big)_{\epsilon; \lambda, \theta}(\rho) := \lambda \big(f^*\big)_\epsilon(\rho) 
 - \theta \sqrt{1+ \sfd^2(\rho,\delta_0)} \in C(\sfX), 
  \quad \lambda >1, \theta >0,
\end{align*}
is a strong point-wise sub-solution to 
\begin{align*}
 (I -\alpha {\bf H}_0 ) \big(f^*\big)_{\epsilon;\lambda, \theta} 
 \leq h_{0;\epsilon, \lambda, \theta}
\end{align*}
with 
\begin{align*}
 h_{0;\epsilon,\lambda, \theta}(\rho)& := \sup\Big\{ \lambda h(\gamma) : \gamma \in \sfX\ \text{ s.t. } \sfd^2(\gamma,\rho) \leq 2 \epsilon \big(\sup_\sfX \overline{f} + 
 \beta (\sfd(\rho,\delta_0))\big) \Big\} \\
 & \qquad \qquad + \frac{\lambda \alpha}{2} C_{U,V}^2\epsilon
  - \theta \sqrt{1+\sfd^2(\rho,\delta_0)} + \alpha {\rm Err}_{\lambda, \theta}(\rho),
\end{align*}
where the function ${\rm Err}_{\lambda,\theta}$ is defined in \eqref{Sec8:epsrho}.

In a similar way, we define Yosida approximation of the $\underline{f} \in \LSC(\sfX)$ by
\begin{align*}
\big( \underline{f}\big)_\epsilon(\gamma):= \inf_{\rho \in \sfX} \Big( \underline{f}(\rho) +\frac{\sfd^2(\rho,\gamma)}{2\epsilon}\Big).
\end{align*}
Then, by Lemma~\ref{Sec6:H1SupUlt}, $\underline{f}$ is a super-solution in the point-wise strong sense to \eqref{Sec6:H1sup} given by Hamiltonian operator $H_1$. By Lemma~\ref{Sec6:bbH1Eqn}, it is a strong point-wise super-solution to \eqref{Sec6:bbH1res} given by Hamiltonian ${\mathbb H}_1$. By Lemma~\ref{Sec8:bbH1bfH1}, the $(\underline{f})_\epsilon$ is a super-solution in the strong point-wise sense to \eqref{Sec8:bH1Sup} with operator ${\bf H}_1$.
By Lemma~\ref{Sec8:bfH1Pert},
\begin{align*}
\big(\underline{f}\big)_{\epsilon,\lambda,\theta}
:= \big(\underline{f}\big)_{\epsilon,\lambda,\theta}(\gamma) :=
\lambda^{-1} \big(\underline{f}\big)_\epsilon(\gamma) 
  + \theta \sqrt{1+\sfd^2(\gamma,\delta_0)} \in C(\sfX), 
   \quad \lambda>1, \theta>0,
\end{align*}
is a strong point-wise viscosity solution to 
\begin{align*}
(I-\alpha {\bf H}_1) \big(\underline{f}\big)_{\epsilon,\lambda,\theta} 
 \geq h_{1; \epsilon,\lambda,\theta},
\end{align*}
with 
\begin{align*}
 h_{1; \epsilon,\lambda,\theta}(\gamma) := \inf \Big\{ \lambda^{-1} h(\rho) :
 \rho \in {\mathcal N}_\epsilon(\gamma) \Big\} - \frac{\alpha C_{U,V}^2}{2 \lambda}\epsilon
  +\theta \sqrt{1+\sfd^2(\gamma, \delta_0)} 
   -\alpha {\rm Err}_{\lambda,\theta}(\gamma),
\end{align*}
where the $\sfd$-bounded neighborhood ${\mathcal N}_\epsilon(\gamma)$ is 
defined in \eqref{Sec8:NeDef}.

Next, we apply the comparison principle established in Theorem~\ref{Sec7:CMP}  to arrive at
\begin{align*}
 \text{LHS} & := \limsup_{\lambda \to 1^+}  \limsup_{\theta \to 0^+} \limsup_{\epsilon \to 0^+}
 \sup_\sfX  \Big( \big(f^*\big)_{\epsilon; \lambda, \theta}
  - \big(\underline{f}\big)_{\epsilon,\lambda,\theta}\Big) \\
  & \qquad \leq 
 \limsup_{\lambda \to 1^+}   \limsup_{\theta \to 0^+} \limsup_{\epsilon \to 0^+}
   \sup_\sfX \Big( h_{0; \epsilon,\lambda,\theta} 
  -  h_{1; \epsilon,\lambda,\theta}\Big)=: \text{RHS}.
\end{align*}
For every $\rho \in \sfX$ fixed, by Lemmas~\ref{Sec8:bbH0bfH0} and \ref{Sec8:bbH1bfH1}, we have
\begin{align*}
 (f^*- \underline{f})(\rho) \leq \limsup_{\epsilon \to 0^+}
 \Big( \big(f^*\big)_\epsilon - \big(\underline{f}\big)_\epsilon\Big)(\rho) 
  \leq \text{LHS}.
\end{align*}
To evaluate the right hand side, first, by Conditions~\ref{U0CND},~\ref{VCND}, we have  
\begin{align*}
 -\langle U+V*\rho, \rho \rangle \leq -\inf_\sfX U - \inf_\sfX V<+\infty.
\end{align*}
Consequently, 
\begin{align*}
  \limsup_{\lambda \to 1^+} \limsup_{\theta \to 0^+} \limsup_{\epsilon \to 0^+}
  \sup_{\rho \in \sfX} {\rm Err}_{\lambda, \theta} (\rho) \leq 0.
\end{align*}
Therefore, recall definition of neighborhood ${\mathcal N}_\epsilon(\sigma)$ in \eqref{Sec8:NeDef}, 
we also introduce another neighborhood 
\begin{align*}
 \hat{\mathcal N}_\epsilon(\sigma):=\big\{ \gamma \in \sfX:   \sfd^2(\gamma, \sigma) 
  \leq 2 \epsilon \big(\sup_\sfX \overline{f} + 
 \beta\circ\sfd( \sigma,\delta_0)\big) \big\}.
\end{align*}
Writing
\begin{align*}
 F_{\theta,\epsilon}(\rho, \gamma, \sigma)= 
 \big( h(\gamma) - h(\rho) \big)\vee 0 - 2 \theta \sqrt{1+\sfd^2( \sigma,\delta_0)},
\end{align*}
then,
\begin{align*}
{\rm RHS} & = 
 \limsup_{\lambda \to 1^+} \limsup_{\theta \to 0^+} \limsup_{\epsilon \to 0^+}
    \sup_\sfX \big( h_{0; \epsilon,\lambda,\theta} 
  -  h_{1; \epsilon,\lambda,\theta}\big)  \\
& =  \limsup_{\lambda \to 1^+} \limsup_{\theta \to 0^+} \limsup_{\epsilon \to 0^+}
 \sup_{\sigma \in \sfX} \Big( \sup_{\gamma \in \sfX} \big\{ \lambda h(\gamma) :  
     \gamma \in\hat{\mathcal N}_\epsilon(\sigma) \big\} \\
 & \qquad \qquad \qquad  - \inf_{\rho \in \sfX} \big\{ \lambda^{-1} h(\rho) :
 \rho \in {\mathcal N}_\epsilon(\sigma) \big\} 
 - 2 \theta \sqrt{1+\sfd^2(\sigma,\delta_0)} \Big) \\
 & \leq \limsup_{\theta \to 0^+} \limsup_{\epsilon \to 0^+}
   \Big\{ \sup_{ (\rho,\gamma,\sigma) \in \sfX \times \sfX \times \sfX }
 \big\{   F_{\theta,\epsilon}(\rho, \gamma, \sigma) :  \forall
    \gamma \in\hat{\mathcal N}_\epsilon(\sigma), 
   \forall \rho \in {\mathcal N}_\epsilon(\sigma) \big\}  \Big\}.
\end{align*}

Next, we claim that for each $\theta>0$ fixed, there exists a finite $M_\theta>0$ which is independent of the $\epsilon>0$, such that
\begin{align}\label{Sec8:Fineq}
  \sup_{ (\rho,\gamma,\sigma) \in \sfX \times \sfX \times \sfX }
 \big\{   F_{\theta,\epsilon}(\rho, \gamma, \sigma) & :
   \gamma \in\hat{\mathcal N}_\epsilon(\sigma), 
    \rho \in {\mathcal N}_\epsilon(\sigma),
     \sfd(\rho,\delta_0) >M_\theta \big\}   \leq 0, \quad \forall \epsilon \in (0,1].  
\end{align}
Then the above implies 
\begin{align*}
& \sup_{ (\rho,\gamma,\sigma) \in \sfX \times \sfX \times \sfX }
 \big\{   F_{\theta,\epsilon}(\rho, \gamma, \sigma) :  
    \gamma \in\hat{\mathcal N}_\epsilon(\sigma), 
  \rho \in N_\epsilon(\sigma) \big\} \\
   &\leq  \sup_{ \sigma \in \sfX } 
 \big\{  \big( h(\gamma) - h(\rho) \big) \vee 0   :   
    \gamma \in\hat{\mathcal N}_\epsilon(\sigma), 
    \rho \in {\mathcal N}_\epsilon(\sigma), 
   \sfd(\rho,\delta_0) \leq M_\theta  \big\} \vee 0 \\
   & \leq \sup_{ \sigma \in \sfX } 
 \big\{  \big( h(\gamma) - h(\rho) \big) \vee 0  :   
    \sfd^2(\rho,\gamma) \leq 2\epsilon [\sup_\sfX \underline{f}
     + \beta(M_\theta)],  \sfd(\rho,\delta_0) \leq M_\theta  \big\} \vee 0 \\
   & \leq \omega_{h;M_\theta}\big(\sqrt{2\epsilon [\sup_\sfX \underline{f} + \beta(M_\theta)]}\big);
\end{align*}
where the last step above follows from assumption of $h$ having a modulus of continuity $\omega_{h;M}$ in bounded $\sfd$-balls with finite radius $M>0$ (See Definition~\ref{Sec8:ClassC}.\ref{Sec8:CC6} about class $\mathcal C$). Consequently
\begin{align*}
 {\rm RHS} \leq 0.
\end{align*}

We prove \eqref{Sec8:Fineq} next.  
First, by sub-linear growth at $+\infty$ assumption on the $\beta$, there exists 
finite $C_\beta >0$ such that $\beta(r) \leq C_\beta + r$. Therefore, 
from $\rho \in {\mathcal N}_\epsilon(\sigma)$, we can find finite $C_{\underline{f},\beta}>0$ such that $\sfd^2(\rho,\sigma) \leq 2 \epsilon 
[C_{\underline{f},\beta} + \sfd(\rho,\sigma) 
 + \sfd(\sigma,\delta_0)]$, implying 
\begin{align}\label{Sec8:drsEst}
 \sfd(\rho,\sigma) \leq \epsilon + \sqrt{\epsilon^2 
  + 2\epsilon\big( C_{\underline{f},\beta} + \sfd(\sigma,\delta_0) \big)}.
\end{align}
Second, it follows then the following holds for every 
$\rho \in {\mathcal N}_\epsilon(\sigma)$, 
$\gamma \in \hat{\mathcal N}_\epsilon(\sigma)$, and $\epsilon \in (0,1]$:
\begin{align*}
F_{\theta,\epsilon}(\rho,\gamma,\sigma)
&\leq \sup_\sfX \overline{f}  + \beta\circ \sfd(\rho,\delta_0) 
-2 \theta \sqrt{1+\sfd^2(\sigma,\delta_0)} \\
& \leq \sup_\sfX \overline{f}  + \beta\big( \sfd(\rho,\sigma)+\sfd(\sigma,\delta_0) \big) -2 \theta \sqrt{1+\sfd^2(\sigma,\delta_0)} \\
& \leq \sup_\sfX \overline{f}  + \beta\big(1 + \sqrt{ 1
  + 2[ C_{\underline{f},\beta} + \sfd(\sigma,\delta_0) ]}
   +\sfd(\sigma,\delta_0) \big) -2 \theta \sqrt{1+\sfd^2(\sigma,\delta_0)}.
\end{align*}
By the sub-linear growth assumption of $s \mapsto \beta(s)$ as $s  \to +\infty$,
and the linear growth of $r \mapsto \sqrt{1+r^2}$, there exists a finite $N_\theta>0$ which is independent of the $\epsilon>0$, such that
right hand side of the above becomes negative when 
$\sfd(\sigma,\delta_0)>N_\theta$.
Third, from \eqref{Sec8:drsEst} and $\sfd(\rho,\delta_0) \leq \sfd(\rho,\sigma) + \sfd(\sigma,\delta_0)$, we have existence of $M_\theta$ such that $\sfd(\rho,\delta_0)>M_\theta$ implies $\sfd(\sigma,\delta_0)>N_\theta$. Combine the above three steps, we verified \eqref{Sec8:Fineq}.
 
In summary, we have 
\begin{align*}
 f^* - \underline{f} \leq {\rm LHS} \leq {\rm RHS} \leq 0.
\end{align*}
But by construction, and noting Lemma~\ref{Sec6:Lfstar}, $\underline{f} \leq f^*$. Consequently $f^* =\underline{f} \in C(\sfX)$, and \eqref{Sec6:fNCvgf} follows. 

The sequence $(\{{\mathfrak f}_N\}_{N\in \N},f)$ satisfies various properties of being in 
class ${\mathcal C}$: property~\ref{Sec8:CC1} in Definition~\ref{Sec8:ClassC} follows from Lemma~\ref{SubQuotHJ}; property~\ref{Sec8:CC2} from Lemma~\ref{Sec6:Lfstar}, the definition of $\underline{f}$ and the fact that $f^*=f$ as proved above; property~\ref{Sec8:CC3} from Lemma~\ref{parHJprop} and estimate \eqref{Sec5:frafUbd}; property~\ref{Sec8:CC4} from Lemma~\ref{Sec6:Lfstar} (the part regarding ${\mathscr P_N}$ property).
\end{proof}

In the next section (Theorem~\ref{Sec9:MainThm2}), we construct an explicit variational representation for the limiting $f$ in above theorem.

\newpage

\section{Lagrangian dynamics in space of probability measures}
\label{Sec9}
Ambrosio-Gigli-Savar\'e~\cite{AGS08} discussed concept and properties of absolute continuous curves in metric spaces. The space $\sfX:={\mathcal P}_2(\R^d)$, with Wasserstein order-2 metric $\sfd$, is an Alexandrov metric space. Following Definition~\ref{veloCurv}, we we introduce velocity of a curve. We also recall identification of tangent cones in  Lemma~\ref{Sec2:TanIden}. 
\begin{lemma}\label{Sec9:IntP0}
Let $\sigma(\cdot) \in AC([0,\infty); \sfX)$ be defined as in Chapter 1 of \cite{AGS08}.
 Then, for Lebesgue a.e. $t \geq 0$ the following holds:
\begin{enumerate}
\item there exists $v(t):=\Pi_{\sigma(t)}\big(v(t)\big) \in L^2_{\nabla, \sigma(t)} \subset \Tan_{\sigma(t)}$ such that in the sense of distribution
\begin{align}\label{Sec9:ConEq}
 \partial_t \sigma + \div (v \sigma) =0  \text{ in }  \R^d \times [0,\infty);
\end{align}
\item the following derivative exists in the sense of Definition~\ref{veloCurv} 
\begin{align*}
\dot{\sigma}(t):=\frac{d}{dt_+} \sigma(t) =  \bnu(t) = \bnu(t;dx,d\xi)= \delta_{v(t,x)}(d\xi)\sigma(t,dx) \in \Tan_{\sigma(t)}. 
\end{align*}
\end{enumerate}
 \end{lemma}
\begin{proof}
See Theorem 8.3.1, Propositions 8.4.5 and 8.4.6 of Ambrosio, Gigli and Savar\'e~\cite{AGS08}. 
\end{proof}

Let $L: {\mathcal P}_2(\R^d \times \R^d) \mapsto \R \cup \{+\infty\}$ be defined as in \eqref{Sec1:Lbnu} with the $\bar{\sfL}_{U,V}$ given by \eqref{Sec1:LUV} and $\bar{\sfL}$ by \eqref{Sec1:barsfL}. We also introduce action functional $A$ for continuous curve $\sigma(\cdot) \in C([0,\infty); \sfX)$ as in \eqref{Sec1:ActMV}:
\begin{align*}
A:=A_T[\sigma(\cdot)]:= 
\begin{cases}
 \int_0^T L(\bnu(t)) dt, & \text{ when } \sigma\in AC([0,T];\sfX), \\
 +\infty, & \text{ otherwise.}
\end{cases}
\end{align*} 
Let $\alpha>0$ and $h \in C(\sfX)$ with $\sup_\sfX h <+\infty$,
we define value function $f:= {\bf R}_\alpha h: \sfX \mapsto \R$ by
\begin{align}\label{Sec9:R}
{\bf R}_\alpha h (\rho)& := \sup \Big\{ \int_0^\infty e^{-\alpha^{-1}s} \Big( \frac{h\big(\sigma(s)\big)}{\alpha}  - L\big(\bnu(s)\big) \Big) ds : \sigma \in AC, \sigma(0) = \rho\Big\}\\
 & = \sup \Big\{ \int_0^\infty \frac{e^{-\alpha^{-1}s}}{\alpha} \Big( h\big(\sigma(s)\big)   - \int_0^s L\big(\bnu(r)\big)dr \Big) ds : \sigma \in AC, \sigma(0) = \rho\Big\}, \nonumber
\end{align} 
where the last equality follows by Fubini theorem.

\subsection{Convexity}
Following Definition~\ref{BarPDef} and Lemma~\ref{Bproj}, each $\bnu \in \Tan_\sigma$ for some $\sigma \in \sfX$ induces a Barrycentric projected curve onto special subsets $L^2_{\nabla, \sigma}$ (see \eqref{defLnab}) of the tangent cones:
\begin{align*}
 u_\bnu:=u_\bnu(x):= \int_{\R^d} v \bnu(dv|x) \in L^2_{\nabla, \sigma} \subset \Tan_{\sigma}.
\end{align*}
The space $L^2_{\nabla,\sigma}$ has a linear structure. Since the following is convex
\begin{align*}
 v \mapsto \bar{\sfL}_{U,V}(x,v;\rho),
\end{align*}
 by Jensen's inequality and Lemma~\ref{JensenLike},
\begin{align*}
& {\bf R}_\alpha h (\rho)\\
& = \sup_{\sigma \in AC, \sigma(0) = \rho}
 \Big\{ \int_{s=0}^\infty e^{-\alpha^{-1}s} \Big( \frac{h\big(\sigma(s)\big)}{\alpha} 
  - \int_{x\in \R^d}\bar{\sfL}_{U,V}\big(x, u_\bnu(s,x); \sigma(s)\big) \sigma(s, dx)\Big) ds   \Big\}\\
 & = \sup_{\sigma \in AC, \sigma(0) = \rho} \Big\{ \int_{s=0}^\infty \frac{e^{-\alpha^{-1}s}}{\alpha} \Big( h\big(\sigma(s)\big)   - \int_{(r,x) \in[0,s] \times \R^d} \bar{\sfL}_{U,V}\big(x, u_\bnu(r,x); \sigma(r)\big) \sigma(r, dx) dr \Big) ds\Big\}.
\end{align*} 

\subsection{Some properties of the value function and viscosity solutions}

\subsubsection{Some useful estimates}
We begin this subsection by recalling those mass transport theory notations in Section~\ref{appmass}.
Our main goal is to prove Lemma~\ref{Sec9:Pmax}. However, for such purpose, we need some preparatory results first.

 \begin{lemma}\label{Sec9:Hf1LSC} 
 For each $f_1 \in {\mathcal S}^{-,\infty}$,
${\mathbb H}_1 f_1 \in \LSC(\sfX; \R)$. In fact, assume that $\gamma_n, \gamma_0, \rho_0 \in \sfX$ satisfy $\gamma_n \Rightarrow \gamma_0$ in the narrow (i.e. weak) convergence of probability sense; also assume that   
\begin{align*}
 \sup_n \sfd(\gamma_n, \rho_0) <+\infty
\end{align*}
(recall that $\sfd$ is the 2-Wasserstein metric); then we have
\begin{align*}
 \liminf_{n \to \infty} {\mathbb H}_1 f_1(\gamma_n) \geq {\mathbb H}_1 f_1(\gamma_0).
\end{align*}

 \end{lemma}
 \begin{proof}
 The ${\mathbb H}_1 f_1$ is defined in \eqref{Sec6:bbH1} with an equivalent expression in \eqref{Sec6:bbH1A}. 
 
 Let $\rho_1, \ldots, \rho_K$ be those in the expression of $f_1$ (see \eqref{mtrang}). We note the following property:   
$\bmu_{k,n} \in \exp^{-1}_{\gamma_n}(\rho_k)$ implies that $\{ \bmu_{k,n} : n \in \N\}$ 
is relatively compact in $\mathcal P(\R^d \times \R^d)$ in the narrow topology. Moreover, at least along subsequences, 
$\bmu_{k,n} \Rightarrow \bmu_k \in \exp^{-1}_{\gamma_0}(\rho_k)$ as $n \to \infty$.
The above observation, together with Fatou's lemma, imply conclusion of the above lemma.
 \end{proof}
 
 \begin{lemma}\label{Sec9:IntP1}
 Let $\rho_0, \gamma_0 \in \sfX$ and $\bnu \in G(\gamma_0) \subset \Tan_{\gamma_0} $. We define a curve
 \begin{align*}
 \sigma(t):= \big( \pi^1 + t \pi^2 \big)_\# \bnu, \quad t \in [0,1].
\end{align*}
Then
\begin{enumerate}
\item the map 
\begin{align*}
t \mapsto \sfd^2\big(\sigma(t), \rho_0\big) - t^2 \int_{\R^{2d}} |v|^2 \bnu(dy, dv)
\end{align*}
is concave in $t \in [0,1]$;
\item  for Lebesgue a.e. $ t \in [0,1]$ where the $\dot{\sigma}(t)$ exits (Lemma~\ref{Sec9:IntP0}),   
\begin{align*}
  \frac{d}{dt_+} \sfd^2(\sigma(t), \rho_0)  =    \inf_{\bM} (-2)\int_{\R^{3d}} \big( ( x - y ) \cdot v \big) \bM(dy; dx, dv),
 \end{align*}
where the $\bM \in \mathcal P_2(\R^{3d})$ is over all those with 
$\pi^{1,2}_\# \bM = \Gamma^{\opt}(\sigma(t),\rho_0)$ and $\pi^{1, 3}_\# \bM = \dot{\sigma}(t)$. We note that $\Vert \dot{\sigma}(t) \Vert_{\sigma(t)} \leq \Vert \bnu\Vert_{\sigma(0)}$.
 \end{enumerate}
\end{lemma}
 \begin{remark}
 Note that, because of $\bnu \in G(\gamma_0)$, $\sigma$ is a geodesic curve for $t \in [0, \delta]$ for some $\delta>0$. However,  this may fail to hold for $t \geq \delta$.
 \end{remark}
\begin{proof}
The first claim is a part of Theorems 7.3.2 of Ambrosio, Gigli and Savar\'e~\cite{AGS08}. 
The second calim is just a special case of Remark~\ref{dggamma}.
\end{proof}

\begin{lemma}\label{Sec9:IntP2}
Let $f_1 \in  {\mathcal S}^{-,\infty}$ and the $\beta_k$s be defined according to \eqref{betaDef} (the $g$ there is the $f_1$ here). We define $\bnu$ and $\sigma(t)$ according to Lemma~\ref{Sec9:IntP1}. Then  the following holds:
\begin{enumerate}
\item $\dot{\sigma}(0) = \bnu$;
\item $\dot{\sigma}(t)$ exists a.e. $t \in [0,1]$ (see Lemma~\ref{Sec9:IntP0});
\item there exists a finite constant $C_{f_1}>0$ which only depends on $f_1$,
and a modulus of continuity $\omega_{f_1}$ which depends on $f_1$ and 
$ C:=\sup_{t \in [0,1]} \int |y|^2 \sigma(t;dy)<\infty$, such that
for those $t \in [0,1]$ that $\dot{\sigma}(t)$ exists, we have
\begin{align*}
& \big(d_{\sigma(0)}f_1\big)\big(\dot{\sigma}(0)\big) \leq 
   \big(d_{\sigma(t)}f_1\big)\big(\dot{\sigma}(t)\big)
      + C_{f_1} t \Vert \bnu\Vert^2_{\sigma(0)} 
  +  \omega_{f_1} \big(\sfd(\sigma(t), \sigma(0)\big) \Vert \bnu\Vert_{\sigma(0)}.
\end{align*}
\end{enumerate}
 \end{lemma}
\begin{proof}
The first claim follows because that $\sigma$ is a constant speed geodesic for short time. 

We prove the third claim next. By the concavity result in Lemma~\ref{Sec9:IntP1}, 
\begin{align*}
 [0,1] \ni t \mapsto \frac{d}{dt_+}  \Big( \sfd^2\big(\sigma(t), \rho_0\big) - t^2 \int_{\R^{2d}} |v|^2 \bnu(dy, dv) \Big)  = \frac{d}{dt_+}  \sfd^2\big(\sigma(t), \rho_0\big) - 2t \Vert \bnu \Vert_{\sigma(0)}^2
\end{align*}
is a non-increasing function. Note that $\beta_k \geq 0$.
Consequently, for every $t>0$ such that $\dot{\sigma}(t)$ exists, we have
\begin{align*}
 \big(d_{\sigma(0)}f_1\big)\big(\dot{\sigma}(0)\big)  
 &=  \sum_{k=1}^K - \beta_k\big(\sigma(0); \rho_1, \ldots, \rho_K\big)
  \frac{d}{dt}\big|_{t=0^+} \sfd^2(\sigma(0), \rho_k)  \\
  & \leq  
   \sum_{k=1}^K - \beta_k\big(\sigma(0); \rho_1, \ldots, \rho_K\big)
 \big( \frac{d}{dt_+} \sfd^2(\sigma(t), \rho_k) + 2t \Vert \bnu\Vert^2_{\sigma(0)} \big) \\
  & \leq \sum_{k=1}^K - \beta_k\big(\sigma(t); \rho_1, \ldots, \rho_K\big)
  \frac{d}{dt_+} \sfd^2(\sigma(t), \rho_k)  \\
  & \qquad \quad  + C_{f_1}  t \Vert \bnu\Vert^2_{\sigma(0)} 
  +  \omega_{f_1} \big(\sfd(\sigma(t), \sigma(0)\big)\Vert \bnu\Vert_{\sigma(0)}  \\
   &= \big(d_{\sigma(t)}f_1\big)\big(\dot{\sigma}(t)\big) + C_{f_1} t \Vert \bnu\Vert^2_{\sigma(0)} 
  +  \omega_{f_1} \big(\sfd(\sigma(t), \sigma(0)\big) \Vert \bnu\Vert_{\sigma(0)} .
 \end{align*}
The last inequality above follows from estimate on $\frac{d}{dt_+} \sfd^2(\sigma(t), \rho_k)$, which can be obtained from the second part of Lemma~\ref{Sec9:IntP1}.
\end{proof}

With the above preparations, we give the main result of this subsection.
 \begin{lemma}\label{Sec9:Pmax}
 For each $f_1 \in {\mathcal S}^{-,\infty}$ and $\gamma \in \sfX$, there exists a 
 $\sigma (\cdot) \in AC([0,\infty);\sfX)$ with $\sigma(0) =\gamma$ and (in sense of Definition~\ref{veloCurv})
\begin{align*}
\bnu(t):=\dot{\sigma}(t)   \in \Tan_{\sigma(t)}, \quad t \text{ a.e.},
\end{align*}
such that  
\begin{align*}
 \int_0^t ({\mathbb H}_1 f_1)\big(\sigma(r)\big) dr 
 &\leq  f_1\big(\sigma(t)\big) - f_1\big(\sigma(0)\big) - \int_0^t   L\big(\bnu(r)\big)  dr \\
 & = \int_0^t \Big( \big(d_{\sigma(r)} f_1\big)\big(\bnu(r)\big) - L\big(\bnu(r)\big)\Big) dr, \qquad \forall t >0.
 \end{align*}
 \end{lemma}
 \begin{proof}
It is sufficient to construct such a curve $\sigma(\cdot) \in AC([0,1];\sfX)$. 

{\bf Step one: Constructing approximate curves.}
For each $n \in \N$, we partition the $[0,1]$ into equally sized intervals $0:=t_0< t_1< \ldots, t_n:=1$. We define a curve $\sigma_n(\cdot) \in C([0,1];\sfX)$ through iteration: Let $\sigma_n(t_0)=\gamma$. For each $t_i$, let $\bnu_n(t_i) \in G(\sigma_n(t_i))$ be such that (see Lemma~\ref{Sec7:HfG})
\begin{align}\label{Sec9:aOptMu}
 {\bf H} f_1(\sigma_n(t_i)) \leq \frac1n + (d_{\sigma_n(t_i)}f_1)(\bnu_n(t_i)) - L(\bnu_n(t_i)).
\end{align}
We construct curves:
\begin{align*}
 \bnu_n(t) : = \big(\pi^1+ (t-t_i) \pi^2, \pi^2\big)_\# \bnu_n(t_i),
  \quad \forall t \in [t_i, t_{i+1}]; \qquad \sigma_n(t) : =  \pi^1_\# \bnu_n(t).
\end{align*}
Recall that such $\sigma_n$ is geodesic only for $t \in [t_i, t_i+\delta_i]$ for some $\delta_i>0$, which may not be big enough to cover $[t_i, t_{i+1})$.  Nevertheless, the above construction gives estimate
\begin{align}\label{Sec9:modCur}
 \sfd\big(\sigma_n(t), \sigma_n(s)\big) \leq \Vert \bnu_n(t_i) \Vert_{\sigma_n(t_i)} (t-s), \quad s,t \in [t_i, t_{i+1}].
\end{align}
In particular,
\begin{align*}
 \Vert \dot{\sigma}_n(t) \Vert_{\sigma_n(t)} \leq \Vert \bnu_n(t_i) \Vert_{\sigma_n(t_i)}, 
 \quad  \text{a.e. } t \in [t_i, t_{i+1}].
\end{align*}
See Section~\ref{appmass} for notations regarding mass transport theory.

Next, we verify that 
\begin{align}\label{Sec9:nuEst}
 \sup_{n \in \N} \sup_{i=1,\ldots, n}\Vert \bnu_n(t_i)\Vert_{\sigma_n(t_i)} < \infty.
\end{align}
On one hand, by the estimate in Lemma~\eqref{PbarH} about $\bar{\sfL}$, there exists finite constants $c, \tilde{C}, C>0$ such that for $n$ sufficiently large, there exists an $\epsilon_0>0$,
\begin{align*}
{\bf H} f_1(\sigma_n(t_i)) & \leq  \tilde{C} + |D^+_{\sigma_n(t_i)} f_1| \Vert \bnu_n(t_i)\Vert_{\sigma_n(t_i)} -c \Vert \bnu_n(t_i)\Vert_{\sigma_n(t_i)}^2 - \langle U + V*\sigma_n(t_i), \sigma_n(t_i)\rangle \\
& \leq C- \epsilon_0 \Vert \bnu_n(t_i)\Vert_{\sigma_n(t_i)}^2 - \langle U + V*\sigma_n(t_i), \sigma_n(t_i)\rangle.
\end{align*}
On the other hand, taking $\bnu_{n,0}(t_i):= \bnu_{n,0}(dy,dv;t_i):= \delta_0(dv) \sigma_n(dy;t_i)$, we also have
\begin{align*}
 {\bf H} f_1\big(\sigma_n(t_i)\big)  
  \geq (d_{\gamma_n(t_i)} f_1)\big(\bnu_{n,0}(t_i)\big) - L\big(\bnu_{n,0}(t_i)\big) = 0 - 
  \bar{\sfL}(0) - \langle U + V*\sigma_n(t_i), \sigma_n(t_i)\rangle.
\end{align*}
Consequently, \eqref{Sec9:nuEst} holds.

 From \eqref{Sec9:modCur} and \eqref{Sec9:nuEst}, we conclude
\begin{align*}
 \sup_{n\in N} \sup_{t \in [0,1]} \int_{\R^d} |y|^2 \sigma_n(t; dy) <\infty.
\end{align*}
Choosing any metric that gives the narrow convergence on ${\mathcal P}(\R^d)$, 
by a version of Arzel\`a-Ascoli theorem, $\sigma_n(\cdot)$ converges to a limiting trajectory $\sigma(\cdot) \in C\big([0,1];{\mathcal P}(\R^d)\big)$.
By Fatou's lemma and the estimate \eqref{Sec9:modCur}, one further conclude that $\sigma(\cdot) \in AC([0,1];\sfX)$. 

Since
\begin{align*}
\sup_{n \in \N} \int_{(t,x,v) \in [0,1] \times \R^d \times \R^d} 
\big( |t|^2+ |x|^2 + |v|^2 \big) \Big(\bnu_n(t; dy, dv) dt\Big)<\infty,
\end{align*}
the measure $\bnu_n(t; dy, dv) dt$ is tight in ${\mathcal P}([0,1] \times \R^d \times \R^d)$, hence relatively compact in the narrow convergence topology. Since the marginal measure of the time-variable is always $dt$, we have that any limiting measure has to be of the form $\bnu(t;dy,dv) dt \in {\mathcal P}([0,1] \times \R^d \times \R^d)$. By Fatou's lemma, $\bnu(t) \in \mathcal P_2(\R^d \times \R^d)$ for Lebesgue a.e. $t \in [0,1]$.
Choose a convergent subsequence and relabel if necessary, we write
\begin{align*}
\bnu_n(t;dy,dv) dt \Rightarrow \bnu(t;dy,dv) dt.
\end{align*}
We note that $\pi^1_\# \bnu(t)=\sigma(t)$ Lebesgue a.e. $t\in [0,1]$.

{\bf Step two: Limiting curve satisfies continuity equation.}
We want to show that  
\begin{align*}
\bnu(t)=\dot{\sigma}(t), \quad \text{ Lebesgue a.e. in } t \in [0,1], 
\end{align*}
We write, for each $t \in [0,1]$,
\begin{align*}
 u(t,x) :=  \int_{v \in \R^d}  v \bnu(t;dv|x), \quad \text{ a.e. } \sigma(t;dx).
\end{align*}
In view of Lemma~\ref{Sec9:IntP0}, it is sufficient to show that 
\begin{align*}
 \partial_t \sigma + \div_x( \sigma u) =0 \in {\mathcal D}^\prime((0,1) \times \R^d).
\end{align*}

Take an arbitrary $\varphi:=\varphi(t,x) \in C_c^\infty((0,1) \times \R^d)$. 
By construction of the $\bnu_n$,  
\begin{align*}
0&=\lim_{n \to \infty} \int_{(t,x,v)\in (0,1) \times \R^d \times \R^d} \big((\partial_t  - v \cdot \nabla_x )\varphi(t,x)\big) \bnu_n(t, dx,dv)dt \\
& =   \int_{(t,x)\in (0,1) \times \R^d \times \R^d} \big( (\partial_t - v\cdot \nabla_x )\varphi(t,x)\big)  \bnu(t,dx,dv)dt \\
 & = \int_{(t,x)\in (0,1) \times \R^d}  \big(\partial_t\varphi(t,x)
  - \nabla_x \varphi(t,x) \cdot u(t,x)\big) \sigma(t; dx) dt.
  \end{align*}

{\bf Step three: Another class of approximating curves and some limiting inequalities.}
We also consider measure-valued piece-wise constant curves  
\begin{align*}
\hat{\bnu}_n(t):= \bnu_n(t_i), \quad \hat{\sigma}_n(t) 
:= \pi^1_\# \hat{\bnu}_n(t) =\sigma_n(t_i),
\quad \forall  t \in [t_i, t_{i+1}).
\end{align*}
Two Wasserstein distance estimates follow
\begin{align}\label{Sec9:nuhnu}
 \sfd(\sigma_n(t), \hat{\sigma}_n(t))\leq \sfd(\bnu_n(t), \hat{\bnu}_n(t)) \leq |t-t_i| \Vert \bnu_n(t_i)\Vert_{\sigma_n(t_i)}, 
 \quad t \in [t_i, t_{i+1}).
\end{align}

Using similar arguments for the measures $\bnu_n$, we have that $\{ \hat{\bnu}_n(t; dy, dv) dt\}_n$ is tight in ${\mathcal P}([0,1] \times \R^d \times \R^d)$, hence relatively compact in the narrow convergence topology.  In view of \eqref{Sec9:nuhnu}, a
limiting measure can be chosen to be of the form $\hat{\bnu}(r;dy,dv) dt \in {\mathcal P}([0,1] \times \R^d \times \R^d)$ with $\hat{\bnu}(t)= \bnu(t)$ a.e. in $t$. That is
\begin{align*}
\hat{\bnu}_n(t;dy,dv) dt \Rightarrow  \bnu(t; dy,dv) dt .
\end{align*}
By Fatou's lemma,
\begin{align*}
 & \liminf_{n \to \infty}   \sum_{s \leq \ldots < t_i < \ldots \leq t} L\big(\bnu_n(t_i)\big)({t_{i+1}}-{t_i}) \\
 & = \liminf_{n \to \infty} \int_s^t \int_{\R^{2d}} 
 \bar{\sfL}_{U,V}\big(x,v;\hat{\sigma}_n(r)\big) \hat{\bnu}_n(r; dx, dv) dr \\
 & \geq  \int_s^t \int_{\R^{2d}} 
 \bar{\sfL}_{U,V}\big(x,v;\hat{\sigma}_n(r)\big) \bnu(r; dx, dv) dr
   =\int_s^t  L\big(\bnu(r)\big) dr. 
\end{align*}

Noting ${\bH} f_1 \geq {\mathbb H}_1 f_1$ (see~\eqref{Sec7:IneqH1s}),  by Lemma~\ref{Sec9:Hf1LSC} and by Fatou's lemma (in view of \eqref{Sec9:nuhnu}),
\begin{align*}
\liminf_{n \to \infty}   \sum_{s \leq \ldots < t_i < \ldots \leq t}  {\bH}f_1\big(\sigma_n(t_i)\big) ({t_{i+1}}-{t_i})  
   \geq \liminf_{n \to \infty} \int_s^t {\mathbb H}_1f_1\big(\hat{\sigma}_n(r)\big) dr 
    \geq \int_s^t {\mathbb H}_1f_1\big(\sigma(r)\big) dr.
\end{align*}

By Lemma~\ref{Sec9:IntP2} and estimates \eqref{Sec9:modCur} and \eqref{Sec9:nuEst}, there exists a modulus $\omega$ (uniform with respect to $n$) such that 
\begin{align*}
  (d_{\sigma_n(t_i)}f_1)(\bnu_n(t_i)) - (d_{\sigma_n(t)}f_1)(\dot{\sigma}_n(t))  
  \leq   \omega(|t_{i+1} - t_i|), \quad t \in [t_i, t_{i+1}).
\end{align*}
Consequently,  
\begin{align*}
& \limsup_{n \to \infty}   \sum_{s \leq \ldots < t_i < \ldots \leq t}  \big(d_{\sigma_n(t_i)}f_1\big)\big(\bnu_n(t_i)\big) ({t_{i+1}}-{t_i}) \\
&  \leq  \limsup_{n \to \infty} 
\int_s^t \big(d_{\sigma_n(r)}f_1\big)\big(\dot{\sigma}_n(r)\big)dr
\leq  \limsup_{n \to \infty} f_1(\sigma_n(t)) -  \liminf_{n \to \infty} f_1(\sigma_n(s)) .
\end{align*}
 
 {\bf Step four: Conclusion.}
Combine the above estimates together (and in view of \eqref{Sec9:aOptMu}),
\begin{align*} 
\int_0^t {\mathbb H}_1f_1\big(\sigma(r)\big) dr \leq  \limsup_{n \to \infty} f_1(\sigma_n(t)) -  \liminf_{n \to \infty} f_1(\sigma_n(0)) - \int_0^t L \big(\bnu(r)\big) dr. 
\end{align*}
 Noting $\limsup_{n \to \infty} f_1(\sigma_n(t)) \leq f_1(\sigma(t))$ and that $\sigma_n(0) =\gamma$ is fixed, we conclude.   \end{proof}
%
%
 
 \subsubsection{Resolvent estimates lead to viscosity solution property}
\begin{lemma}\label{Sec9:gen}
\begin{align*}
{\bf R}_\alpha (f_0- \alpha {\bf H} f_0) & \leq f_0, 
 \quad f_0 \in {\mathcal S}^{+,\infty}; \\
{\bf R}_\alpha (f_1 - \alpha {\mathbb H}_1 f_1) & \geq f_1, 
 \quad f_1 \in {\mathcal S}^{-,\infty}.
\end{align*}
\end{lemma}
\begin{proof}
The proof of Lemma 8.19 in Feng and Kurtz~\cite{FK06} works here: Lemma~\ref{Sec9:Pmax} verifies the required Condition 8.11 in \cite{FK06}. We also note that for every $\sigma(\cdot) \in AC([0,\infty);\sfX)$, 
\begin{align*}
f(\sigma(t)) - f(\sigma(s))  = \int_s^t \big(d_{\sigma(r)} f\big) (\dot{\sigma}(r)) dr, 
 \quad \forall f \in {\mathcal S}^{+,\infty} \cup {\mathcal S}^{-,\infty}.
 \end{align*}
\end{proof}
\begin{remark}
We point out that equation (8.15) in Condition 8.11 in \cite{FK06} involves time integrals of the form $\int_{t_1}^{t_2}\ldots$ for every $0\leq t_1 <t_2$. However, we indeed only need a slightly weaker version of that condition involving integrals of the form $\int_0^t \ldots$ for $t \geq 0$. See the proof of Lemma 8.19 in middle of page 147 in \cite{FK06}, which is the only place that condition is used.  Lemma~\ref{Sec9:Pmax} in this paper verified this weaker version, which is good enough for  Lemma~\ref{Sec9:gen}. 
\end{remark}

\begin{lemma}\label{Sec9:fGrow}[Growth and modulus estimate]
Suppose that $h: \sfX \mapsto \R$  is such that $\sup_\sfX h<+\infty$.
Then 
\begin{enumerate}
\item the $f:={\bf R}_\alpha h$ is bounded from above $\sup_\sfX f<+\infty$;
\item there exists a non-decreasing sub-linear function $\beta: \R_+ \mapsto \R$ such that
\begin{align*}
   f (\rho) - h(\rho) \geq - \beta \circ \sfd(\rho, \delta_0);
\end{align*}
\item if the $h$ is bounded below in $\sfd$-balls of finite radius
\begin{align*}
 \inf_{\substack{\sigma \in \sfX \\ \sfd(\sigma,\delta_0)\leq R}} h(\sigma) > -\infty, \quad \R \in \R_+,
\end{align*}
then for each $R \in \R_+$, there exists a modulus of continuity $\omega_R \in C(\R_+;\R_+)$ such that 
\begin{align*}
  f(\rho) - f(\gamma) \leq \omega_R\big(\sfd(\rho,\gamma)\big), 
   \quad \forall \rho, \gamma \in \sfX, \text{ with } \sfd(\rho,\delta_0) +\sfd(\gamma,\delta_0) \leq R.
\end{align*}
\end{enumerate}
\end{lemma}
\begin{proof}
From the definition in \eqref{Sec9:R} and Conditions~\ref{U0CND}, \ref{VCND}, we know that $\sup_\sfX f<+\infty$. The existence of sub-linear function $\beta$ can be proved using same method as in the first part of Lemma~\ref{parHJprop}.

The modulus of continuity part  follows from essentially the same proof of Lemma~\ref{Sec5:fnmest}.
\end{proof}

\begin{lemma}\label{Sec9:Resol}
For each $h: \sfX \mapsto \R$ with $\sup_\sfX h <+\infty$, we have
\begin{align*}
{\bf R}_\alpha h = {\bf R}_\beta \Big( {\bf R}_\alpha h - \beta \frac{{\bf R}_\alpha h - h}{\alpha} \Big), \qquad \forall \alpha >\beta >0.
\end{align*}
\end{lemma}
\begin{proof}
The proof of Lemma 8.20 in \cite{FK06} works here.
\end{proof}

\begin{lemma}\label{Sec9:Rcntr}
Suppose $h_i: \sfX \mapsto \R$ for $i=1,2$ is such that
$\sup_\sfX h_1 <+\infty$ and $h_2 \geq -\beta \circ\sfd(\cdot, \delta_0)$ for some non-decreasing, sub-linear function $\beta: \R_+ \mapsto \R_+$.
Then
\begin{align*}
 \sup_\sfX \big( {\bf R}_\alpha h_1 - {\bf R}_\alpha h_2 \big) \leq \sup_\sfX (h_1 - h_2), \qquad \forall \alpha >0.
\end{align*}
\end{lemma}
\begin{proof}
The same proof of Lemma 8.21 in \cite{FK06} works here. 
\end{proof}

\begin{lemma}\label{Sec9:bfHH1Cmp}
Let $\alpha>0$, $h \in C(\sfX)$ with $\sup_\sfX h <+\infty$ and 
$h \geq -\beta \circ \sfd(\cdot,\delta_0)$ for some non-decreasing sub-linear function $\beta: \R_+ \mapsto \R_+$.   Then $f:={\bf R}_\alpha h$ is a sub-solution in the sequential viscosity solution sense to 
\begin{align}\label{Sec9:bfHeqn}
 f - \alpha {\bf H} f \leq h;
\end{align}
and a super-solution in the sequential viscosity solution sense to 
\begin{align}\label{Sec9:supbbH1}
 f - \alpha {\mathbb H}_1 f \geq h.
\end{align}
\end{lemma}
\begin{proof}
Combining Lemmas~\ref{Sec9:gen} \ref{Sec9:Resol} and \ref{Sec9:Rcntr}, the same method of proof in Theorem 8.27 in \cite{FK06} gives the sequential viscosity sub- and super- solution properties.
\end{proof}

\subsection{Continuity of the $f = {\bf R}_\alpha h$}
We can obtain such continuity through direct estimates as in the proof of Lemma~\ref{Sec9:fGrow}.
We can also obtain the continuity indirectly through the following comparison arguments.

Let  
\begin{align*}
f_*:=f_*(\gamma):= \lim_{\epsilon \to 0^+} \inf \{ f(\rho) : \sfd(\rho,\gamma)<\epsilon\}  
\end{align*}
be a lower semicontinuous regularization of the $f$, with respect to the metric $\sfd$. 
\begin{lemma} \label{Sec9:fReg}
Let $h$ satisfy the same condition as in Lemma~\ref{Sec9:bfHH1Cmp}. Moreover, we assume that the $h$ has modulus of continuity in every bounded $\sfd$-metric balls as assumed in Theorem~\ref{Sec7:CMP}.  Then
\begin{enumerate}
\item $f \in \USC(\sfX)$, and is a point-wise strong viscosity sub-solution for \eqref{Sec9:bfHeqn}, and for equation
\begin{align}\label{Sec9:subbfH0}
 f - \alpha {\bf H}_0 f \leq h.
\end{align}
\item $f_* \in \LSC(\sfX)$ is a point-wise strong viscosity super-solution for 
\eqref{Sec9:supbbH1}.
\item Indeed, $f =f_* \in C(\sfX)$.
\end{enumerate}
\end{lemma}
\begin{proof}
First, since $\bnu \mapsto L(\bnu)$ is lower-semicontinuous in the weak convergence of probability measure (i.e. narrow convergence) topology, using the method of proof in Lemma 8.17 in \cite{FK06}, we can conclude that $f \in \USC(\sfX)$. It can be verified that ${\bf H}$ is a local operator satisfying the property described in \eqref{Sec3:gapp} in Lemma~\ref{Sec3:seq2spw}. In view of the sequential sub-solution result in Lemma~\ref{Sec9:bfHH1Cmp}, apply Lemma~\ref{Sec3:seq2spw}, the $f$ is also a point-wise strong viscosity sub-solution for \eqref{Sec9:bfHeqn}, and for \eqref{Sec9:subbfH0} (see  Lemma~\ref{Sec7:Hcom}).
 
Second, by Lemma~\ref{Sec9:bfHH1Cmp}, we know that the $f$ is sequential viscosity super-solution to \eqref{Sec9:supbbH1}. We note that domain $D({\mathbb H}_1) \subset C(\sfX)$ 
and each ${\mathbb H}_1 f_1 \in \LSC(\sfX)$ (Lemma~\ref{Sec9:Hf1LSC}). Using the method of proof in Theorem 8.27 in the last line on page 153 of \cite{FK06} (which uses Lemma~\ref{App:dissip} -- see the first line on page 154 of \cite{FK06}),  we can verify that the $f_*$ is a sequential viscosity super-solution for \eqref{Sec9:supbbH1} as well. The ${\mathbb H}_1$ verifies also a super-solution version of \eqref{Sec3:gapp}. By Remark~\ref{Sec3:seq2SUP}, the $f_*$ is a point-wise strong viscosity super-solution to \eqref{Sec9:supbbH1}. 

Third, we define Yosida regularization of the $f_*$ as 
\begin{align*}
f_{*,\epsilon}(\gamma) := \inf_{\rho \in \sfX}\Big( f_*(\rho) 
 + \frac{\sfd^2(\rho,\gamma)}{2\epsilon}\Big) \leq f_*(\gamma).
\end{align*}
From Lemma~\ref{Sec8:H1bfH1Alt}, it follows that the $f_{*,\epsilon}$ is a strong point-wise super-solution to 
\begin{align*}
 f_{*,\epsilon} - \alpha {\bf H}_1 f_{*,\epsilon} \geq h_\epsilon,
\end{align*}
where the $h_\epsilon$ is defined as in \eqref{Sec8:h1espD} with the $h_1$ replaced by $h$.
Next, we introduce perturbations
\begin{align*}
f_{*;\epsilon,\lambda,\theta}(\gamma) & := \lambda^{-1} f_{*,\epsilon}(\gamma)
 + \theta \sqrt{1+\sfd^2(\gamma,\delta_0)}, \\
h_{\epsilon,\lambda,\theta}(\gamma) & :=  \lambda^{-1} h_\epsilon(\gamma) 
 +\theta \sqrt{1+\sfd^2(\gamma,\delta_0)} - \alpha {\rm Err}_{\lambda,\theta}(\gamma),
\quad \forall \theta>0, \lambda>1. 
\end{align*}
See \eqref{Sec8:epsrho} for the definition of ${\rm Err}_{\lambda,\theta}$ term.
Then, according to Lemma~\ref{Sec8:bfH1Pert}, the $f_{*;\epsilon,\lambda,\theta}$ is a strong point-wise super-solution to 
\begin{align*}
 (I -\alpha {\bf H}_1) f_{*;\epsilon,\lambda,\theta} \geq h_{\epsilon,\lambda,\theta}.
\end{align*}

Finally, we now are in a position to apply the comparison principle 
in Theorem~\ref{Sec7:CMP} to obtain  
  \begin{align*}
  f(\gamma) - \lambda^{-1} f_*(\gamma)
 - \theta \sqrt{1+\sfd^2(\gamma,\delta_0)} \leq 
\sup_\sfX ( f - f_{*;\epsilon,\lambda,\theta}) 
\leq  \sup_\sfX (h - h_{\epsilon,\lambda,\theta}), \quad \forall \gamma \in \sfX.
\end{align*}
We note that $\sup_\sfX h<\infty$ and that, for each $\theta>0$ fixed,
$\gamma \mapsto h_{\epsilon,\lambda,\theta}(\gamma)$ 
grows to $+\infty$ at a rate which is linear with respect to size of $\sfd$-metric balls. 
There exists finite constant $C_\theta>0$ such that 
\begin{align*}
\sup_\sfX (h - h_{\epsilon,\lambda,\theta}) \leq  
\sup_{\substack{\gamma \in \sfX \\ \sfd(\gamma,\delta_0) \leq C_\theta}}
 \big(h(\gamma) - h_{\epsilon,\lambda,\theta}(\gamma)\big) \vee 0.
\end{align*}
The above gives
\begin{align*}
 \limsup_{\lambda \to 1^+} \limsup_{\theta \to 0^+} \limsup_{\epsilon \to 0^+} 
  \sup_\sfX (h - h_{\epsilon,\lambda,\theta}) \leq 0.
\end{align*}
Consequently
\begin{align*}
  f(\gamma) -   f_*(\gamma)  \leq 0, \quad \forall \gamma \in \sfX.
\end{align*}
Hence $f = f_* \in C(\sfX)$.
\end{proof}
 
 \subsection{Weak upper semicontinuity of the $f={\bf R}_\alpha h$ in $\sfX$}
 \begin{lemma}\label{Sec9:fweakUSC}
Let $h : \sfX \mapsto \R$ be such that $\sup_\sfX h<+\infty$, and that $\sfd_{p=1}$-upper semi-continuous in $\sfX$ (See Definition~\ref{Sec6:weakUSC}), 
then the $f:={\bf R}_\alpha h$ is $\sfd_{p=1}$-upper semi-continuous in $\sfX$.
 \end{lemma}
 \begin{proof}
 The proof in Lemma 8.17 on page 145 of Feng and Kurtz~\cite{FK06} (which also uses the proof of Proposition 8.13 in \cite{FK06}) can be adapted here. 
 \end{proof}
 
 \subsection{Lagrangian representation}
 \begin{lemma}\label{Sec9:RhSol}
The $f={\bf R}_\alpha h$ in Lemma~\ref{Sec9:bfHH1Cmp} is a point-wise strong viscosity sub-solution to \eqref{Sec6:bbH0Sub}, and a point-wise strong viscosity super-solution to \eqref{Sec6:bbH1res}, with the $h_0=h_1=h$.
\end{lemma}
\begin{proof}
We recall the inequalities in Lemma~\ref{Sec7:Hcom}, when considering 
operator ${\mathbb H}_0$ in place of ${\bf H}_0$. 
The sub-solution property in Lemma~\ref{Sec9:fReg} implies that the $f$ is also a point-wise strong sub-solution to \eqref{Sec6:bbH0Sub}.
The case of super-solution is just the super-solution part of Lemma~\ref{Sec9:fReg}. 
\end{proof}
 
\begin{theorem}\label{Sec9:MainThm2}
There is a unique $f \in C(\sfX)$, which has at most sub-linear growth with respect to the $2$-Wasserstein metric $\sfd$, such that it is a sub-solution to \eqref{Sec6:bbH0Sub} and super-solution to \eqref{Sec6:bbH1res}, both in the point-wise strong viscosity sense. 

Moreover,
\begin{enumerate}
\item such $f={\bf R}_\alpha h$;
\item such $f$ is the same one as arising from limit \eqref{Sec6:fNCvgf} in Theorem~\ref{Sec8:MainThm1}.
\item \label{Sec9:Thm2Tri} in the context of Theorem~\ref{Sec8:MainThm1}, the convergent sequence
$(\{ {\mathfrak f}_N \}_{N\in \N}, f) \in {\mathcal C}$. 
\end{enumerate}
\end{theorem}
\begin{proof}
As in the proof of Theorem~\ref{Sec8:MainThm1}, through upper- and lower- Yosida approximations and proper perturbation arguments, we can apply the comparison principle in Theorem~\ref{Sec7:CMP} to conclude uniqueness for a function which is both sub-solution to 
\eqref{Sec6:bbH0Sub} and super-solution to \eqref{Sec6:bbH1res}, both in the point-wise strong viscosity sense.

The existence (hence representation of the solution) follows from Lemma~\ref{Sec9:RhSol}.

The rest of the conclusion follows by combining the above result with that of Theorem~\ref{Sec8:MainThm1}, and the properties we proved for $R_\alpha h$ in Lemma~\ref{Sec9:fGrow} and in Lemma~\ref{Sec9:fweakUSC}.
\end{proof}

\begin{remark}\label{Sec9:SGRmk}
Theorem~\ref{Sec9:MainThm2}.~\ref{Sec9:Thm2Tri} is a result on convergence for viscosity-solutions of ``resolvent" type problems. Indeed, such result also implies convergence of associated Cauchy (or nonlinear operator semigroup) type problems -- namely, convergence of $S_N(t) \mapsto S(t)$ with the $S_N, S$ defined respectively by \eqref{Sec1:SNSG} and \eqref{Sec1:SSG}, solving \eqref{Sec1:MCauN} and \eqref{Sec1:MCau} . 

In 1958, Trotter~\cite{Tro58} introduced an interesting method on semigroup convergence.
Subsequently, Kurtz~\cite{Kurtz69,Kurtz70,Kurtz73} generalized the method to more applicable settings. Through this type of techniques, convergence of semigroups follows from semigroup generation theorems on a sequence space. See Proposition (1-8) in \cite{Kurtz70} or Section 2 in \cite{Kurtz73} for quick introductions.   Although developed with linear operator semigroup setting in mind at the beginning, this method is readily adapted to nonlinear semigroup settings after Crandall and Liggett~\cite{CLigg71} discovered a nonlinear semigroup generation theorem.
In fact, the result of \cite{Kurtz73} is formulated on (possibly-) nonlinear semigroups.  Using modern viscosity solution language and techniques, the Crandall-Liggett semigroup generation theorem can be replaced by existence and uniqueness (through the comparison principle) and convergence of viscosity solutions. Assemble all these steps together,  Feng and Kurtz~\cite{FK06} adapted the above strategy to develop a viscosity solution convergence approach to the theory of large deviation for Markov processes in metric spaces. See Proposition 5.5 in \cite{FK06} for convergence posed in  nonlinear semigroup language, and then Theorem 7.17 there for a translation in viscosity solution language, in that book. Here, we can re-adapt the procedure to extend the resolvent convergence result in Theorem~\ref{Sec9:MainThm2} to semigroup convergence of the $S_N$s. Since such development is expected to be lengthy but relatively routine, we do not pursue details anymore.

We informally summarize the ingredients for showing semigroup convergence: We introduce operator (see VI.3. of Crandall and Lions~\cite{CL85} and the proof of Theorem 7.17 in \cite{FK06}) 
\begin{align*}
\hat{H} & := \cup_{\alpha>0} \Big\{\big( {\bf R}_\alpha h, \frac{{\bf R}_\alpha h -h}{\alpha} \big) :
  h \in C(\sfX), \sup_\sfX h<+\infty, \text{ and $h$ satisfies } \\
& \qquad \qquad \qquad \text{ requirements \eqref{Sec8:CC5}
  and \eqref{Sec8:CC6}  in Condition $\mathcal C$
   in Definition~\ref{Sec8:ClassC}}\Big\}.
\end{align*}
Such $\hat{H}$ satisfies range condition in semigroup theory:
\begin{align*}
 D(\hat{H}) \subset R(I - \alpha \hat{H}).
\end{align*}
With the semigroups $S_N(t)$ and $S(t)$ defined in \eqref{Sec1:SNSG} and \eqref{Sec1:SSG}, resolvents ${\bf R}_{N;\alpha}$ and ${\bf R}_\alpha$ 
in \eqref{Sec1:fiDVal} and \eqref{Sec1:RaMV},
we also introduce 
\begin{align*}
 {\mathcal S}(t) \big( \{ {\mathfrak f}_N \}_{N\in \N}, f\big)
& : = \big(\{ S_N(t) {\mathfrak f}_N \}_{N\in \N}, S(t) f\big), \\
{\mathcal R}_\alpha \big( \{ {\mathfrak f}_N \}_{N\in \N}, f\big)
& : = \big(\{ {\bf R}_{N;\alpha}   {\mathfrak f}_N \}_{N\in \N}, {\bf R}_\alpha f\big).  
 \end{align*}
The results (Theorems~\ref{Sec8:MainThm1}, \ref{Sec9:MainThm2}) in this paper allow us to apply the convergence method of Trotter-Kurtz, we obtain 
\begin{align*}
 S_N(t) {\mathfrak f}_N   \to S(t) f, \quad \text{ whenever } {\mathfrak f}_N \to f;
\end{align*}
with a notion of convergence properly defined. Moreover,
\begin{align*}
 {\mathcal S}(t) = \lim_{n \to \infty} {\mathcal R}_{n^{-1}}^{[nt]}. 
\end{align*}
In particular, from Theorem~\ref{Sec9:MainThm2}, we see that, in context of Theorem~\ref{Sec8:MainThm1}, 
$(\{ {\mathfrak h}_N \}_{N\in \N}, h) \in {\mathcal C}$ implies that 
$(\{ {\mathfrak f}_N \}_{N\in \N}, f) \in {\mathcal C}$. This implies that the $\mathcal C$ is an invariant set under the map ${\mathcal R}_\alpha$; 
hence the ${\mathcal S}(t)$ for every $t\geq 0$.
 \end{remark}
\newpage

 \appendix
 \section{Miscellaneous results on metric space}
We list some abstract concepts and results that we invoked in the main text regarding analysis in metric spaces.

\subsection{Semi-continuity}
Let $(\sfX, \sfd)$ be a metric space and $\Lambda$ be an index set.
\begin{lemma}\label{infsupSC}
If $f_\alpha \in \LSC(\sfX; \bar{\R})$ for every $\alpha \in \Lambda$, then $\sup_{\alpha \in \Lambda} f_\alpha \in \LSC(\sfX; \bar{\R})$. 
Suppose additionally that $\Lambda$ is a finite set, then $\min_{\alpha \in \Lambda} f_\alpha \in \LSC(\sfX;\bar{\R})$.
More generally, suppose that $(\Lambda,r) $ is a compact metric space and $(x,\alpha) \mapsto f_\alpha(x) \in \LSC(\sfX \times \Lambda; \bar{\R})$. Then $F:=\inf_{\alpha \in \Lambda} f_\alpha \in \LSC(\sfX; \bar{\R})$.

Similarly, if $f_\alpha \in \USC(\sfX; \bar{\R})$ for every $\alpha \in \Lambda$, then $\inf_{\alpha \in \Lambda} f_\alpha \in \USC(\sfX; \bar{\R})$. 
Suppose additionally that $\Lambda$ is a finite set, then $\max_{\alpha \in \Lambda} f_\alpha \in \USC(\sfX;\bar{\R})$.
More generally, suppose that $(\Lambda,r) $ is a compact metric space and $(x,\alpha) \mapsto f_\alpha(x) \in \USC(\sfX \times \Lambda; \bar{\R})$. Then $F:=\sup_{\alpha \in \Lambda} f_\alpha \in \USC(\sfX; \bar{\R})$.
\end{lemma}
\begin{proof}
We only verify the lower semi-continuous properties. The upper semi-continuous situation follows by replacing the $f_\alpha$s by $-f_\alpha$s and applying the lower semi-continuous results. 

The first two claims follow by definition. We verify the last one which assume that $\Lambda$ is compact.  Let $x_n, x_0\in \sfX$ be such that $\lim_{n \to \infty}\sfd(x_n,x_0)=0$. Then there exists $\alpha_n:=\alpha_n(x_n) \in \Lambda$ such that
\begin{align*}
 F(x_n) \geq f_{\alpha_n}(x_n) - \frac1n.
\end{align*}
 By compactness of $\Lambda$ and through extracting subsequence $\{ n(k) :k =1,2, \ldots\}$ if necessary, we have $\alpha_{n(k)} \to \alpha_0 \in \Lambda$ for some $\alpha_0$ and 
 \begin{align*}
 \liminf_{n \to \infty} F(x_n) \geq \liminf_{k \to \infty} f_{\alpha_{n(k)}}(x_{n(k)}) \geq f_{\alpha_0}(x_0) \geq F(x_0). 
\end{align*}
\end{proof}

\subsection{A slope estimate}
The following is a direct consequence of the definition of slopes in Definition~\ref{slopDef}.
\begin{lemma}
Let $(\sfX, \sfd)$ be a metric space and $f, f_0: \sfX \mapsto \R$. Suppose $x_0 \in \sfX$ is such that 
\begin{align*}
 f (x_0) - f_0(x_0) = \sup_\sfX( f - f_0).
\end{align*}
Then the following estimate for downward slopes hold
\begin{align*}
 |D_{x_0}^- f_0 | \leq |D^-_{x_0}f|.
\end{align*}
\end{lemma}
 
\subsection{Dissipativity in function spaces}\label{App:Diss}
Let $(\sfX,\sfd)$ be a metric space. The following is Lemma 7.8 of Feng and Kurtz~\cite{FK06}. The original proof contains an error because it implicitly used a condition which was not assumed.  However, the results remain true in the way originally stated. Below, we provide a new proof taken from Errata of \cite{FK06} for completeness. 
\begin{lemma}\label{App:dissip}
Let $f, g : \sfX \mapsto \bar{\R}$ and $f - \epsilon g \in M(\sfX; \bar{\R})$
\footnote{This means in particular that  $\infty - \infty$ or $-\infty + \infty$ won't occur for the $f - \epsilon g$.}
for every $\epsilon \in (0,\epsilon_0)$. Suppose that
\begin{align*}
-\infty< \sup_\sfX f \leq \sup_{\sfX} (f - \epsilon g)< \infty, \quad \forall \epsilon \in (0,\epsilon_0).
\end{align*}
Then there exists $x_n \in \sfX$ such that 
\begin{align*}
\lim_{n \to \infty} f(x_n) = \sup_\sfX f, \text{ and } 
\limsup_{n \to \infty} g(x_n) \leq 0.
\end{align*}
\end{lemma}
\begin{proof}
Let $(0,\epsilon_0) \in \epsilon_n \to 0$, we can choose $x_n \in \sfX$ such that 
\begin{align}\label{fgdiss}
 \sup_\sfX f  \leq \sup_\sfX (f - \epsilon_n g) < f(x_n) - \epsilon_n g(x_n) + \epsilon_n^2. 
\end{align}
From the above, we have
\begin{align*}
  g(x_n) < \epsilon_n^{-1} \big( f(x_n) - \sup_\sfX f \big) + \epsilon_n \leq \epsilon_n.
\end{align*}
To conclude the lemma, we only need to show $\limsup_{n \to \infty} f(x_n) \geq \sup_\sfX f$, which also follows from \eqref{fgdiss} provided we can establish estimate $\liminf_{n \to \infty} g(x_n) >-\infty$.

Let $\varphi_n(\epsilon):=  f(x_n) - \epsilon g(x_n) + \epsilon^2 - \sup_\sfX f$.
We observe that $\varphi_n(0) <0$ and $\varphi_n(\epsilon_n) >0$. By continuity of $\varphi_n$, there exists 
  $\epsilon_n^\prime \in (0,\epsilon_n)$ such that $\varphi_n(\epsilon_n^\prime) =0$, which gives
  \begin{align*}
 \lim_{n \to \infty} \big( f(x_n) - \epsilon_n^\prime g(x_n) \big) = \sup_\sfX f.
\end{align*}
Take a fixed $\epsilon \in (0, \epsilon_0)$, when $\epsilon_n^\prime<\epsilon$, 
\begin{align*}
 \sup_\sfX (f- \epsilon g) \geq f(x_n) -\epsilon g(x_n) 
     = \big( f(x_n) - \epsilon_n^\prime g(x_n)\big) - (\epsilon - \epsilon_n^\prime) g(x_n).        
    \end{align*}
Taking $n \to \infty$ gives the estimate 
 \begin{align*}
  \liminf_{n \to \infty} g(x_n) \geq \epsilon^{-1} \big(\sup_\sfX f - \sup_\sfX (f - \epsilon g)\big) > -\infty.
\end{align*}
\end{proof}

\subsection{Perturbed optimization principle}\label{App:BorPre}
Let $(\sfX, \sfd)$ be a complete metric space. Let $F \in \USC( \sfX; \R \cup\{-\infty\})$, $F \not \equiv -\infty$ and $\sup_\sfX F< +\infty$. 
We state a special version of the Borwein-Preiss~\cite{BorPre87} generalization on the Ekeland's perturbed optimization principle~\cite{Eke79}.
 \begin{lemma}[Borwein-Preiss]\label{BorPre}
Let $\epsilon >0$ and $x_0 \in \sfX$ be such that 
\begin{align*}
 F(x_0) > \sup_\sfX F - \epsilon.
\end{align*}
Then there exists a convergence sequence of  $\{ x_{\epsilon,k} \}_{k \in \N}\subset \sfX$ with limit point 
$x_\epsilon \in \sfX$ that has the following properties: 
By introducing a barrier function $\Delta: \sfX \mapsto \R$ given by
\begin{align}\label{App:DelDef}
\Delta (x):= \Delta_{\epsilon, x_0}(x):=\sum_{k=0}^\infty \frac{1}{2^{k+1}} \sfd^2(x, x_{\epsilon,k}),
\end{align}
and a perturbed function 
\begin{align*}
 F_\epsilon :=  F - \sqrt{\epsilon} \Delta,
\end{align*}
we have
\begin{enumerate}
\item $F_\epsilon (x_\epsilon) = \sup_\sfX F_\epsilon$;
\item $F(x_\epsilon) > \sup_\sfX F_\epsilon - \epsilon$;
\item $ \lim_{k \to \infty} \sfd(x_{\epsilon, k}, x_\epsilon) =0$, 
$\sup_{k\in \N} \sfd(x_{\epsilon,k},x_\epsilon) < \epsilon^{1/4}$ 
and $\sfd(x_\epsilon, x_0) < \epsilon^{1/4}$;
\item $|\Delta (x_\epsilon)| < \sqrt{\epsilon}$;
\item  the following estimate on local Lipschitz constant holds
\begin{align*}
|D_{x_\epsilon} \Delta | := \limsup_{y \to x_\epsilon} \frac{|\Delta(y) - \Delta(x_\epsilon)|}{\sfd(y, x_\epsilon)} \leq 2 \epsilon^{1/4}.
\end{align*}
\end{enumerate}
\end{lemma}
\begin{proof}
Following the proof of Theorem 2.6 in \cite{BorPre87}, we take 
$g:= -F: \sfX \mapsto \R \cup\{+\infty\}$, $p=2$, $\lambda=\epsilon^{1/4}$ and $\epsilon_1=\frac{\epsilon}{3}$, and select $\mu =\frac12$, then the conclusions follow.
\end{proof}

\subsection{Submetry and metric foliations}\label{App:SubM}
We denote $(\sfX, \sfd_X)$ and $(\sfY, \sfd_\sfY)$ two metric spaces. 
\begin{definition}[Submetry and Strong Submetry]\label{subMDef}
A map $\sfp: \sfY \mapsto \sfX$ is called a {\em submetry}, if 
\begin{align*}
\sfp\big(B_\sfY(y,r)\big) = B_\sfX\big(\sfp(y),r\big), \quad \forall y \in \sfY, r \geq 0.
\end{align*}
 In the above, $B(y,r)$ is an open ball with radius $r$. 
 We call $\sfp$ a {\em strong submetry}, if the open balls above are replaced by closed balls
\begin{align*}
\sfp\big(\bar{B}_\sfY(y,r)\big) = \bar{B}_\sfX\big(\sfp(y),r\big), \quad \forall y \in \sfY, r \geq 0.
\end{align*}
 \end{definition}
The closed ball formulation was the original one that Berestovskii used, when first introducing the concept of submetry.
It follows from definition that, a submetry is a continuous and open, surjective map.
In particular, $\sfp^{-1}(x)$ is a closed subset in $\sfY$ for every $x \in \sfX$.
 It also follows that, if $\sfp^{-1}(x)$ is proper in $\sfY$ (i.e. ball compact for every finite radius balls) for every $x \in \sfX$, then $\sfp$ being a submetry implies that it is a strong submetry.

Submetry is a generalization of submersion to metric space setting. Therefore, we expect the structure of submetry can be viewed from a different perspective using foliations.  A result from Galaz-Garc\'ia, Kell, Mondino and Sosa~\cite{GKMS18} confirms this. 

\begin{definition}[Foliation in metric spaces]\label{MFoliate}
 A partition $\mathcal F$ of a metric space $(\sfY, \sfd_Y)$ into a family  of closed disjoint subsets
\begin{align*}
 \sfY:= \bigsqcup_{\mathcal F_\alpha \in \mathcal F} \mathcal F_\alpha,
\end{align*}
is called a {\em foliation}. Each $\mathcal F_\alpha$ is called a leaf. 

If, in addition, the foliation $\mathcal F$ satisfies the following {\em equi-distant property}
\begin{align*}
 \sfd_\sfY(\mathcal F_\alpha,  \mathcal F_{\alpha^\prime}) 
 = \sfd_\sfY(y, \mathcal F_{\alpha^\prime}), \quad \forall y \in \mathcal F_\alpha,  
\end{align*} 
then we call the foliation $\mathcal F$ as a {\em metric foliation}, and the $\sfY$ is {\em metrically foliated} by $\mathcal F$.
 \end{definition}
 Note that distance between two subsets is defined, as always, as $\sfd(A,B):=\inf_{x \in A, y \in B} \sfd(x,y)$. 
 
Let $(\sfY, \sfd_\sfY; \mathcal F)$ be a metric foliation. If we denote $\mathcal F_{[y]}$ the leaf containing $y$, then this induces an equivalent relation $\sim$  
 \begin{align*}
 y_1 \sim y_2 \text{ if and only if } \mathcal F_{[y_1]} = \mathcal F_{[y_2]},
\end{align*}
Quotient space $\sfX:=\sfY/{\sim}$ is the set of equivalence classes. 
If we denote $\sfp: \sfY \mapsto \sfX$ the projection onto the quotient space.
Then for each $x \in \sfX$,  there is a canonical association of leaf $\mathcal F_{x}:=\sfp^{-1}(x) \in \mathcal F$,
and it follows that $\mathcal F=\{ \mathcal F_{x} : x \in \sfX \}$. Moreover, it can be directly verified that the following defines a metric
\begin{align}\label{QMetDef}
\sfd_\sfX(x_1,x_2) := \inf \{ \sfd_\sfY(y_1 ,y_2 ) : y_1 \in \sfp^{-1}(x_1), y_2 \in \sfp^{-1}(x_2)  \};
\end{align}
and that $(\sfX, \sfd_\sfX)$ is a metric space.
 
 There is a 1-1 correspondence, up to an isometry, between submetry and metric foliations.
\begin{lemma}\label{subMFol}
Suppose that metric space $\sfY$ is metrically foliated into $\sfX:=\sfY/\sim$ with the natural projection map $\sfp$. Then the $\sfp : \sfY \mapsto \sfX$ is a submetry.

Suppose that $\sfX, \sfY$ are two metric spaces and $f: \sfY \mapsto \sfX$ is a submetry. 
Then the foliation given by $\bigsqcup_{x \in \sfX} f^{-1}(x)$ is a metric foliation. Moreover, let $\sfY^*:=\sfY/\sim$ denote the quotient space induced by the foliation and $\sfp: \sfY \mapsto \sfY^*$ the natural projection. Then there is an isometry 
$\iota_f : \sfX \mapsto \sfY^*$ such that 
\begin{align*}
  \iota_f \circ f =\sfp.
\end{align*}
\end{lemma}
\begin{proof}
This is Lemma 8.4 of Galaz-Garc\'ia, Kell, Mondino and Sosa~\cite{GKMS18}. 
\end{proof}

The notion of strong submetry can be equivalently viewed through the following 2-point property.
\begin{definition}[2-point lifting property]
A map $\sfp: \sfY \mapsto \sfX$ is said to have {\em $2$-point lifting property}, if for each $x_1, x_2 \in \sfX$ and $y_1 \in \sfp^{-1}(x_1) \subset \sfY$, there exists $y_2 \in \sfp^{-1}(x_2)$ such that $\sfd_\sfY(y_1, y_2) = \sfd_\sfX(x_1, x_2)$.
\end{definition} 
\begin{lemma}\label{2pts}
A strong submetry $\sfp: \sfY \mapsto \sfX$ has the $2$-point lifting property.
In addition, within the class of $1$-Lipschitz maps, $2$-point lifting property implies strong submetry. 
\end{lemma}
\begin{proof}
First, we assume the $2$-point lifting property. Let $x_i \in \sfX$, $i=1,2$ and $R>0$ be such that 
$\sfd_\sfX(x_1, x_2) \leq R$. By the $2$-point lifting property, there exists
$y_i\in \sfY$ with $\sfp(y_i) = x_i$, $i=1,2$, such that $\sfd_\sfY(y_1,y_2) = \sfd_\sfX(x_1, x_2) \leq R$. Therefore 
\begin{align*}
\bar{B}_\sfX(\sfp(y_1);R) \subset  \sfp(\bar{B}_\sfY(y_1,R)).
\end{align*}
In addition, the $1$-Lipschitz property implies a reversed inclusion relation holds in the above as well. Consequently, $\sfp$ is a strong submetry.

Second, we assume that $\sfp$ is a strong submetry. 
For $x_i \in \sfX$ and $y_1 \in \sfY$ with 
$\sfp(y_1) = x_1$, let $\sfd_\sfX(x_1, x_2) =R$. One one hand, from $x_2 \in \bar{B}(\sfp(y_1);R) \subset \sfp(\bar{B}_\sfY(y_1,R))$, we can find $y_2 \in  \sfY$ with $x_2 = \sfp(y_2)$ and $\sfd_\sfY(y_1, y_2) \leq R = \sfd_\sfX(x_1, x_2)$. On the other hand, $\sfp: \sfY \mapsto \sfX$ being $1$-Lipschitz map means $\sfd_\sfX(x_1, x_2) \leq \sfd_\sfY(y_1, y_2)$. Hence the two are equal, giving the $2$-point lifting property.
\end{proof}

%

\subsection{Quotient given by isometric actions of groups}\label{App:MQ}
 A large class of metric foliations/submetries are given by isometric group actions on metric spaces. 
 
Let $\sfG$ be a group and denote $\sfG \times \sfY \mapsto \sfY$ by $ (g,y) \mapsto gy$ an action by isometry of the group $\sfG$ on the metric space $(\sfY, \sfd_\sfY)$. 
 We assume that the group orbit $\sfG(y):=\{ gy : g \in \sfG\}$ is closed.  We summarize the above requirements into the following condition.
 \begin{condition}\label{condG}\
 
 \begin{enumerate}
 \item $(gh)y =g(hy)$ for every $y \in \sfY$ and $g,h \in \sfG$;
 \item $e y = y$ for every $y \in \sfY$ and where $e$ is the unit element of the group $\sfG$;
 \item for every $g \in \sfG$, the map $\tau_g : \sfY \mapsto \sfY$ by $y \mapsto \tau_g(y):= gy$ is an isometry
 \begin{align*}
\sfd_\sfY \big(\tau_g(y_1), \tau_g(y_2)\big) = \sfd_\sfY(y_1,y_2).
\end{align*}
\item for every $y \in \sfY$, the orbit $\sfG(y)$ is a closed subset of $\sfY$.
 \end{enumerate}
 \end{condition}
Being in the same orbit defines an equivalence relation $\sim$.
We define $\sfX:=\sfY/\sim := \sfY / \sfG$ and   
\begin{align*}
\sfd_\sfX(x_1,x_2) & :=\inf\{ \sum_{i=1}^k \sfd_\sfY(p_i,q_i) : \forall p_i, q_i \in \sfY \\& \quad \qquad \text{ such that } p_1 \in x_1, q_k\in x_2, q_i \in \sfG(p_{i+1}) \text{ and } k \in \N\},
\end{align*}
and denote $\sfp: \sfY \mapsto \sfX$ the quotient projection. Then the following holds.
\begin{lemma}\label{App:GsubMe}
The $(\sfX, \sfd_\sfX)$ is a metric quotient space, $\sfp$ is a submetry, and  
the $\sfY$ is metrically foliated by
\begin{align*}
\sfY= \bigsqcup_{x \in \sfX} \mathcal F_{x}, \quad \mathcal F_{x} := \sfp^{-1}(x).
\end{align*}
If $(\sfY, \sfd_\sfY)$ is complete (respectively, length space), then $(\sfX, \sfd_\sfX)$ is complete (respectively, length space).
\end{lemma}
We mention that, while quotients by isometric group actions give metric foliation, the concept of metric foliation can be more general than that.

\section{Variational formulae for an effective Hamiltonian $\bar{\sfH}(P)$}\label{varbarH}
Our approach to hydrodynamic limit relies upon equation \eqref{cell}. To recapitulate, 
let $\sfH: \R^d \times \R^d \mapsto \R$ and write 
\begin{align*}
 \sfH^P(q,p):= \sfH(q,P+p), \quad \forall (q,p) \in \R^d \times \R^d, P \in \R^d,
\end{align*}
we are concerned with solution $(\varphi, c)$  to the following (cell) PDE problem in the viscosity solution sense  
\begin{align}\label{viscCell}
 \sfH^P(q, \nabla_q \varphi \big) = c, \quad \forall q \in \R^d,
\end{align} 
where $\varphi:= \varphi(q)$ is a function and $c$ is a finite constant. We call $c:=c_P:= \bar{\sfH}(P)$ the one particle level effective Hamiltonian. In this section, we presents its variational representations and a few regularity estimates for $\bar{\sfH}(P)$ as a function.  To simply presentation and highlight our main concern about hydrodynamic limits in this paper, we only work under the assumption that the $\sfH$ has a periodic structure in $q$ (see Condition~\ref{PerCND}). 
With exception of Section~\ref{supinfrep}, results in this appendix can be found in exiting literature on nonlinear-homogenization and weak KAM theory. For references, see Lions, Papanicolaou and Varadhan~\cite{LPV87}, Fathi~\cites{Fa97a, Fa97b, Fa98a, Fa98b}, \cite{FathiBook},  E~\cites{E91,E99}, unpublished works of Ma\~n\'e (see Contreras-Iturriaga-Paternain-Paternain~\cite{CIPP98} for summary and references), as well as Evans and Gomez~\cites{EG01, EG02a,EG02b}.  For various weak KAM results without periodic (or more generally without compact state state space) assumption, we mention Ishii~\cite{Ishii08}, Barles and Roquejoffre~\cite{BR06}, and Ishii and Siconolfi~\cite{IshiiSi20}. At least one approach to extend our hydrodynamic limit problem to such setting seems possible. It involves additional technical steps by introducing space of probability measures for $q$-variable with a weakened topology. We don't pursue it in this paper. Finally, there is an interesting parallel  between results here and those arising from homogenization and averaging on large deviation of Markov processes. See Chapters 11, 12 and Appendix B of Feng and Kurtz~\cite{FK06}. This should be not be a surprise, since these seemingly different topics are indeed identical in nature once formulated using two-scale Hamiltonian convergence.

As in \eqref{sfLdef}, we define $\sfL$ the Legendre transform of $\sfH$ in the $p$-variable. We also define $\bar{\sfL}$ according to \eqref{effsfLdef} and introduce its Legendre transform
\begin{align}\label{defBsfH}
  \bar{\sfH}(P):= \sup_{v \in \R^d} \big( vP - \bar{\sfL}(v) \big).
\end{align} 
The main purpose of this part of the Appendix is to establish the following.   
\begin{proposition}\label{avgHrep}
Assume that Condition~\ref{PerCND} holds. Then there is a unique $c:=c_P \in \R$ such that \eqref{viscCell} admits a viscosity solution $\varphi \in C_\per(\R^d) \cap \Lip(\R^d)$ in the sense of Definition~\ref{sfHvisc}. Furthermore,  
\begin{align*}
 c_P =\inf_{\varphi \in C^\infty_{\per}(\R^d)} \sup_{q \in \R^d} \sfH^P(q, \nabla_q \varphi)   
  = \sup_{\varphi \in C^\infty_{\per}(\R^d)} \inf_{q \in \R^d} \sfH^P(q,  \nabla_q \varphi) 
 = \bar{\sfH}(P).
\end{align*}
 \end{proposition}

 \subsection{Definition of viscosity solution in current context}
  
For $u \in \USC(\R^d)$ and $v \in \LSC(\R^d)$, we define
\begin{align*}
D^+ u & := \Big\{(q,p)  : p = \nabla \phi(q), \exists (q, \phi)  \in \R^d \times C^1(\R^d), 
  \text{ s.t.} (u - \phi)(q) = \sup_{\R^d} (u -\phi) \Big\},   \\
D^- v & := \Big\{ (q,p) : p = \nabla \phi(q), \exists (q, \phi) \in \R^d \times C^1(\R^d), 
 \text{ s.t.} (\phi - v)(q) = \sup_{\R^d} (\phi -v) \Big\}.
\end{align*}

\begin{definition}[Viscosity solution]\label{sfHvisc}
We say that $u \in \USC(\R^d)$ is a viscosity sub-solution to \eqref{viscCell} (formally written $\sfH(q,  \nabla_q u) \leq c$),  if it holds that
\begin{align*}
 \sfH(q,   p) \leq c, \quad \forall (q,p) \in D^+u. 
\end{align*}
Similarly, we say that $v \in \LSC(\R^d)$ is a viscosity super-solution to \eqref{viscCell} (formally written $\sfH(q,   \nabla_q v) \geq c$),  if it holds that 
\begin{align*}
 \sfH(q,  p) \geq c, \quad \forall (q,p) \in D^-v. 
\end{align*}
If a function is both a sub-solution as well as super-solution, then it is called a solution.
\end{definition}
In the context of equation \eqref{viscCell}, there are a number of equivalent definitions of viscosity solution, we will use them interchangeably without further mentioning. For their relations and properties, see expository text such as Crandall, Ishii and Lions~\cite{CIL92}, Bardi and Capuzzo-Dolcetta~\cite{BC97}, Cannarsa and Sinestrari~\cite{CS04}. In particular, we recall that locally Lipschitz viscosity solution are almost everywhere solutions when the gradient is interpreted in the sense of Rademacher theorem (e.g.  Proposition 1.9 of~\cite{BC97}). 
\subsection{A few concepts in Lagrangian dynamic of Hamiltonian systems}
\label{App:wKAMCnp}
Let
\begin{align*}
 \sfL^P(q, \xi):= \sup_{p \in \R^d} \big( p \cdot \xi - \sfH^P(q,p)\big) = \sfL(q, \xi) - P \cdot \xi.
\end{align*}
We define a two-fixed-time-point action by
\begin{align*}
 A^P_T[q^\prime, q] := \inf\big\{ \int_0^T \sfL^P\big(\zeta(s), \dot{\zeta}(s)\big) ds : \zeta \in AC ([0,T]; \R^d) \text{ with } \zeta(0) =q^\prime, \zeta(T) =q \big\},
\end{align*}
where the AC stands for absolute continuous curves. Let 
$c^+_P \in \R$ be the largest constant that admits a viscosity sub-solution to \eqref{viscCell}. We define critical Ma\~n\'e potential as
\begin{align*}
\sfd_{\sfH^P}(q^\prime, q):= \inf \big\{ A^P_t[q^\prime, q] + c_P^+ t : t >0 \big\};
\end{align*}
and projected Aubry set 
\begin{align}\label{App:Aub}
 {\mathbb A}_{\sfH^P} := \big\{ q^\prime \in \R^d : q \mapsto \sfd_{\sfH^P}(q^\prime, q) \text{ is a viscosity solution to }   \eqref{viscCell} \big\},
\end{align}
and Peierls' barrier
\begin{align}\label{App:Pei}
 \sfP_{\sfH^P}(q^\prime, q):= \inf_{q^{\prime \prime} \in {\mathbb A}_{\sfH^P}} 
     \big\{  \sfd_{\sfH^P}(q^\prime, q^{\prime \prime}) +  \sfd_{\sfH^P}(q^{\prime\prime}, q) \big\}, 
     \quad \forall q, q^\prime \in \R^d.
\end{align}
Recall the notion of closed probability measure in Definition~\ref{Sec1:closedM}.
We define set of Mather measures (where the $c^+_P$ is defined a few lines below),
\begin{align}\label{App:Mather}
  {\mathscr M}_{\sfH^P}:= \Big\{ \mu:= \mu(dq, d \xi) \in {\mathcal P}(\R^{2d}):  \langle \mu, \sfL^P \rangle + c^+_P =0,   \mu  \text{ is closed} \Big\}.
\end{align} 
By set of projected Mather measures, we mean
\begin{align*}
 {\mathcal M}_{\sfH^P} := \Big\{ \sfm:= \pi^1_\# \mu : \mu \in {\mathscr M}_{\sfH^P}  \Big\}.
\end{align*}
We also define Mather set 
\begin{align*}
   M_{\sfH^P}:= \bigcup_{\mu \in {\mathscr M}_{\sfH^P}} \supp[ \mu],
\end{align*}
and projected Mather set 
\begin{align*}
  {\mathbb M}_{\sfH^P}:= \bigcup_{\sigma \in {\mathcal M}_{\sfH^P}} \supp[ \sigma].
\end{align*}

\begin{lemma}\label{App:wKAMSets}
Suppose that Condition~\ref{PerCND} holds and $P \in \R^d$. Then 
 \begin{enumerate}
 \item  ${\mathbb A}_{\sfH^P}$ is non-empty,
 \item   ${\mathscr M}_{\sfH^P}$ is non-empty,
 \item ${\mathbb M}_{\sfH^P} \subset {\mathbb A}_{\sfH^P}$.
 \end{enumerate}
\end{lemma}
\begin{proof}
See Proposition 3.6 of \cite{DFIZ16}, and Theorem 5.2.8 of \cite{FathiBook}.
\end{proof}

\subsection{Variational representations, the $\inf_\varphi \sup_q$ case}
Let
\begin{align*}
c^+  := c^+_P&:= \inf \Big\{ a \in \R : \exists  \phi \in C_{\per}(\R^d), \sfH^P\big(q,  \nabla_q \phi \big) \leq a \text{ in viscosity sense} \Big\} \\
 &=  \inf_{u \in  \USC_{\per}(\R^d)} \sup_{(q,p) \in D^+ u } \sfH^P(q,   p), \\
c^- := c^-_P&:= \sup\Big\{ a \in \R : \exists \phi \in C_{\per} (\R^d), \sfH^P\big(q,  \nabla_q \phi\big) \geq a \text{ in viscosity sense}\Big\} \\
& = \sup_{v \in \LSC_{\per}(\R^d)} \inf_{(q,p) \in D^-v}\sfH^P(q,p).
\end{align*}
Suppose that there is a viscosity solution $(\varphi, c) \in C_{\per}(\R^d) \times \R$ of \eqref{viscCell}, then $c^+ \leq c \leq c^-$ by definition. 
 
 \begin{lemma}\label{Cell:exist}
Under Condition~\ref{PerCND}, for each $P \in \R^d$, there exists a viscosity solution $(c, \varphi) :=(c_P, \varphi_P)$
for \eqref{viscCell} with $c =c^+ \in \R$ and $ \varphi \in \Lip_\per(\R^d)$. In particular, 
the following is a special solution  
\begin{align}\label{DFIZSol}
\varphi:= \varphi_P(q):= \min \big\{ \int_{\R^d} \sfP_{\sfH^P}(q^\prime, q) \sfm(dq^\prime) : \sfm \in {\mathcal M}_{\sfH^P} \big\},
\end{align}
where the $ \sfP_{\sfH^P}$ is Peierls' barrier and ${\mathcal M}_{\sfH^P}$ is the set of projected Mather measures for Hamiltonian $\sfH^P$.
\end{lemma}
\begin{proof} Existence of a solution $(\varphi,c)$ has been constructed by Lions, Papanicolaou and Varadhan\cite{LPV87}.  The particular solution $\varphi$  in \eqref{DFIZSol} is constructed as a limit problem in Theorem 4.3 of Davini, Fathi, Iturriaga and Zavidovique~\cite{DFIZ16}.   
\end{proof}

\begin{lemma}\label{Cell-Ishii}
We assume Condition~\ref{PerCND}. Then, for each $P \in \R^d$,
\begin{enumerate}
\item the constant $c$ is unique in the sense that if $(c^\prime, \varphi^\prime) \in \R \times C_{\per}(\R^d)$ is another solution, then $c=c^\prime$.
\item  the following holds
\begin{align}\label{cpmeq}
c= c^+ = c^- = \inf_{\varphi \in  \Lip_\per} \sup_{(q,p) \in D^+ \varphi } \sfH^P(q, p) =
  \sup_{\varphi \in \Lip_\per} \inf_{(q,p) \in D^-\varphi}\sfH^P(q, p).
\end{align}
\item it also holds that
\begin{align}\label{cstar}
c^+= c^*:= c^*_P:=  \inf_{\varphi \in  C^\infty_\per(\R^d)} \sup_{q \in \R^d} \sfH^P \big(q,   \nabla_q \varphi \big).
\end{align}
\end{enumerate}
\end{lemma}
\begin{proof}
Since a Lipschitz viscosity solution for \eqref{viscCell} exists,
\begin{align*}
 c^+\leq \inf_{\varphi \in  \Lip_{\per}} \sup_{(q,p) \in D^+ \varphi } \sfH^P(q, p) \leq c \leq 
  \sup_{\varphi \in \Lip_{\per}} \inf_{(q,p) \in D^-\varphi}\sfH^P(q,p) \leq c^-.
\end{align*}
The reverse inequality  $c^- \leq c^+$ and uniqueness of $c$  follows from well-known comparison arguments for ergodic type Hamilton-Jacobi equation first appeared in \cite{LPV87}.  See also comparison principle Theorem 8.2.4 of Fathi~\cite{FathiBook}.
\eqref{cstar} is a well-known result in the weak KAM literature and can be found in, for instance, Theorem 2.5 (proof follows from Proposition 3.3) of Nakayasu~\cite{Naka19}.
\end{proof}
   
The $c_P$ also admits another variational representation from a Lagrangian perspective. Instead of studying minimal orbits of Hamiltonian systems,  Mather~\cite{Mather91} focused on occupation measures associated with these orbits. He gave a minimizing invariant measure interpretation of the variational constants $c^+$ that we studied earlier. See also Chapter 3 of Ma\~n\'e~\cite{Mane91} and Evans and Gomes~\cite{EG02b}. In control theory literature, the idea of using measure-based linear programming to study trajectory-based optimal controls has an  even earlier history. See Manne~\cite{Man60}, Vinter and Lewis~\cites{VL78a, VL78b}, Fleming and Vermes~\cite{FV89} and Stockbridge~\cite{Stock90}, etc. 
 
\begin{lemma}
Under Condition~\ref{PerCND},  $c^*_P= \bar{\sfH}(P)$.
 \end{lemma}
\begin{proof}
First of all, the following sequence of relations hold by definition
 \begin{align*}
 c^*_P & = \inf_{\varphi \in C^\infty_\per(\R^d)} \sup_{q \in \R^d} \sfH^P(q, \nabla_q \varphi) \\
 & = \inf_{\varphi \in C^\infty_\per(\R^d)} \sup_{(q,\xi) \in \R^{2d}} 
  \Big( \xi \cdot (P+ \nabla_q \varphi) - \sfL(q, \xi) \Big) \\
 & = \inf_{\varphi \in C^\infty_\per(\R^d)} \sup_{\mu \in \mathcal P(\R^{2d})} 
 \int_{\R^{2d}} \Big( \xi \cdot (P+ \nabla_q \varphi) - \sfL(q, \xi) \Big) \mu(dq, d\xi) \\
 & \geq  \sup_{\mu \in \mathcal P(\R^{2d})} \inf_{\varphi \in C^\infty_\per(\R^d)}
 \int_{\R^{2d}} \Big( \xi \cdot (P+ \nabla_q \varphi) - \sfL(q, \xi) \Big) \mu(dq, d\xi)   \\
 & = \sup_{\mu \in \mathcal P(\R^{2d})} \Big\{
 \int_{\R^{2d}} \big( \xi P - \sfL(q, \xi) \big) \mu(dq, d\xi)  : 
 \int_{\R^{2d}} (\xi \nabla_q \varphi) \mu(dq, d\xi) =0, \forall \varphi \in C^\infty_\per(\R^d) \Big\} \\
 & = \sup_{\mu   \in \mathcal P(\R^{2d})} \Big\{
  v P    -\int_{\R^{2d}}   \sfL(q, \xi) \mu(dq, d\xi)  : \text{ where }  v= \int_{(q, \xi) \in \R^{2d}} \xi \mu(dq, d\xi),
   \mu \text{ is closed} \Big\}\\
 &= \sup_{v \in \R^d} \Big\{ vP - \bar{\sfL}(v) \Big\}  .
\end{align*}
Therefore, our conclusion follows if the inequality above is an equality by a minimax theorem type argument. In the following, we present a more streamlined alternative proof.  

Secondly, since $\mathscr M_{\sfH^P}$ is non-empty, for each $\mu_0 \in \mathscr M_{\sfH^P}$,    
\begin{align*}
  \sup_{\mu \in \mathcal P(\R^{2d})} \Big\{
 \int_{\R^{2d}} \big( \xi P - \sfL(q, \xi) \big) \mu(dq, d\xi)  : &
 \int_{\R^{2d}} (\xi \nabla_q \varphi) \mu(dq, d\xi) =0, \forall \varphi \in C^\infty_c(\R^d) \Big\} \\
 & \geq    \langle - \sfL^P, \mu_0 \rangle = c^+_P. 
\end{align*}

Since $c^*_P = c^+_P$ by Lemma~\ref{cstar}, we conclude.
 \end{proof}

\subsection{Variational representations,  the $\sup_\varphi \inf_q$ case}\label{supinfrep}
In the presence of convexity of $p \mapsto \sfH(q,p)$ and Lipschitz regularity on a viscosity solution, the sub-solution property is approximately stable under the usual mollification by convolution technique. This is how \eqref{cstar} is proved. However,  such approximation procedure becomes unstable for the super-solution property. In Lasry and Lions~\cite{LL86}, the authors introduced a nonlinear Moreau-Yosida type regularization procedure, for approximating a continuous function in $\R^d$ by $C^{1,1}_{\loc}(\R^d)$ functions. Both sub- and super-solution viscosity solution properties are approximately stable. Next, we adapt such technique to our context for yet another variational representation of the critical constant $c_P$.

Let $w \in C(\R^d)$ be such that
\begin{align}\label{wgro}
 |w(q)| \leq C_w(1+ |q|^2), \quad \exists C_w \in \R_+.
\end{align}
For $0<\epsilon < \frac{1}{2C_w}$, we introduce non-linear mollifications  
\begin{align}
w_\epsilon(q^\prime) &:= \inf_{q^{\prime \prime} \in \R^d}  \Big( w(q^{\prime \prime}) 
+\frac{1}{2 \epsilon} |q^\prime - q^{\prime \prime}|^2\Big), \label{ueps} \\
v_\epsilon(q)&:= \sup_{q^\prime\in \R^d} \Big(  w_\epsilon(q^\prime)  - \frac{1}{\epsilon} |q-q^\prime|^2\Big). \label{veps}
\end{align}
We note that if $w$ is periodic, then $w_\epsilon$ is periodic:
\begin{align*}
 w_\epsilon(q^\prime +k)  & = \inf_{q^{\prime \prime}}\Big( w(q^{\prime \prime}) 
+\frac{1}{2 \epsilon} |q^\prime +k - q^{\prime \prime}|^2\Big)  \\
& = \inf_{q^{\prime \prime}}\Big( w(q^{\prime \prime} -k) 
+\frac{1}{2 \epsilon} |q^\prime - (q^{\prime \prime}-k)|^2\Big)  
= w_\epsilon(q^\prime).
\end{align*}
In the same way, $v_\epsilon$ becomes periodic too.
\begin{lemma}\label{MY1}
We have
\begin{enumerate}
\item $w_\epsilon \leq w$,  $w_\epsilon(q^\prime) - \frac{1}{2\epsilon}|q^\prime|^2$ is concave, and 
\begin{align}\label{ugrth}
 w_\epsilon(q^\prime) \leq \frac{1}{2\epsilon} |q^\prime|^2 + w(0), \quad \forall q \in \R^d.
\end{align}
\item For every $q_0^\prime \in \R^d$,  minimizer $q_0^{\prime \prime} \in \R^d$ in 
the definition of $w_\epsilon(q_0^\prime)$ exists. Moreover, any such minimizer satisfies
\begin{align}\label{argw}
 |q_0^\prime - q_0^{\prime \prime}|^2 \leq 2 \epsilon \big(w(q_0^\prime) - w(q_0^{\prime \prime}) \big),
\end{align}
and
\begin{align}\label{DLwest}
\frac{q_0^\prime - q_0^{\prime \prime}}{\epsilon}  \in D_{q_0^{\prime \prime}}^- w.
\end{align}
\item If $(q_0^\prime, p_0^\prime) \in D^- w_\epsilon$, 
then $\nabla_{q_0^\prime}w_\epsilon$ exists in the classical sense,  with
\begin{align}\label{Dmwe}
  \frac{1}{\epsilon} (q_0^\prime - q_0^{\prime \prime})  = p_0^\prime = \nabla_{q_0^\prime} w_\epsilon \in  D^-_{q_0^{\prime \prime}} w .
\end{align}
\end{enumerate}
\end{lemma}
\begin{proof}
The concavity of $q^\prime \mapsto w_\epsilon(q^\prime) - \frac{1}{2 \epsilon}|q^\prime|^2$ follows from representation
\begin{align*}
 w_\epsilon (q^\prime) - \frac{1}{2\epsilon} |q^\prime|^2
  = - \sup_{q^{\prime \prime} \in \R^d} \Big( \frac{1}{\epsilon} q^\prime q^{\prime \prime} 
  - \big( w(q^{\prime \prime}) + \frac{1}{2 \epsilon} |q^{\prime \prime}|^2 \big)   \Big).
\end{align*}

The existence of minimizer $q_0^{\prime \prime}$, in the definition of $w_\epsilon(q_0^\prime)$, 
follows from \eqref{wgro}. 
The minimizing property in the definition of $w_\epsilon(q_0^\prime)$ gives
\begin{align*}
 w(q_0^{\prime \prime}) + \frac{1}{2\epsilon} |q_0^\prime - q_0^{\prime \prime}|^2
  \leq w(q^{\prime \prime}) + \frac{1}{2\epsilon} |q_0^\prime - q^{\prime \prime}|^2, \quad \forall q^{\prime \prime} \in \R^d.
\end{align*}
Consequently, 
\begin{align*}
 \frac{q_0^{\prime} - q_0^{\prime \prime}}{\epsilon} \in D_{q_0^{\prime \prime}}^- w.
\end{align*}

We already proved that $w_\epsilon \in \SCC_{\loc}(\R^d)$. Hence for $(q_0^\prime, p_0^\prime) \in D^- w_\epsilon$, $\nabla_{q_0^\prime} w_\epsilon$ exists (e.g. part b of Proposition 4.7 in \cite{BC97}). Therefore, there exists a $\varphi \in C^1(\R^d)$ with 
 $p_0^\prime=\nabla_{q_0^\prime} w_\epsilon=\nabla_{q_0^\prime} \varphi$ such that $w_\epsilon - \varphi$
  attains a local minimum at $q_0^\prime$. That is, by definitions of the $w_\epsilon$ and the $q_0^{\prime \prime}$,
  \begin{align*}
 w(q_0^{\prime \prime}) + \frac{1}{2\epsilon}|q_0^\prime -q_0^{\prime \prime}|^2 - \varphi(q_0^\prime) 
& =   w_\epsilon(q^\prime_0) - \varphi(q^\prime_0) \\
& \leq w_\epsilon(q^\prime) - \varphi(q^\prime)  
 \leq w(q^{\prime \prime}) + \frac{1}{2\epsilon}|q^\prime -q^{\prime \prime}|^2 - \varphi(q^\prime), 
  \quad \forall q^\prime , q^{\prime \prime} \in \R^d.
\end{align*}
Take $q^{\prime \prime}= q_0^{\prime \prime}$, then
\begin{align*}
 \frac{1}{\epsilon} (q_0^\prime - q_0^{\prime \prime}) = \nabla_{q_0^\prime} \varphi.
\end{align*}
Summarizing all the above, we arrive at \eqref{Dmwe}. 
 \end{proof}

\begin{lemma}\label{MY2}
We have
\begin{enumerate}
\item $v_\epsilon \geq w_\epsilon$ and $v_\epsilon(q) + \frac{1}{\epsilon} |q|^2$ is convex. 
\item  $v_\epsilon(q) - \frac{1}{t} |q|^2$ is concave for every $0<t<\epsilon$.
\item $v_\epsilon \in C^{1,1}_{\loc}(\R^d)$. 
\item For each $q_0 \in \R^d$, there exists a unique maximizer $q_0^\prime \in \R^d$ 
in the variational definition of $v_\epsilon(q_0)$. Moreover,
\begin{align}\label{mxrDv}
    \nabla_{q_0} v_\epsilon = \frac{2(q_0^\prime - q_0)}{\epsilon},
\end{align}
and for any minimizer $\tilde{q}_0$ in the definition of $w_\epsilon(q_0)$,
\begin{align}\label{argwe}
|q_0^\prime-q_0|^2 & \leq 8 \epsilon \big( w(q_0) - w(\tilde{q}_0)\big).
\end{align}
\item for every $\tilde{q}_0 \in \R^d$ which is a minimizer in the definition of $w_\epsilon(q_0)$, we have
\begin{align}\label{DvDw}
 |\nabla_{q_0} v_\epsilon| \leq 4 \inf  \big\{ | \tilde{p}_0 | :  
  \tilde{p}_0 \in D_{\tilde{q}_0}^-w   \big\}.
\end{align}
\end{enumerate}
\end{lemma}
 \begin{proof}
Similar to the arguments verifying concavity of the 
$q^\prime \mapsto w_\epsilon(q^\prime) - \frac{1}{2\epsilon} |q^\prime|^2$ in Lemma~\ref{MY1}, 
we have that $v_\epsilon(q) + \frac{1}{\epsilon} |q|^2$ is convex.  

We note that
\begin{align*}
 v_\epsilon(q) -\frac{|q|^2}{t} = \sup_{q^\prime \in \R^d} 
  \big( w_\epsilon(q^\prime) -\frac{|q^\prime|^2}{2\epsilon} 
    + \frac{|q^\prime|^2}{2\epsilon} -\frac{|q-q^\prime|^2}{\epsilon} -\frac{|q|^2}{t} \big)  
    =: \sup_{q^\prime \in \R^d} G(q,q^\prime).
\end{align*}
Since $w_\epsilon(q^\prime) -\frac{1}{2\epsilon}|q^\prime|^2$ is concave, the above $G: \R^d \times \R^d \mapsto \R$ is concave. By the lemma on page 265 of Lasry and Lions~\cite{LL86}, we conclude  $v_\epsilon(q) - \frac{1}{t}|q|^2$  is concave for every $0<t<\epsilon$.

Hence $v_\epsilon \in C^{1,1}_{\loc}(\R^d)$ by Lemma 3.3.8 of Cannarsa and Sinestrari~\cite{CS04}.
 
Next, we take an arbitrary but fixed $q_0 \in \R^d$. 
Because of estimate \eqref{ugrth},  there exists maximizer $q_0^\prime\in \R^d$ 
of the $v_\epsilon(q_0)$. The uniqueness follows from strict convexity of 
\begin{align*}
 q^\prime \mapsto w_\epsilon(q^\prime) -\frac{|q-q^\prime|^2}{\epsilon}
  = \big( w_\epsilon(q^\prime) -\frac{|q^\prime|^2}{2\epsilon}\big)
   + \big( \frac{|q^\prime|^2}{2\epsilon} -\frac{|q-q^\prime|^2}{\epsilon}\big).
\end{align*} 
Since $v_\epsilon \in C^1(\R^d)$, we can find a $\varphi \in C^1(\R^d)$ such that the $q_0$ is a local maximum for $v_\epsilon - \varphi$ with $\nabla_{q_0} v_\epsilon = \nabla_{q_0} \varphi$. Therefore, there exists a neighborhood of the $q_0$, for every $q$ in this neighborhood and for every $q^\prime \in \R^d$, we have
\begin{align}\label{wphimax}
 w_{\epsilon}(q_0^\prime) - \frac{1}{\epsilon} |q_0-q_0^\prime|^2 - \varphi(q_0) 
 \geq w_{\epsilon}(q^\prime) - \frac{1}{\epsilon} |q-q^\prime|^2 - \varphi(q).
\end{align}
Taking $q^\prime =q_0^\prime$, we obtain $ \frac{2}{\epsilon}(q_0^\prime - q_0) = \nabla_{q_0} \varphi$,
 giving \eqref{mxrDv}. 

Next, we verify \eqref{argwe}. From the maximizing property,  
\begin{align}\label{q0dis}
 |q_0^\prime-q_0|^2 \leq \epsilon \big( w_\epsilon(q_0^\prime) - w_\epsilon(q_0)\big).
 \end{align}
We further estimate right hand side of the above inequality through concavity property (Lemma~\ref{MY1}) of 
$g(q):=g_\epsilon(q):= w_\epsilon(q) - \frac{1}{2\epsilon}|q|^2$. Let $\tilde{q}_0$ be a minimizer in the definition of $w_\epsilon(q_0)$. First,  by concavity,
\begin{align*}
 g(q_0^\prime) - g(q_0) \leq \frac{g(q_t^\prime) - g(q_0)}{t}, \quad \forall q_t^\prime:=(1-t) q_0 + t q_0^\prime, t \in (0,1).
\end{align*}
Second, following the definition of $w_\epsilon$,  
\begin{align*}
w_\epsilon(q^\prime_t) - w_\epsilon(q_0) \leq w_\epsilon(\tilde{q}_0) +\frac{|\tilde{q}_0 - q_t^\prime|^2}{2\epsilon} 
 - w_\epsilon(\tilde{q}_0) -\frac{|\tilde{q}_0 - q_0|^2}{2\epsilon} = t \frac{q_0 + q_t^\prime - 2 \tilde{q}_0}{2\epsilon}  (q_0^\prime - q_0).
\end{align*}
By definition of $g$,
\begin{align*}
g(q_t^\prime) - g(q_0) \leq t   \frac{(q_0 + q_t^\prime - 2 \tilde{q}_0 - 2 q_0)}{2\epsilon}(q_0^\prime - q_0).
\end{align*}
Therefore
\begin{align*}
  g(q_0^\prime) - g(q_0) \leq \limsup_{t \to 0^+}\frac{g(q_t^\prime) - g(q_0)}{t} \leq \frac{-\tilde{q}_0}{\epsilon}(q_0^\prime - q_0),
\end{align*}
implying
\begin{align*}
 w_\epsilon (q_0^\prime) - w_\epsilon (q_0) \leq \frac{(q_0-\tilde{q}_0)}{\epsilon}(q_0^\prime - q_0) +\frac{|q_0^\prime -q_0|^2}{2\epsilon}.
\end{align*}
Combined with \eqref{q0dis}, we arrive at 
\begin{align}\label{qq0est}
  |q_0^\prime-q_0| \leq 2 |q_0-\tilde{q}_0|.
\end{align}
In the current context,  \eqref{argw} in Lemma~\ref{MY1} becomes $|q_0 -\tilde{q}_0|^2 \leq 2 \epsilon \big( w(q_0) - w(\tilde{q}_0)\big)$. Combined with \eqref{qq0est}, we have \eqref{argwe}.
 
Finally, \eqref{DvDw} follows from combing \eqref{mxrDv} with \eqref{qq0est} and \eqref{DLwest}. 
 
 \end{proof}
 
 \begin{lemma}
Let $\sfH \in C(\R^d \times \R^d)$. Suppose that  $w \in \Lip_{\loc}(\R^d)$
is a viscosity super-solution to  
\begin{align}\label{sfHphi0}
\sfH^P(q,  \nabla_q w) \geq c.
 \end{align}
 Then for each $q_0 \in \R^d$ and the unique (Lemma~\ref{MY2}) maximizer $q_0^\prime \in \R^d$ in the variational definition of $v_\epsilon(q_0)$, there exists a minimizer $q_0^{\prime \prime} \in \R^d$ in the variational definition of $w_\epsilon(q_0^\prime)$, such that
\begin{align}\label{sfHcsig}
\sfH^P(q_0^{\prime \prime},  \nabla_{q_0} v_\epsilon)  \geq  c.   
\end{align} 
 \end{lemma}
 \begin{proof}
 Following notations in the proof of Lemma~\ref{MY2}, we take $q=q_0$ in inequality \eqref{wphimax}. Then, for every $q^\prime \in \R^d$ and each $q_{q^\prime}^{\prime \prime}$ 
which is a minimizer in the variational definition of the $w_\epsilon(q_0^\prime)$,
we have 
\begin{align*}
 \frac{1}{\epsilon}|q^\prime - q_0|^2 -\frac{1}{\epsilon} |q_0^\prime - q_0|^2 
 & \geq w_\epsilon(q^\prime) - w_\epsilon(q_0^\prime) \\
 & = \sup_{q^{\prime \prime \prime}\in \R^d} \inf_{q^{\prime \prime}\in \R^d} \big( w(q^{\prime \prime}) 
 + \frac{|q^\prime - q^{\prime \prime} |^2}{2 \epsilon} - w(q^{\prime \prime \prime}) 
 - \frac{|q_0^\prime - q^{\prime \prime \prime} |^2}{2 \epsilon}\big) \\
 & \geq \frac{|q^\prime -q_{q^\prime}^{\prime \prime}|^2}{2 \epsilon} - \frac{|q_0^{\prime} - q_{q^\prime}^{\prime \prime}|^2}{2 \epsilon} 
  \geq (q^\prime - q_0^\prime) \frac{q_0^{\prime} -q^{\prime \prime}_{q^\prime}}{\epsilon}.
\end{align*}
By \eqref{DLwest},
 \begin{align}\label{plotspri}
 p^{\prime \prime}_{q_{q^\prime}^{\prime \prime}}:=
 \frac{q_0^{\prime} -q^{\prime \prime}_{q^\prime}}{\epsilon} 
 \in  D^-_{q_{q^\prime}^{\prime \prime}} w.
 \end{align} 
Combined with \eqref{mxrDv}, we conclude
\begin{align}\label{QcpDw}
\frac{|q^\prime - q_0^\prime|^2}{\epsilon} + (q^\prime - q_0^\prime) \nabla_{q_0} v_\epsilon \geq  (q^\prime - q_0^\prime) p^{\prime \prime}_{q_{q^\prime}},  \qquad \forall q^\prime \in \R^d.
\end{align}
Let $q^\prime \to q_0^\prime$. In view of the growth estimate \eqref{wgro}, we have at least along subsequence that $q_{q^\prime}^{\prime \prime} \to q_0^{\prime \prime}$ where the limiting point $q_0^{\prime \prime}$ is a minimizer of $w_\epsilon(q_0^\prime)$. 
Therefore
\begin{align*}
 \lim_{q^\prime \to q_0^\prime} p^{\prime \prime}_{q_{q^\prime}^{\prime \prime}} 
  = \frac{q_0^{\prime} -q^{\prime \prime}_0}{\epsilon} = : p_0^\prime.
\end{align*}
We claim that  
\begin{align*}
p_0^\prime  = \nabla_{q_0} v_\epsilon.
\end{align*}

To verify the claim, we take a particular choice of $q^\prime = q_0^\prime + tn$ in \eqref{QcpDw}, where $n \in \R^d$ with $|n|=1$ and $t\in (0,1)$ are arbitrary. Then 
\begin{align*}
 | p_0^\prime - \nabla_{q_0} v_\epsilon|
 = \sup_{\substack{n \in \R^d,\\ |n|=1}} n (p_0^\prime - \nabla_{q_0} v_\epsilon) 
 & = \sup_{\substack{n \in \R^d,\\ |n|=1}} n (p_0^\prime - p^{\prime \prime}_{q_{q^\prime}^{\prime \prime}}
    + p^{\prime \prime}_{q_{q^\prime}^{\prime \prime}} - \nabla_{q_0} v_\epsilon) \\
 & \leq \sup_{\substack{q^\prime=tn,\\ n \in \R^d, |n|=1}}
  | p_0^\prime -p^{\prime \prime}_{q_{q^\prime}^{\prime \prime}} |
  +  \frac{t}{\epsilon},
\end{align*}
where the last inequality follows from \eqref{QcpDw}. Taking $\lim_{t \to 0^+} \lim_{q^\prime \to q^\prime_0}$ verifies the claim. 

 Finally, by viscosity super-solution property for \eqref{sfHphi0} and \eqref{plotspri}, for every $q^\prime \in \R^d$,
\begin{align*}
 \sfH^P(q^{\prime \prime}_{q^\prime}, p^{\prime \prime}_{q_{q^\prime}^{\prime \prime}})
     \geq c.
\end{align*}
Consequently, by continuity of the $(q,p) \mapsto \sfH(q,p)$,
\begin{align*}
 \sfH^P(q_0^{\prime \prime},  \nabla_{q_0} v_\epsilon) 
= \lim_{q^\prime \to q_0^\prime} \sfH^P\big( q_{q^\prime}^{\prime \prime}, p^{\prime \prime}_{q_{q^\prime}^{\prime \prime}}   ) \geq  c.
\end{align*}

\end{proof}

\begin{lemma}\label{App:smEigen}
Under Condition~\ref{PerCND},  
\begin{align*}
c_*:= c_{*,P} := \sup_{\varphi \in C^\infty_{\per}(\R^d)} \inf_{q \in \R^d} \sfH^P(q, \nabla_q \varphi) \geq c.
\end{align*}
\end{lemma}
\begin{proof}
By another approximation step, we only need to verify the above when the $C^\infty_\per(\R^d)$ is replaced by $C^{1,1}_\loc(\R^d)\cap C_\per(\R^d)$.

Take $w:= \varphi_P$ as defined in \eqref{DFIZSol}. Then $(w, c)$ solve is a viscosity solution to \eqref{viscCell} and $w \in \Lip_\loc(\R^d) \cap C_\per(\R^d)$. We mollify this $w$ through \eqref{ueps} and \eqref{veps} to arrive at the $v_\epsilon \in C_{\loc}^{1,1}(\R^d) \cap C_\per(\R^d)$ where the periodicity follows from periodicity of the $w$. Moreover, let $L$ be a Lipschitz constant for $w$, by \eqref{DvDw}, then 
$\sup_q |\nabla_q v_\epsilon| \leq 4 L$. By the periodicity assumption on $q \mapsto \sfH(q,p)$ and continuity of $\sfH$, there exists a modulus $\omega_L \in C(\R_+; \R_+)$ with $\omega_L(0)=0$, such that
\begin{align*}
|\sfH(q_0 , \nabla_{q_0} v_\epsilon) - \sfH(q_0^{\prime \prime}, \nabla_{q_0} v_\epsilon)| \leq \omega_L (|q_0 - q_0^{\prime \prime}|) 
\leq \omega_L(6 \sqrt{\epsilon} \Vert w\Vert_{L^\infty}),
\end{align*}
where the last inequality follows from \eqref{argw} and \eqref{argwe}.
 
Combine the above estimate with \eqref{sfHcsig},  the conclusion follows. 
\end{proof}

From definitions, $c_* \leq c^-$. On the other hand, combine the above result with \eqref{cpmeq}, $c_* \geq c=c^+=c^-$. Hence $c^*=c^+=c^-$.

\subsection{Some variational properties on Hamiltonians for infinite particles}
In this section, we first recall the definition of effective Hamiltonian $\bar{\sfH}:= \bar{\sfH}(P)$ in \eqref{defBsfH} and its many equivalent representations (e.g. Lemma~\ref{avgHrep}) under Condition~\ref{PerCND}. 
We also recall notation $\bar{\sfH}(x,P) := \bar{\sfH}(P) - U_0(x)$ defined in \eqref{barHxP}. 

\begin{lemma}\label{PbarH}
Under Condition~\ref{PerCND},  $\bar{\sfH}: \R^d \mapsto \R$ is convex and locally Lipschitz.
Furthermore, if we additionally assume that Condition~\ref{Sec1:Tone2} holds, then 
\begin{align*}
 -c + C^{-1} |P|^2 \leq \bar{\sfH}(P) \leq c+ C|P|^2,
\end{align*}
with $c, C$ the same constants in Condition~\ref{Sec1:Tone2}.
Also, same type estimate for $\bar{\sfL}$ holds by convexity argument.
\end{lemma}
\begin{proof}
Since $\bar{\sfH}$ is Legendre transform of the $\bar{\sfL}$, it is convex. Such $\bar{\sfH}$ is finite everywhere, because that
\begin{align*}
 -\infty < \inf_{q \in \R^d} \sfH(q,P) & \leq  \sup_\varphi \inf_{q \in \R^d} \sfH(q, P+\nabla_q \varphi)  \\
 &= \bar{\sfH}(P)   =\inf_\varphi \sup_{q \in \R^d} \sfH(q, P+\nabla_q \varphi) \leq \sup_{q \in \R^d} \sfH(q,P)<\infty.
\end{align*}
Hence it is locally Lipschitz.

With Condition~\ref{Sec1:Tone2},  
\begin{align*}
-c +C^{-1}|P|^2 \leq \inf_q \sfH(q, P) 
 \leq \bar{\sfH}(P) \leq \sup_q \sfH(q, P) \leq c+ C|P|^2. 
\end{align*}
\end{proof}

 We recall the definition of ${\mathcal F}_0$ in \eqref{defF0}.
  
   \begin{lemma}\label{infintH}
For $\bnu \in \mathcal P_2(\R^d \times \R^d)$, we have
 \begin{align*}
  & \inf_{\phi \in {\mathcal F}_0} \int_{\R^d \times \R^d}  
  \sup_{q \in \R^d}   \sfH \big(q, P + \nabla_q \phi(x,P;q)\big)  \bnu(dx,dP)  \\
  & \qquad \qquad = \sup_{\phi \in {\mathcal F}_0}  \int_{\R^d \times \R^d}  
  \inf_{q \in \R^d}   \sfH\big(q,P+ \nabla_q \phi(x,P;q)\big)   \bnu(dx,dP)   
   = \int_{\R^d \times \R^d}   \bar{\sfH}(P) \bnu(dx, dP). 
  \end{align*}
 \end{lemma}
 \begin{proof}
Part one: We establish identity
\begin{align*}
 \LHS :=&  \sup_{\phi \in {\mathcal F}_0} \int_{\R^{2d}} \inf_{q \in \R^d} 
  \sfH \big(q, P+ \nabla_q \phi(x,P;q)\big)  \bnu(dx,dP) \\
 & = \int_{\R^{2d}} \sup_{\varphi \in C^\infty_\per(\R^d)} \inf_{q \in \R^d} 
    \sfH \big(q, P + \nabla_q \varphi(q)\big)  \bnu(dx,dP) =: \RHS.
\end{align*} 
It is sufficient to verify that $\LHS \geq \RHS$. For notational convenience, we denote 
\begin{align*}
 h(x,P;\varphi):= \inf_{q \in \R^d} \sfH(q,P+\nabla_q \varphi(q)), \quad \varphi \in C^1(\R^d).
\end{align*}
Noting $\inf_{(q, P) \in \R^{2d}} \sfH(q,P) >-\infty$, we assume with no loss of generality that $h \geq 0$ in the following proof. 
By a density argument, we can find a countable set of $\{ \varphi_i \in C_\per^\infty(\R^d) : i=1,2,\ldots\}$
such that 
\begin{align*}
\RHS= \int_{\R^{2d}} \sup_{i \in \N}  h(x,P;\varphi_i)  \bnu(dx,dP).
\end{align*}
By Lemma 2.35 of Ambrosio, Fusco and Pallara~\cite{AFP00}, 
\begin{align}\label{AFPID}
\int_{\R^{2d}} \sup_{i \in \N} h(x,P; \varphi_i) \bnu(dx, dP) = 
\sup \Big( \sum_{i \in I} \int_{A_i} h(x,P; \varphi_i) \bnu(dx,dP) \Big)
\end{align}
where the supremum ranges over all finite sets $I \subset \N$ and all families
$ \{A_i : i \in I \}$ of pairwise disjoint open sets with compact closure in $\R^d \times \R^d$.
Let $\{ \alpha_i : n=1,2,\ldots\}$ be a smooth partition of unity with
\begin{align*}
\alpha_i(x,P)=
\begin{cases}
1 , & \forall (x,P) \in A_i,   \\
0,  & \forall (x,P) \in A_j, \text{ when } j \neq i. 
\end{cases}
 \end{align*}
Then letting 
$\phi(x, P; q) :=  \sum_{j \in I} \alpha_j(x,P) \varphi_j(q) \in {\mathcal F}_0$, 
we have $h(x, P; \varphi_i) =h\big(x, P; \phi(x,P; \cdot)\big)$ for $(x,P) \in A_i$.
Consequently for every finite index set $I \subset \N$, and every family $\{ A_i : i \in I\}$ of pairwise disjoint open sets with compact closure in $\R^d \times \R^d$, we have
\begin{align*}
 \sum_{i \in I} \int_{A_i} h(x, P; \varphi_i) \bnu(dx,dP)  
&  =      \sum_{i \in  I} \int_{A_i} h\big(x, P; \phi(x,P; \cdot)\big) \bnu(dx, dP) \\
 & \leq \int_{\R^{2d}}  h\big(x, P; \phi(x,P; \cdot)\big) \bnu(dx, dP) \\
 & \leq \sup_{\phi \in {\mathcal F}_0} \int_{\R^{2d}} h \big(x, P; \phi(x,P;\cdot)\big) \bnu(dx,dP) =\LHS.
\end{align*}
Combined with \eqref{AFPID}, we conclude that $\LHS \geq \RHS$.

Part two: Denoting
\begin{align*}
 g(x,P;\varphi):= \sup_{q \in \R^d} \sfH(q,P+\nabla_q \varphi), \quad \forall \varphi \in C^1(\R^d),
\end{align*}
we prove that
\begin{align*}
L:= \inf_{\phi \in {\mathcal F}_0} \int_{\R^{2d}} g\big(x, P; \phi(x,P;\cdot)\big) \bnu(dx,dP) =
 \int_{\R^{2d}} \inf_{\varphi \in C^\infty_\per(\R^d)} g(x, P; \varphi) \bnu(dx,dP)=:R.
\end{align*} 
We only need to prove $L \leq R$. It follows from Part one of the proof that 
\begin{align*}
  \inf_{\phi \in {\mathcal F}_0} \int_{\R^{2d}} (g \wedge k)\big(x, P; \phi(x,P;\cdot)\big) \bnu(dx,dP) & =
 \int_{\R^{2d}} \inf_{\varphi \in C^\infty_\per(\R^d)} (g \wedge k)(x, P; \varphi)   \bnu(dx,dP), \\
 & \leq \int_{\R^{2d}} \inf_{\varphi \in C^\infty_\per(\R^d)} g (x, P; \varphi)   \bnu(dx,dP). \quad \forall k \in \R_+.
\end{align*} 
Therefore, denoting
\begin{align*}
 F(k,\phi):=  \int_{\R^{2d}} (g \wedge k)\big(x, P; \phi(x,P;\cdot)\big) \bnu(dx,dP), \quad \forall \phi \in {\mathcal F}_0, k \in \bar{\R}_+:= \bar{\R}_+ \cup \{+\infty\},
\end{align*}
we only need to show that
\begin{align}\label{App:Sion}
  \sup_{k \in \bar{\R}_+}\inf_{\phi \in {\mathcal F}_0} F(k,\phi)
  =\inf_{\phi \in {\mathcal F}_0} \sup_{k \in \bar{\R}_+} F(k,\phi).
\end{align}
The map $\bar{\R}_+ \ni k \mapsto F$ is concave, and $\phi \mapsto F$ is convex. 
The $\bar{\R}_+$ is endowed with one-point compactification topology of $\R_+$. 
$\bar{\R}_+ \ni k \mapsto F$ is continuous by monotone convergence theorem. 
Consequently \eqref{App:Sion} follows from a version of minimax Theorem 4.2 in Sion~\cite{Sion58}.

We now conclude the lemma in view of variational representations of the $\sfH$ by Proposition~\ref{avgHrep}.
\end{proof}

For purpose of proving viscosity extension Lemmas~\ref{Sec6:bbH0Eqn} 
and \ref{Sec6:bbH1Eqn} in the main text, we need strengthened versions of the above result. 
 
 \begin{lemma}\label{App:supintH0}
 Let $\mathcal N \subset \subset {\mathcal P}_2(\R^{2d})$ be a compact subset with respect to the topology given by $2$-Wasserstein metric. Then
  \begin{align*}
    \inf_{\phi \in {\mathcal F}_0} \sup_{\bnu \in {\mathcal N}}\int_{\R^{2d}}  
  \sup_{q \in \R^d}   \sfH \big(q, P + \nabla_q \phi(x,P;q)\big)  \bnu(dx,dP)  
  = \sup_{\bnu \in {\mathcal N}} \int_{\R^{2d}}   \bar{\sfH}(P) \bnu(dx, dP). 
 \end{align*}  
\end{lemma}
\begin{proof}
First, letting
\begin{align*}
f(\bnu,\phi) := \int_{\R^d \times \R^d}  
  \sup_{q \in \R^d}   \sfH \big(q, P + \nabla_q \phi(x,P;q)\big)  \bnu(dx,dP),  
\end{align*}
then $(\bnu, \phi) \mapsto f(\bnu,\phi)$ is concave-convex-like in the sense of Sion~\cite{Sion58}, and $\bnu \mapsto f(\bnu, \phi)$ is continuous with each $\phi$ fixed. Consequently, by Theorem 4.2 in \cite{Sion58}, 
\begin{align*}
    \inf_{\phi \in {\mathcal F}_0} \sup_{\bnu \in {\mathcal N}} f(\bnu,\phi)
  =  \sup_{\bnu \in {\mathcal N}} \inf_{\phi \in {\mathcal F}_0} f(\bnu,\phi).
\end{align*}  
Second, by Lemma~\ref{infintH},
$ \inf_{\phi \in {\mathcal F}_0} f(\bnu,\phi) = \int_{\R^{2d}} \bar{\sfH}(P) \bnu(dy,dP)$.
 Hence we conclude. 
\end{proof}

\begin{lemma}\label{App:infintH1}
 Let $\mathcal N \subset \subset {\mathcal P}_2(\R^{2d})$ be a compact subset with respect to the topology given by $2$-Wasserstein metric. Then
    \begin{align*}
& \sup_{\phi:=\phi(x,P;q) \in {\mathcal F}_0}   \inf_{ \bnu \in \mathcal N} 
   \int_{\R^{2d}} \inf_{q \in \R^d}   \sfH \Big( q, P+\nabla_q \phi\big(y,P; q \big) \Big)  \bnu(dy, dP) \\
   & \qquad \qquad 
   = \inf_{ \bnu \in \mathcal N} \int_{\R^{2d}} \bar{\sfH}(P) \bnu(dy, dP).
\end{align*}
\end{lemma}
\begin{proof}
We denote the left (and right) hand side of the above identity by LHS (respectively RHS). Then by Proposition~\ref{avgHrep},  we only need to show $\LHS \geq \RHS$.

Let $\epsilon>0$ be arbitrarily given by fixed. 

First of all, by Lemma~\ref{App:smEigen}, for each $P \in \R^d$, there exists $P$-parametrized functions $\varphi_\epsilon:=\varphi_\epsilon(P; \cdot) \in C^\infty_\per(\R^d)$ such that 
\begin{align*}
\inf_{q \in \R^d} \sfH\big(q, P+ \nabla_q \varphi_\epsilon(q,P)\big) 
 \geq \bar{\sfH}(P) -\frac{\epsilon}{2}.
\end{align*}
The above implies the existence of $P$-parametrized  vector fields 
$\xi_\epsilon:=\xi_\epsilon(P;\cdot) : \R^d \mapsto \R^d$,
 which is continuous in both $P, q$ variables,  such that 
\begin{align*}
\inf_{q \in \R^d}  \Big(\big( P + \nabla_q\varphi_\epsilon(P;q)\big) \xi_\epsilon(P;q) - \sfL\big(q, \xi_\epsilon(P;q)\big) \Big) 
 \geq \bar{\sfH}(P) -\epsilon.
\end{align*}
Secondly, writting
\begin{align*}
g(\phi, \bnu):= \int_{\R^{2d}} \inf_{q \in \R^d} \Big(\big(P+ \nabla_q \phi (y,P;q)\big) \xi_\epsilon(P;q)
      - \sfL\big(q, \xi_\epsilon(P;q)\big) \Big) \bnu(dy,dP),
\end{align*}
then $(\phi,\bnu) \mapsto g$ is concave-convex-like in the sense of \cite{Sion58}, and $\bnu \mapsto g$ is lower semi-continuous with each $\phi$ fixed. Sion's minimax Theorem $4.2^\prime$ in \cite{Sion58} applies. We arrive at
\begin{align*}
   \LHS   \geq \sup_{\phi \in {\mathcal F}_0} \inf_{\bnu \in {\mathcal N}}
     g(\phi, \bnu)   
      =  \inf_{\bnu \in {\mathcal N}} \sup_{\phi \in {\mathcal F}_0} g(\phi, \bnu),
\end{align*}
where the inequality follows from $\sfH$ being Legendre transform of $\sfL$, and the equality follows from the minimax theorem.
Third, applying the same arguments as in the proof of Lemma~\ref{infintH}, we have that
\begin{align*}
  \sup_{\phi \in {\mathcal F}_0}  g(\phi, \bnu)  
  =  \int_{\R^{2d}} \sup_{\varphi \in C_\per^\infty(\R^d)}\inf_{q \in \R^d} 
 \Big(\big(P+ \nabla_q \varphi (q)\big) \xi_\epsilon(P;q)
      - \sfL\big(q, \xi_\epsilon(P;q)\big) \Big) \bnu(dy,dP). 
\end{align*}
 
Combine the above three steps together,
\begin{align*}
 \LHS & \geq   \inf_{\bnu \in {\mathcal N}}  \int_{\R^{2d}} \bar{\sfH}(P) \bnu(dy,dP) -\epsilon.
\end{align*}
By arbitrariness of the $\epsilon>0$, we conclude.
\end{proof}

\end{document}